%% file: main.tex
\documentclass[en,phd]{ldlthesis}

\usepackage[dvips]{color}
\usepackage{algorithm}
\usepackage{algpseudocode} 
\usepackage{comment} 
\usepackage{url} 
\usepackage{setspace}

\university{Universidade Federal de Pernambuco}

\institute{Centro de Ciências Exatas e da Natureza}

\program{Pós-graduação em Matemática Computacional}

\majorfield{Topologia Computacional}

\title{\textsc{Blink}: a language to view, recognize, \\ classify and manipulate 3D-spaces}

\date{January, 2007}

\author{Lauro Didier Lins}

\adviser{Sóstenes Luiz Soares Lins}

 \include{macros}

\begin{document}

\frontmatter


\thispagestyle{empty}

\phantom{Lauro} 

\vspace{0.5cm}

\begin{center}
\begin{large}

{\Large \textsc{Blink}: a language to view, recognize, \linebreak classify and manipulate 3D-spaces} \\
\vspace{1.4cm}
{\large by} \\
\vspace{1.4cm}
{\large Lauro Didier Lins}

\vspace{5cm}
A thesis presented to the Universidade Federal de Pernambuco \linebreak
in partial fulfillment of the requirements for the 
degree of \linebreak Doutor em \textsc{Matemática Computacional}

\vspace{4.5cm}
Recife, January, 2007

\end{large}
\end{center}
\vfill


\begin{dedicatory}
to  Sofia
\end{dedicatory}

\include{acknowledgements}



\include{resumo}

\include{abstract}


\tableofcontents

\listoffigures


\mainmatter

\include{chapter1}

\include{chapter2}
\include{chapter3}
\include{chapter4}
\include{chapter5}
\include{chapter6}


\backmatter

\appendix

\include{appendix1}
\include{appendix2}
\include{appendix3}
% ...

\setlength{\topmargin}{0cm}
\nocite{*}
\bibliographystyle{alpha}
\bibliography{references}


\end{document}

%% file: macros.tex
\newcommand{\IS}{\mathbb{S}}
\newcommand{\IR}{\mathbb{R}}

\newcommand{\IN}{\mathbb{N}}
\newcommand{\IZ}{\mathbb{Z}}
\newcommand{\IC}{\mathbb{C}}
\newcommand{\B}{\mathcal{B}}

\newcommand{\C}{\mathcal{C}}
\newcommand{\K}{\mathcal{K}}

\newcommand{\red}{\color{red}}

\def\mapright#1#2{\smash{\mathop{\hbox to 0.10cm{\rightarrowfill}}\limits^{#1}_{#2}}}
\def\mapleft#1#2{\smash{\mathop{\hbox to 0.10cm{\leftarrowfill}}\limits^{#1}_{#2}}}
\newcommand{\Aequiv}{\buildrel {\scriptstyle A} \over {\sim}}
\newcommand{\Sequiv}{\buildrel {\scriptstyle S} \over {\sim}}
\newcommand{\Requiv}{\buildrel {\scriptstyle R} \over {\sim}}
\newcommand{\ie}{{\it i.e.}\ }
\newcommand{\eg}{{\it e.g.}\ }

%% file: acknowledgements.tex
\acknowledgements

\begin{center}{\it If I had the right to thank only one person, this
acknowledgement would be $\ldots$}
\end{center}

I would like to thank my supervisor Sóstenes Lins, who is also my
father. I owe to him many things: to exist is obviously the most
important one, but here I want to mention the opportunity that
he gave me to, using my skills, contribute to the yet mysterious
field of 3-dimensional spaces. He presented me
an unexplored and elegant way to present spaces and said:
``Let's use computers to explore this''. And this is what we
did, and I am very happy with the process and with the result we
achieved.

\begin{center}
{\it but this constraint does not exist, so I can continue$\ldots$}
\end{center}

Thank you

Sofia, for filling my life with joy; my mother Bernardete and
my sister Isis, for your unrestricted support always; the new
generation: Pedrinho, Fernandinha, Joãozinhozinho, Arthur and Mariana; Nadja, Jorginho,
Joãozinho, Adelaide, Eneida, Heloiza, Dulce, Niedja,
vovó Lourdes, vovô Luiz, vovó Myriam and vovô Lauro for being a family
that makes me feel beloved; Roberta, Maria, and Fred for growing my family.

I am also grateful to: Silvio Melo, for his correct quantum
invariant implementation that guided me;
Paulo Soares, my best mate in the courses and
also my analysis tutor; Jalila and Donald Pianto, from whom I
learned in our joint works during our common disciplines;
Professors Klaus and Francisco Cribari from whom I learned
some probability and statistics; Valéria, who kept me informed
of all bureaucracy and due dates; Sérgio Santa Cruz and Francisco
Brito, for their help with the hyperbolic plane; UFPE;
CAPES, for financial support. 

%% file: resumo.tex
\resumo

Um {\it blink} é um grafo plano onde cada aresta ou é vermelha ou é verde. Um
{\it espaço 3D} ou, simplesmente, um {\it espaço} é uma variedade 3-dimensional
conexa, fechada e orientada. Neste trabalho exploramos pela primeira vez
em maiores detalhes o fato de que todo blink induz um espaço e todo espaço é
induzido por algum blink (na verdade por infinitos blinks). Qual o espaço de
um triângulo verde? E de um quadrado vermelho? São iguais? Estas perguntas
foram condensadas numa pergunta cuja busca pela resposta guiou em grande parte
o trabalho desenvolvido: quais são todos os espaços induzidos por blinks
pequenos (poucas arestas)? Nesta busca lançamos mão de um conjunto
de ferramentas conhecidas: os {\it blackboard framed links} (BFL),
os {\it grupos de homologia}, o {\it invariante quântico} de Witten-Reshetikhin-Turaev,
as {\it 3-gems} e sua teoria de simplificação. Combinamos a estas ferramentas
uma teoria nova de decomposição/composição de blinks e, com isso, conseguimos
identificar todos os espaços induzidos por blinks de até 9 arestas (ou BFLs
de até 9 cruzamentos). Além disso, o nosso esforço resultou
também num programa interativo de computador chamado \textsc{Blink}.
Esperamos que ele se mostre útil no estudo de espaços e, em particular,
na descoberta de novos invariantes que complementem o invariante quântico
resolvendo as duas incertezas deixadas em aberto neste trabalho.

\vspace{1cm}

\par\vskip\baselineskip\noindent{{\bf Palavras-chave}:
topologia, 3-variedades fechadas conexas e orientadas, 
grafos planos, espaços, {\it graph encoded manifolds}.
}

%% file: abstract.tex
\abstract

A {\it blink} is a plane graph with its edges being red or green. A
{\it 3D-space} or, simply, a {\it space} is a connected, closed and oriented
3-manifold. In this work we explore in details, for the first time,
the fact that every blink induces a space and any space is induced by
some blink (actually infinitely many blinks). What is the space of a green triangle?
And of a red square? Are they the same? These questions were condensed into
a single one that guided a great part of the developed work: what are all
spaces induced by small blinks (few edges)? In this search we used a known
set of tools: the {\it blackboard framed links} (BFL),
the {\it homology groups}, the {\it quantum invariant} of Witten-Reshetikhin-Turaev,
the {\it 3-gems} and its simplification theory. Combining these tools with a
new theory of decomposition/composition of blinks we could identify all
spaces induced by blinks with up to 9 edges (or BFLs with up to 9 crossings).
Besides that, our effort resulted in an interactive computer program named
\textsc{Blink}. We hope that this program becomes useful in the study
of spaces, in particular, in the discovery of new invariants that
complement the quantum invariant and homology group solving the two
uncertainties that we left open in this work.

\vspace{1cm}

\par\vskip\baselineskip\noindent{{\bf Keywords}: 
topology, closed connected oriented 3-manifolds, plane graphs, 
spaces, graph encoded manifolds.
}

%% file: chapter1.tex
\chapter{Introduction}

\vspace{-1cm}
\section{Initial motivation}
\label{sec:initialMotivation}

Unexplored simplicity. This was the reason for the birth of this work.
Let me explain. Topology deals with, among other objects, the so
called {\em 3-manifolds}. A basic type of 3-manifold is a {\em
closed, connected, oriented 3-manifold} and, in this work, the term
{\em space} will be used as a synonym for it.\footnote{The term
``space'' is also used as a synonym for a 3-manifold (any
3-manifold), but here we will use it as a synonym only for a closed,
connected, oriented 3-manifold.} In 1994, Kauffman and Lins
\cite{KauffmanAndLins1994} introduced a new way to present spaces.
They named their new presentation {\em blinks}. In fact, on their
work, blinks were appreciated as being a new concise presentation
for spaces, but its main importance there was its use as an
intermediate step to convert a {\em blackboard framed link}
presentation of a space into a {\em 3-gem} presentation of the same
space. This conversion, among other results, is central to this work,
but we will get to it later. What I want to say now is what
attracted me to this object named blink. Here it is.
\vspace{-0.3cm}
{\begin{center} \it
Blink states that a triangle is a space. That a square is a space.
That any plane graph (\ie any drawing in a plane made of
points connected by curves that do not intersect) is a space.
\end{center}}
\vspace{-0.3cm}
\noindent What space has the form of a
triangle? And the form of a square? I found this connection
between plane graphs and spaces very elegant. This elegance
becomes even more special when we see that other space
presentations have more ``complicated'' drawings than
blink presentation does. For example, the following three
drawings are different presentations for space $\IS^3\!/\!<3,3,2>$. \footnote{That is, the quotient of the $3$-sphere
under the action of the binary tetrahedral subgroup, which is a non-abelian group of order 12.}
\begin{center}
\includegraphics[width=8cm]{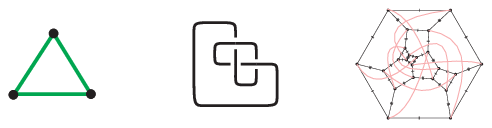}
\end{center}
The first is a blink presentation, the second is a blackboard
framed link presentation and the third is a 3-gem presentation.
In this example it is clear that the blink presentation is
simpler. Its perceptual complexity is smaller. Indeed this is
always the case. No other space presentation has simpler drawings
than blink does. Blackboard framed links don't, 3-gems don't.
{\it special spines} don't, {\it Heegaard diagrams} don't.

This simplicity aspect of blinks allied to the fact that they were
not studied before, except as a by product on the books
\cite{KauffmanAndLins1994, Lins1995} made me decide to explore it.
To see the greatest number of spaces through blink
drawings was the initial intent, but the hope was also that
so much elegance and simplicity is not in vain and it could,
by the fact that it was not explored before, hold
some yet unknown secret of spaces.

\section{Historical overview}
\label{sec:historicalOverview}

In 1962 Lickorish \cite{Lickorish1962} revolutionized the area of
spaces by proving in a purely combinatorial way a result first
published by Wallace \cite{Wallace1960} two years before by means of
differential topology. The result has, as a corollary, the following
fact:

{\em Given any space $M$, there exists a finite number $k$ of
disjoint solid tori $T_i\subset M$ such that $M \backslash
\bigcup_{i=1}^k T_i$ is homeomorphic to
$\IS^3\backslash\bigcup_{i=1}^k T'_i$, for a corresponding set of
disjoint solid tori $T'_i \subset \IS^3$.}

Another consequence of Lickorish's Theorem permits the
presentation of an arbitrary space as a link diagram in the plane
with an integer attached to each one of the components, the so
called {\em framed links}. The revolution was completed in 1978 when
Kirby \cite{Kirby1978} published his now famous calculus on framed links.
This paper spurred an enormous activity in the
area paving the way to prove deep theorems by a specific kind of
diagrammatic calculus with well defined rules.

Framed links can be freed from the integers associated to their
components if one introduces curls in their projections so as to
have the frame of a component equal to its writhe, thus producing a
{\em blackboard framed link} or a {\em BFL}. In this way, every
space presents itself as a blackboard framed link. A reformulation
of Kirby's calculus into blackboard framed link language is presented
in \cite{Kauffman1991}: a formal BFL calculus. The importance of the
BFL presentation is
testified by the fact that we can obtain from
it the Witten-Reshetikhin-Turaev quantum invariant (or WRT-invariant) of
the induced space \cite{Witten1989,ReshetikhinAndTuraev1991}. This is
not possible, at least it is unknown at the present, from a
triangulation \cite{Hempel1976}, from a Heegaard diagram
\cite{Rolfsen1976} or from a {\em special spine} \cite{Matveev2003}
of the space. In this work the {\em QI} or quantum invariant of a
space mean their WRT-invariants computed at $A=e^{i\pi/2r}$
\cite{KauffmanAndLins1994}.

A BFL can be reformulated into the object named {\it blink} also
introduced in \cite{KauffmanAndLins1994}. The presentation by blink
is concise and permits, as BFLs does, the computation of best invariants
available. But how do we prove that two blinks yielding the same
invariants are manifestations of the same space? Kirby's moves, although a
theoretical masterpiece, are in this practical context of no use,
except in very limited circumstances. For this task, the TS-moves
and U-move on gems of Lins \cite{Lins1995} are much better. The
applicability of the gem simplifying dynamics is available because
from a blink it is straightforward to produce a gem inducing the
same space \cite{KauffmanAndLins1994}. Summarizing, the task of
classifying blinks has to rely on gem theory which is, at present,
an indispensable complement to perform the task. The topological
classification of gems with 30 vertices (extending the classification
of \cite{Lins1995} using again TS- and U-moves) has been recently
accomplished~\cite{CasaliCrist2006}.

\newpage

\section{What we did}

In Section~\ref{sec:initialMotivation} we said our intent was to see
the greatest number of distinct spaces through blinks (\ie blink drawings).
To be able to accomplish this task, the most important question we must know
how to answer is whether two blinks $A$ and $B$ are
actually presentations of the same space. A fact we did
not mention before is that any space has infinitely
many distinct blink presentations. So, to answer this question is
not just checking that $A = B$ or $A \neq B$. The relation
space-blink is not a one to one relation.

As mentioned in Section~\ref{sec:historicalOverview}, blink is a direct
reformulation of blackboard framed link (BFL). It is easy to go from a blink
presentation of a space to a BFL presentation of the same space and
vice-versa. We also mentioned that a formal calculus for BFLs is well
known. This
means that a set of BFL operations or moves is known such that any pair
of BFLs induce the same space if and only if there
is a path, \ie a finite sequence of operations or moves in this set, that
transforms one BFL into the other. So, we could answer our question like
this: obtain a BFL presentation for $A$, obtain a BFL presentation for
$B$, then show a path in the BFL calculus transforming $A$ into $B$, proving
that they induce the same space; or show that such a path does not exist, and
conclude that they induce different spaces. Although this approach is correct and theoretically
possible, it is not practical. How to show a proof that a path does not
exist? Despite of this practical gap, one of the contributions of
this thesis is a reformulation of the BFL calculus into a purely blink
calculus.

Section~\ref{sec:historicalOverview} also mentions that from a blink we
are able to calculate some space invariants. For example, it is possible
to obtain the homology group or the Witten-Reshetikhin-Turaev
quantum invariant of its space. Indeed, this is the first thing we do
to answer whether blink $A$ induces the same space as blink $B$. We
calculate these two invariants and if any of them are different
we can answer for sure that the blinks induce different spaces. But,
if they are both the same, we cannot say their spaces are the same.
No complete invariant for spaces is known. So, any known space
invariant may fail to distinguish different spaces.

When the space of blink $A$ is not distinguished from the space
of blink $B$ by the space invariants, we must use another tool to answer
our question. This other tool is 3-gem theory. The book
\cite{KauffmanAndLins1994} shows a way to obtain a 3-gem presentation
from a blink presentation inducing the same space. We improved this
algorithm in Chapter~\ref{chap:gems}. In \cite{Lins1995}
a nice algorithm to simplify 3-gems is presented. So the last thing
we do to answer our question is to check whether the
blinks $A$ and $B$ not distinguished by the space invariants have
their 3-gems simplified to a common 3-gem. If this is the case then
we are sure that $A$ and $B$ induce the same space. If not, then we
are not sure. It is a hint that they are different but this cannot
be said. For small blinks, as we will see later, there are only two
uncertainties left out of $\approx 500$.

This approach of testing the homology group and the WRT-quantum invariant
to distinguish blinks and then, if not distinguished,
applying the 3-gem simplifying algorithm to show that they induce
the same space was very successful in our experiments as we will
see in Chapter~\ref{chap:census}. Its only constraint is that
it works only for small blinks and 3-gems. The computational
effort to calculate quantum invariants or to simplify 3-gems is exponential
in the sizes of the blinks and of the 3-gems, respectively.

Let's return to our initial intent: isolating the largest number of
spaces through blinks. We already know how to test if two (small)
blinks do induce the same space or not. The next important thing
to define is for what blinks we are going to ask these questions.
To try all possible blinks is prohibitive and unnecessary. As we
will see, we can search for spaces in only a small fraction of all possible
blinks and yet not lose anything. To get to this optimization
we first developed a useful decomposition/composition theory
of blinks in Chapter~\ref{chap:blinks} (actually the theory was developed
for its combinatorial counterpart: the g-blink). Then, using the
results and operations of blink calculus and BFL calculus we filtered
some redundant blinks. This resulted in a small set of blinks for which we
could identify all spaces. We also could isolate interesting sets
of small blinks inducing the same spaces such that a path in blink
calculus or BFL calculus is not trivial to identify. These sets may
lead to new ideas for theorems or space invariants.

A contribution of this thesis is a computer program named \textsc{Blink}.
It was responsible for most of the figures in this document. It also
supports the most important concepts discussed in the following
chapters and we hope that it will become popular and help researchers and
students to learn, do research or just appreciate spaces through
the language of blinks, blackboard framed links or 3-gems which are
the objects that it supports at the moment.

\section{The structure of this thesis}

In Chapter~\ref{chap:topologyEtc} we begin with a
review of the basic topological language and concepts needed in the remaining of
the thesis. We then introduce {\it knots and links} and their diagrams.
After that we introduce {\it framed links} and
{\it blackboard framed links} (BFL) and show how they encode spaces.
A calculus for blackboard framed links is presented.

In Chapter~\ref{chap:blinks} we define the motivating object of this
work: a {\it blink}. We show that blinks are a simple
reformulation of blackboard framed links with the advantage of having
simple drawings. We then reformulate the BFL calculus shown in
Chapter~\ref{chap:topologyEtc} in blink language. From
blinks we define a new combinatorial object named {\it g-blinks}. We show
how to obtain the homology group and quantum invariant of
Witten-Reshetikhin-Turaev invariants from a g-blink. The code
of a g-blink is then presented. Then some involutions
of g-blinks are defined: dual, reflection and refDual. The concept
of a representative g-blink is introduced based on the previous
results shown for g-blinks. We end the chapter showing how to
identify all spaces that are induced by small blinks and what
is the missing piece to do that: a way to prove homeomorphisms.

In Chapter~\ref{chap:gems} we define 3-gems: another way to present
spaces. We then show some moves that can be done in
3-gems without changing the induced space. These moves yields a
viable computational way to prove homeomorphism of spaces: a
combinatorial simplification dynamics of 3-gems. To connect blinks
and 3-gems we show an improved way to obtain, from a blink,
a 3-gem inducing the same space. We finish this chapter with the
proof, via 3-gems, of a theorem on g-blinks stated in
Chapter~\ref{chap:blinks}: the partial reflection theorem.

In Chapter~\ref{chap:census} we present the computational
experiments and results that we have obtained. We define formally
what we are searching: a census of prime spaces. Then we construct
a small set named $U$ that has the {\it 9-prime-unavoidable} property.
We then show how we topologically identified the spaces of every
g-blink in $U$. We finish the chapter exploring another set of
g-blinks: simple 3-connected monochromatic blinks.

In Chapter~\ref{chap:conclusion} we review the main contributions
of this work, talk a little about the program \textsc{Blink} and
about a theoretical contribution that we did not finish on
time to this thesis: a polynomial algorithm to obtain the blink
of a 3-gem. Some research problems that can be explored as a
continuation of this work are also shown.

%% file: chapter2.tex
\chapter{Topology, manifolds, links and blackboard framed links}
\label{chap:topologyEtc}

\section{Topology, manifolds and what we call here ``spaces''}
\label{sec:topology}

In topology the shape of a cup of coffee is equivalent to the shape
of a doughnut. Everybody knows that if we try to put coffee on a
doughnut the result is not the same as if we try to put coffee on a
cup of coffee. So, our first conclusion is: what topology states as
``equivalent shapes'' is definitively not aligned with our practical
understanding of equivalent shapes. One of the main problems that
topology deals is classifying shapes as equivalent or not and, at
the end, describing what are all possible shapes. In this section,
based on the clean and direct approach of \cite{elementaryTopology},
we present an introduction to elementary topology to settle the
vocabulary and the basic concepts needed. At the end we define what
are {\it manifolds} and the specific class of closed, connected,
oriented manifolds which are the ``shapes'' that we are interested
in this work.

\begin{figure}[htp]
   \begin{center}
      \leavevmode
      \includegraphics[width=5cm]{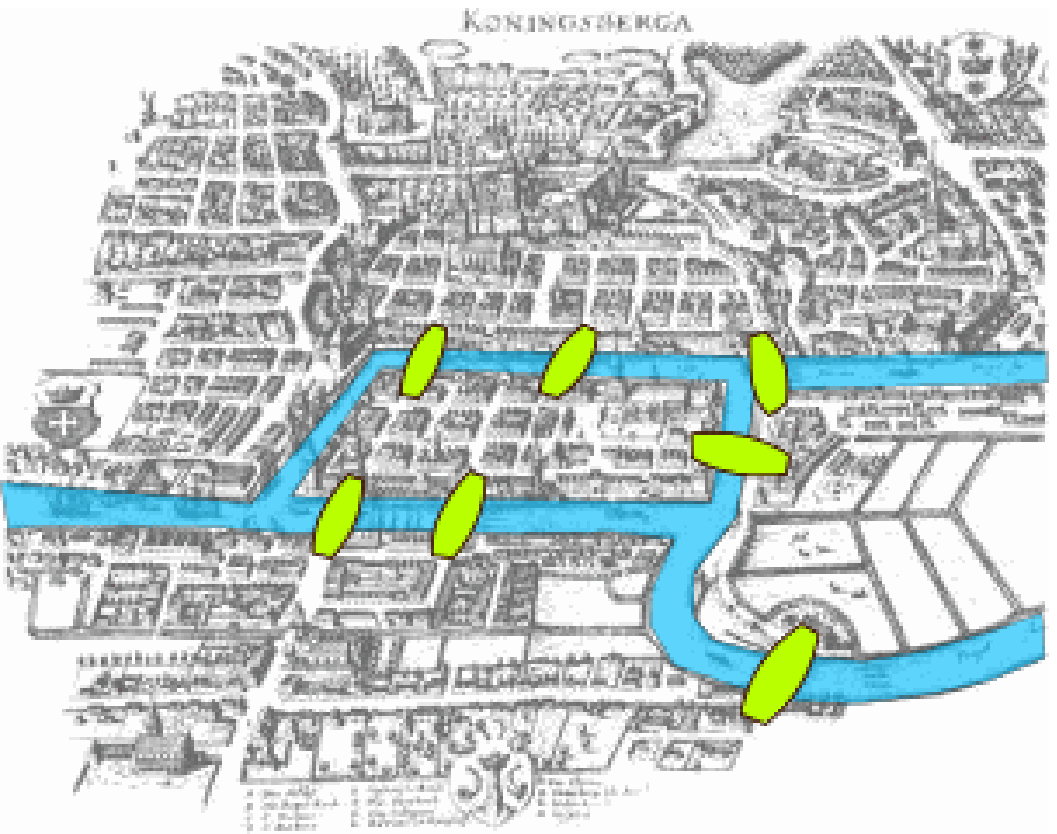}
      \hspace{1cm}
      \includegraphics[width=5cm]{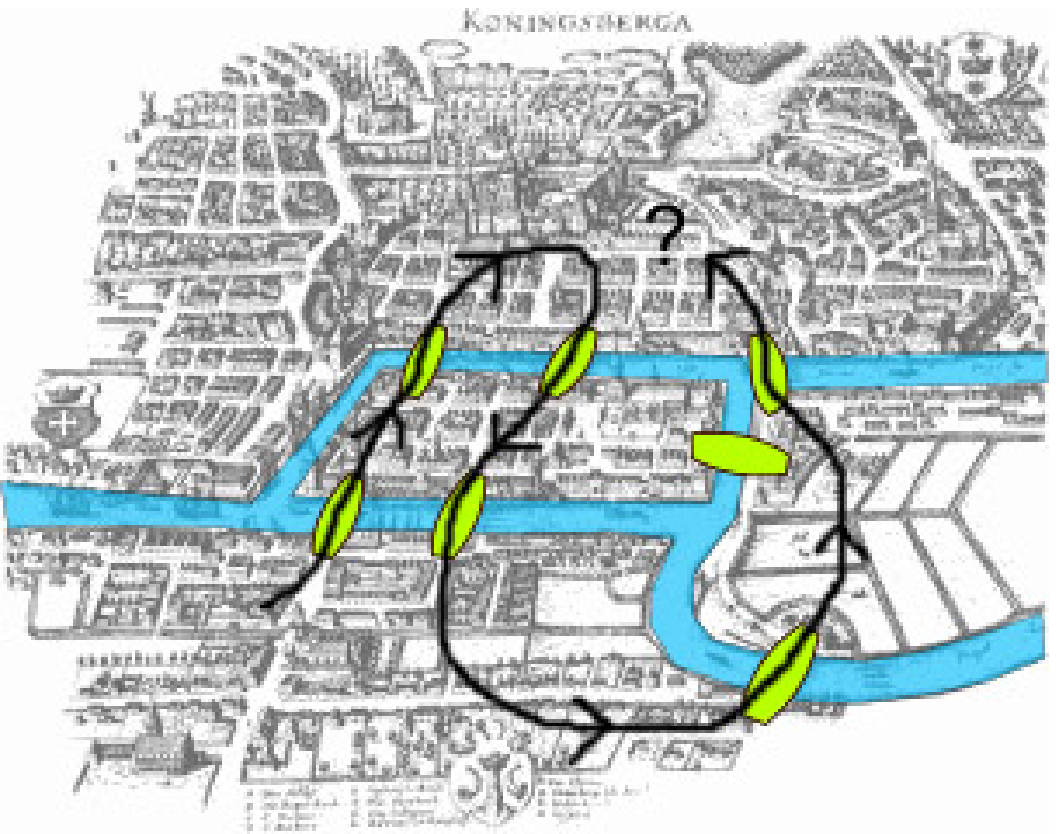} \\
   (A) \hspace{5.5cm} (B)\\
   \end{center}
   \vspace{-0.7cm}
   \caption{Seven bridges of Königsberg}
   \label{fig:sevenBridges}
\end{figure}

In the XVIII century, the city of Königsberg, Prussia (now
Kaliningrad, Russia) had seven bridges over the Pregel river
connecting two islands and other parts of the city as is shown in
Figure~\ref{fig:sevenBridges}A. A famous problem concerning
Königsberg was whether it was possible to take a walk through the
town in such a way as to cross over every bridge only one time.
Figure~\ref{fig:sevenBridges}B shows a wrong walk attempt: by the
time the sixth bridge is crossed the only uncrossed bridge is
unreachable.

No one was able to do this walk, and yet nobody knew how to prove
that it could not be done. In 1735, some college students sent this
problem to Leonhard Euler, one of the greatest mathematician of all
time. Euler was able to prove mathematically that this walk was
impossible. This result did not depend on the lengths of the
bridges, nor on their distance from one another, but only on
connectivity properties: which bridges are connected to which
islands or riverbanks. What Euler captured with the ``Problem of the
Seven Bridges of Königsberg'' is the motivating insight behind
topology:
\begin{center}\it
some geometric problems depend not on the exact shape of the objects
\linebreak involved, but rather on the ``way they are connected
together''.\end{center} Leonhard Euler's 1736 paper on Seven Bridges
of Königsberg is regarded as one of the first topological results
and also led to graph theory, a branch of mathematics with
``infinite'' applications
\cite{wiki:topology,wiki:sevenBridgesOfKonigsberg}.

Topology, in its present form, long after Euler, uses the term {\it
topological space} for what we called a ``shape'' on the beginning
of this section. Before defining what are {\it topological spaces}
we define {\it metric spaces}, as they are the source for the
concrete ``shapes'' or {\it topological spaces} that we are
interested.

\bigskip \centerline{\bf \textsc{Metric Spaces}} \bigskip

A {\it metric} or a {\it distance} in a set $X$ is a function $\rho:
X \times X \rightarrow \IR_+ = \{x \in \IR \, | \, x \geq 0 \}$ that
satisfies $$\begin{array}{l}
  \hbox{(1) \,} \rho(x,y) = 0, \hbox{ iff }x = y,\\[0.2cm]
  \hbox{(2) \,}\rho(x,y) = \rho(y,x), \hbox{ for every }x,y \in X,\\[0.2cm]
  \hbox{(3) \,}\rho(x,y) \leq \rho(x,z) + \rho(z,y), \hbox{ for every }x,y,z \in X. \hbox{ ({\it triangle inequality})}\\
\end{array}$$
The pair $(X,\rho)$, where $\rho$ is a metric in $X$, is called a
{\it metric space}. The function
$$\IR^n \times \IR^n \rightarrow \IR_{+} : (x,y) \mapsto
\sqrt{\sum_{i=1}^{n}{(x_i-y_i)^2}}$$ is a metric in $\IR^n$ and is
called {\it euclidean} metric.

Let $(X,\rho)$ be a metric space, let $a$ be a point in $X$, and let
$r$ be a positive real number. The sets
$$\begin{array}{l}
  B_r(a) = \{ x \in X \, | \, \rho(a,x) < r \},\\[0.2cm]
  D_r(a) = \{ x \in X \, | \, \rho(a,x) \leq r \},\\[0.2cm]
  S_r(a) = \{ x \in X \, | \, \rho(a,x) = r \}\\[0.2cm]
\end{array}$$
are called, respectively, {\it open ball}, {\it closed ball}, and
{\it sphere} of the space $(X,\rho)$ with center $a$ and radius $r$.
If $(X,\rho)$ is a metric space and $A \subset X$, then the
restriction of metric $\rho$ to $A \times A$ is a metric in $A$, and
$(A,\rho|_{A \times A})$ is a metric space. It is called a {\it
subspace} of $(X,\rho)$. The ball $D_1(0)$ and the sphere $S_1(0)$
in $\IR^n$ with the euclidean metric are denoted by symbols $D^n$
and $S^{n-1}$ and called {\it $n$-dimensional ball} and {\it
$(n-1)$-dimensional sphere}. They are considered as metric spaces
with the metric restricted to $\IR^n$. Note that: $D^1$ is the
segment $[-1,1]$; $D^2$ is a disk; $S^0$ is the pair of points
\{-1,1\}; $\IS^1$ is a circle; $\IS^2$ is a sphere; $D^3$ is a ball. The
words disk, circle, sphere and ball were used, in last sentence,
appealing to their common sense meaning. Now, for this work,
they have a formal meaning: a {\it disk} is $D^2$, a {\it circle} is
$\IS^1$, a {\it sphere} is $\IS^2$ and a {\it ball} is $B^3$.

\newpage

\bigskip \centerline{\bf \textsc{Topological Spaces}} \bigskip

A {\it topological space} is a set $X$ with a collection $\Omega$ of
subsets of $X$ satisfying the following three axioms:

\smallskip

\begin{tabular}{l}
(1) \,\, the empty set $\varnothing$ and $X$ are in $\Omega$, \\[0.2cm]
(2) \,\, the union of any collection of sets in $\Omega$ is in $\Omega$, \\[0.2cm]
(3) \,\, the intersection of any pair of sets in $\Omega$ is in $\Omega$. \\[0.2cm]
\end{tabular}

The collection $\Omega$ is called a {\it topological structure} or a
{\it topology} in $X$. The sets in $\Omega$ are called {\it open}.
The elements of $X$ are called {\it points}. A set $F \in X$ is said
{\it closed} in the space $(X,\Omega)$ if its complement $X
\backslash F$ is open ({\it i.e. $X \backslash F \in \Omega$}). Note
that $\varnothing$ and $X$ are both open and closed. A {\it
neighborhood} of a point is any open set containing that point. A
collection $\Sigma$ of open sets is called a {\it base} for a
topology ({\it i.e.} topological structure), if each nonempty open
set is a union of sets belonging to $\Sigma$.

The following result connects metric spaces and topological spaces:
\begin{center}
\it the collection of all open balls in a metric space $(X,\rho)$
\linebreak is a base for some topology in $X$.
\end{center}
For example, consider $\IR^2$ with the euclidean metric. Then, a
topology for $\IR^2$ is the set of all unions of open balls (open
disks in the plane). This topology is the default topology when
nothing else is mentioned.

Let $(X,\Omega)$ be a topological space, and $A \subset X$. Denote
by $\Omega_A$ the collection of sets $A \cap V$, where $V \in
\Omega$. Then,
\begin{center} \it
$\Omega_A$ is a topological structure in $A$.
\end{center}
The pair $(A,\Omega_A)$ is called a {\it subspace} of the space
$(X,\Omega)$. The collection $\Omega_A$ is called the {\it subspace
topology} or the topology {\it induced} on $A$ by $\Omega$, and its
elements are the open sets in $A$.

At this point, we can think, for instance, of $\IS^2$ as a topological
space. We know that the collection of open balls of $\IR^3$ (as a
metric space with the euclidean metric) is a base for a topology in
$\IR^3$. Consider this topology to view $\IR^3$ as a topological
space. Restrict this topology of $\IR^3$ to $\IS^2$ to obtain a
topology for $\IS^2$: $\IS^2$ is now a topological space. In this work
this logical sequence to obtain a topology for a subset of $\IR^n$
is always the one considered. So, from now on, every subset of
$\IR^n$ may also be viewed as a topological space. For example the
surface of doughnut and of the coffee cup considered in the
beginning of this section may now be viewed as subsets of $\IR^3$
and, consequently, as topological spaces.

\bigskip \centerline{\bf \textsc{Maps}} \bigskip

In the context of topology, the terms {\it map} and {\it mapping}
are synonyms of function. A mapping $f: X \rightarrow Y$ is called a
{\it surjective map}, or just a {\it surjection} if every element of
$Y$ is an image of at least one element of $X$. It is called an {\it
injective map}, {\it injection} or {\it one-to-one map} if every
element of $Y$ is an image of, at most, one element of $X$. A
mapping is called a {\it bijective map}, {\it bijection}, or {\it
invertible} if it is surjective and injective.

The {\it image} of a set $A \subset X$ under a map $f: X \rightarrow
Y$ is the set of images of all points of $A$. It is denoted by
$f(A)$. Thus,
$$f(A) = \{ f(x) : x \in A \}.$$
The image of the entire set $X$ ({\it i.e.} $f(X)$) is called the
{\it image} of $f$. The {\it preimage} of a subset of $B \subset Y$
under map $f: X \rightarrow Y$ is the set of elements of $X$ whose
images belong to $B$. It is denoted by $f^{-1}(B)$. Thus,
$$f^{-1}(B) = \{ x : f(x) \in B \}.$$

\bigskip \centerline{\bf \textsc{Continuous Maps}} \bigskip

Let $X$, $Y$ be topological spaces. A map $f:X \rightarrow Y$ is
said to be {\it continuous} if the preimage of any open subset of
$Y$ is an open subset of $X$. A map $f:X \rightarrow Y$ is said to
be {\it continuous at point }$a \in X$ if for every neighborhood $U$
of $f(a)$ there exists a neighborhood $V$ of $a$ such that $f(V)
\subset U$. One result about continuous maps is that: a map $f:X
\rightarrow Y$ is continuous iff it is continuous at each point of
$X$. Another result is that this notion of continuity coincides with
the one that is usually studied in calculus:
\begin{center} {\it
Let $X, Y$ be metric spaces, and $a \in X$. A map $f:X \rightarrow
Y$ is \linebreak continuous at the point $a$, iff for every
$\epsilon
> 0$ there exists a $\delta > 0$ \linebreak such that for every point $x \in X$
inequality  $\rho(x,a) < \delta$ implies \linebreak
$\rho(f(x),f(a))<\epsilon$.}
\end{center}

\bigskip \centerline{\bf \textsc{Homeomorphism}} \bigskip

Now we are able to formally define the ``topologically equivalence''
concept. An invertible mapping is called a {\it homeomorphism} if it
is continuous and its inverse is also continuous. A topological
space $X$ is said to be {\it homeomorphic} to space $Y$ if there is
a homeomorphism $X \rightarrow Y$. Being homeomorphic is {\it an
equivalence relation}. Let $X$, $Y$ and $Z$ be topological spaces
then: (1) $X$ is homeomorphic to $X$; (2) if $X$ is homeomorphic to
$Y$ then $Y$ is homeomorphic to $X$; and (3) if $X$ is homeomorphic to
$Y$ and $Y$ is homeomorphic to $Z$ then $X$ is homeomorphic to $Z$.

Some examples of homeomorphic topological spaces: $[0,1]$ and
$[a,b]$ for any $a < b$; $(-1,1)$ and $\IR$; an open disk and the
plane $\IR^2$; $\IS^n \backslash \{\hbox{point in }\IS^n\}$ and $\IR^n$.
Some examples of non-homeomorphic topological spaces: balls $D^p,
D^q$ with $p \neq q$; spheres $S^p, S^q$ with $p \neq q$; punctured
plane $\IR^2 \backslash \{\hbox{point}\}$ and a plane with a hole
$\IR^2 \backslash \{(x,y):x^2+y^2<1\}$.

From the topological point of view homeomorphic spaces are
completely identical: a homeomorphism $X \rightarrow Y$ establishes
one-to-one correspondence between all phenomena in $X$ and $Y$ which
can be expressed in terms of topological structures. Thus, two
spaces are {\it topologically equivalent} or {\it the same for the
purposes of topology} if there is a homeomorphism between them.
There is a homeomorphism between the surface of a doughnut and the
surface of a coffee cup, so they are topologically equivalent.

As we pointed out on the first paragraph of this section, not yet in
the correct language, the topological equivalence problem or {\it
homeomorphism problem} is one of the classic and important problems
of topology:
\begin{center}
\it Given two topological spaces, are they homeomorphic?
\end{center}
To prove that topological spaces are homeomorphic, it is enough to
present a homeomorphism between them. Essentially this is what is
done in this case. However, to prove that topological spaces are not
homeomorphic, it does not suffice to consider any special mapping,
and is usually impossible to review all mappings. Therefore for proving
non-existence of a homeomorphism one uses indirect arguments. In
particular, one finds a property or a characteristic shared by
homeomorphic spaces such that one of the spaces has it, while the
other does not. Properties and characteristics shared by homeomorphic
spaces are called {\it topological properties} or {\it invariants}. For instance, the
cardinality of the set of points and of the set of open sets is a
topological invariant.

\bigskip \centerline{\bf \textsc{Embedding}} \bigskip

A continuous mapping $f: X \rightarrow Y$ is called a {\it
(topological) embedding} if the mapping \hbox{$f':X \rightarrow
f(X)$} is a homeomorphism, where $f'(x) = f(x)$ for all $x \in X$.
Embeddings $f_1, f_2: X \rightarrow Y$ are said to be {\it
equivalent} if there exist homeomorphisms $h_X:X \rightarrow X$ and
$h_Y:Y \rightarrow Y$ such that $f_2 \circ h_X = h_Y \circ f_1$.

Note that homeomorphisms are special kind of embeddings, where the
mapping is surjective.

\bigskip \centerline{\bf \textsc{Cover}} \bigskip

A collection $\Gamma$ of subsets of a set $X$ is called a {\it
cover} or a {\it covering} if the union of the elements of
$\Gamma$ contains $X$, {\it i.e.}, $X \subset \cup_{A \in \Gamma}A$.
A cover $\Gamma$ of a topological space $X$ is said to be an {\it
open cover} if every element of $\Gamma$ is an open set. A cover
$\Gamma$ of a topological space $X$ is said to be a {\it closed
cover} if every element of $\Gamma$ is a closed set. If $\Sigma$
covers $X$ and $\Gamma$ covers $X$ and $\Sigma \subset \Gamma$, then
$\Sigma$ is a {\it subcover} or {\it subcovering} of $\Gamma$.

\bigskip \centerline{\bf \textsc{Connectedness}} \bigskip

A topological space $X$ is said to be {\it connected} if it has only
two subsets which are both open and closed: $\varnothing$ and the
entire $X$. Although this definition is clear, at first, it is not
intuitive. Let's get a more intuitive definition. A {\it partition}
of a set is a cover of this set with pairwise disjoint sets. To {\it
partition} a set means to construct such a cover. Now the equivalent
notion of connectedness of a topological space:
\begin{center}
\it A topological space is connected iff it cannot be \linebreak
partitioned into two nonempty open sets iff it cannot be \linebreak
partitioned into two nonempty closed sets.
\end{center}
For instance, consider the topological space obtained as a subspace
of the plane that consists of two disjoint open disks (open balls)
({\it e.g}: one open ball $B_1(-1,-1)$ and $B_1(1,1)$). This
topological space is not connected, because the two open disks, that
are open sets, form a partition of the entire space.

A {\it connected component} of a space $X$ is a maximal connected
subset of $X$ ({\it i.e.} a connected subset, that is not contained
strictly in other larger subset of $X$). Some properties of
connected components: every point belongs to some connected
component; connected components are closed; two connected
components are disjoint or coincident. The image of a connected
space under a continuous mapping is connected, so connectedness is a
topological property. The number of connected components is a
topological invariant.

\bigskip \centerline{\bf \textsc{Compactness}} \bigskip

A topological space $X$ is said to be {\it compact} if any of its
open covers has a finite subcover, {\it i.e.} if $\Gamma$ is a cover
for $X$ then exists a finite $\Sigma \subset \Gamma$ that also
covers $X$. The image of a compact space by a continuous mapping is
also compact, so compactness is a topological property.

Compactness is a sort of topological counter-part for the property
of being finite in the context of set theory. Example of s non-compact
space: $\IR^n$. Example of a compact space: $\IS^n$. Indeed a subset of $\IR^n$ is compact
if and only if it is closed and bounded (\ie contained in an open ball).

\bigskip \centerline{\bf \textsc{Homotopy}} \bigskip

Let $f, g$ be continuous maps of a topological space $X$ to a
topological space $Y$, and $H: X \times [0,1] \rightarrow Y$ a
continuous map such that $H(x,0) = f(x)$ and $H(x,1) = g(x)$ for any
$x \in X$. Then $f$ and $g$ are said to be {\it homotopic} and $H$
is called a {\it homotopy} between $f$ and $g$. Homotopy of maps is
an equivalence relation: (1) if $f: X \rightarrow Y$ is a continuous
map then $H:X \times [0,1] \rightarrow Y$ defined by $H(x,t) = f(x)$
is a homotopy between $f$ and $f$; (2) if $H$ is a homotopy between
$f$ and $g$ then $H'$ defined by $H'(x,t) = H(x,1-t)$ is a homotopy
between $g$ and $f$; (3) if $H$ is a homotopy between $f$ and $f'$
and $H'$ is a homotopy between $f'$ and $f''$ then $H''$ defined by
$$ H''(x,t) =
 \left\{
 \begin{array}{ll}
 H(x,2t) & \hbox{for } t \leq 1/2,\\
 H'(x,2t-1) & \hbox{for } t \geq 1/2,\\
 \end{array}\right.$$ is a homotopy between $f$ and $f''$.

\bigskip \centerline{\bf \textsc{Isotopy}} \bigskip

Let $X$, $Y$ be topological spaces, $h, h': X \rightarrow Y$
homeomorphisms. A homotopy $h_t:X \rightarrow Y$, $t \in [0,1]$
connecting $h$ and $h'$ ({\it i.e.}, with $h_0 = h$, $h_1 = h'$) is
called a {\it isotopy} between $h$ and $h'$ if $h_t$ is a
homeomorphism for each $t \in [0,1]$. Homeomorphisms $h, h'$ are
said to be {\it isotopic} if there exists an isotopy between $h$ and
$h'$. Being isotopic is an equivalence relation on the set of
homeomorphisms $X \rightarrow Y$.

The concept of isotopy may also be applied to embeddings. Let $X$,
$Y$ be topological spaces, $h, h': X \rightarrow Y$ topological
embeddings. A homotopy $h_t:X \rightarrow Y$, $t \in [0,1]$
connecting $h$ and $h'$ ({\it i.e.}, with $h_0 = h$, $h_1 = h'$) is
called an {\it (embedding) isotopy} between $h$ and $h'$ if $h_t$ is
an embedding for each $t \in [0,1]$. Embeddings $h, h'$ are said to
be {\it isotopic} if there exists an isotopy between $h$ and $h'$.
Being isotopic is an equivalence relation on the set of embeddings
$X \rightarrow Y$.

A family $A_t$, $t \in I = [0,1]$ of subsets of a topological space
is called an {\it isotopy of the set $A = A_0$} if the graph $\Gamma
= \{(x,t) \in X \times I \, | \, x \in A_t\}$ of the family is {\it
fibrewise homeomorphic} to the cylinder $A \times I$, {\it i.e.}
there exists homeomorphisms $A \times I \rightarrow \Gamma$ mapping
$A \times \{t\}$ to $\Gamma \cap X \times \{t\}$ for any $t \in I$.
Such a homeomorphism gives rise to an isotopy of embeddings $ \Phi_t
: A \rightarrow X$, $t \in I$ where $\Phi_0$ is the identity mapping
and $\Phi_t(A) = A_t$. An isotopy of a subset is also called a {\it
subset isotopy}. Subsets $A$ and $A'$ of the same topological space
are said to be {\it isotopic in $X$}, if there exists a subset
isotopy $A_t$ of $A$ with $A' = A_1$. The isotopic relation over the
set of subsets of a topological space $X$ is an equivalence
relation.

An isotopy of a subset $A \in X$ is said to be {\it ambient}, if it
may be accompanied with an embedding isotopy $\Phi_t : A \rightarrow
X$ extendible to an isotopy $\tilde{\Phi_t}: X \rightarrow X$ of the
identity homeomorphism of space $X$. The isotopy $\tilde{\Phi_t}$ is
said to be {\it ambient} for $\Phi_t$. Two isotopic subsets of a
topological space may not be ambient isotopic. Any pair of circles
$\IS^1$ embedded in $\IS^3$ is isotopic, but a circle
(Figures~\ref{fig:links3d}A) and a trefoil
(Figures~\ref{fig:links3d}B) are not ambient isotopic.

\bigskip \centerline{\bf \textsc{Manifolds}} \bigskip

Let $n$ be a non-negative integer. A topological space $X$ is called
{\it locally euclidean space of dimension $n$} if each point of $X$
has a neighborhood homeomorphic either to $\IR^n$ or $\IR^n_{+}$
({\it i.e.} $\IR^n_+ = \{ x \in \IR^n : x_1 \geq 0 \}$, defined for
$n \geq 1$). Examples of locally euclidean spaces: $\IR^n$; $\IS^n$,
$D^n$.

A topological space is called {\it Hausdorff space} or just {\it
Hausdorff} if any two distinct points possess disjoint
neighborhoods. Result: any metric space is Hausdorff ({\it i.e.} the
topological space with topology induced from the metric space is
Hausdorff).

A space is said to satisfy the {\it second axiom of
countability} or to be {\it second countable} if it has a countable
base. Result: any metric space is Hausdorff ({\it i.e.} the
topological space with topology induced from the metric space has
a countable base).

A {\it manifold of dimension $n$} or $n$-manifold is a topological
space that satisfies:
\smallskip

\begin{tabular}{l}
(1) \,\, it is a locally euclidean space of dimension $n$, \\[0.2cm]
(2) \,\, it is Hausdorff, \\[0.2cm]
(3) \,\, it is second countable. \\[0.2cm]
\end{tabular}

\smallskip

Examples of {\it $n$-manifolds}: $\IR^n$; $\IS^n$, $D^n$.

The definitions until now were very formal, but this one will not be
formal. A manifold of dimension $n$ is called {\it non-orientable}
if it is possible to take the homeomorphic image of an
\hbox{$n$-dimensional} ball in the manifold and move it through the
manifold and back to itself, so that at the end of the path, the
ball has been reflected. The Möbius band and Klein bottle are the
most famous examples of non-orientable manifolds. A manifold which
is not non-orientable is {\it orientable}. An orientable space has
two orientations and the choice of one of them makes the space an
{\it oriented} space.

\bigskip \centerline{\bf \textsc{Spaces for us$\ldots...$}} \bigskip

In this work we are interested in studying a specific kind of
``shape'' or, as we learned, topological space. This is it:
\begin{center}
\it \large connected, closed, oriented 3-manifolds.
\end{center}
The adjective {\em closed} applied to a 3-manifold means that it is boundary free and compact.
We use, from now on, the word {\it space} to denote these
topological spaces. This is not a perfect choice, as ``space'', even
in mathematics, has a lot of meanings. It could be a metric space. A
vector space. A topological space not matching these constraints of
being compact, connected, oriented. A common use for space, for
instance, is any 3-dimensional topological space or any 3-manifold.
These include our spaces but others too. In spite of all these,
here, space will be exactly this: a connected, compact, oriented
3-manifold. So let's put it big:
\begin{center}
\it \large A {\huge space} in this work is the same as \linebreak a
connected, closed, oriented 3-manifold.
\end{center}

\newpage

\section{Knots, links and diagrams}
\label{sec:knotsAndLinks}

In general terms, {\it knot theory} studies the {\it placement
problem}. As stated in \cite{Kauffman1987}, this problem is
\begin{center}
\it Given topological spaces $X$ and $Y$, classify how $X$ may be
\linebreak placed within $Y$. Here the ``how'' is usually an
 embedding, and classify often \linebreak means up to some form of
movement of $X$ in $Y$ (isotopy, for example).
\end{center} When $X$ is the circle $\IS^1$ and $Y$ is
the 3-dimensional space $\IR^3$ or $\IS^3$, we have {\it classical
knot theory}. In this section we see, for this classical knot
theory, a characterization of what are the equivalent embeddings.
Things shown here form the basis for the approach to the problem of
characterizing homeomorphic spaces (connected, compact, oriented
3-manifolds) that we show later.

\begin{figure}[htp]
   \begin{center}
\begin{tabular}{ccccc}
      \includegraphics[width=2.5cm]{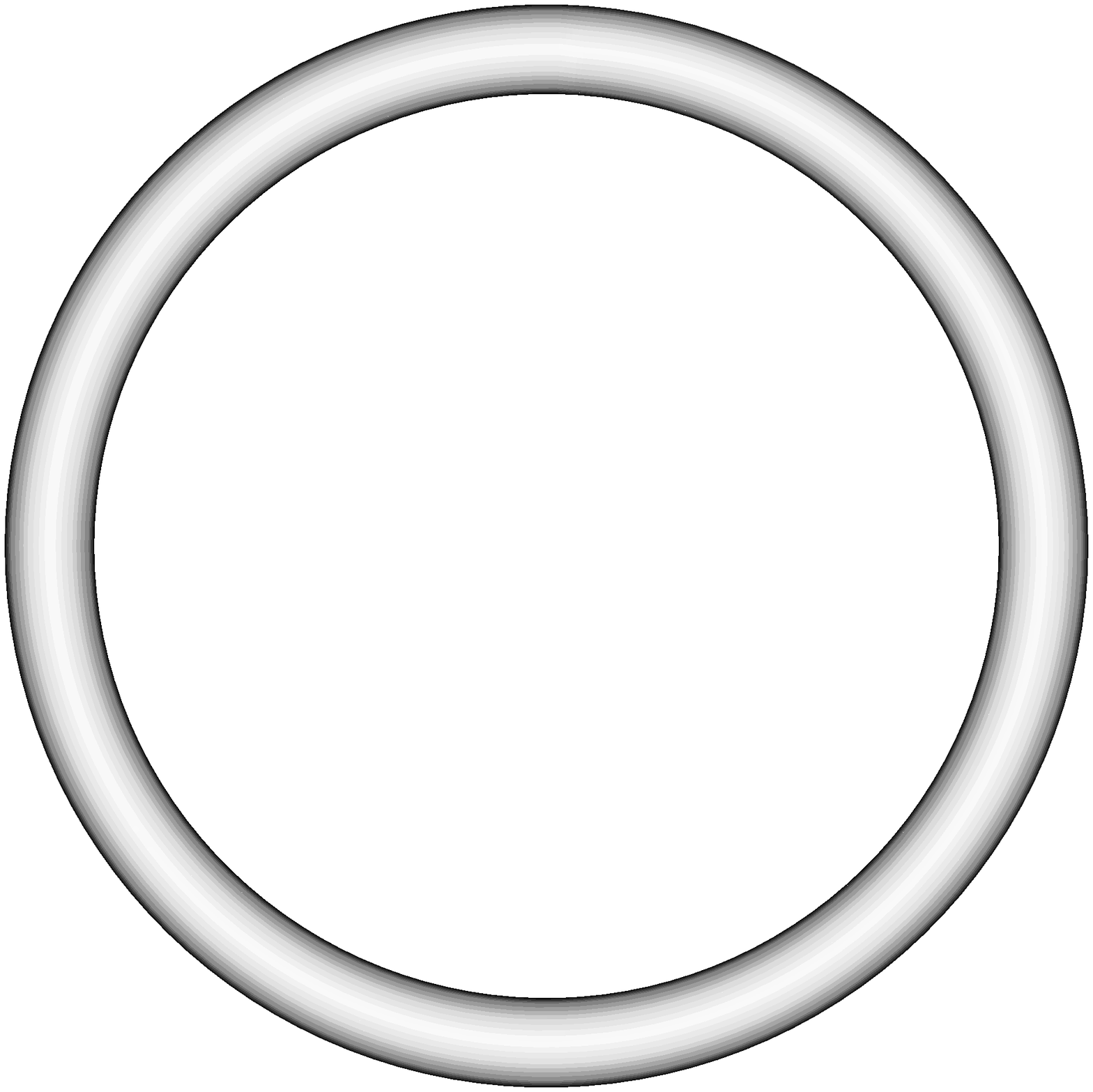} &
      \includegraphics[width=2.5cm]{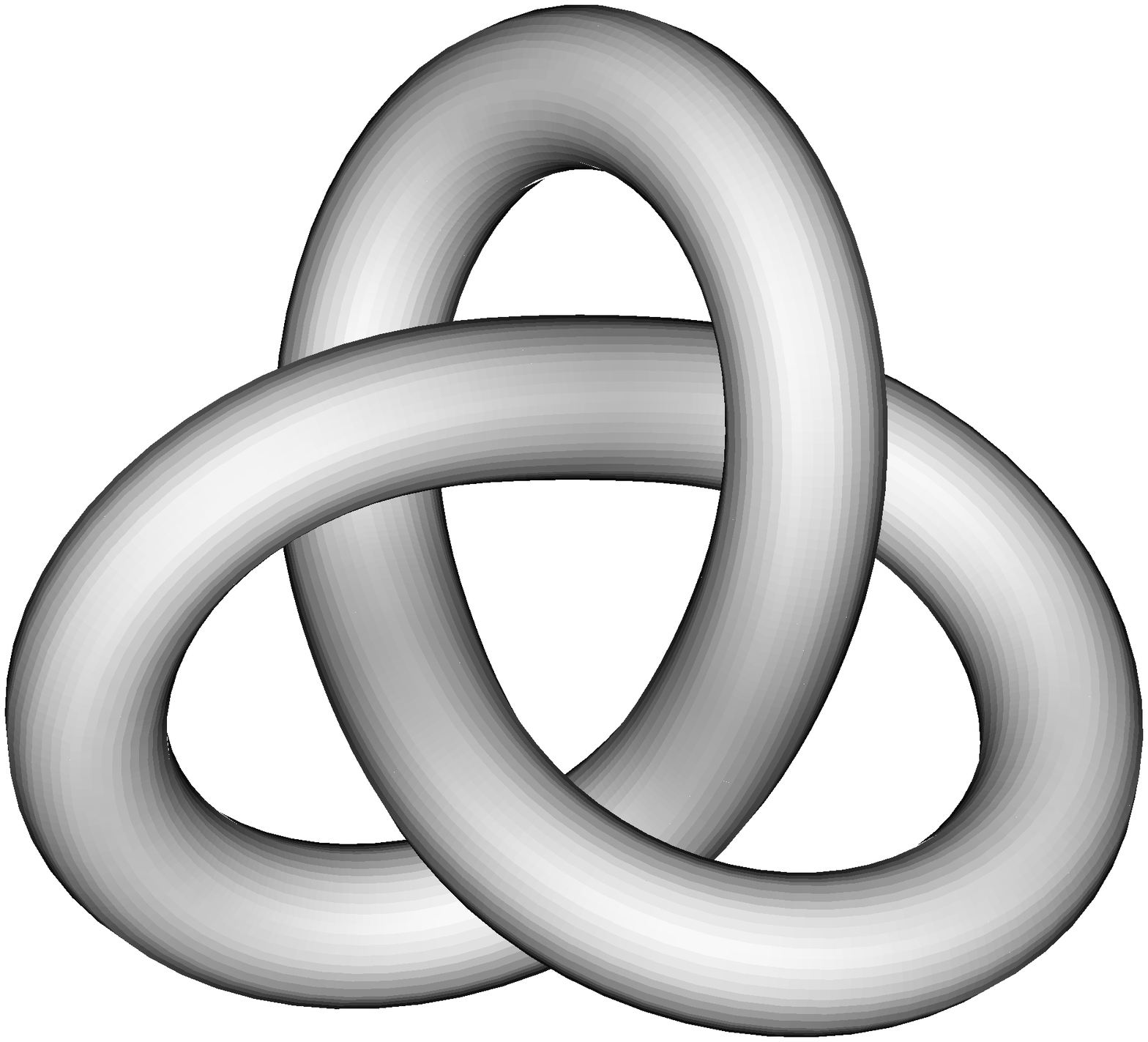} &
      \includegraphics[width=2.5cm]{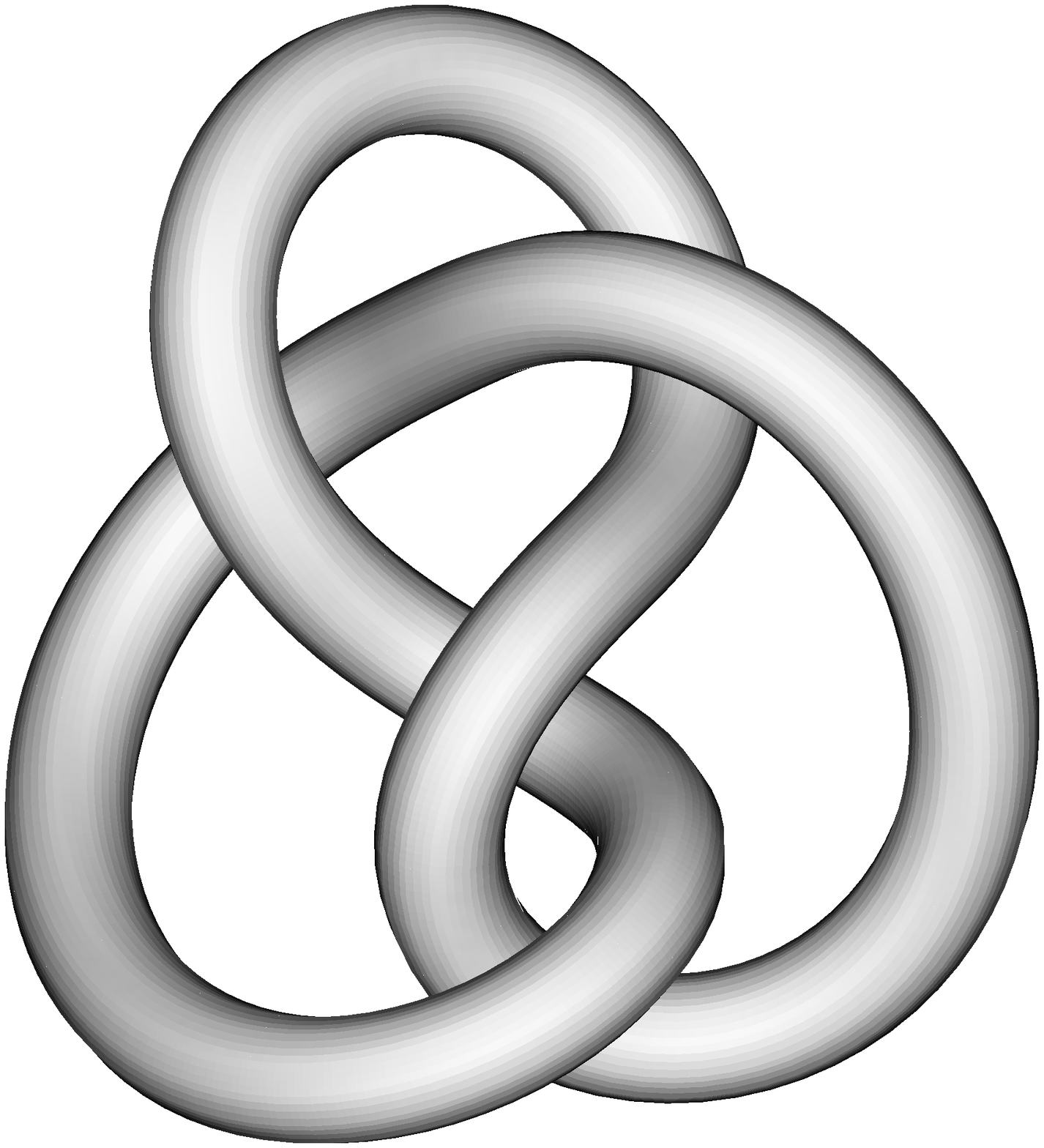} &
      \includegraphics[width=2.5cm]{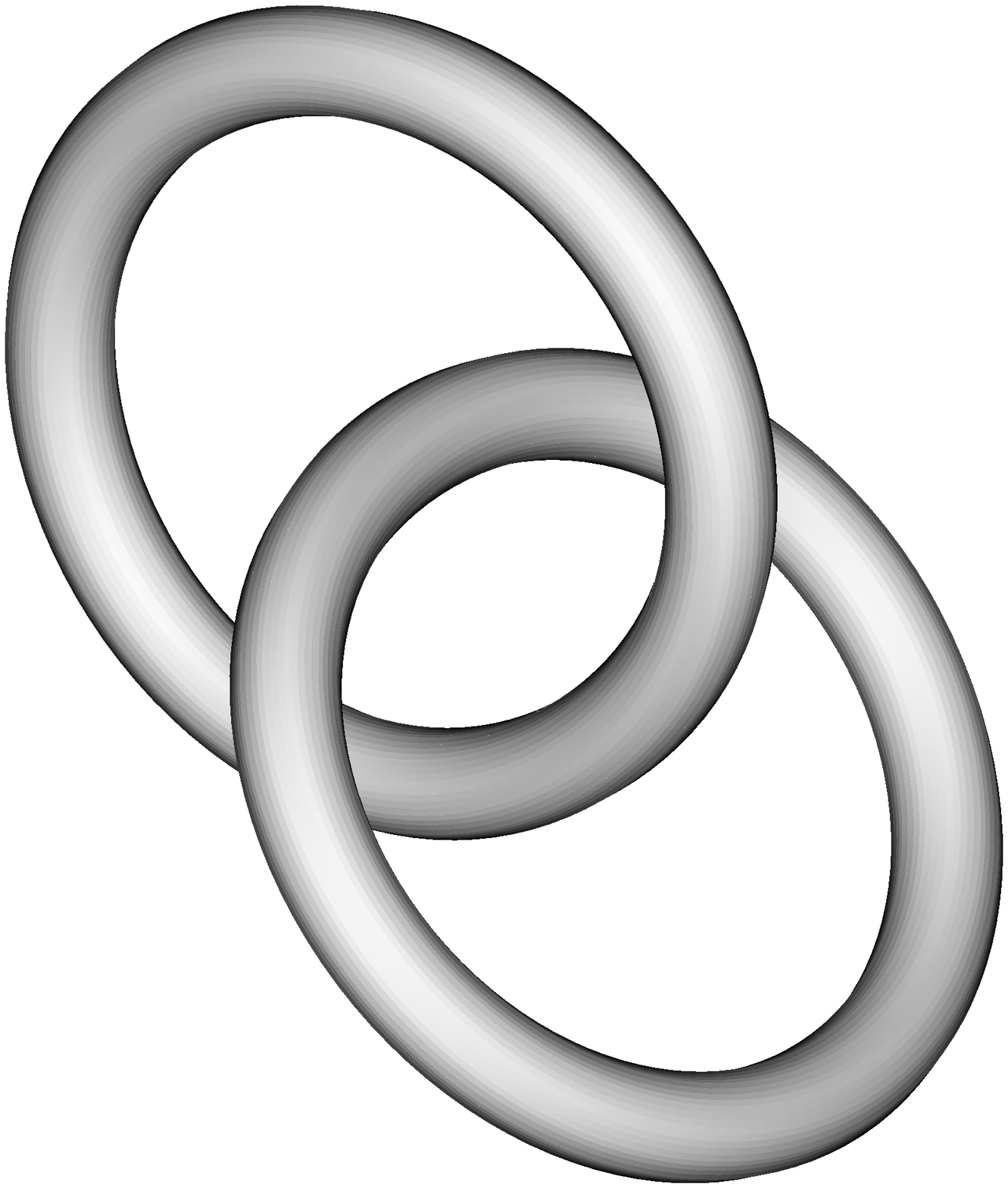} &
      \includegraphics[width=2.5cm]{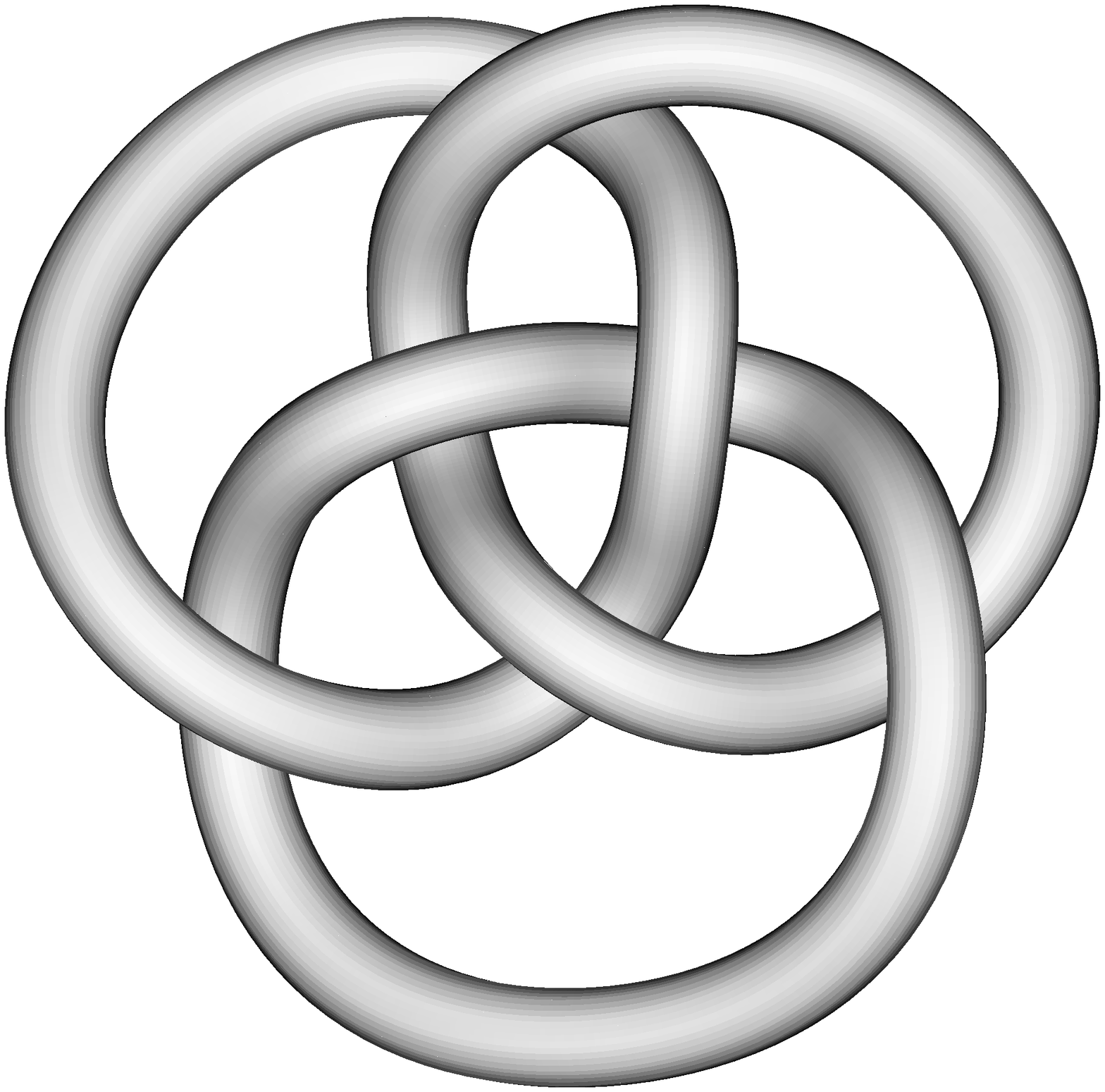} \\
      (A) & (B) & (C) & (D) & (E) \\
\end{tabular}
   \end{center}
   \vspace{-0.7cm}
   \caption{ Knots and links}
   \label{fig:links3d}
\end{figure}

An embedding of a circle $\IS^1$ in the 3-dimensional space $\IR^3$ or
in the 3-sphere $\IS^3$ is called a {\it knot}. An embedding of a
collection of circles in the 3-dimensional space $\IR^3$ or in the
3-sphere $\IS^3$ is called a {\it link}. Each circle (or the image of
one) in a link is called a {\it component} of the link. So, a knot
is a link with only one component. Figure \ref{fig:links3d} shows
some links\footnote{These figures were created using the beautiful
tool called \textsc{KnotPlot} that was part of the phd thesis of
Robert Scharein \cite{Scharein1998}.}. The link of Figure
\ref{fig:links3d}A is also a knot (one component) and is
suggestively called {\it unknot}. Figure \ref{fig:links3d}B is the
knot called {\it trefoil}. Figure \ref{fig:links3d}C is the knot
called {\it figure eight knot}. Figure \ref{fig:links3d}D is a
link with two components. Figure \ref{fig:links3d}E is a link with
three components and it is called the {\it borromean rings}. This link has
an interesting property: we cannot separate the three rings without
breaking one of them, {\it i.e.} the three rings are {\it linked},
even though, any two rings are separable without breaking ({\it
i.e.} any pair of rings is unlinked).

\begin{figure}[htp]
   \begin{center}
\begin{tabular}{ccccc}
      \includegraphics[width=2.5cm]{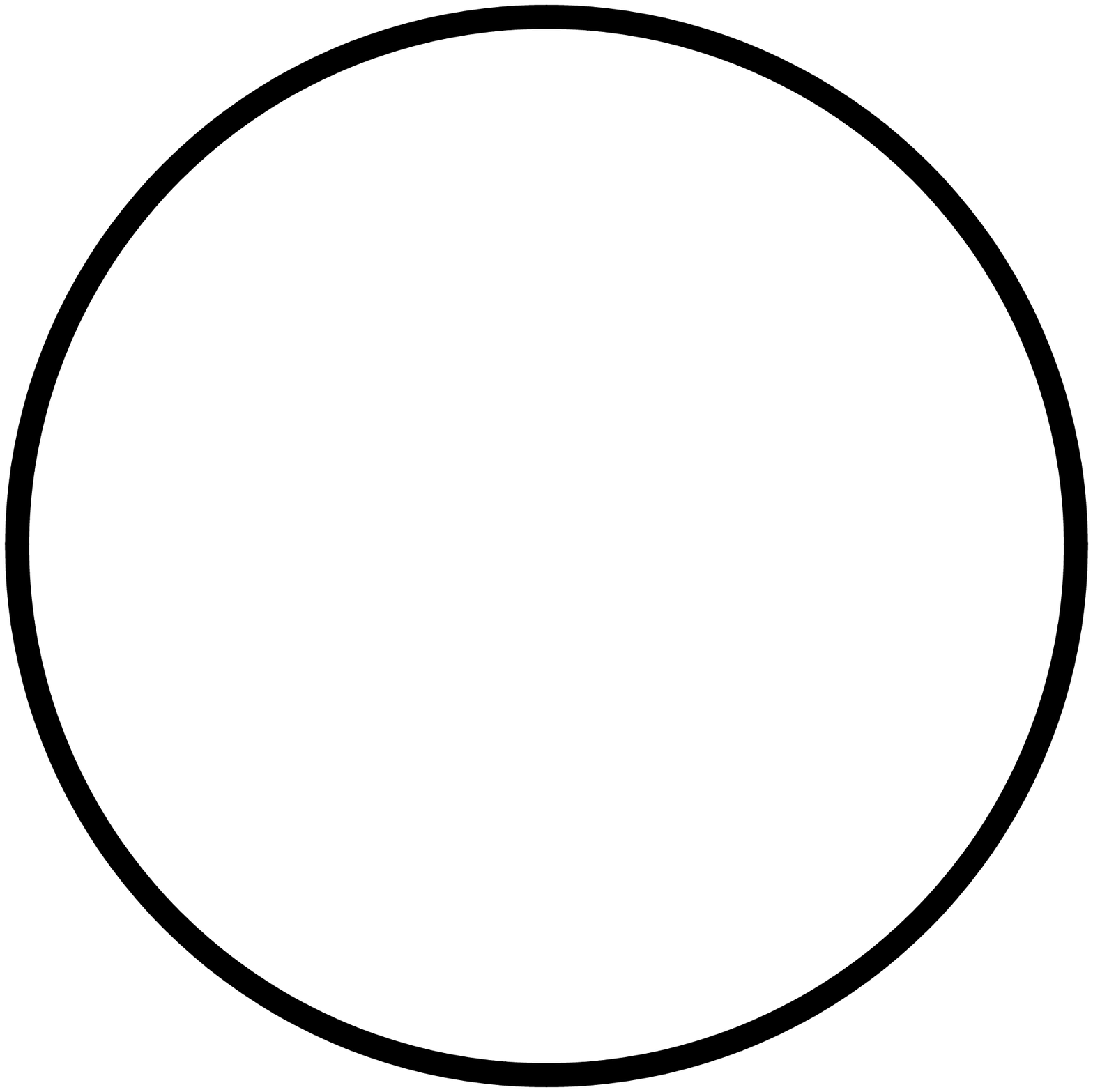} &
      \includegraphics[width=2.5cm]{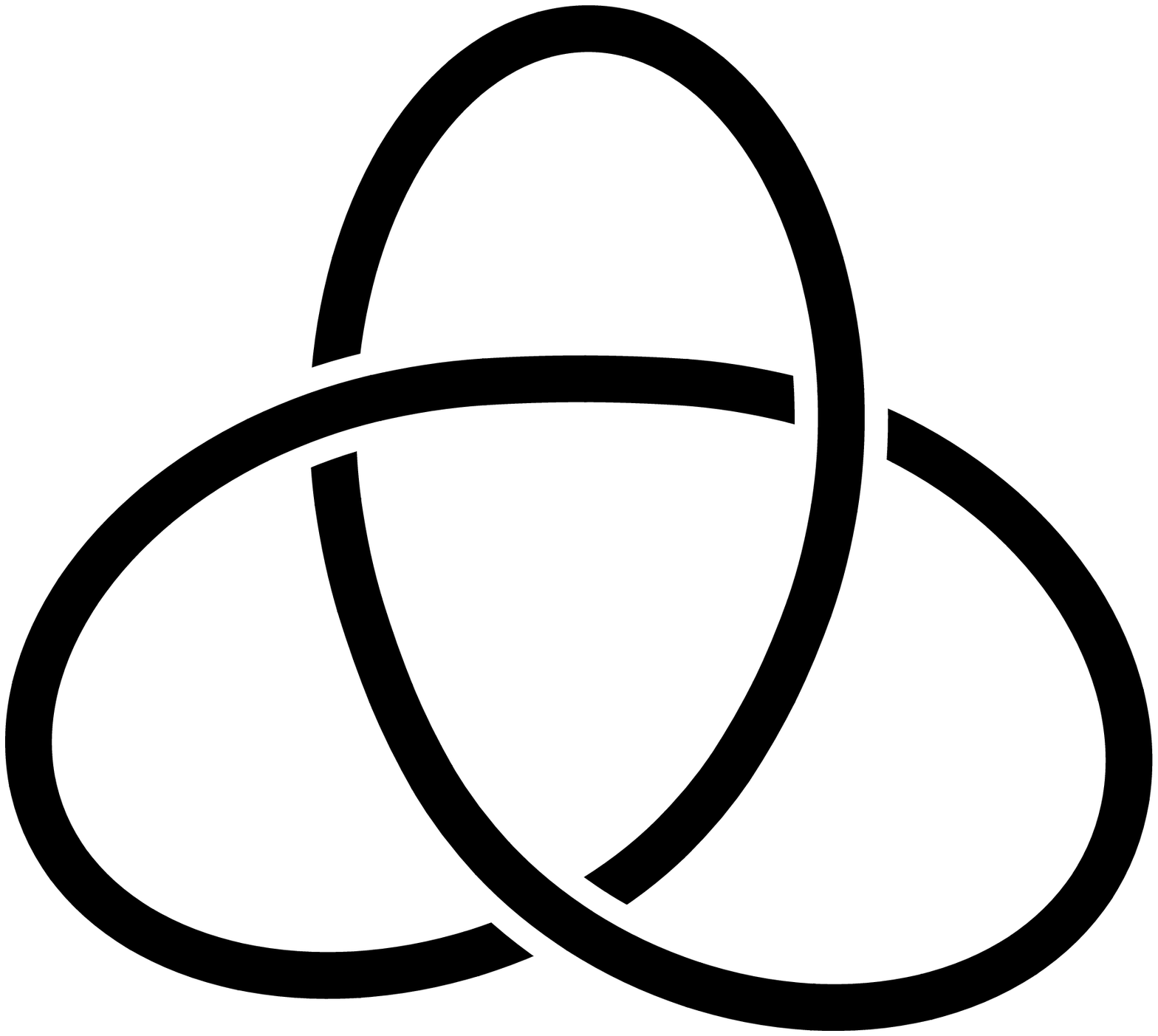} &
      \includegraphics[width=2.5cm]{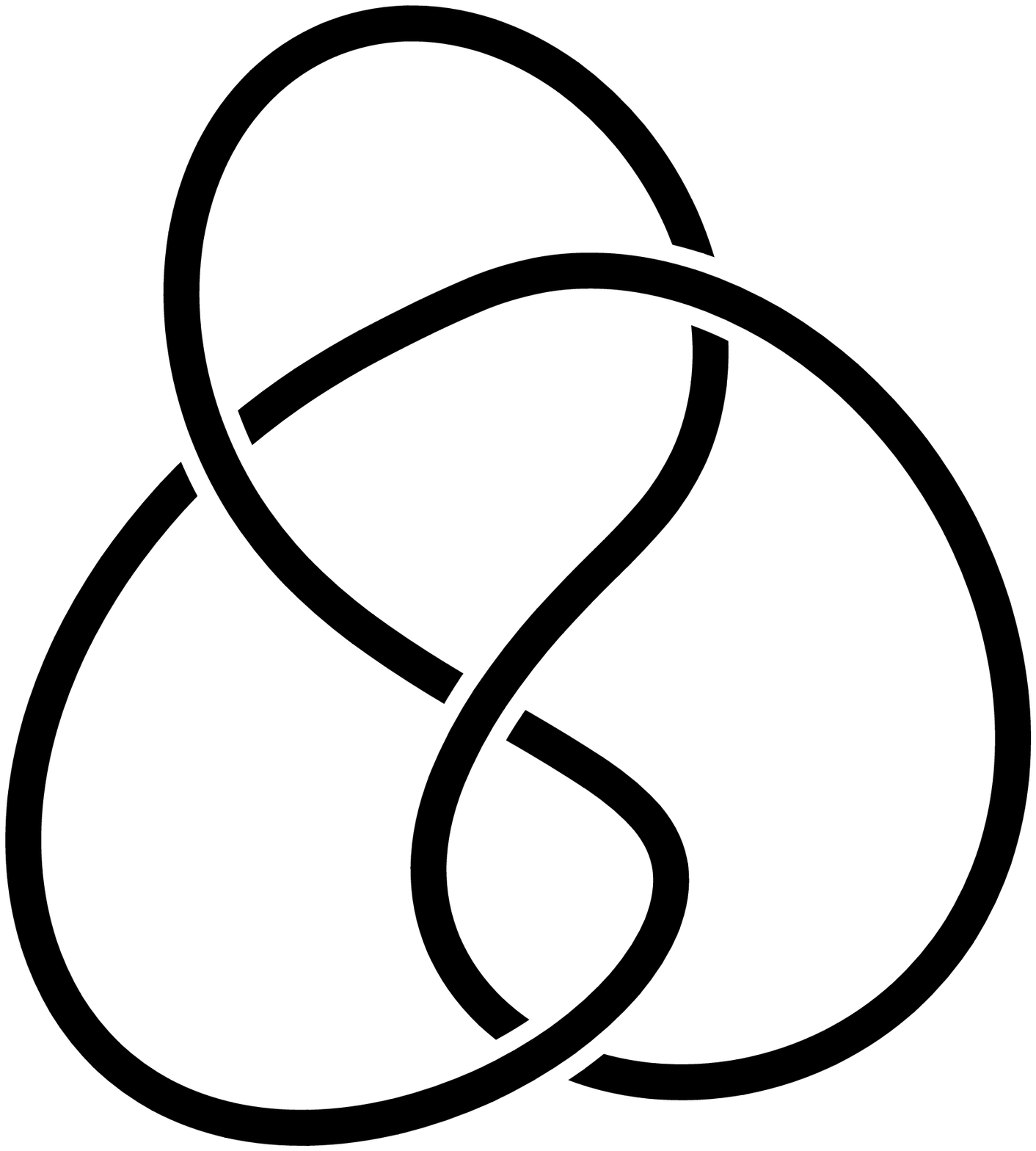} &
      \includegraphics[width=2.5cm]{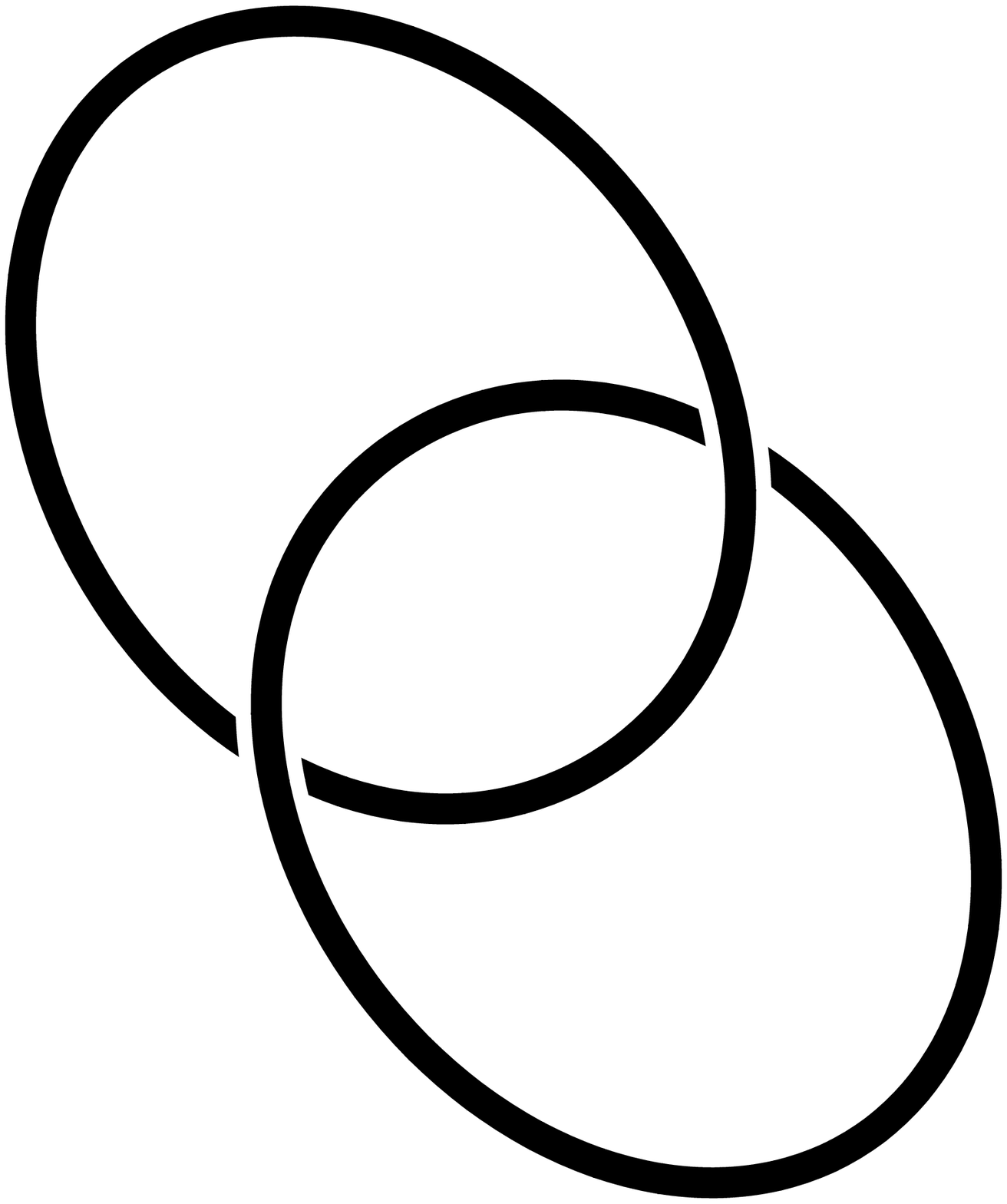}  &
      \includegraphics[width=2.5cm]{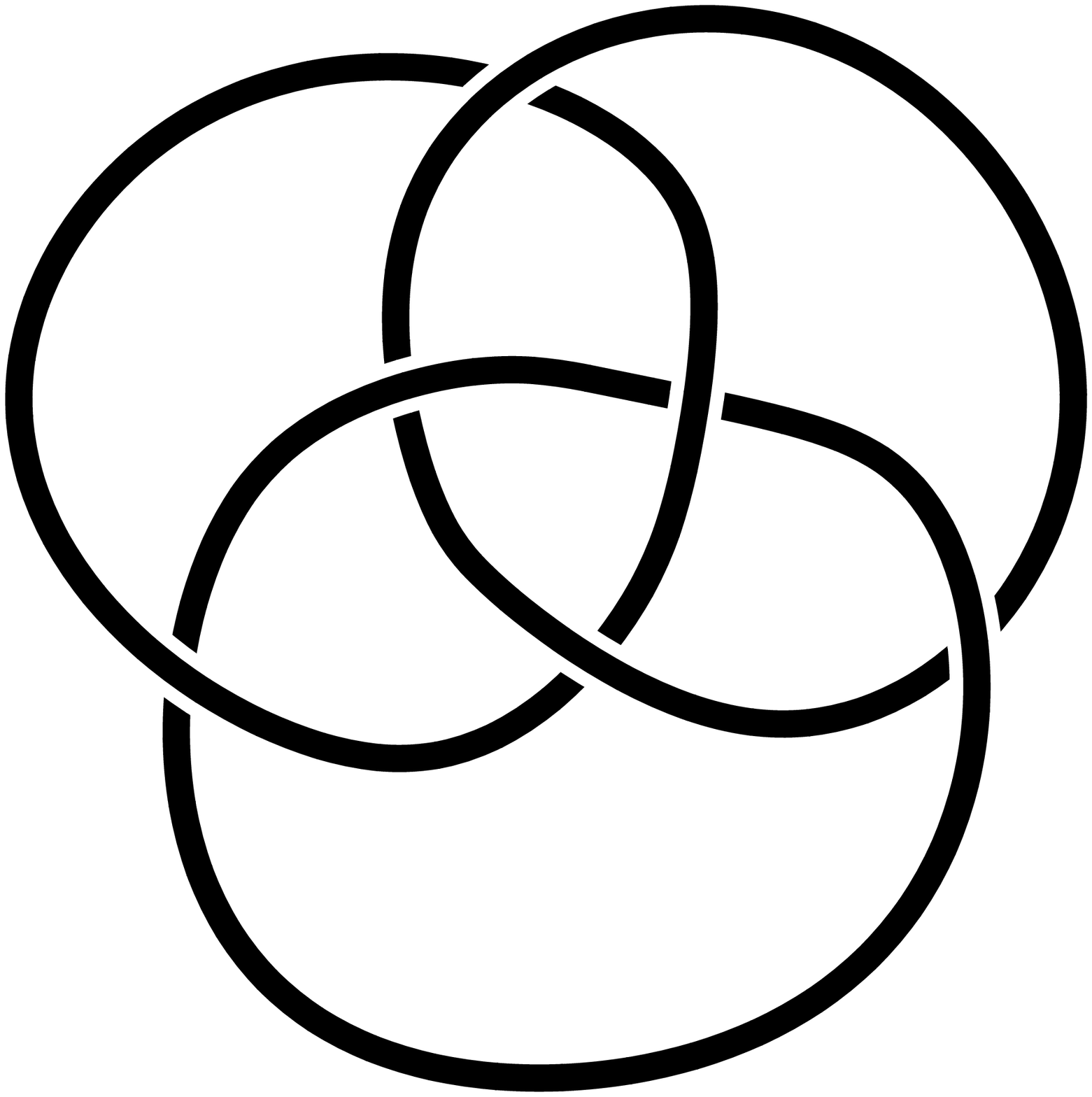} \\
      (A) & (B) & (C) & (D) & (E) \\
\end{tabular}
   \end{center}
   \vspace{-0.7cm}
   \caption{ Knots and links diagrams}
   \label{fig:links3ddiagrams}
\end{figure}

Actually, Figure \ref{fig:links3d} presents the projection on the
``plane of this paper'' of thin cylinders centered and following the
1-dimensional strings that are the 3-dimensional image of the
circles through the embeddings or links. It happens that we could
replace Figure \ref{fig:links3d} by Figure \ref{fig:links3ddiagrams}
without loosing any important information. Each of this drawings is
called a {\it knot diagram} (if only one component) or {\it link
diagram} (any number of components). On each crossing that appears
in the plane projection of the cylinders, there is one cylinder
segment on top of another cylinder segment. This is represented by a
continuous curve (top segment) and a broken curve (bottom segment).

\begin{figure}[htp]
   \begin{center}
   \includegraphics[width=10cm]{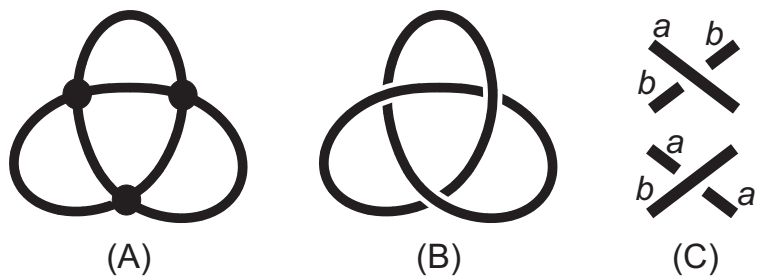}
   \end{center}
   \vspace{-0.7cm}
   \caption{ Link diagram as plane graphs}
   \label{fig:linksDiagramElements}
\end{figure}

So, a link diagram can be seen as a 4-regular {\it plane graph} with
an extra information on each vertex. For example, the trefoil may be
seen as the plane graph of Figure \ref{fig:linksDiagramElements}A.
The extra information of the vertices is shown on Figure
\ref{fig:linksDiagramElements}B and it encodes, in an intuitive way,
exactly the information of which ``cylinder segment'' is on top and
which ``cylinder segment'' is below. Note that there are two
possibilities for this ``extra information''. They are shown in
Figure~\ref{fig:linksDiagramElements}C. The $a$ curve ($b$ curve) in
this figure is said to be the {\it overcrossing} ({\it
undercrossing}) line in the top case and the {\it undercrossing}
({\it overcrossing}) line in the bottom case.

\begin{figure}[htp]
   \begin{center}
   \includegraphics[width=10cm]{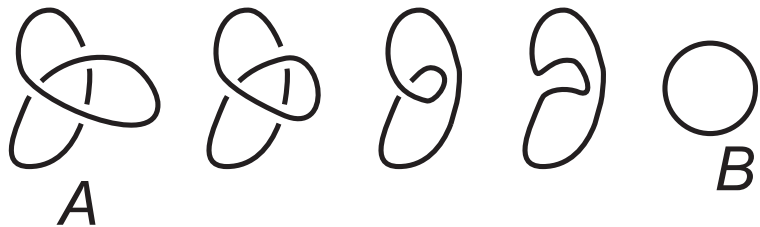}
   \end{center}
   \vspace{-0.7cm}
   \caption{ Ambient isotopic knots}
   \label{fig:exampleOfAmbientIsotopicKnots}
\end{figure}

Given two links, an interesting question to answer is whether these
links can be aligned without tearing any of the strings. For example
the links $A$ and $B$ given by their diagrams on
Figure~\ref{fig:exampleOfAmbientIsotopicKnots} can be aligned as is
shown. Imagine this sequence of ``moves'' transforming $A$ and $B$
occurring on the 3-dimensional space. It is intuitive that we need
no tearing.
On the other hand, the circle and the trefoil (note the crossings on
$A$ to see that it is not a trefoil) cannot be aligned without
tearing. These are examples of the placement problem for links. We
say two links are placed the same way in 3-dimensional space if this
alignment can be done. The formal term for this alignment, defined
in Section~\ref{sec:topology}, is: ambient isotopy. Ambient isotopy
is an equivalence relation and when we say that links are equivalent
we are referring to the ambient isotopy relation. So $A$ and $B$ in
Figure~\ref{fig:exampleOfAmbientIsotopicKnots} are equivalent, but
are not equivalent to the trefoil.
\begin{figure}[htp]
   \begin{center}
   \includegraphics[width=15cm]{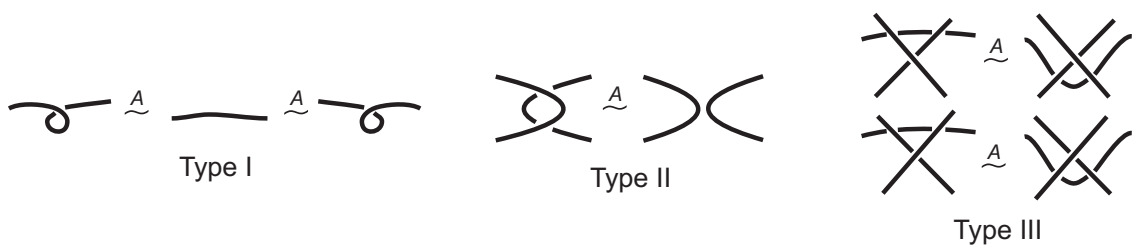}
   \end{center}
   \vspace{-0.7cm}
   \caption{ Reidemeister moves}
   \label{fig:ReidemeisterMoves}
\end{figure}

Reidemeister \cite{Reidemeister1932} proved the following result
about link equivalence in the diagrammatic language:
\begin{center}\it two links are equivalent ({\it i.e.} ambient isotopic) if
and only if  \linebreak any diagram of one link can be transformed
into a diagram for the other link \linebreak  via a sequence of
Reidemeister Moves (Figure~\ref{fig:ReidemeisterMoves}).
\end{center}
We use the symbol $\Aequiv$ between two link diagrams or detached
pieces of link diagrams (where the boundaries of these pieces have a
correspondence that should be clear) to denote that they are ambient
isotopic. Note that the Reidemeister moves we used in the
transformation of Figure~\ref{fig:exampleOfAmbientIsotopicKnots}
were type II move, type I move and alignments. The three
Reidemeister moves will be also called $RE_1$, $RE_2$ and $RE_3$ for
moves of type I, type II and type III, respectively. Two link
diagrams that differ by a finite sequence of Reidemeister moves
$RE_2$ and $RE_3$ are said to be {\it regular isotopic}. The
notation $A \Requiv B$, where $A$ and $B$ are link diagrams, is used
to say that $A$ and $B$ are regular isotopic. Note that regular
isotopic diagrams are always ambient isotopic,
$$ A \Requiv B \implies A \Aequiv B,$$
while the converse is not always true. Observe that the regular
isotopic relation between link diagrams defines an equivalence
relation on the set of link diagrams. This relation is called {\it
regular isotopy}.

Link diagrams interpreted as  {\it blackboard framed links}, as we
will see later, is a presentation for spaces ({\it i.e.} connected,
compact, oriented 3-manifolds). This connection is
essential to the contribution of this work: a prime space catalog of
small BFLs or blinks.

\newpage

\section{Linking number, writhe and linking matrix}

\begin{figure}[htp]
   \begin{center}
   \includegraphics[width=9cm]{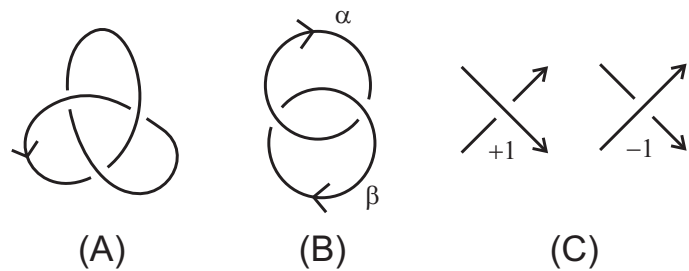}
   \end{center}
   \vspace{-0.7cm}
   \caption{ Oriented links}
   \label{fig:orientedLinks}
\end{figure}

A link is said {\it oriented} if all its components have an {\it
orientation}. There are two possible orientations for each
component. So, a link with $k$ components can be oriented in $2^k$
different ways. To present an oriented link we can draw the link
diagram with one arrow on each component indicating its orientation.
For example, Figure~\ref{fig:orientedLinks}A shows an oriented
trefoil. A crossing on an oriented link diagram may be of two types
as Figure~\ref{fig:orientedLinks}C shows. Each of these types has a
number +1 or -1 which is called the {\it sign of the (oriented)
crossing}. When the undercrossing line, on its orientation sees the
overcrossing line go from left to right then the sign is +1, else,
if it sees the overcrossing line going right to left, the sign is
-1.

Let $D$ be an oriented link diagram of link $L$. Let $\alpha$ be a
component of $L$. The sum of the signs of auto-crossings of $\alpha$
(crossings of $\alpha$ with $\alpha$) on $D$ is said to be its {\it
writhe} and is denoted by $w(\alpha)$. For instance, the writhe of
the only component of the oriented trefoil of
Figure~\ref{fig:orientedLinks}A is +3, because all 3 crossing are
auto-crossings and with sign +1. Note that changing the orientation
of a component does not change its writhe. Now, let $\alpha$ and
$\beta$ be two components of a link $L$. The sum of the signs of the
crossings on $D$ of components $\alpha$ and $\beta$ is said to be
its {\it linking number} and is denoted by $\ell k(\alpha,\beta)$.
For instance, on Figure~\ref{fig:orientedLinks}B, the linking number
of $\alpha$ and $\beta$ is -2 as the two crossings are -1.

Let $D$ be an oriented link diagram with components $\alpha_1$,$\alpha_2$,$\ldots$,
$\alpha_n$. The {\it linking matrix} of $D$ is given by
$$\left(\begin{array}{cccc}
    w(\alpha_1) & \ell k(\alpha_1,\alpha2) & \cdots & \ell k(\alpha_1,\alpha_n) \\
    \ell k(\alpha_2,\alpha_1) & w(\alpha_2) & \cdots & \ell k(\alpha_2,\alpha_n) \\
    \vdots & \vdots & \ddots & \vdots \\
    \ell k(\alpha_n,\alpha_1) & \ell k(\alpha_n,\alpha_2) & \cdots & w(\alpha_n)
  \end{array}\right)
$$ From this matrix, as we will see in Section~\ref{sec:homologyGroup}, it is
possible to obtain a space invariant: the homology group. But to
understand this we must first understand what a link diagram has to
do with spaces. This is the theme of next section.

\section{Framed links and blackboard framed links: encoding spaces}
\label{sec:bfl}

This section presents a fundamental result for this work. To get
deep into this result's justification ideas would demand a lot of
work not needed for our aim. But to get a good image of this
result's meaning is easier. Let's get it.

Consider the unknot on $\IS^3$, {\it i.e.} a ring floating inside the
3-dimensional sphere. Now imagine a small tubular volume $T$,
centered on this unknot. In this situation one could ask: is there a
way to replace the interior of this tubular volume $T$ with
something different? Of course there is. We could replace $T$ by
``nothing'', leading to the ``shape'' of $\IS^3$ with a toroidal hole in it.
Although this is a correct thought, it is not what we are imagining
here. We would like to replace $T$ with something different, but not
leaving a hole. In this case, is there something different of
``replacing $T$ by $T$''? The answer is also yes \footnote{For
theory and examples of these replacements see \cite{Rolfsen1976}.}.
We can replace $T$ by another volume that fills in the hole and
leads to a closed 3-manifold different from $\IS^3$. In fact, this
idea generalizes.

Let's call a replacement like the one we mentioned above by {\it
surgery}. Think of a link on $\IS^3$ and a thin tubular volume $T_i$
centered on each of its components. By analogy with the simple
unknot case above, it is easy to note the possibility of obtaining
different closed 3-manifolds by different surgeries ({\it i.e.}
replacement of the thin tubular volumes $T_i$). In fact, as
Lickorish \cite{Lickorish1962} and Wallace \cite{Wallace1960} proved
independently, any closed, connected, oriented 3-manifold may be
obtained by surgeries (the technical name is Denn surgeries) of a
link on $\IS^3$. So, by doing valid replacements of the thin tubular
volumes centered on the components of a link, one can obtain any
closed, connected, oriented 3-manifold. This result is very
important once shows an intrinsic connection between links ({\it
i.e.} embeddings of circles into the 3-dimensional sphere $\IS^1
\rightarrow \IS^3$) with spaces ({\it i.e.} closed, connected,
oriented 3-manifold).

The information that defines how to do the surgery on a component
(replacement of the tubular volume centered on that component) is
called {\it framing}. So, with a link and a framing for each of its
component, a space is defined. A framing as is justified in Chapter
9 of \cite{Rolfsen1976} may be just an integer number. This leads to
the definition we use of {\it framed link}: a link in $\IS^3$ with an
integer associated to each component. So, from a framed link it is
possible to obtain a space. Start with the link, define the thin
tubular volumes $T_i$ centered on each component and, finally, apply
the surgery on each $T_i$ defined by the framing of component $i$.

\begin{figure}[htp]
   \begin{center}
   \includegraphics[width=12cm]{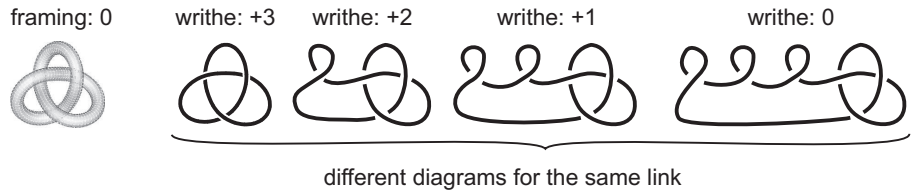}
   \end{center}
   \vspace{-0.7cm}
   \caption{ Aligning framing with writhe}
   \label{fig:blackboardFramedLink}
\end{figure}

In Section~\ref{sec:knotsAndLinks} we saw that a link diagram is a
way to present a link. So, we can present framed links (spaces) by
drawing a link diagram and writing an integer near each component:
its framing. Although this is a nice way to present spaces, there is
an even more concise way to do it based on two facts: (1) the writhe
of a component, as is the framing of a component on a framed link,
is an integer number associated to it on a link diagram; (2) by
adding or removing one curl on a link diagram, application of one
Reidemeister move I, one is able to, without modifying the link,
increment by 1 or -1 the writhe of a component. Using these ideas we
have the following result:
\begin{center}
\it given any framed link, it is possible to draw a link
diagram\linebreak
 for it where the writhe of each component on the link diagram \linebreak
is exactly the framing of the same component.
\end{center}
For example, suppose we want a link diagram for the trefoil with
framing zero. See Figure~\ref{fig:blackboardFramedLink}. We first
draw a link diagram for the trefoil. While the framing of the
component is not equal to its writhe on the diagram, we add curls
until they match. Note that adding these curls do not change the
underlying link.

A {\it blackboard framed link} or {\it BFL} is a link diagram
presentation of a space. The space is the space induced by a framed
link defined by: (1) the link of the framed link is the link on the
diagram; (2) the framing of each component equals to the writhe of
the same component on the diagram. So, any link diagram may be seen
as a blackboard framed link inducing a space and, also, any space
has blackboard framed link presentation for it. One of the main
aims of this work is to identify all prime spaces that have a
``small blackboard framed link'' inducing it.

\bigskip
\bigskip
\centerline{\large \textsc{How to fill the toroidal holes?}}
\bigskip

\begin{figure}[htp]
   \begin{center}
   \includegraphics[width=7cm]{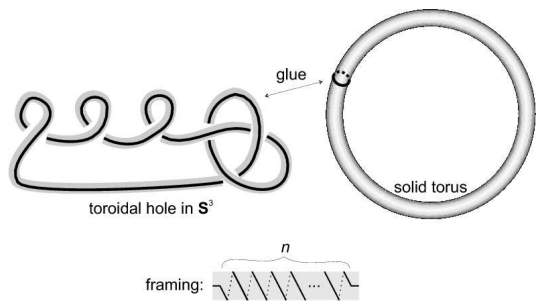}
   \end{center}
   \vspace{-0.7cm}
   \caption{ Gluing a solid torus
   to a toroidal hole: BFL-component and meridian become the same}
   \label{fig:unknotWithMeridian}
\end{figure}

As we said, the framing tells how to close a toroidal
hole in $\IS^3$. But how is that done? Here is how.
A hole is a solid torus embedded in $\IS^3$.
If the framing in zero, then define $c$ as a curve on the surface of the hole
(\ie a solid torus)
parallel to the curve that follows the center of
the hole (see the black line in the left drawing of
Figure~\ref{fig:unknotWithMeridian}). If the framing
is $n \neq 0$ then define $c$ the same way except
that it does $n$ twists on the surface of the hole
before completing a loop (see the bottom drawing
in Figure~\ref{fig:unknotWithMeridian}). To close the
toroidal hole is a matter of doing an abstract gluing:
identify curve $c$ with the meridian curve of a torus
as is shown by the ``glue'' arrow in
Figure~\ref{fig:unknotWithMeridian}. With this
identification a complete homeomorphism is defined
between ``the hole'' and ``the shape'' that replaces
it in a different way.


\newpage

\section{A calculus on blackboard framed links}
\label{sec:BFLcalculus}

\begin{figure}[htp]
   \begin{center}
      \leavevmode
      \includegraphics[width=15cm]{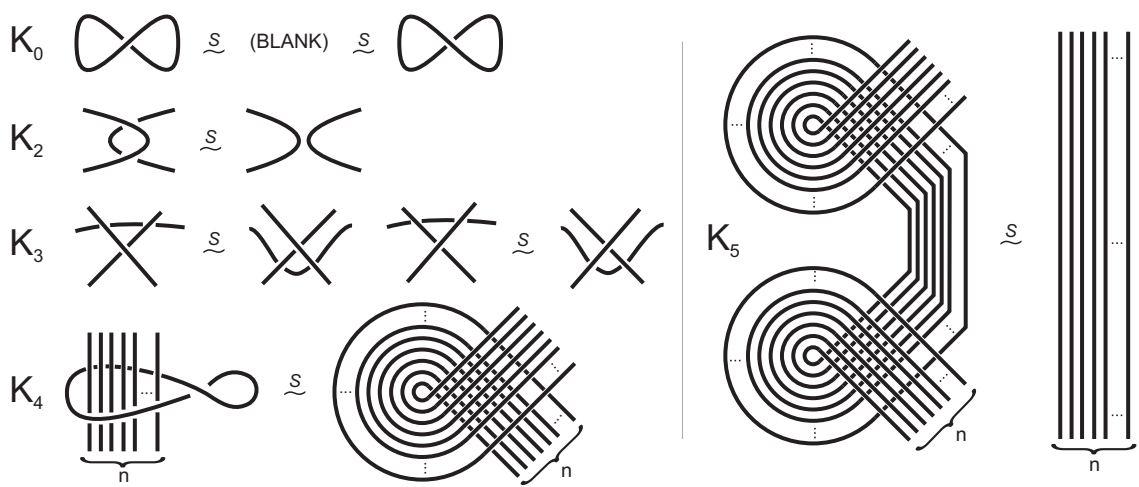}
   \end{center}
   \vspace{-0.7cm}
   \caption{ Kauffman's blackboard framed link formulation of Kirby's
   calculus}
   \label{fig:kauffmancalculus}
\end{figure}

When do two framed links induce the same space? Kirby, in
\cite{Kirby1978}, showed when. Fenn and Rourke
\cite{FennAndRourke1979} reformulated Kirby's ideas, and, from that
point, Kauffman brought Kirby's result to the diagrammatic language
of blackboard framed links. Figure~\ref{fig:kauffmancalculus} shows
Kauffman's blackboard framed link formulation of Kirby's calculus
(page 260 of \cite{Kauffman1991}).

Some notes about Figure~\ref{fig:kauffmancalculus}. The symbol
$\Sequiv$ between two blackboard framed links denotes that both BFLs
induce the same space. When the symbol $\Sequiv$ is used between two
detached pieces of blackboard framed links (the correspondence on
the boundary of these pieces must exist and should be easily
identifiable as in Figure~\ref{fig:kauffmancalculus}) it means that
exchanging these pieces on any blackboard framed link do not change
the induced space. Move $K_0$ states that we can create or eliminate
disjoint knots in form of $\infty$ as we wish, that the induced
space does not change. Note that there is no $K_1$. We reserved this
label for a move shown later. Moves $K_2$ and $K_3$ are,
respectively, Reidemeister moves $RE_2$ and $RE_3$. So, regular
isotopic blackboard framed links (this relation is also defined for
BFL, once BFLs are link diagrams) induce the same space,
$$A \Requiv B \implies A \Sequiv B,$$
because there is a finite sequence of moves in $\{K_2, K_3\}$
connecting them. Moves $K_4$ and $K_5$ are actually a family of
infinite moves indexed by a parameter $n \in \IN$.

\begin{figure}[htp]
   \begin{center}
      \leavevmode
      \includegraphics[width=13cm]{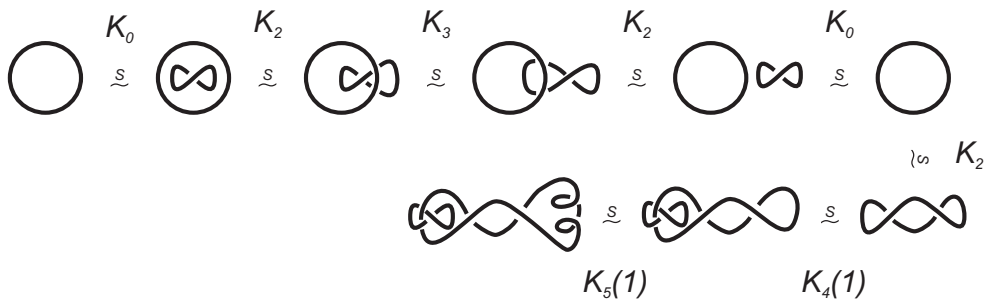}
   \end{center}
   \vspace{-0.7cm}
   \caption{ Example of BFL's inducing the same
   space: $\IS^2 \times \IS^1$}
   \label{fig:exampleBFLCalculus}
\end{figure}

Kauffman's reformulation of Kirby's result states
that
\vspace{-0.2cm}
\begin{center} \it
 if $A$ and $B$ are blackboard framed links, then $A$ and $B$ induce the same
 space \linebreak
 if and only if applying a finite sequence of moves in $\{K_0, K_2, K_3, K_4(n),
 K_5(n)\}$ \linebreak
 one can transform $A$ into $B$.
\end{center}
As an application of this result, see
Figure~\ref{fig:exampleBFLCalculus}. All drawings are
blackboard framed link versions of the same space as they are all
connected by a finite sequence of moves in  $\{K_0, K_2, K_3,
K_4(n), K_5(n)\}$. One can verify, by applying the surgeries on
$\IS^3$ defined by each of these BFL's, that the resulting space is
$\IS^2 \times \IS^1$ (See \cite{Rolfsen1976}).

\begin{figure}[htp]
   \begin{center}
      \leavevmode
      \includegraphics[width=8.5cm]{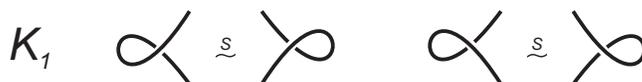}
   \end{center}
   \vspace{-0.7cm}
   \caption{ The {\it ribbon move} or $K_1$}
   \label{fig:ribbonMove}
\end{figure}

We denote by ${\cal K}^0$ Kauffman's set of moves or axioms on
Figure~\ref{fig:kauffmancalculus}: $${\cal K}^0 = \{ K_0, K_2, K_3,
K_4(n), K_5(n) \}.$$ We reserve the remainder of this section to
show that the move defined on Figure~\ref{fig:ribbonMove}, called
the {\it ribbon move}, and denoted by $K_1$, can replace the
infinite class of moves $K_5(n)$ on ${\cal K}^0$ leading to a
simpler and equivalent calculus ${\cal K}^1$.
Let's start by showing that the ribbon~move is a consequence of
${\cal K}^0$. But, before, we need a simple lemma.

\begin{Lem}[Whitney trick] Blackboard framed links that differ by
the pieces below are regular isotopic, so they induce the same
space.
\begin{center}
      \includegraphics[width=8cm]{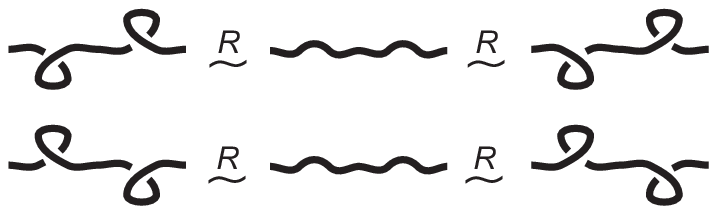}
\end{center}
 \label{lem:whitneytrick}
\end{Lem}

\vspace{-0.6cm}

\begin{proof} The four forms of this lemma are
obtained by combined reflections on the $x$ and $y$ axis of
transformation
\begin{center}
\includegraphics[width=12cm]{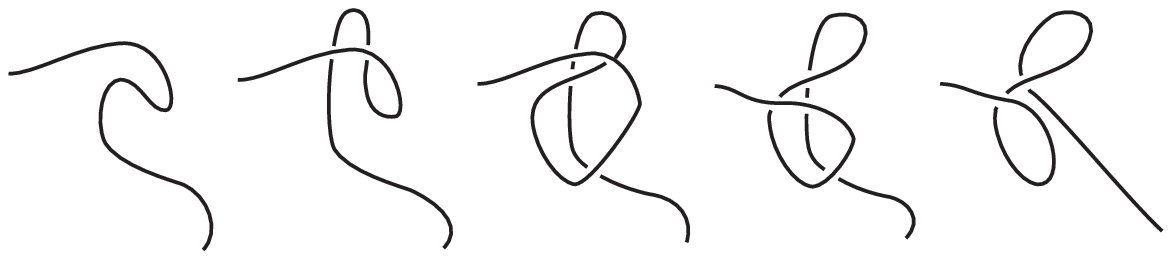}.
\end{center}
Note that each passage is a regular isotopy move in $\{K_2, K_3\}$.
\end{proof}

\begin{Prop} The ribbon move follows from ${\cal K}^0$. More
specifically, from Whitney trick and the $n=1$ version of axiom
$K_5$.
\end{Prop}

\begin{proof} The following picture speaks by itself.
\begin{center}
\includegraphics{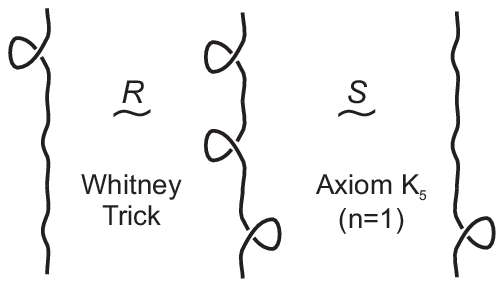}
\end{center}
\vspace{-0.8cm}
\end{proof}

To show that $K_1$ actually can replace $K_5(n)$ it remains to prove
that with the remaining moves and $K_1$ ({\it i.e.} moves in$({\cal
K}^0 \backslash \{K_5(n)\}) \cup \{K_1\})$, we can reproduce
$K_5(n)$, for any $n$. Before doing this, we define some notation
and show some necessary results.

\begin{figure}[htp]
   \begin{center}
      \leavevmode
      \includegraphics[width=13cm]{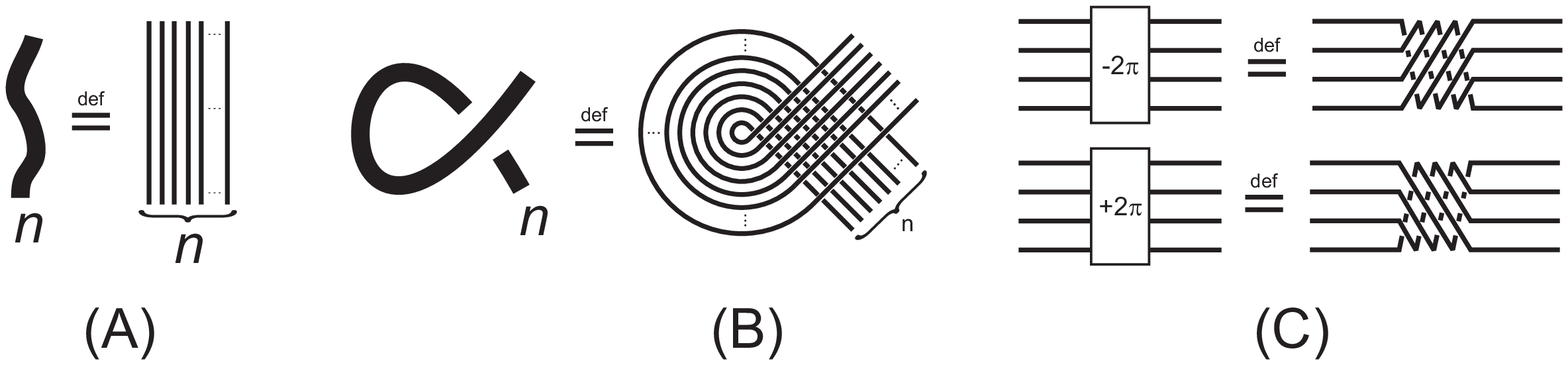}
   \end{center}
   \vspace{-0.7cm}
   \caption{ Some notation}
   \label{fig:someNotation}
\end{figure}

Figure~\ref{fig:someNotation}A shows a thick cable with an $n$ near
it. This notation is a shortcut for $n$ parallel thin lines (the
ones we have been using). Figure~\ref{fig:someNotation}B shows a
thick cable with an $n$ near it doing a curl. This notation is a
shortcut for $n$ parallel lines doing a curl and respecting the
crossings as is shown. When a thicker cable appears without an $n$
and thin cables appear on the same link diagram, the $n$ is implicit
for the thicker cable. Figure~\ref{fig:someNotation}C shows the
definition of a {\em $+2\pi$ twist box} and of a {\em $-2\pi$ twist
box} both with {\em size} equals 4 (number of ``inputs''). The
extension of this definition for any size $\geq 2$ is immediate.

Now we show the last result before proving that $K_1$ may replace
$K_5(n)$. This result uses the $\pm2\pi$ twist boxes notation that
we defined earlier.

\begin{Lem} Regular isotopy leads to
\begin{center}
\includegraphics{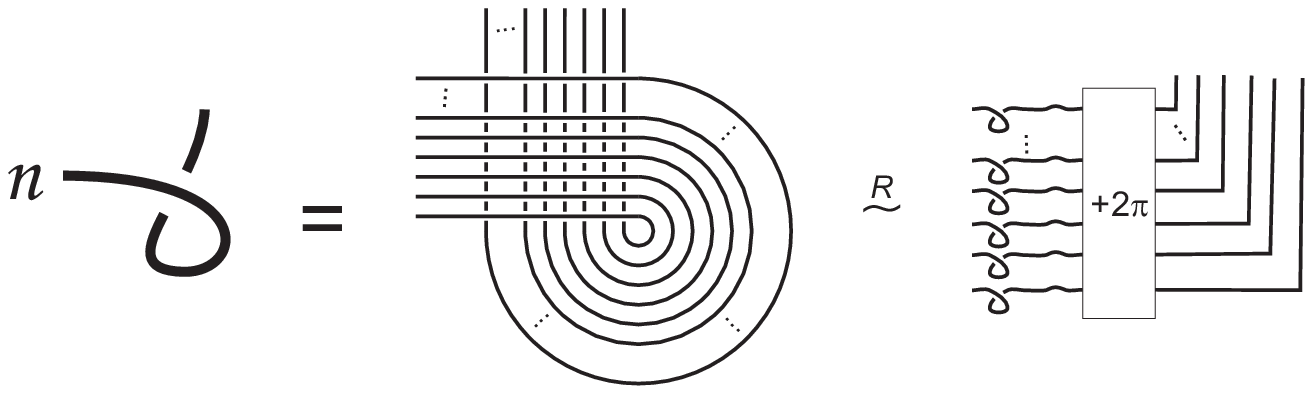}
\end{center}
\end{Lem}
\begin{proof}
Generalization of the following case where $n = 3$.
\begin{center}
\includegraphics[width=15cm]{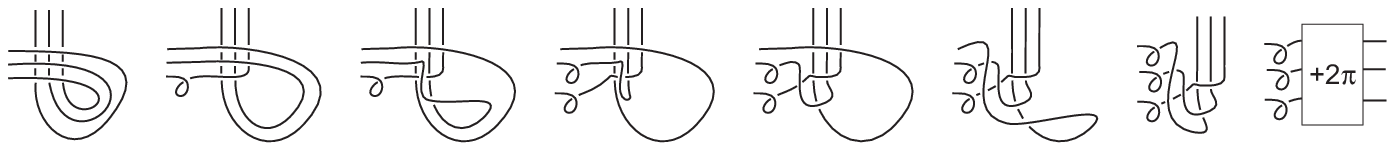}
\end{center}
\vspace{-0.8cm}
\end{proof}

\begin{Lem} Regular isotopy alone is capable of simplifying the left
configuration below with $2i^2$ crossings down to the right one with $2i$
crossings. \label{lem:istranistrands}
\begin{center}
\includegraphics[width=9cm]{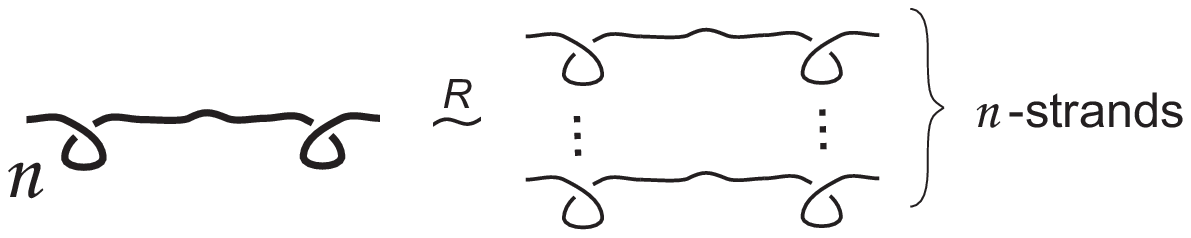}
\end{center}
\end{Lem}
\begin{proof}
This proof is taken from \cite{Kauffman1991}.
\begin{center}
\includegraphics[width=10cm]{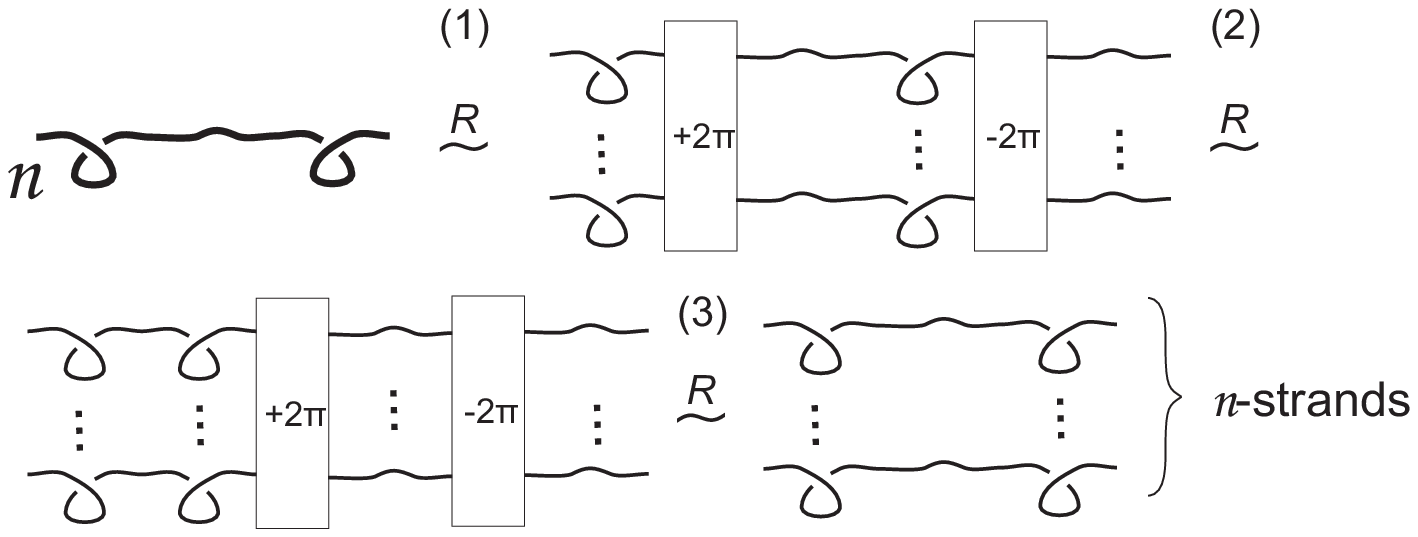}
\end{center}

\vspace{-1cm}
\end{proof}

\begin{Theo} The ribbon move $K_1$ together with regular isotopy
moves $K_2$ and $K_3$ implies move $K_5(n)$.
\end{Theo}

\vspace{-0.3cm}
\begin{proof} Follow this text and the figure below. We begin with
the left side of $K_5$ move. The passage (1) is the application of
Lemma \ref{lem:istranistrands}. The passage number (2) is the
application of ribbon moves on the bottom curl of each strand. The
passage number (3) is the application of Lemma
\ref{lem:whitneytrick} (Whitney trick) on each strand. The rightmost
image is the right side of move $K_5$, so the theorem is proved.
\begin{center}
\includegraphics[width=9cm]{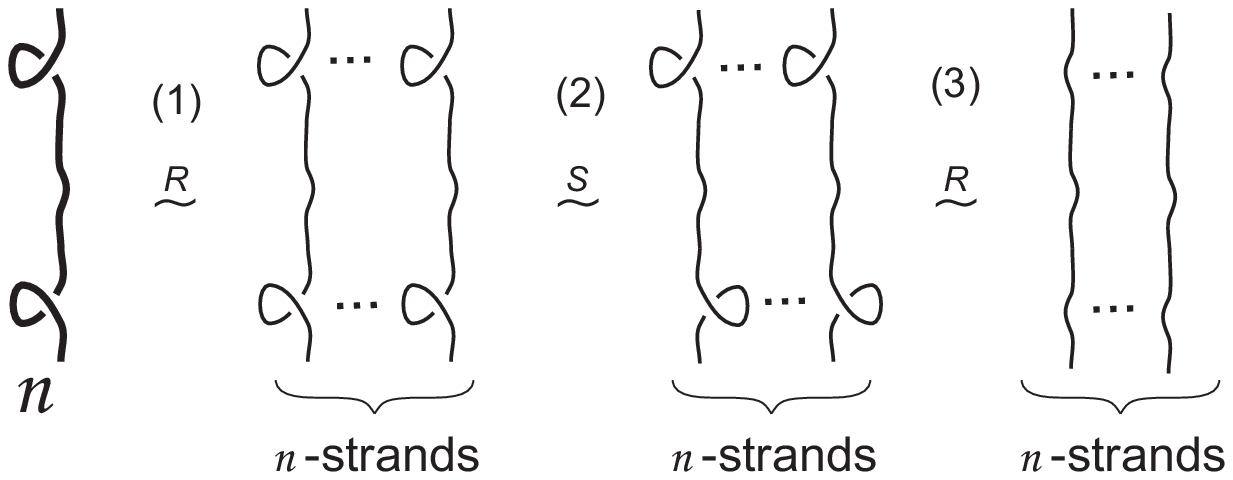}
\end{center}

\vspace{-1cm}
\end{proof}

\begin{figure}[htp]
   \begin{center}
      \leavevmode
      \includegraphics[width=12.5cm]{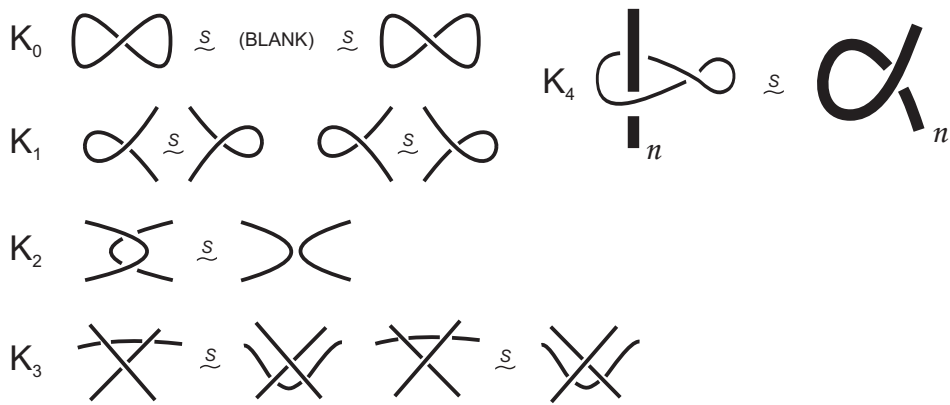}
   \end{center}
   \vspace{-0.7cm}
   \caption{ BFL calculus ${\cal K}^1$, obtained by replacing
   $K_5(n)$ by $K_1$ (ribbon move)}
   \label{fig:BFLcalculusK1}
\end{figure}

Now we present Figure~\ref{fig:BFLcalculusK1} that shows
together all moves of ${\cal K}^1$ calculus: $${\cal K}^1 = \{ K_0, K_1, K_2, K_3, K_4(n)
\}.$$ Two BFLs induce the same space if and only if there is a
finite sequence of moves in $\{ K_0, K_1, K_2, K_3, K_4(n) \}$
transforming one BFL into the other.

\pagebreak

We end this section with some results that are consequence of BFL calculus.

\begin{Lem}[Passing Wall Lemma] These patterns are all regular isotopic

\begin{center}
\includegraphics{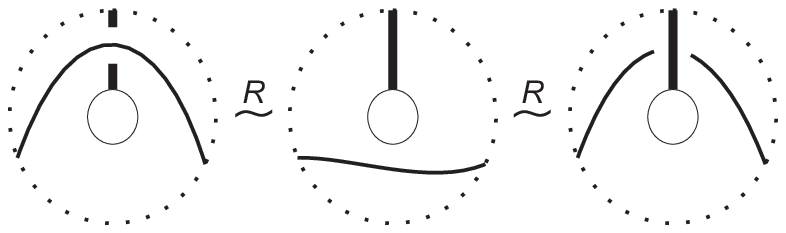}
\end{center}
\end{Lem}
\begin{proof}
This is just the idea: start passing
the horizontal curve under xor (\ie exclusive or) over all
the crossings of the white ball using the moves $K_2$ and $K_3$.
\end{proof}

\begin{Lem}[Passing Cross Lemma] \label{lem:passingCrossLemma} The first
two and last two patterns are all regular isotopic
\begin{center}
\includegraphics{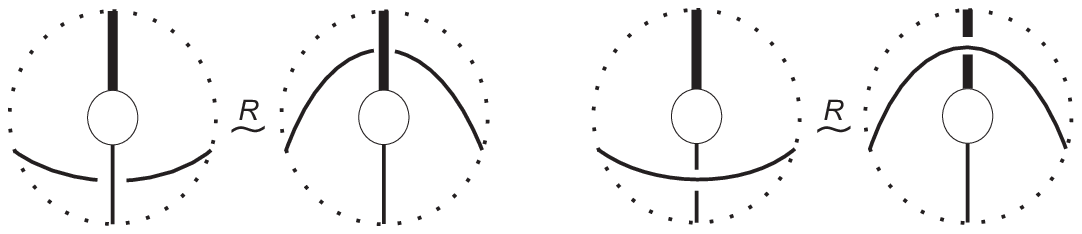}
\end{center}
\end{Lem}
\begin{proof} Follow the picture below. It proves that the first two
patterns of this lemma are regular isotopic. On the passage (1) the
pattern is rearranged to show the structure of the Passing Wall
Lemma. On passage (2) this lemma is applied. On passage (3) we use
the regular isotopy basic move $K_3$ and rearrange its result on
passage (4), arriving at the pattern wanted. The proof for the last
two pattern is analogous to this one.
\begin{center}
\includegraphics{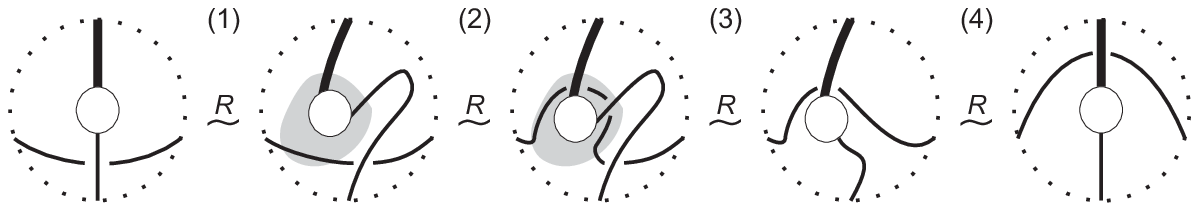}
\end{center}
\end{proof}

\begin{Lem}[Jumping Rope Lemma] \label{lem:jumpingRopeLemma}
The following patterns induces the same space.
\begin{center}
\includegraphics{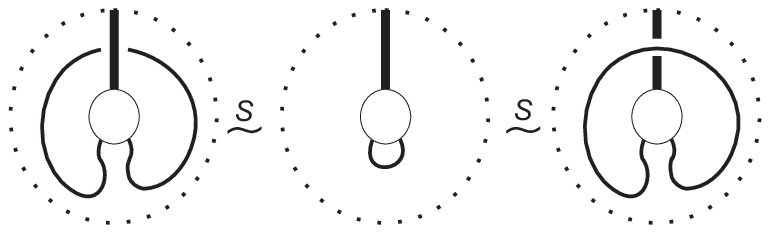}
\end{center}
\end{Lem}
\begin{proof} Follow the picture below. It proves that the first two
patterns of this lemma induce the same space. On the passage (1) the
pattern is rearranged to show the structure of the Passing Wall
Lemma. On passage (2) this lemma is applied. On passage (3) we just
rearrange the pattern to stress the two curls on the different sides
of the cable. On passage (4) the Passing Cross Lemma is repeatedly
applied until the right curl traverse all the cable. On passage (5)
the ribbon move is applied. Finally, the Whitney trick is used on
passage (6). We have thus proved that the first two patterns of this
lemma indeed induce the same space. From these steps it is clear how
the last pattern of this lemma is also proved to induce the same.
\begin{center}
\includegraphics{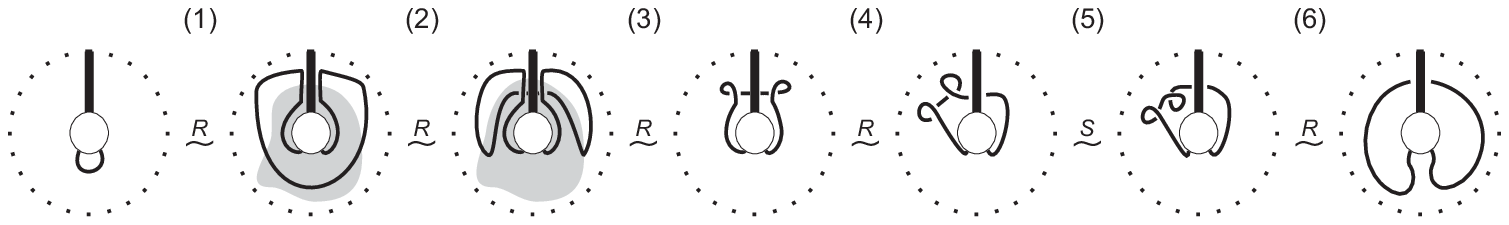}
\end{center}
\vspace{-0.8cm}
\end{proof}


%% file: chapter3.tex
\chapter{Blinks} \label{chap:blinks}

\section{From blackboard framed links to blinks}
\label{sec:BFL2Blinks}

In Section \ref{sec:bfl} we saw that a blackboard framed link or BFL
is a link diagram that induces a space. We now describe a procedure
to build a new object from a blackboard framed link.

\begin{figure}[htp]
   \begin{center}
      \leavevmode
      \includegraphics[width=12.5cm]{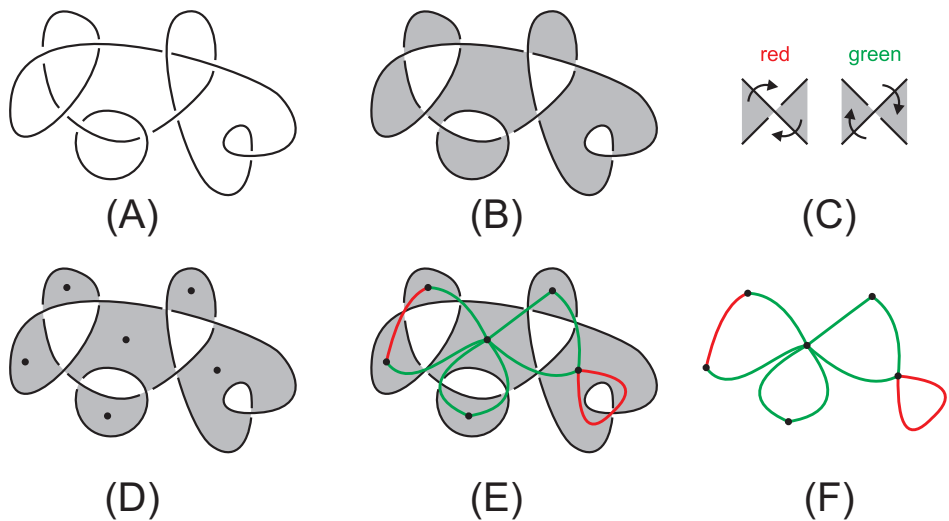}\\
   \end{center}
   \vspace{-0.7cm}
  \caption{Procedure \textsc{BFL2Blink}}
  \label{fig:BFL2Blink}
\end{figure}

Follow the steps described in this paragraph on the example of
Figure~{\ref{fig:BFL2Blink}}. Start with a BFL
(Figure~{\ref{fig:BFL2Blink}A}). We say that two faces on a BFL are
{\it adjacent} if they share a curve (not just a point) that
separates them. The faces of a BFL can be colored black or white
such that no two adjacent faces have the same color. To do this
first define all faces as {\it unassigned}: no color. Then assign
white to the external face of the BFL. Then repeat this: assign
white a face that is adjacent to a black face or assign black a face
that is adjacent to a white face, until all faces are assigned white
or black. This procedure always leads to a unique coloring.
Figure~{\ref{fig:BFL2Blink}B} shows the resulting color assignment
of the BFL of Figure~{\ref{fig:BFL2Blink}A} with the black faces
painted in gray and white faces painted white. The next step is to
classify each crossing of the BFL as {\it red} or {\it green}
(Figure~{\ref{fig:BFL2Blink}C}). A crossing is {\it red} if the
overcrossing line, on the clockwise direction, separates a black
face from a white face. A crossing is {\it green} if the
overcrossing line, on the clockwise direction, separates a white
face from a black face. Now choose one interior point on each black
face as shown in Figure~{\ref{fig:BFL2Blink}D}. For each crossing
$c$, let $A$ and $B$ be the chosen interior points of the two black
faces involved in $c$. Draw a simple curve from $A$ to $B$ such
that: (1) it passes through the crossing point of $c$; (2) all of its
points are black region points or the crossing point of $c$; (3) its
points that are not end-points do not intersect any other crossing
curve. Note that $A$ and $B$ can be the same point. In this case the
curve is a {\it loop}. Figure~{\ref{fig:BFL2Blink}E} shows the
result after drawing all such curves. Figure~{\ref{fig:BFL2Blink}F}
shows the new object after erasing the underlying BFL that guided
its construction. This resulting object is named a {\it blink} and
its general definition is a plane graph with each edge colored either red
or green. Note that a blink may have loops and multiple edges. Each
chosen point on each black face is called a {\it blink vertex} and
each simple curve is called a {\it blink edge}. The {\it size of a blink}
is its number of edges. The blink on Figure~{\ref{fig:BFL2Blink}F}
has size equal to 9.

We name the procedure described in last paragraph as
\textsc{BFL2Blink}. It is always possible to apply it in backwards
and obtain a blackboard framed link from a blink. So the
\textsc{BFL2Blink} when applied in backwards becomes the
\textsc{Blink2BFL} procedure. A blink and a BFL related by the
\textsc{BFL2Blink} or the \textsc{Blink2BFL} are said to be {\it
associated}. So the BFL on Figure~{\ref{fig:BFL2Blink}A} and the
blink on Figure~{\ref{fig:BFL2Blink}F} are associated.

\begin{figure}[htp]
   \begin{center}
      \leavevmode
      \includegraphics[width=9cm]{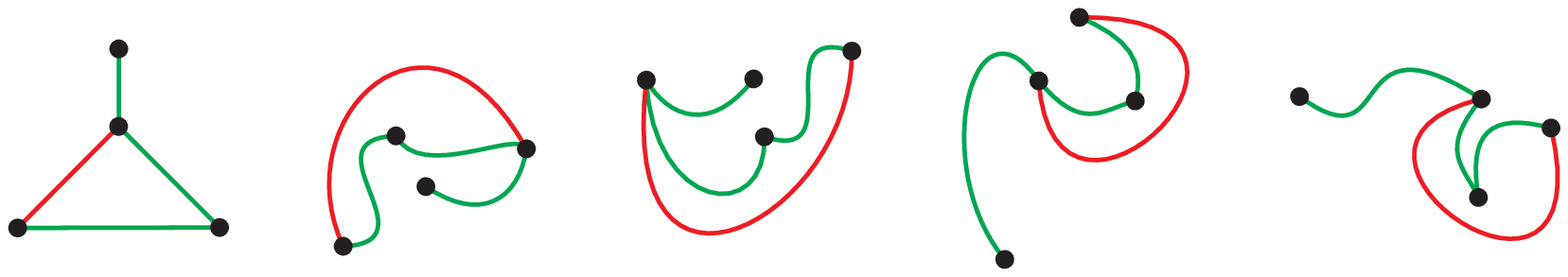}
   \end{center}
   \vspace{-0.7cm}
   \caption{Blinks}
   \label{fig:BlinksPlaneIsotopy}
\end{figure}
In a strict sense, all blinks on Figure~\ref{fig:BlinksPlaneIsotopy}
are different. Their edges are different curves, so, as plane
graphs, they are different. But something connects all these blinks:
there exists a plane isotopy from any of these blinks to any other.
From now on, we say two {\it blinks are equal} if they are
``connected'' by a plane isotopy, otherwise they are {\it
different}. We use the same
convention with blackboard framed links: we say two {\it BFLs are
equal} if they are ``connected'' by a plane isotopy, otherwise they
are {\it different}.

We now claim that all associated BFLs of a class of equivalent
blinks (blinks connected by a plane isotopy) are also connected by
a plane isotopy and vice-versa. So everything fits together and we may
think of a blink as a class of equal blinks and a BFL as a class of
equivalent BFLs. In this sense, a blink (the whole class of
equivalence) is associated to only one BFL (the whole class of
equivalence). By Proposition
\ref{prop:BFLinduceSameSpaceIfPlaneIsotopy} we know that a BFL (the
whole class) induces only one space. This allows us to define the
{\it space of a blink} (the whole class) as the space induced by the
associated BFL.

\begin{Prop}
\label{prop:BFLinduceSameSpaceIfPlaneIsotopy} If there is a plane
isotopy between two blackboard framed links then they induce the
same space.
\end{Prop}
\begin{proof}
Every element involved in the space construction from a BFL is preserved
under plane isotopy.
\end{proof}

\newpage

\section{A calculus for blinks}
\label{sec:blinkCalculus}

\begin{figure}[htp]
   \begin{center}
      \leavevmode
      \includegraphics[width=10cm]{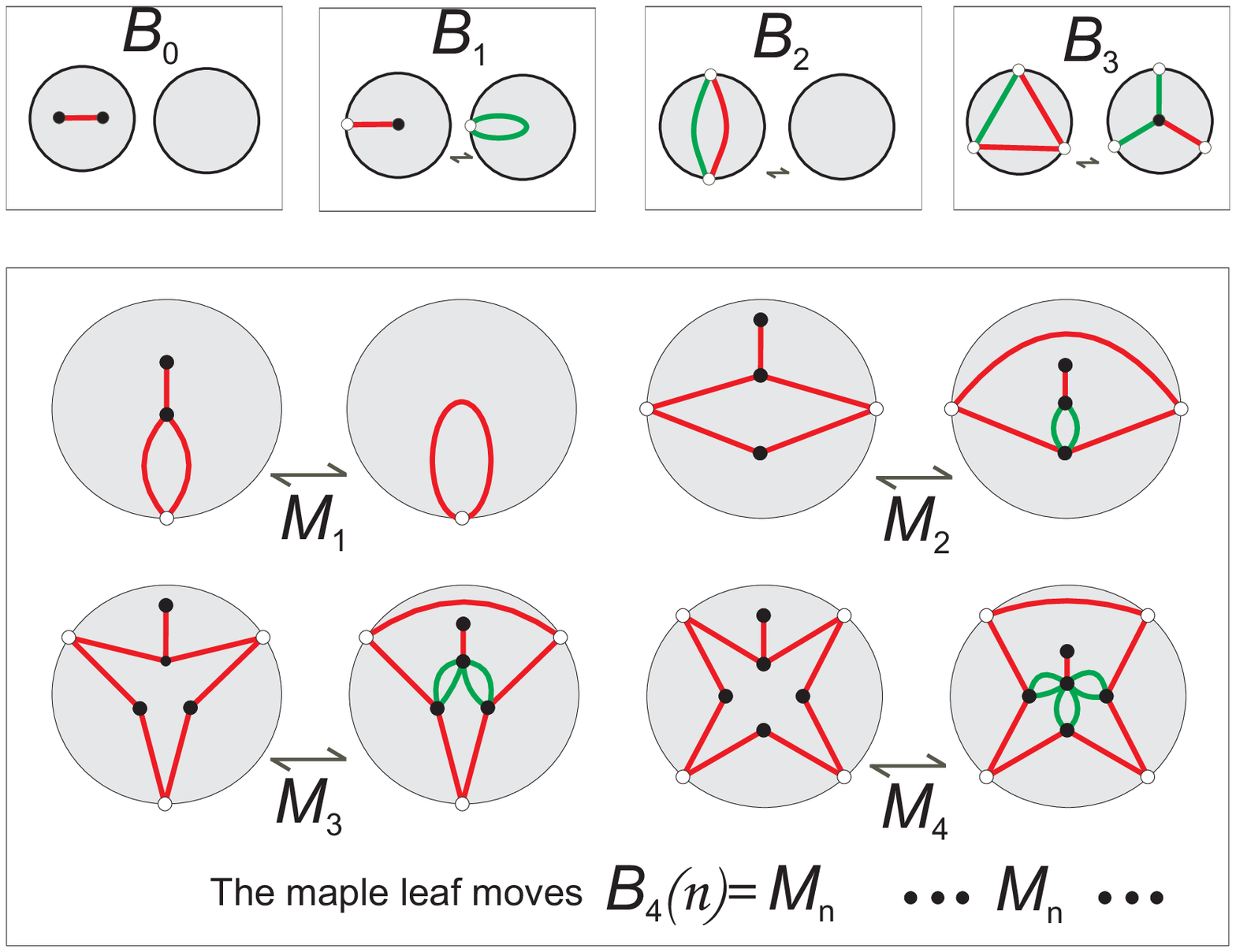}
   \end{center}
   \vspace{-0.7cm}
   \caption{Blink formal calculus by local coins replacements}
   \label{fig:blinkCalculusOnCoins}
\end{figure}

In Section \ref{sec:BFLcalculus} we presented two sets of blackboard
framed link moves: ${\cal K}^0$ and ${\cal K}^1$. These sets have
the strong property of connecting BFLs if and only if they induce
the same space. Here we present a blink version of these sets: a set
of blink moves named ${\cal B}$. Two blinks induce the same space if
and only if there is a finite sequence of moves in ${\mathcal B}$
transforming one blink into the other. The main result of this
section is the following Theorem:
\begin{Theo} \label{theo:blinkCalculus} Two blinks induce the same space if and only if they
are connected by a finite sequence of moves, where each one of them
is one of the ones displayed in Figure
\ref{fig:blinkCalculusOnCoins}, or its red/green twin.
\end{Theo}

Some explanation on Figure \ref{fig:blinkCalculusOnCoins} is in
order. The portion of the blinks which are altered is depicted in
an open disk named a {\em coin}. The interior of the coins modifies
precisely as indicated. The vertices interior to a coin are
displayed as small black circles. The intersection of the blink with
the complement of the coin is a subset of vertices, the {\em
attachment} vertices displayed as small white circles. In this way a
point in the interior of an edge of the blink is either inside or
else outside the coin. We allow arbitrary identifications in the
attachment vertices via deformations of the coins so as to preserve
their interiors (as long as they preserve planarity).

In ${\cal B}$ there are four simple moves twins: $B_0$, $B_1$,
$B_2$, $B_3$, and an infinite family $B_4(1)={M}_1$,\
$B_4(2)={M}_2$,\ $B_4(3)={M}_3$, $\ldots$, named the {\em maple leaf
moves}, $B_4(n) = {M}_n$. By an abuse of notation, each move
$B_i,$ ($i=0,1,2,3$) or $M_n,$ $n \in \IN$, denotes either the move
depicted in Figure \ref{fig:blinkCalculusOnCoins} or its red/green
twin.

The maple leaf move $M_n$ is the manifestation in the blink of the
move $\mu'_n$ on BFLs treated in the subsection which follows the
next one. Move $\mu'_n$ will replace move $\alpha_n$ which
is another name for move $K_4(n)$ shown on Figure~\ref{fig:BFLcalculusK1}.
We stress the point that the set of axioms in the above formal ${\cal
L}$-calculus is a minimal one. For instance, we anticipate the fact
that a move obtained from a move in $\B$ by taking planar duals of
the blinks is a consequence of $\B$.

\bigskip \bigskip \centerline{\bf \large In BFLs, $\mu_n$ is equivalent to $\alpha_n$} \bigskip

\begin{figure}[htp]
   \begin{center}
      \leavevmode
      \includegraphics[width=10cm]{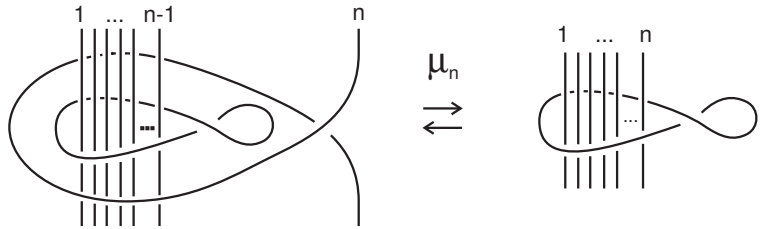}
   \end{center}
   \vspace{-0.7cm}
   \caption{Definition of $\mu_n$, $n \geq 2$.}
   \label{fig:mun}
\end{figure}

We now show that the $\alpha_n$ axiom on BFLs can be replaced by a
new axiom: $\mu_n$. This is useful because the number of crossings
involved in $\mu_n$ is linear on $n$ while in $\alpha_n$ is
quadratic. The axiom $\mu_1$ is defined to coincide with $\alpha_1$.
For $n >1$, $\mu_n$ is defined by Figure \ref{fig:mun}.

\begin{Lem} \label {lem:heartSmoothing}  The heart-shape smoothing
move depicted below is obtained regular isotopies, and a single
ribbon move.
\begin{center}
\includegraphics[width=6cm]{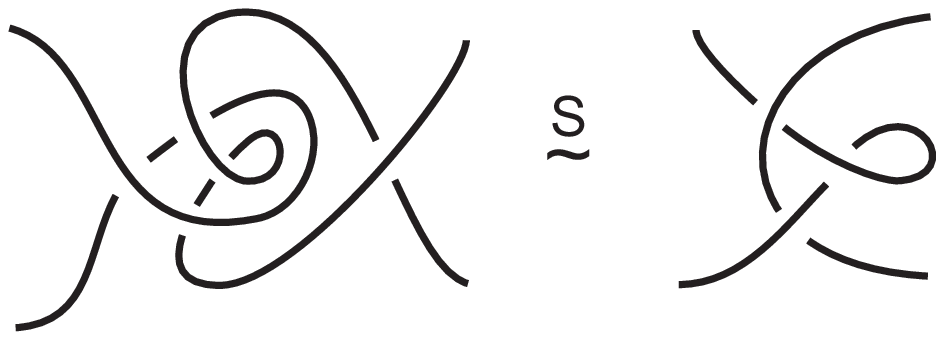}
\end{center}
\end{Lem}
\begin{proof}Follow the proof in the figure below. The ribbon move
is used in the second configuration to prepare for the application
of Whitney´s trick. After this is done we obtain the third
configuration. All other moves are regular isotopies.
\begin{center}
\includegraphics[width=12cm]{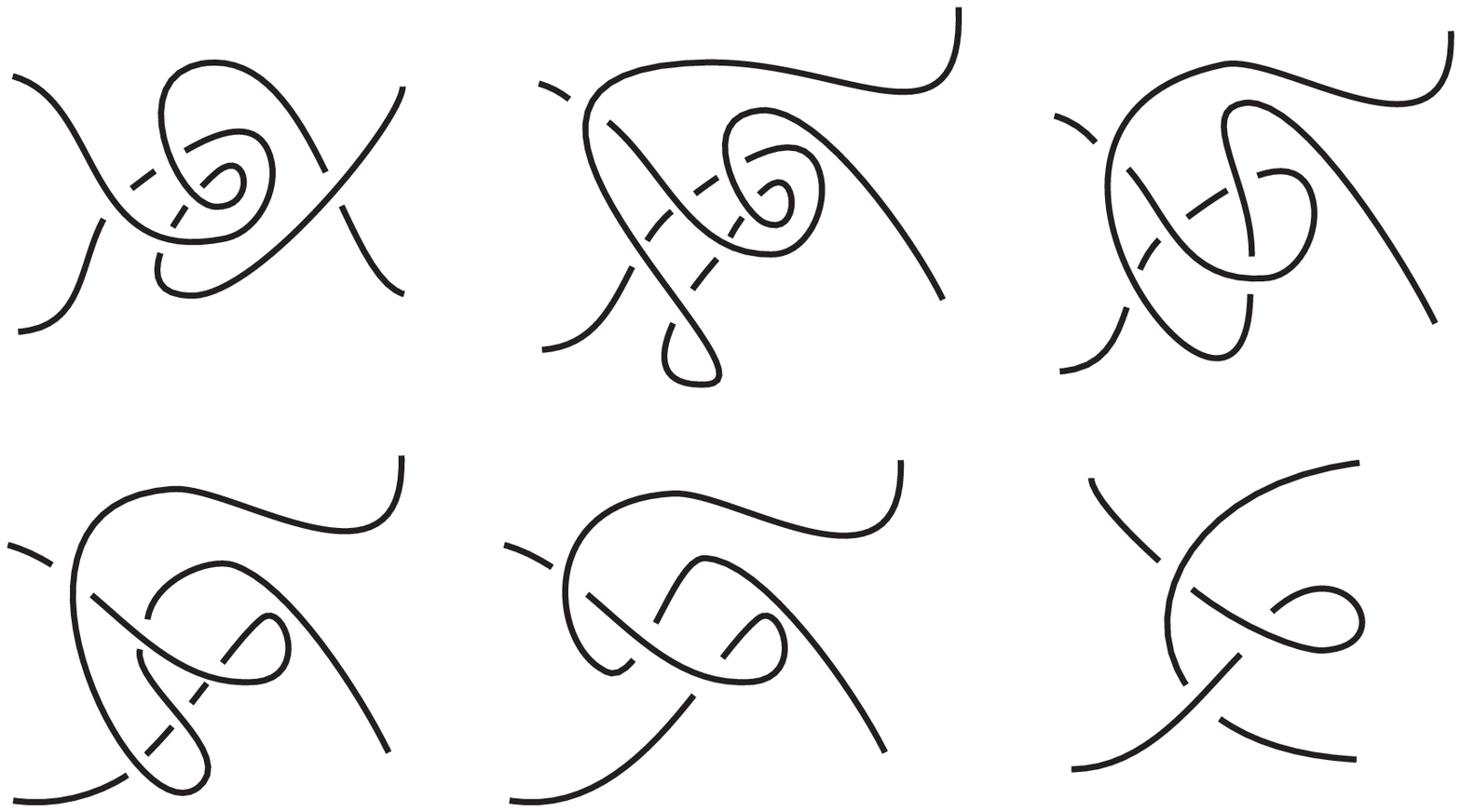}
\end{center}
\vspace{-0.8cm}
\end{proof}

\begin{Lem} The move $\mu_n$ does not change the induced space.
\end{Lem}
\begin{proof} The proof is done for a class of moves that
generalizes $\mu_n$, depicted in Fig. \ref{fig:muNproof}.
\begin{figure}[htp]
   \begin{center}
      \leavevmode
      \includegraphics[width=9cm]{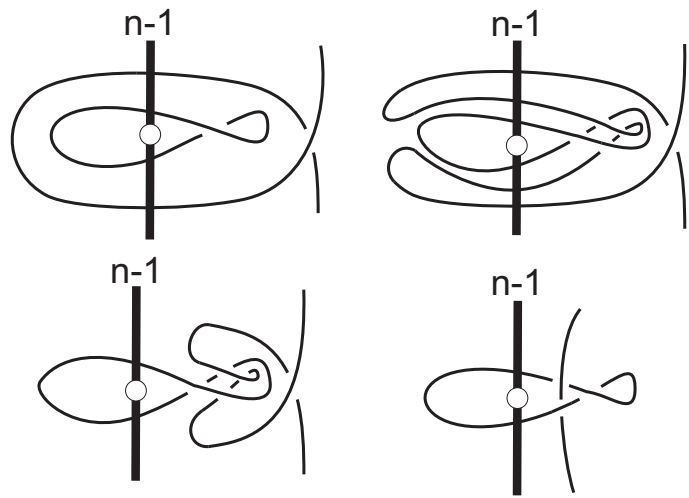}
   \end{center}
   \vspace{-0.7cm}
   \caption{Moves that generalize $\mu_n$.}
   \label{fig:muNproof}
\end{figure}
The white circle separating the cable of $n-1$ parallel strands
means that the $2(n-1)$ individual strands in its boundary are
paired arbitrarily (maintaining planarity, of course). The precise
undercrossings and overcrossings of the individual strands in the
cable are also arbitrary and are left undisplayed: we indicate this
by a real crossing between the thick line and the thinner ones. The
first passage from the first to the second configuration, is a Kirby
handle slide (page 122 of \cite{KauffmanAndLins1994}) obtained by
doubling the $\infty$-shaped component and performing the connected
sum at the external encircling component. Note that (irrespectively
of the individual crossings not shown) the third configuration is
reachable from the second by Reidemeister moves of type II because
consecutive crossings along the individual
strands inside the cable are both over or both under. The third
passage is a consequence of the heart-smoothing move of Lemma
\ref{lem:heartSmoothing}
\end{proof}

\begin{Lem} $\mu_1, \mu_2, \mu_3, \ldots \Rightarrow \alpha_n$ for all $n\geq 1$. In words:
if you have the infinite sequence of moves $\mu_1,\mu_2,\ldots$ then
you can reproduce $\alpha_n$ for any $n \geq 1$.
\end{Lem}
\begin{proof} By induction on $n$. It is obvious that we
have $\alpha_1$ from $\mu_1, \mu_2, \mu_3, \ldots $ once, by
definition, $\alpha_1 = \mu_1$. Suppose we have how to reproduce
$\alpha_i$ from $\mu_1,\mu_2,\ldots$ for all $i < n$. Then, for $n$,
as can be seen on the Figure below, we can apply the induction
hypothesis on the internal $n-1$ strands of the curl and then apply
the $\mu_n$, thus obtaining $\alpha_n$.
\begin{center}
\includegraphics[width=12cm]{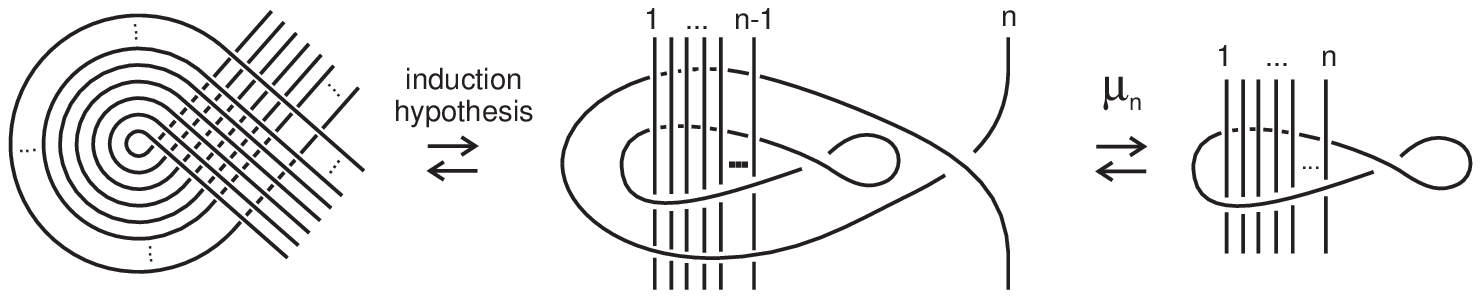}
\end{center}
\end{proof}

\bigskip \bigskip \centerline{\large \bf In BFLs,  $\mu_n$ is equivalent to $\mu'_n$} \bigskip

\begin{figure}[htp]
   \begin{center}
      \leavevmode
      \includegraphics[width=10cm]{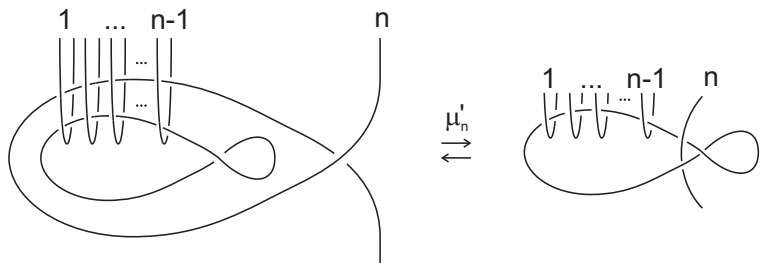}
   \end{center}
   \vspace{-0.7cm}
   \caption{The axiom $\mu'_n$ ($n \geq 2$) : ``better'' blink translation than $\mu_n$}
   \label{fig:MuNLinha}
\end{figure}

By replacing $\alpha_n$ with $\mu_n$ we have simplified our axioms
in the sense that $\mu_n$ has fewer crossings than $\alpha_n$. But,
before translating our axiom system on blackboard framed links to
the blink language, we define the move $\mu'_n$ that is equivalent to
$\mu_n$ but has a ``better'' translation to blinks. The axiom
$\mu'_1$ is equal to $\mu_1$. For $n \geq 2$, $\mu'_n$ is defined by
the schema on Figure \ref{fig:MuNLinha}.

\newpage

\begin{Prop} (Regular isotopy and $\mu'_n$) $\Rightarrow
\mu_n,\,\,\,$ for $n\geq 1$.
\end{Prop}
\begin{proof}
For $n=1$ it is obvious because $\mu_1=\mu_1'$. The figure below
shows the proof for $n=2$. Beginning with the right side of $\mu_2$
we apply regular isotopy ({\it i.e} moves $K_2$ and $K_3$) until we
get to a pattern where $\mu_2'$ can be applied (the second pattern
on the second line). We apply it and then use again regular isotopy
to get to the pattern of the left side of the $\mu_2$ axiom. As all
these transformations are both ways, we have proved the case for
$n=2$. The proof for $n>2$ is analogous to the $n=2$ case and will
not be shown.
\begin{center}
\includegraphics[width=13cm]{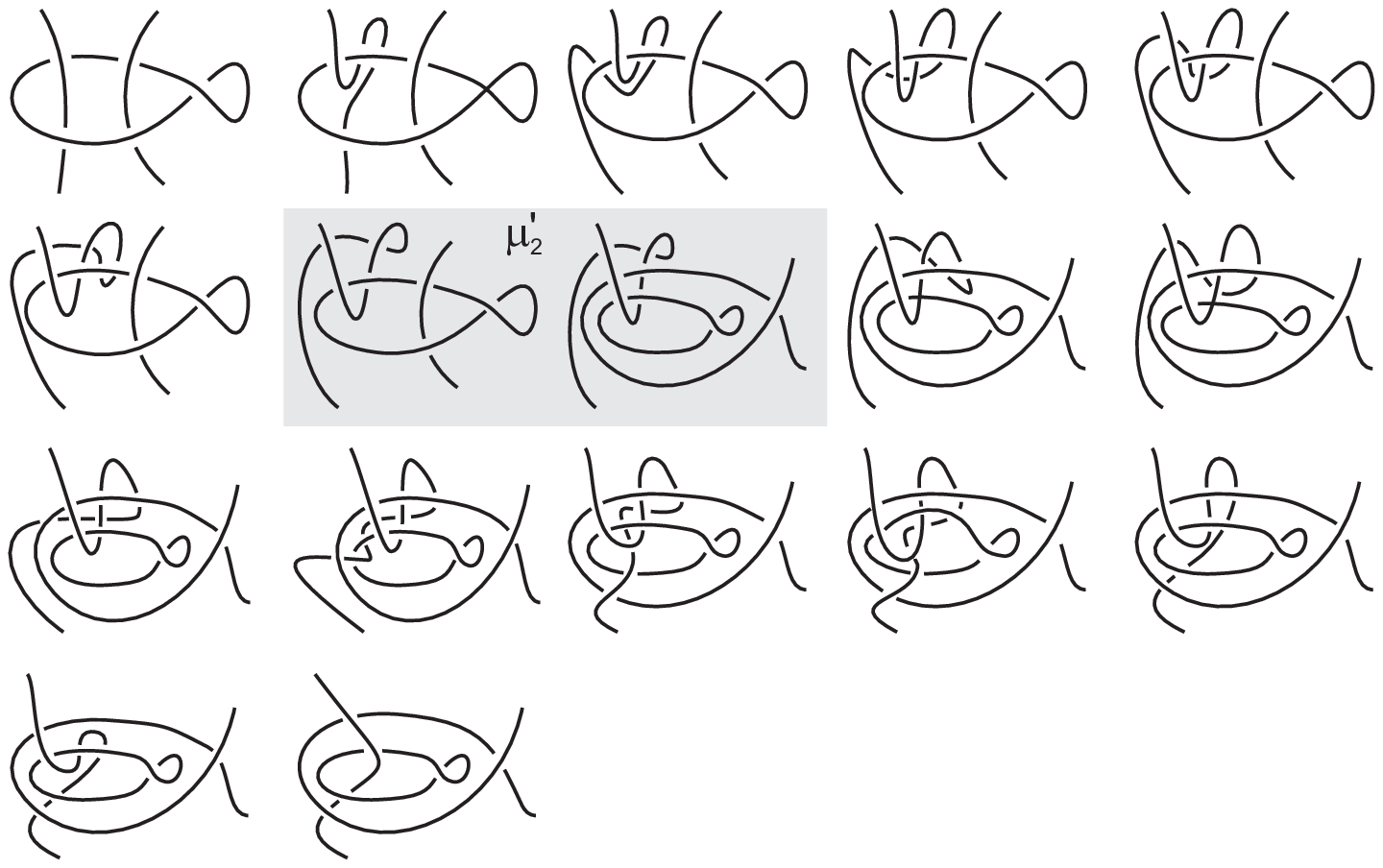}
\end{center}
\vspace{-1cm}
\end{proof}

\newpage

\bigskip \centerline{\bf \textsc{Translation of Blackboard Framed Link Calculus to Blink Calculus}} \bigskip

\begin{figure}[htp]

\begin{center}
\leavevmode
\begin{footnotesize}
\begin{tabular}{|c|c|} \hline
\rule{0cm}{2.7cm}
\includegraphics[width=5cm]{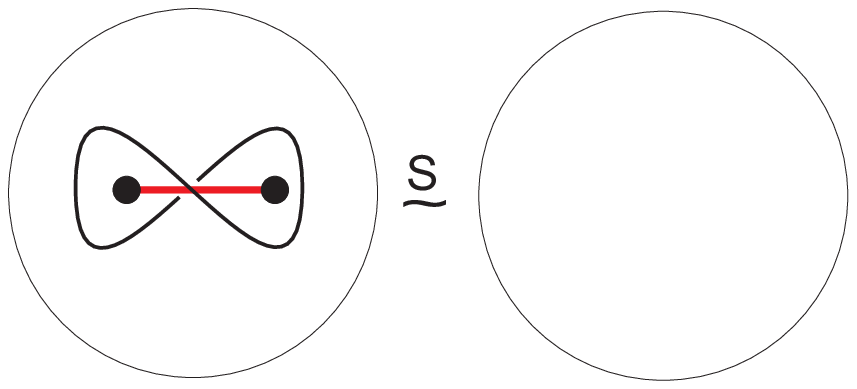} &
\includegraphics[width=5cm]{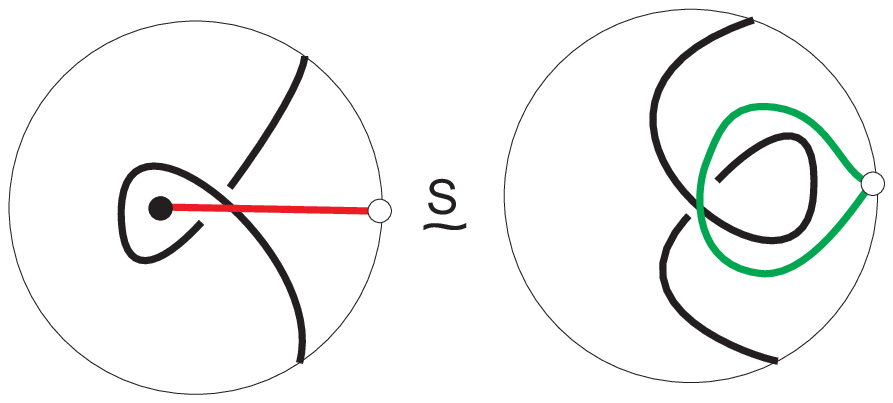} \\
Translation of $K_0$ into $B_0$ &
Translation of the ribbon move $K_1$ into $B_1$ \\ \hline
\rule{0cm}{2.7cm}
\includegraphics[width=5cm]{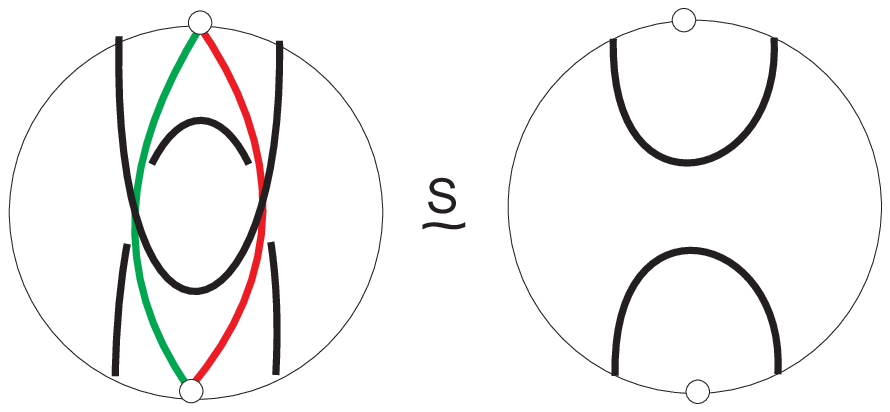} &
\includegraphics[width=5cm]{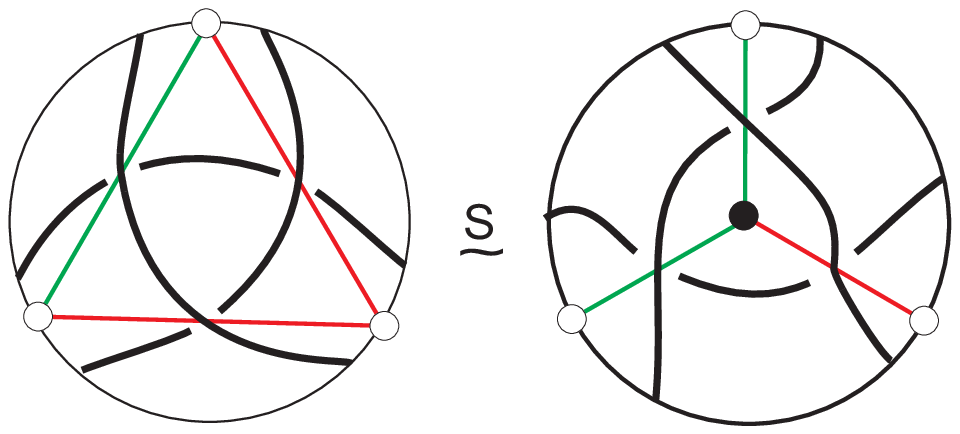} \\
Translation of $K_2$ into $B_2$ &
Translation of $K_3$ into $B_3$ \\ \hline
\multicolumn{2}{|c|}{
      \rule{0cm}{7.4cm} \includegraphics[width=9cm]{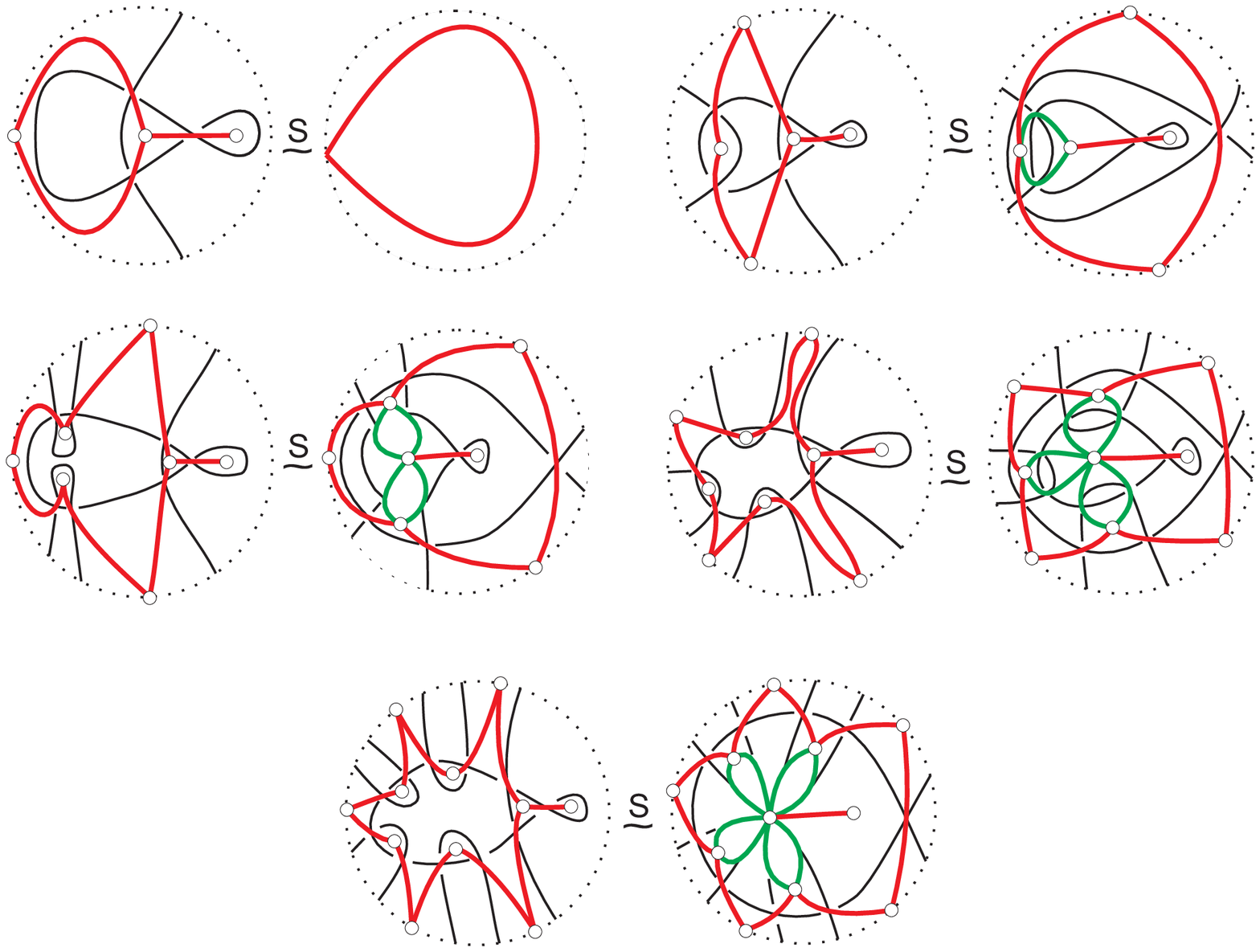}
} \\
\multicolumn{2}{|c|}{
      Translation of $\mu'_1, \ldots, \mu'_5, \ldots$,
   into $M_1, \ldots, M_5, \ldots$
} \\ \hline
\end{tabular}
\end{footnotesize}
\end{center}
   \vspace{-0.5cm}
   \caption{Translation of BFL calculus to blink calculus}
   \label{fig:translationBFLCalculusToBlinkCalculus}
\end{figure}

The translation from $\K$ into $\B$ is depicted in Figure~\ref{fig:translationBFLCalculusToBlinkCalculus}.

\begin{proof}(of Theorem \ref{theo:blinkCalculus}.) The proof is a direct translation
of the moves $K_0, \ldots, K_3$ and $\mu'(1),\ldots \mu'_n,\ldots,$
into the moves $B_0, \ldots B_3$, and $M_1, \ldots K_n, \ldots$
\end{proof}

\newpage

\section{g-blinks}
\label{sec:gblinks}

\begin{figure}[htp]
   \begin{center}
      \leavevmode
      \includegraphics{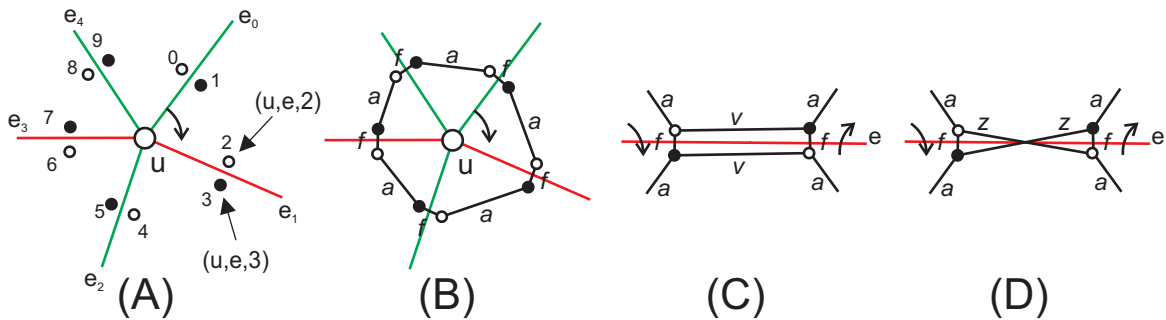}\\
   \end{center}
   \vspace{-0.7cm}
  \caption{Elements on the definition of a g-blink from a blink}
  \label{fig:gblinkFromBlinkElements}
\end{figure}

Let $B$ be a blink. We now describe a procedure to define, from $B$,
a 4-regular graph $G_B$ named a \hbox{\it g-blink}. This procedure is
called \textsc{Blink2GBlink} and associates to any blink
(topological object) a unique g-blink (combinatorial object). Let
$u$ be a vertex of $B$ and $e_0,\ldots,e_{\delta_u-1}$ be the edges
incident to $u$ ordered in clockwise direction ($e_0$ may be any
edge). For each edge $e_i$ with $i \in \{0,\ldots,\delta_u-1\}$ we
define two vertices in $G_B$: one labeled $(u,e_i,2i)$ positioned
close to $e_i$ but before it in clockwise direction; the other is
labeled $(u,e_i,2i+1)$ positioned close to $e_i$ but after it in
clockwise direction (see Figure \ref{fig:gblinkFromBlinkElements}A).
If $(u,e,2j)$ and $(u,e,2j+1)$ are vertices of $G_B$ then they are
the ends of a {\it face-edge} of $G_B$ (Figure
\ref{fig:gblinkFromBlinkElements}B). If $(u,e,2j+1)$ and $(u,f,2j+2
\hbox{ \rm mod } 2\delta_u)$ are vertices in $G_B$ then they are the
ends of a {\it angle-edge} in $G_B$ (Figure
\ref{fig:gblinkFromBlinkElements}B). If $(u,e,j)$ and $(v,e,k)$ are
vertices in $G_B$ and the parity of $j$ is different from the parity
of $k$ then they are the ends of a {\it vertex-edge} of $G_B$
(Figure \ref{fig:gblinkFromBlinkElements}C). If $(u,e,j)$ and
$(v,e,k)$ are vertices in $G_B$ and the parity of $j$ is equal to
the parity of $k$ then they are the ends of a {\em zigzag-edge} of
$G_B$ (Figure \ref{fig:gblinkFromBlinkElements}D).

\begin{figure}[htp]
   \begin{center}
      \leavevmode
      \includegraphics[width=14cm]{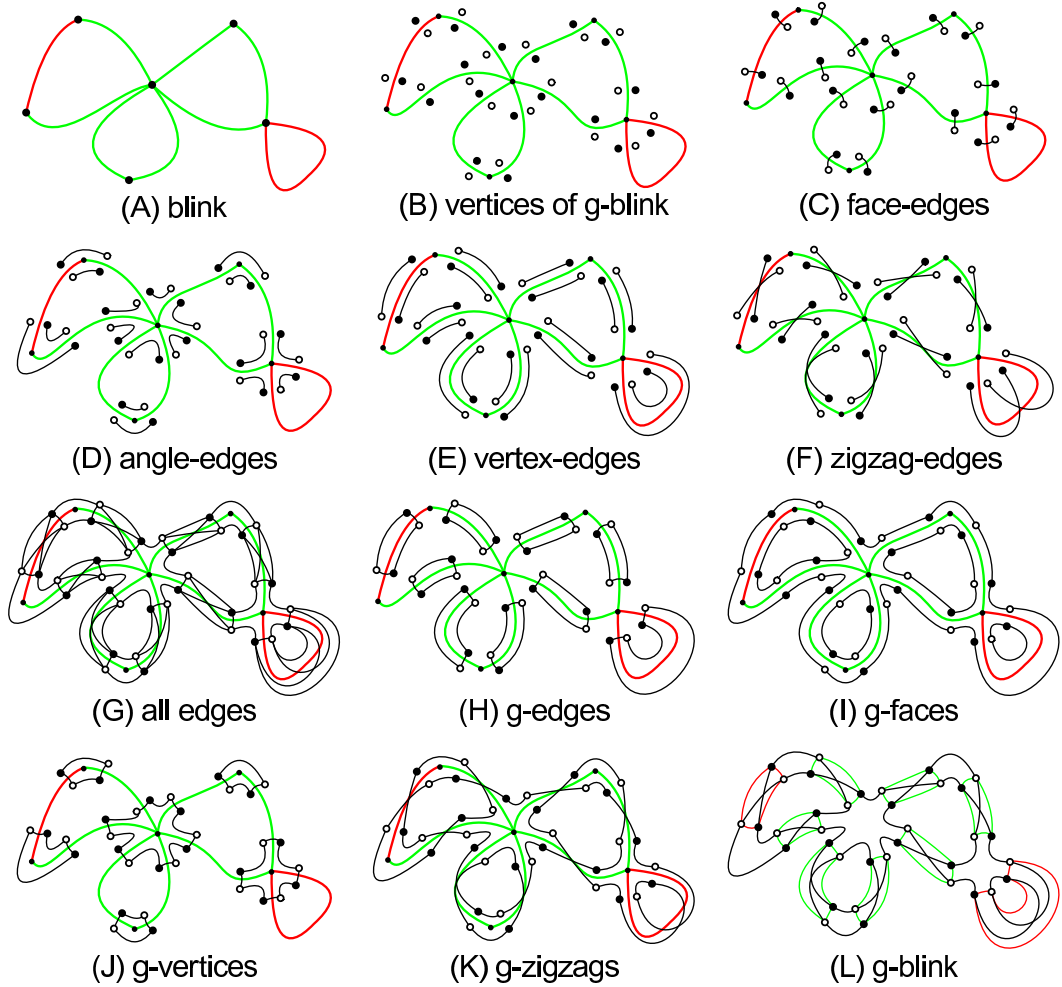}\\
   \end{center}
   \vspace{-0.7cm}
  \caption{Blink, g-blink and attributes: an example}
  \label{fig:Blink2GBlinkExample}
\end{figure}

We define a bipartition $V_0$ and $V_1$ of the vertices of $G_B$
like this: a vertex $v$ labeled with $(\_\, , \_\, ,2j)$, for some integer
$j$, is said to be a {\em parity zero vertex} and it is in $V_0$; a
vertex $v$ labeled with $(\_\, ,\_\, ,2j+1)$, for some integer $j$, is said
to be a {\em parity one vertex} and it is in $V_1$. On the example
of Figure \ref{fig:Blink2GBlinkExample} $V_0$ are the white vertices
and $V_1$ are the black vertices.

If $B$ has $n$ edges, then $G_B$ has $4n$ vertices and $8n$ edges.
Each vertex of $G_B$ has degree 4 and is incident to a face-edge, an
angle-edge, a vertex-edge and a zigzag-edge. If $v$ is a vertex in
$G_B$ we denote by ${\rm adj}_v(v)$, ${\rm adj}_f(v)$, ${\rm
adj}_a(v)$ and ${\rm adj}_z(v)$ the vertices adjacent to $v$ by
vertex-edge, face-edge, angle-edge and zigzag-edge respectively.

A {\it g-edge} is a polygon on a g-blink whose edges alternate
between face-edges and vertex-edges (Figure
\ref{fig:Blink2GBlinkExample}H). A g-edge has always 4 edges and 4
vertices and is associated to an edge on a blink (note that the
vertices of a g-edge are of the form $(\_\, ,e,\_)$). If the
corresponding blink edge of a g-edge is red then this g-edge is also
red. If the corresponding blink edge of a g-edge is green then this
g-edge is also green. A {\it g-face} of a g-blink is any polygon
with vertex-edge alternated with angle-edge (Figure
\ref{fig:Blink2GBlinkExample}I). Each of these polygons corresponds
to a face of the blink. A {\it g-vertex} of a g-blink is any polygon
with face-edge alternated with angle-edge (Figure
\ref{fig:Blink2GBlinkExample}J). Each of these polygons corresponds
to a vertex of the blink. A {\it g-zigzag} of a g-blink is any
polygon with angle-edge alternated with zigzag-edge (Figure
\ref{fig:Blink2GBlinkExample}K). Each of these polygons corresponds
to a component on the blackboard framed link associated with the
blink.

\begin{figure}[htp]
   \begin{center}
      \leavevmode
      \includegraphics[width=6cm]{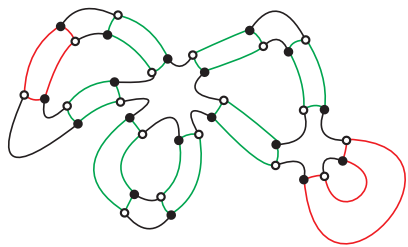}\\
   \end{center}
   \vspace{-0.7cm}
  \caption{g-blink of Figure \ref{fig:Blink2GBlinkExample}L without zigzag-edges: a planar graph}
  \label{fig:gBlinkWithoutZigzags}
\end{figure}

Now, using the notation defined above, we state a definition for
{\it g-blink}. A {\it g-blink} is a graph that satisfies the
following six conditions:

\vspace{-0.1cm}

\begin{enumerate}
\setstretch{1.5}
\setlength{\parskip}{-2pt} \sl

\item[(1)] Its vertices are partitioned in $V_0$ and $V_1$ (white and black
vertices of Figures \ref{fig:gblinkFromBlinkElements} and
\ref{fig:Blink2GBlinkExample});

\item[(2)] vertices in $V_0$ are adjacent by face-edge, vertex-edge and
angle-edge to vertices in $V_1$ and by zigzag-edges to vertices in
$V_0$; vertices in $V_1$ are adjacent by face-edge, vertex-edge and
angle-edge to vertices in $V_0$ and by zigzag-edges to vertices in
$V_1$;

\item[(3)] each vertex is incident to exactly one face-edge, one
vertex-edge, one angle-edge and one zigzag-edge;

\item[(4)] each polygon of alternating face-edge and vertex-edge has
4 edges (is a g-edge) and is assigned color red or green
(see the g-edges on Figure~\ref{fig:Blink2GBlinkExample}L) or,
equivalently, a pair of zigzag edges of the same g-edge is labeled
one edge as {\it overcrossing} and the other edge as {\it
undercrossing};

\item[(5)] the zigzag-edges are the diagonals of the g-edges connecting
the vertices with the same parity. Observe that this implies that
zigzag-edges are redundant when we know the g-edges. They may be
omitted when presenting a g-blink and calculated or shown when
needed. For example Figure~\ref{fig:Blink2GBlinkExample}L can be
easily restored from Figure~\ref{fig:gBlinkWithoutZigzags};

\item[(6)] the 3-regular graph obtained by not considering the zigzag
edges is a planar graph (see Figure \ref{fig:gBlinkWithoutZigzags}).

\end{enumerate}

It is important to note that a g-blink is a combinatorial object,
although, to visualize the connection with blinks, we show drawings
of g-blinks where edges are curves and vertices are points. These
drawings are just to help visualization. The g-blink relevant
information is combinatorial: a set of vertices, the neighbor of
each vertex by each type of edge, the parity of the vertices and the
color of the g-edges.

\begin{figure}[htp]
   \begin{center}
      \leavevmode
      \includegraphics[width=12cm]{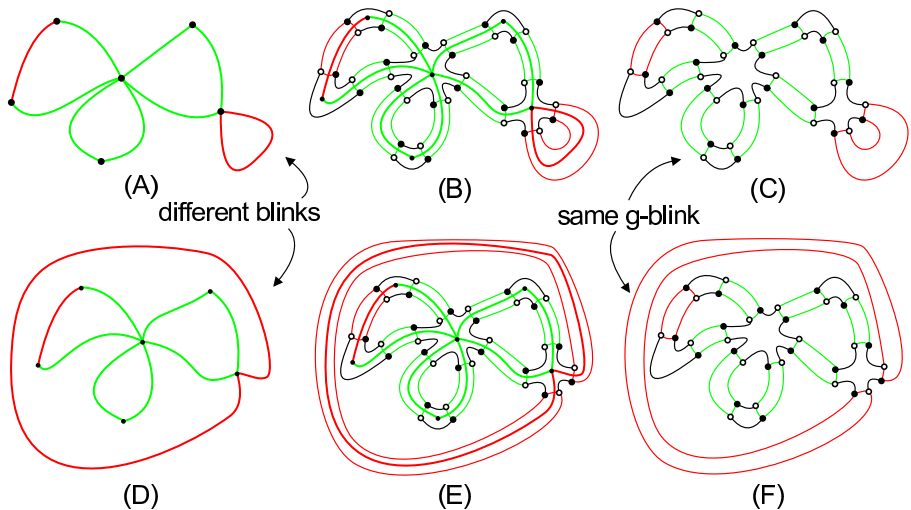}\\
   \end{center}
   \vspace{-0.7cm}
  \caption{Different blinks with the same g-blink}
  \label{fig:DifferentBlinksSameGBlink}
\end{figure}

Note that the definition of g-blinks is independent of blinks. As we
saw, it is a 4-regular graph with some additional structures and
constraints. Now, with this observation in mind, consider the
situation shown on Figure~\ref{fig:DifferentBlinksSameGBlink}. The
blinks of Figure~\ref{fig:DifferentBlinksSameGBlink}A and
Figure~\ref{fig:DifferentBlinksSameGBlink}D are different in the
strict sense (their plane graphs are different) but are
different\footnote{When referring to blinks, this looser sense
concept of ``difference'' is the one we adopted as our convention on
Section \ref{sec:BFL2Blinks}. We could just say the blink of
Figure~\ref{fig:DifferentBlinksSameGBlink}A and the blink of
Figure~\ref{fig:DifferentBlinksSameGBlink}D are different.} in a
looser sense also: there is no plane isotopy between these two
blinks (we would have to tear the red loop on
Figure~\ref{fig:DifferentBlinksSameGBlink}A). On the other hand,
their g-blinks are the same as can be seen on
Figure~\ref{fig:DifferentBlinksSameGBlink}C and
Figure~\ref{fig:DifferentBlinksSameGBlink}F (remember that the edges
on g-blinks presented as drawings are important only to define who
is the neighbor of who, their curve shape is not important).

To obtain a blink from a g-blink we must first embed (forgetting
zigzag-edges) the g-blink on a plane respecting the convention we
used on the procedure \textsc{Blink2GBlink}: the orientation of all
g-vertices induced by orienting the face-edges from white (parity 0
vertex) to black (parity 1 vertex) is always clockwise. Let's name
this convention as {\em convention} $\circlearrowright$. The
embedding part is always possible once a g-blink without
zigzag-edges is a planar graph. It is always possible to
satisfy convention $\circlearrowright$. For example, if all
\hbox{g-vertices} are counterclockwise we may reflect horizontally or
vertically all the embedding correcting the situation. If all
g-vertices are correct except for the external g-vertex (which is
the external face in this case) then we may redraw the curve of an
external angle-edge making it go around all the embedding (see edge
$e$ for an example of this on Figure~\ref{fig:dualGBlinksInduceSameSpace}). If a
blink $B$ is obtained from a g-blink $G$ then we say that $G$
induces $B$.

\begin{figure}[htp]
   \begin{center}
      \leavevmode
      \includegraphics[width=11cm]{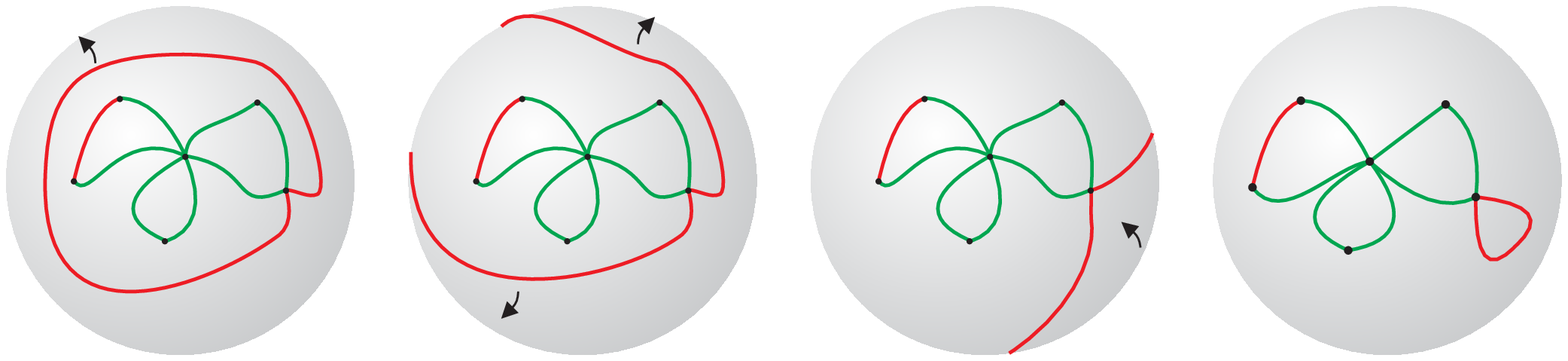}
   \end{center}
   \vspace{-0.7cm}
   \caption{Isotopy on the sphere $\IS^2$}
   \label{fig:SphereIsotopy}
\end{figure}

What are the blinks induced by a g-blink? We must answer this to
continue. Name $A$ the blink of
Figure~\ref{fig:DifferentBlinksSameGBlink}A and $B$ the blink on
Figure~\ref{fig:DifferentBlinksSameGBlink}D. We know there is no
plane isotopy between $A$ and $B$. But $A$ and $B$ are both
obtainable from the same g-blink as
Figure~\ref{fig:DifferentBlinksSameGBlink} shows. How could we
connect $A$ and $B$? The answer is shown on
Figure~\ref{fig:SphereIsotopy}. On the sphere $\IS^2$ there is an
isotopy between $A$ and $B$. One can check that blinks obtainable
from a g-blink are blinks that when embedded on a sphere (draw it on
the plane and then use stereographic projection to get this
embedding) may be transformed one into the other by an isotopy of
the sphere.

Do the blinks obtainable from a g-blink induce the same space? Once
we are, at the end, interested in spaces, g-blink would not be a
useful object if its blinks induced different spaces. The answer
is yes. All blinks of a g-blink induce the same space. The reason is
the \hbox{Blink Jumping Rope Lemma \ref{lem:BlinkJumpingRopeLemma}}.
Note how this Lemma is exactly what is needed to prove that the
blinks of Figure~\ref{fig:DifferentBlinksSameGBlink}A and
Figure~\ref{fig:DifferentBlinksSameGBlink}D induce the same space.

\pagebreak

\begin{Lem}[Blink Jumping Rope Lemma] \label{lem:BlinkJumpingRopeLemma}
The (meta-)blinks shown below induce the same space.
    \begin{center}
        \includegraphics[width=5cm]{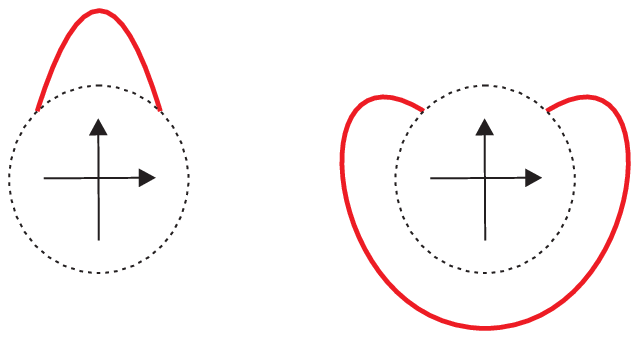}
    \end{center}
\end{Lem}
\begin{proof} Follow the figure below. First we show the BFL
associated with the left blink (the crossing correspondent to red
edge). Then we apply the Jumping Rope Lemma
\ref{lem:jumpingRopeLemma} for BFLs and regular isotopy to get to
our target blink. As our moves preserve the space, we have the
result.
    \begin{center}
        \includegraphics[width=14cm]{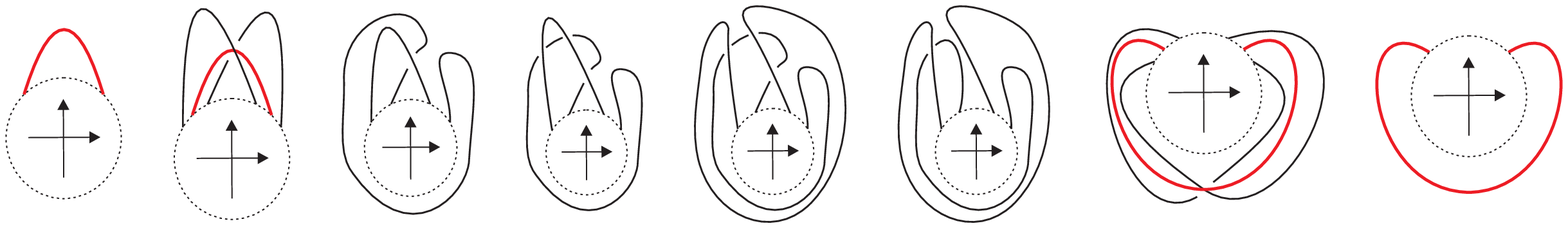}
    \end{center}
\vspace{-1cm}
\end{proof}

With this last result we now define the {\it space of a g-blink} as
the space induced by any blink induced by the g-blink. As we saw,
this space is unique. Observe that the blinks induced by a g-blink
are divided into $|F|$ plane isotopy classes, where $F$ is the set
of g-faces of the g-blink. For each g-face $f \in F$ there is a
blink which has the face corresponding to $f$ as its external face.
So, not considering symmetries that may occur, each g-blink
corresponds to $|F|$ distinct blinks.

Although our initial motivation was to work with blinks, in practice
we did this indirectly through g-blinks. It turned out that this was more
adequate once g-blinks are simpler (\ie to encode a single blink we would
have to have the current g-blink information plus an extra one: which
g-face is the external one), more expressive (\eg one g-blink actually
encodes $|F|$ blinks that induce the same space) and we could prove a
set of g-blink interesting properties that enabled us to do the
experiments we wanted (\eg find all distinct spaces that had a small
blink/BFL/g-blink presentation).

Before ending this section, one last observation: note that a single g-blink
also encodes $|F|$ BFLs: the ones obtained from the $|F|$ blinks by the
\textsc{Blink2BFL} procedure. So we may see a g-blink through $|F|$ blink
views and through $|F|$ BFL views.

\newpage

\section{Homology group from g-blink}
\label{sec:homologyGroup}

The homology group is a topological invariant obtained from the
abelianization of the fundamental group. It is easy to obtain a
presentation of the fundamental group from a blackboard framed link.
However, the problem of deciding if two presentations of a group are
isomorphic is an undecidable problem. This does not occur with the
homology group. It is presented as a pair $(b,t)$, where $b$ is the
{\it Betti number} and $t=(t_1,\ldots,t_p)$ is a sequence with
$p\geq 0$. Each $t_i$ in $t$ is called the $i$-th torsion
coefficient. This sequence also satisfies: $t_1 \geq 2$, if $p > 0$
and $t_i$ divides $t_{i+1}$ for $i<p$. The homology group $(b,t)$
may be obtained from the Smith Normal Form of the linking matrix
of a BFL (see \cite{NemhauserAndWolsey1999} for definition and how to obtain
this normal form). This is so because the linking matrix is a relation matrix for
the homology group. The number of zeros in this diagonal
is the Betti number $b$ and appear all at the end. Throw away the
entries equal to 1. The torsion coefficients $t=(t_1,\ldots,t_p)$
are the other entries on the diagonal. The remainder of this section
shows how to calculate the linking matrix from a g-blink.

\begin{figure}[htp]
   \begin{center}
      \leavevmode
      \includegraphics[width=8cm]{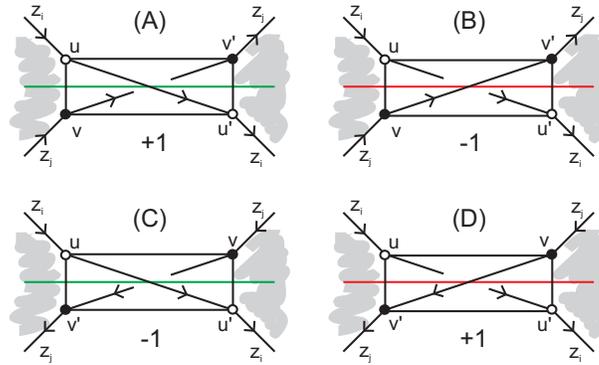}
   \end{center}
   \vspace{-0.7cm}
   \caption{ Signs of a g-edge $a$ for the linking matrix}
   \label{fig:linkingMatrix}
\end{figure}

Let $Z = \{ z_1, \ldots, z_k \}$ be the set of g-zigzags of the g-blink
$G$. So every $z$ in $Z$ is a polygon with alternating zigzag-edges
and angle-edges. We want to define a matrix $N$ of dimension $k
\times k$. First we orient each g-zigzag $z$ in $Z$. This can be
done by mounting a list $v_1, \ldots, v_m$ of the vertices of z such
that $v_i$ is adjacent to $v_{i+1}$ by an edge in $z$, $v_m$ is
adjacent to $v_1$ by an edge in $z$ and the orientation of the edges
in $z$ is defined by the way its end vertices appear in the list: the
edge of $z$ between $v_i$ and $v_{i+1}$ is oriented from $v_i$ to
$v_{i+1}$ for $1\leq i \leq m-1$ and the edge of $z$ whose ends are
$v_1$ and $v_m$ is oriented from $v_m$ to $v_1$. Initialize all
entries of $N$ with zero. For each g-edge $a$, let $u$ and $v$ be
vertices in $a$ such that: $u$ has parity zero (in $V_0$ or white);
$v$ has parity one (in $V_1$ or black); $z_i$ is the g-zigzag
incident to $u$; $z_j$ is the g-zigzag incident to $v$; the
zigzag-edge in $z_i$ incident to $u$ and $u'$ is oriented this way
from $u$ to $u'$; the zigzag-edge in $z_j$ incident to $v$ and $v'$
is oriented this way from $v$ to $v'$. Aligning each g-edge $a$ to
this standard leads to one of the situation shown in
Figure~\ref{fig:linkingMatrix} where the sign $s_a$ of $a$ is also
shown. If $a$ is green and $u$ is adjacent to $v$ by a face-edge
then $s_a=+1$ (Figure~\ref{fig:linkingMatrix}A). If $a$ is red and
$u$ is adjacent to $v$ by a face-edge then $s_a=-1$
(Figure~\ref{fig:linkingMatrix}B). If $a$ is green and $u$ is
adjacent to $v$ by a vertex-edge then $s_a=-1$
(Figure~\ref{fig:linkingMatrix}C). If $a$ is red and $u$ is adjacent
to $v$ by a vertex-edge then $s_a=+1$
(Figure~\ref{fig:linkingMatrix}D).

Knowing the sign of $a$ we update $N$ by $$N_{i,j} \leftarrow
N_{i,j} + s_a$$ and, if $i \neq j$, we also do $$N_{j,i} \leftarrow
N_{j,i} + s_a.$$ Note that $N$ is symmetric. Once $N$ is defined, to
calculate the homology group is to calculate the Smith Normal Form
of $N$ and then collect the pair $(b,t)$ as was already described
above.

\newpage

\section{Quantum invariant from g-blink}
\label{sec:quantumInvariant}

 In this section we show how to calculate
the Witten-Reshetikhin-Turaev quantum invariant for a space from a
g-blink inducing it. This calculation is a translation to g-blinks
of the one showed on \cite{Lins1995} that operates over blackboard
framed links. For further details for this invariant see
\cite{KauffmanAndLins1994}.

The Witten-Reshetikhin-Turaev invariant for a space $M$ is a
function ${\rm wrt}_M:\{3,4,\ldots\} \rightarrow \IC$. This function
maps every integer $i\geq 3$ into a complex number ${\rm wrt_M}(i)
\in \IC.$ If two spaces $A$ and $B$ satisfy ${\rm wrt}_{A}(i) \neq
{\rm wrt}_{B}(i)$ for some $i$ then $A$ and $B$ are different
spaces.

Let $M$ be a space and $r \geq 3$ an integer for which we want to
obtain ${\rm wrt}_M(r).$ Let ${\cal I} = \{0,1,\ldots,r-2\}$. Let
$A$ be a $(4r)$-primitive-root of 1. For $n \in {\cal I}$ define
$$\Delta_n = (-1)^n \, \frac{A^{2n+2}-A^{-2n-2}}{A^{2}-A^{-2}},$$
$${[n]} = \frac{A^{2n}-A^{-2n}}{A^2-A^{-2}} =
(-1)^{n-1}{\Delta_{n-1}}\,.$$ Define $q = A^2$ and, for reasons
inherited from physics, call $[n]$ by {\em q-deformed quantum
integer } and
$$ [n]! = \prod_{1\leq m \leq n}{[m]}$$ by {\em q-deformed
quantum factorial.} Note that although $A$ is a complex number,
$\Delta_n$ and $[n]$ are real numbers. Three numbers $a,b,c \in {\cal
I}$ are said to be an {\em $r$-admissible triple} if $a+b+c \leq
2r-4$ and the numbers $a+b-c,$ $b+c-a,$ $c+a-b$ are non-negative
even numbers.

Let $F$ be the set of g-faces of $G_B,$ $V$ the set of g-vertices of
$G_B$ and $Z$ the set of g-zigzags of $G_B.$ Let $E_a[G_B]$ denote
the angle-edges of $G_B.$ Let $x:F\cup V\cup Z \rightarrow {\cal I}$
be a function that maps an integer in ${\cal I}$ for each g-face,
g-vertex and g-zigzag of $G_B$. We define $x_i = x(i)$ for $i$ in
the domain of $x$. We say that function $x$ is a {\em state.} Denote
by ${\cal X}$ all possible states. Note that ${\cal X}$ is finite.
For every state $x$ exists a complex number $c_x$ defined by
($\alpha$ and $\beta$ are defined after):
$$c_x = \left(\prod_{f \in F}{x_f}\right)
        \left(\prod_{v \in V}{x_v}\right)
        \left(\prod_{z \in Z}{x_z}\right)
        \left(\prod_{a \in E_a[G_B]}\alpha(a,x)\right)
        \left(\prod_{e \in E[B]}\beta(e,x)\right). $$
The value of function ${\rm raw}$ for space $M$ at integer $r$ is the sum
of $c_x$ for every possible state $x$
$${\rm raw}_M(r) = \sum_{x \in {\cal X}}{c_x}.$$

\begin{figure}[htp]
   \begin{center}
      \leavevmode
      \includegraphics{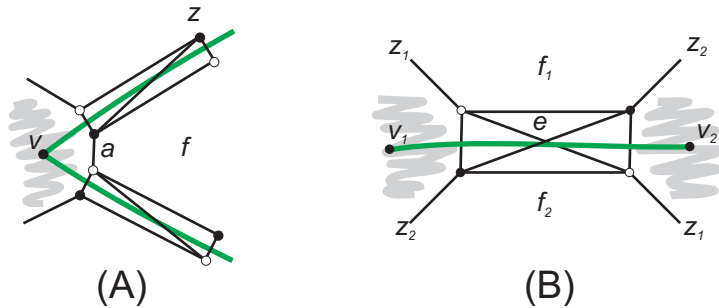}
   \end{center}
   \caption{ Elements for the quantum invariant}
   \label{fig:qialphabeta}
\end{figure}

Now the missing elements: $\alpha$ and $\beta.$ Starting with
$\alpha.$ An angle-edge $a$ may have a drawing like the one shown in
Figure~\ref{fig:qialphabeta}A. Note that the angle-edge $a$ belongs
to one g-face $f$, one g-vertex $v$ and one g-zigzag $z.$ Then we
define
$$\alpha(a,x) = \frac{1}{\theta(x_f,x_v,x_z)}.$$ The function $\theta$
is defined as
$$ \theta(a,b,c) = \left\{%
\begin{array}{ll}
    {\displaystyle \frac{(-1)^{m+n+p}[m+n+p+1]![n]![m]![p]!}{[m+n]![n+p]![p+m]!}}, & \hbox{\footnotesize if $(a, b, c)$ is $r$-admissible;} \\[0.1cm]
    0, & \hbox{\footnotesize otherwise;}\\[0.1cm]
\end{array}%
\right. $$ where $m = (a + b -c)/2$, $n=(b+c-a)/2$, $p=(c+a-b)/2.$

An edge $e$ of the blink $B$ corresponds in $G_B$ to a schema like
the one on Figure~\ref{fig:qialphabeta}B. In this situation, the
elements involved are the g-vertices $v_1$ and $v_2,$ the g-faces
$f_1$ and $f_2$ and the g-zigzags $z_1$ and $z_2.$ It is always
possible, for every edge $e$, to draw a schema like this and
follow this standard: the angle-edges of $z_1$ that appear on the
schema fall between $v_1$ and $f_1$ in one side and between $v_2$
and $f_2$ on the other side. We now define
$$ \beta(e,x) = \left\{
\begin{array}{ll}
    {\displaystyle \frac{{\rm Tet}(x_{f_1},x_{v_1},x_{f_2},x_{v_2},x_{z_2},x_{z_1}) \,\,
                          \lambda(x_{f_1},x_{z_1},x_{v_1}) }
                        { \lambda(x_{v_2},x_{z_1},x_{f_2}) }}, & \hbox{\footnotesize if $e$ is green}
    \\[0.3cm]
    {\displaystyle \frac{{\rm Tet}(x_{f_1},x_{v_1},x_{f_2},x_{v_2},x_{z_2},x_{z_1}) \,\,
                          \lambda(x_{v_2},x_{z_1},x_{f_2}) }
                        { \lambda(x_{f_1},x_{z_1},x_{v_1}) }}, & \hbox{\footnotesize if $e$ is red}
    \\[0.1cm]
\end{array}\right.,$$
where ${\rm Tet}:{\cal I}^6 \rightarrow \mathbb{R}$ is defined as
$${\rm Tet}(a,b,c,d,e,f) = \frac{{\rm Int}!}{{\rm Ext}!} \sum_{m \leq s \leq
M}{\frac{(-1)^s[s+1]!}{\prod_{1 \leq i \leq 4}{[s-a_i]!\prod_{1\leq
j \leq 3}{[b_j-s]!}}}},
$$
in case the triples $(a,b,f)$, $(b,c,e)$, $(c,d,f)$, $(a,d,e)$ are
$r$-admissible and considering
$$\begin{array}{rcl}
{\rm Int}!&=&\prod_{1 \leq i \leq 4 \atop 1 \leq j \leq 3}[b_j-a_i]!\\ [0.2cm] {\rm Ext}!&=&[a]![b]![c]![d]![e]![f]!\\
[0.2cm] a_1&=&\frac{1}{2}(a+b+f)\qquad b_1=\frac{1}{2}(b+d+e+f)\\
[0.2cm] a_2&=&\frac{1}{2}(b+c+e)\qquad b_2=\frac{1}{2}(a+c+e+f)\\
[0.2cm] a_3&=&\frac{1}{2}(c+d+f)\qquad b_3=\frac{1}{2}(a+b+c+d)\\
[0.2cm]
a_4&=&\frac{1}{2}(a+d+e)\qquad m\!=\!\max\{a_i\}\quad M\!=\!\min\{b_j\}.\\
\end{array}.$$
In case any of the triples is not $r$-admissible, the value of ${\rm
Tet}$ is zero. The function $\lambda: {\cal I}^3 \rightarrow \IC$ is
defined by
$$ \lambda(a,b,c) =
\left\{%
\begin{array}{ll}
    (-1)^{(a+b-c)/2}A^{[a(a+2)+b(b+2)-c(c+2)]/2}, & \hbox{\footnotesize if $a,b,c$ is $r$-admissible;} \\
    0, & \hbox{\footnotesize otherwise.} \\
\end{array}%
\right.$$ Finally, the function ${\rm wrt}_M$ is defined as
$$ {\rm wrt}_M(r) = \frac{{\rm raw}_M(r)}{{\rm raw}_{S_1 \times
S_2}(r)}$$

Note that ${\rm wrt}$ is normalized by the raw values of the space
$S_1 \times S_2$. The Figure~\ref{fig:qinvariante} presents the
values of the quantum invariant to the Poincaré Sphere, $E$, for $3
\leq r \leq 30$.
\begin{figure}[htp]
   \begin{center}
      \leavevmode

\begin{footnotesize}
\begin{tabular}{|rrclr|} \hline
 $r$ & \multicolumn{3}{c}{${\rm wrt}_E(r)$} & ev \\[0.1cm] \hline
 3 & 0.7071067811 & + & 0.0000000000$i$ & 2 \\
 4 & -0.5000000000 & + & 0.0000000000$i$ & 4 \\
 5 & -0.3007504775 & - & 0.9256147934$i$ & 6 \\
 6 & 0.2886751346 & + & 0.0000000000$i$ & 9 \\
 7 & -0.8460344491 & - & 0.0447830425$i$ & 12 \\
 8 & 0.0000000000 & - & 0.7325378163$i$ & 16 \\
 9 & -0.1761268770 & + & 0.4020460816$i$ & 20 \\
 10 & -0.7663118960 & - & 0.5567581822$i$ & 25 \\
 11 & 0.2998611170 & - & 0.1557368892$i$ & 30 \\
 12 & -0.7886751345 & + & 0.1830127018$i$ & 36 \\
 13 & -0.1148609711 & - & 0.7426524382$i$ & 42 \\
 14 & -0.1074423864 & + & 0.3977522621$i$ & 49 \\
 15 & -0.7770955704 & - & 0.5344039501$i$ & 56 \\
 16 & 0.3141711649 & - & 0.1762214752$i$ & 64 \\ \hline
\end{tabular}
\begin{tabular}{|rrclr|} \hline
 $r$ & \multicolumn{3}{c}{${\rm wrt}_E(r)$} & ev \\[0.1cm] \hline
 17 & -0.7804263387 & + & 0.1428530500$i$ & 72 \\
 18 & -0.0590950525 & - & 0.7636697702$i$ & 81 \\
 19 & -0.1301847177 & + & 0.3730119013$i$ & 90 \\
 20 & -0.7085827791 & - & 0.6254313947$i$ & 100 \\
 21 & 0.3410488374 & - & 0.1495290291$i$ & 110 \\
 22 & -0.7854601781 & + & 0.0248114386$i$ & 121 \\
 23 & 0.0600389356 & - & 0.7749612722$i$ & 132 \\
 24 & -0.1814470028 & + & 0.3376768599$i$ & 144 \\
 25 & -0.5895059790 & - & 0.7441570346$i$ & 156 \\
 26 & 0.3666499557 & - & 0.0969412734$i$ & 169 \\
 27 & -0.7726037705 & - & 0.1263662241$i$ & 182 \\
 28 & 0.2079977942 & - & 0.7581679950$i$ & 196 \\
 29 & -0.2393556663 & + & 0.2887208942$i$ & 210 \\
 30 & -0.4276587373 & - & 0.8531721152$i$ & 225 \\ \hline
\end{tabular}
\end{footnotesize}
   \end{center}
   \vspace{-0.3cm}
   \caption{ Example of quantum invariant: Poincarè's sphere}
   \label{fig:qinvariante}
\end{figure}

\begin{figure}[h!]
        \begin{center}
            \includegraphics[width=4cm]{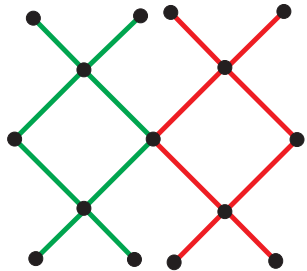}
        \end{center}
        \caption{A peculiar space: its quantum invariants are integers and every integer appears}
        \label{fig:strange}
    \end{figure}

By playing with computation of the quantum invariants from various
blinks we discovered a rather peculiar space.
\begin{Conj}
The quantum invariants
of the space induced by the blink of Figure~\ref{fig:strange} are:
$q_r = \frac{2}{3}\,r$ if $r \equiv 0$ mod $3$, $q_r =
\frac{1}{3}\,(r+1)$ if $r \equiv 2$ mod $3$, $q_r =
\frac{1}{3}\,(r-1)$ if $r \equiv 1$ mod $3$.
\end{Conj}
We have checked this
result to a precision of 10 decimal places and up to $r=45$. The
fact that the quantum invariants are real is evident since the blink
is red-green symmetric. The fact that they are all integer values
and that every integer appears is rather pleasing.

\newpage

\section{Code of a g-blink}
\label{sec:gblinkcode}

\newcommand{\pack}{\textsc{Pack}}
\newcommand{\reds}{\textsc{Reds}}
\newcommand{\GBlinkLabel}{\textsc{GBlinkLabel}}

In this section we are interested in defining a ``word'' with all
information of a g-blink, one that from it we are able to rebuild
the g-blink. This word is said to be the {\em code of the g-blink}.
Let $G$ be a g-blink. One of the ingredients to define this
``code'' is the \GBlinkLabel\ algorithm that labels the vertices of
a g-blink from an initial vertex $v$ (this initial vertex will be
labeled~1).
\begin{algorithm}
\begin{small}
\caption{$\GBlinkLabel(G,v)$}
\label{alg:maplabel}
\begin{algorithmic}[1]
 \State $S \leftarrow$ empty stack; $i \leftarrow 1$; $\forall u, L_u \leftarrow
 \bot$ \Comment{$\bot$ = not defined}
 \State push $v$ into $S$
 \While{$S$ not empty}
    \State $a \leftarrow$ pop $S$
    \If{$L_a = \bot$}
       \State $b \leftarrow {\rm adj}_f(a)$; $c \leftarrow {\rm adj}_v(b)$; $d \leftarrow {\rm adj}_v(a)$
       \State $L_a \leftarrow i; L_b \leftarrow i+1; L_c \leftarrow i+2; L_d \leftarrow i+3$
       \State push ${\rm adj}_a(b)$ into $S$; push ${\rm adj}_a(d)$ into $S$
       \State $i \leftarrow i+4$
    \EndIf
 \EndWhile
 \State {\bf return} $L$
\end{algorithmic}
\end{small}
\end{algorithm}

When we talk about a {\em labeling of a g-blink or of a blink}, we
are referring to a labeling of the vertices of the g-blink given by
{\GBlinkLabel} with a starting vertex being some vertex with parity
1 in $G$. With this constraint, the set of vertices with even label
defined by {\GBlinkLabel}  is exactly the set $V_0$ of $G$ and the
set of vertices with odd label is exactly the set $V_1$ of $G$.
Other important properties of a labeling are: adjacent vertices by
face, vertex or angle edges in $G$ have labels with different
parity; the vertices of the same g-edge have labels $4k-3$, $4k-2$,
$4k-1$, $4k$ for some $k \geq 1$. From the label of a vertex it is
possible to know the label of its neighbor by face, vertex and
zigzag edge. For instance, if $u$ has label $4k-2$ (for some integer
$k\geq 1$) then its neighbor by face edge has label $4k-3$, by
vertex edge has label $4k-1$ and by zigzag edge has label $4k$. One
consequence of this fact is that it is possible to rebuild all
edges of $G$ annotating only the angle edge's neighbors, once the
face edge, vertex edge and zigzag edge are all known from the vertex
label.

Let $L$ be a labeling for $G$. Let $a_1,a_2 \ldots, a_{4n}$ be the
labels of the adjacent vertices by angle edges of the vertices
$1,2,\ldots,4n$ under the $L$ labeling. As we saw, this list is
sufficient to restore the vertices and the edges of $G$. Note also
that, by the property that adjacent vertices have labels with
different parity, this list is made of even labels followed by odd
labels and that if $a_i = j$, then $a_j = i$. From these two
observations it follows that from $\frac{a_1}{2}, \frac{a_3}{2},
\ldots, \frac{a_{4n-1}}{2}$ it is possible to restore $a_1,a_2
\ldots, a_{4n}$ and, consequently, the vertices and edges of $G$. We
denote the list (with labels divided by 2) as the
{\em packed representation of $L$}. Note that the packed
representation is a permutation of $1,\ldots,2n$. If $L$ is a
labeling, we denote by $\pack(L)$ the packed representation of $L$.
\begin{figure}[htp]
   \begin{center}
      \leavevmode
      \includegraphics{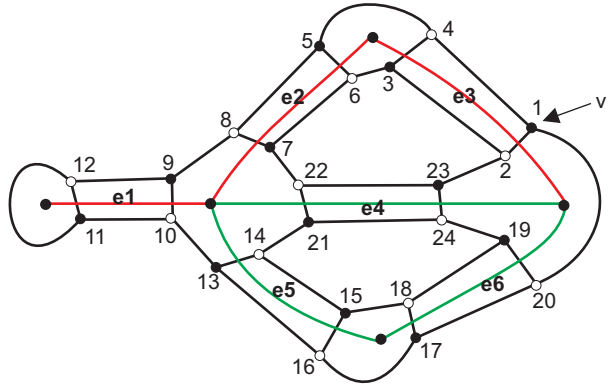}
   \end{center}
   \vspace{-0.7cm}
   \caption{ blink $B$, g-blink $G_B$ and labeling $\textsc{GBlinkLabel}(G_B,v)$}
   \label{fig:gBlinkLabeling}
\end{figure}

The Figure \ref{fig:gBlinkLabeling} presents a blink $B$, its
induced g-blink $G_B$ and the labeling resulted of
${\GBlinkLabel}(G_B,v)$. In this case, the label of the adjacent
vertices by angle edge of $1,\ldots,24$ are $a_1,a_2,\ldots,a_{24}
=$ 20, 23, 6, 5, 4, 3, 22, 9, 8, 13, 12, 11, 10, 21, 18, 17, 16, 15,
24, 1, 14, 7, 2, 19. Its packed representation is $\frac{a_1}{2},
\frac{a_3}{2}, \ldots, \frac{a_{23}}{2} = $ 10, 3, 2, 11, 4, 6, 5,
9, 8, 12, 7, 1.

We represent the bicoloration of the g-edges of a g-blink under the
labeling $L$ by the set of integers $\reds(G,L)$ defined this way:
$k$ is in $\reds(G,L)$ if g-edge with vertices
$4k-3$, $4k-2$, $4k-1$ and $4k$ is red, otherwise $k$
is not in $\reds(G,L)$.

Let $L$ be a labeling for the g-blink $G$, then the {\em pre-code}
of $G$ for labeling $L$ is the pair $$(\pack(L),\reds(G,L)).$$
In the example of Figure \ref{fig:gBlinkLabeling}, the edges $e_1,$
$e_2,$ $e_3,$ $e_4,$ $e_5,$ $e_6$ are labeled with 3, 2, 1, 6, 4, 5
respectively. It follows that the {\em pre-code} for this blink
under the presented labeling (starting at $v$, Figure
\ref{fig:gBlinkLabeling}) is: $$((10, 3, 2, 11, 4, 6, 5, 9, 8, 12,
7, 1),\{1,2,3\}).$$

It is easy to see that different labelings define different
pre-codes and that the same blink may have different labelings
by changing the vertex~1 on the procedure {\GBlinkLabel}. This creates a
difficulty: two different pre-codes for the same g-blink. To resolve
this we define the order relation $\preceq$  on the set of pre-codes. Let
$(\pi_1,R_1)$ and $(\pi_2,R_2)$ be two pre-codes, then
$$
(\pi_1,R_1) \preceq (\pi_2,R_2) \hbox{ if }
\left\{
\begin{array}{l}
  |\pi_1| < |\pi_2| \,\,\, \hbox{ or }\\
  |\pi_1| = |\pi_2| \hbox{ and } \pi_1 < \pi_2 \,\,\, \hbox{ or }\\
  \pi_1 = \pi_2 \hbox{ and } |R_1| < |R_2| \,\,\, \hbox{ or }\\
  \pi_1 = \pi_2 \hbox{ and } R_1 = R_2 \,\,\, \hbox{ or }\\
  \pi_1 = \pi_2 \hbox{ and } |R_1| = |R_2| \hbox{ and } \min(R_1 \backslash R_2) < \min(R_2 \backslash R_1),\\
\end{array}
\right.
$$
where $|\pi|$ is the length of the permutation $\pi$ and $|R|$ is
the size of set $R$. The {\em code of the g-blink} $G$ is its greatest pre-code under
the relation $\preceq$:

\newcommand{\MaxPreCode}{\lower1.5ex\hbox{$\buildrel {\displaystyle max} \over {\scriptstyle \preceq}$}}

$$ \kappa(G) =
   \MaxPreCode \left\{
(\, \pack(L_v),\reds(G,L_v)\,)  \, \left|
\begin{array}{l}
      L_v = \GBlinkLabel(G,v), \\
      v \in V_1[G]
\end{array} \right. \right\}.$$
The {\em code of a blink $B$} is defined as
$\kappa(B) = \kappa(G),$ where $G$ is the induced g-blink of $B$.
A labeling $L$ of a g-blink is said to be a {\it code labeling of $G$}
if $(\,\pack(L),\reds(G,L)\,) = \kappa(G)$. We extend the relation
$\preceq$ on pre-codes to g-blinks and blinks in this natural way:
g-blink $G_1$ {\it is smaller or equal to} g-blink $G_2$, $G_1 \preceq G_2$,
if $\kappa(G_1) \preceq \kappa(G_2)$; blink $B_1$ {\it is smaller
or equal to} blink $B_2$, $B_1 \preceq B_2$, if their induced g-blinks
satisfy $G_1 \preceq G_2$.

\newpage

\section{\textsc{Dual}, \textsc{Reflection} and \textsc{RefDual}
of a g-blink} \label{sec:dualReflectionRefDual}

In this section we study the effects of simple changes on the structure
of a g-blink. For instance, what happens to the induced space of a g-blink
if we swap the parity of its vertices? And what happens to its blink
presentations if we do this? We are interested in studying three types of
modifications in the structure of a g-blink. One of them is swapping the
parity of the vertices and we denote it by (P). Before naming the
other two we establish the convention we use to encode g-blinks.

We saw in the definition of g-blinks that we may encode the red-green
coloring of the \hbox{g-edges} directly or, alternatively, we may encode it by
registering the overcross/undercross status of the zigzag-edges of each
g-edge. In this section we assume that we are using this second
alternative. So, here, the color of the g-edges is a consequence of the
overcross/undercross state of the zigzag-edges of the g-blink.

Besides (P), the other two modifications in the structure of a g-blink that
we study are: swapping the role of face-edges with vertex-edges, denoted
by (FV); and swapping the overcrossing/undercrossing state of each
zigzag-edge, denoted by (C).

\begin{figure}[htp]
   \begin{center}
      \leavevmode
      \includegraphics[width=10cm]{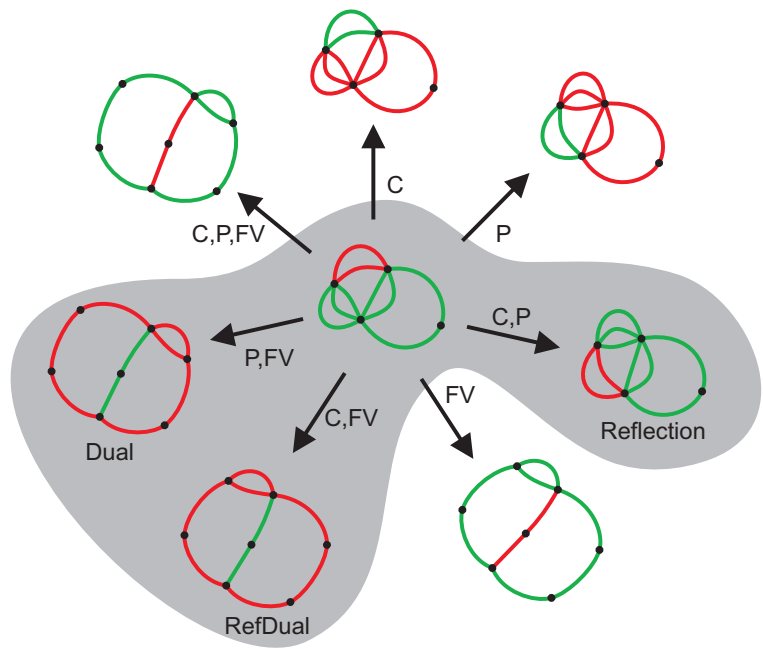}
   \end{center}
   \vspace{-0.7cm}
   \caption{The effect on blinks of applying all combinations of (C), (P)
   and (FV) on its g-blink}
   \label{fig:blinkWith8Operations}
\end{figure}

The central blink in Figure~\ref{fig:blinkWith8Operations} is a blink
presentation for our reference g-blink. By applying all combinations of
(P), (C) and (FV) on this g-blink we obtain new g-blinks
inducing the blinks shown. We can learn from this figure the effects
on blinks of these g-blink modifications. Applying (C), \ie changing the
undercross/overcross status on the zigzag-edges, the only effect is
to swap the colors of the edges of the blinks; applying (P), \ie
changing the parity of the vertices, the effect is to do
one reflection of the blink drawing and change the color of
its edges; applying (C) and (P), \ie changing the undercross/overcross
state of each zigzag-edge and the parity of the vertices,  the effect
is just a reflection of the blink drawing; applying (FV), \ie swapping
the roles of face-edge and vertex-edge, the blink becomes
the dual of the original blink (the dual is in the sense of a dual map or
dual plane graph) whose dual edges preserve
the same color as the original edges, followed by a
reflection; applying (C) and (FV), \ie swapping the roles of face-edge and vertex-edge and
changing the overcross/undercross status of the zigzag-edges, the effect is
a dual blink with the dual edges having changed color (\eg a dual edge that
``crosses'' a red edge becomes a green one) followed by a reflection;
applying (P) and (FV) the effect is the dual blink with dual edges
having the changed colors (\eg a dual edge that
``crosses'' a red edge becomes a green one); applying (C), (P) and (FV),
\ie all three modifications, the effect is the dual of the blink with
the dual edges having the same color as the original ones (\eg a dual edge
that crosses a green edge is itself a green edge).
Note that the blinks (and g-blinks) differing by the application of
two distinct modifications have special names. If $G$ is a g-blink,
the g-blink obtained from $G$ by applying (C) and (P) is said to be
the {\it reflection} of $G$ and is denoted by \textsc{Reflection}($G$);
the g-blink obtained from $G$ by applying (C) and (FV) is said to be
the {\it refdual} of $G$ and is denoted by \textsc{RefDual}($G$); the
g-blink obtained from $G$ by applying (P) and (FV) is said to be
the {\it dual} of $G$ and is denoted by \textsc{Dual}($G$). They
were shown in Figure~\ref{fig:blinkWith8Operations} over a gray
region because they all share one important property as the next
proposition shows.

\newpage

\begin{Prop}\label{prop:dualGblink}
The spaces induced by g-blinks $G$, \textsc{Reflection($G$)},
\textsc{RefDual($G$)} and \textsc{Dual($G$)} are the same.
\end{Prop}

\vspace{-0.5cm}

\begin{proof}
$(G \Sequiv \textsc{Dual}(G))$ -- Consider the dual g-blinks on
Figure~\ref{fig:dualGBlinksInduceSameSpace}A and
Figure~\ref{fig:dualGBlinksInduceSameSpace}D. Following the drawings
in each row of this figure we see how to obtain one induced BFL from
a g-blink. Now observe that the BFLs on
Figure~\ref{fig:dualGBlinksInduceSameSpace}C and
Figure~\ref{fig:dualGBlinksInduceSameSpace}F induce the same space
because we can get from one to the other by applying, on $e$, the
space preserving move for BFLs of Lemma \ref{lem:jumpingRopeLemma}
(Jumping Rope Lemma for BFLs). It is easy to see that this argument
generalizes to any pair $G$ and  $\textsc{Dual}(G)$, so
$G \Sequiv \textsc{Dual}(G)$.
\begin{figure}[htp]
   \begin{center}
      \leavevmode
      \includegraphics[width=13cm]{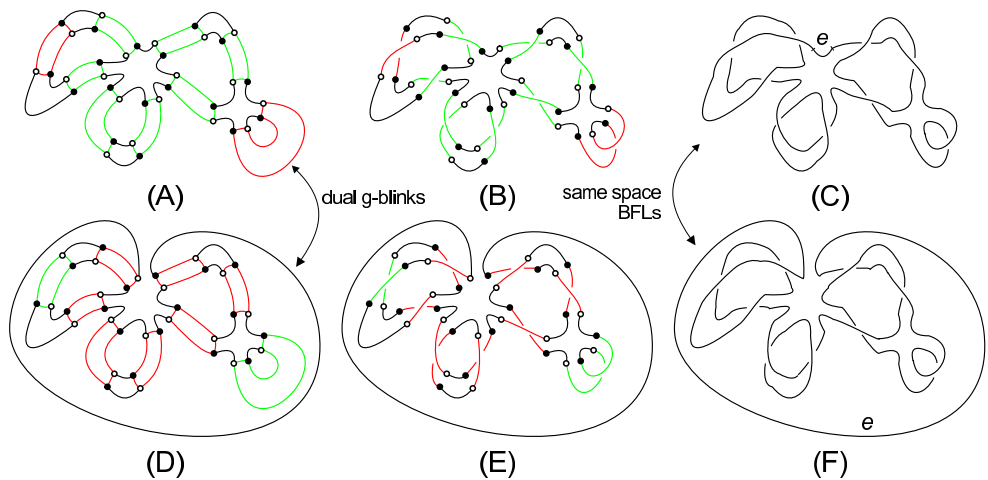}
   \end{center}
   \vspace{-0.7cm}
   \caption{Dual g-blinks induce the same space}
   \label{fig:dualGBlinksInduceSameSpace}
\end{figure}

$(G \Sequiv \textsc{Reflection}(G))$ --
We prove the result in BFL language. Consider a plane disk $D^2$
containing the BFL. Do a $3D$-flip of $D^2$ carrying the BFL along.
Clearly this maintains the ambient isotopy link associated to the
BFL. The writhe of the components do not change: the blink has been
reflected but also its crossings are switched, thus maintaining all
crossing signs. The Proposition is established.

$(G \Sequiv \textsc{RefDual}(G))$ -- Note that
$\textsc{RefDual}(G) = \textsc{Reflection}(\textsc{Dual}(G))$, once applying
(P) and (FV) followed by (P) and (C) is the same as applying only (FV) and (C).
So, by the transitivity of the $\Sequiv$ relation, using the previous two results,
$G \Sequiv \textsc{RefDual}(G)$.
\end{proof}

What about the blinks that did not fall in the gray region of
Figure~\ref{fig:blinkWith8Operations}? What spaces do they induce?
First, it is easy to see that they all induce the same space once,
taking as reference the top most blink (north), the northeast
blink is its reflection (\ie to get there we must apply (C) and (P)),
the southeast blink is its refdual and the northwest blink is its dual.
So, by Proposition~\ref{prop:dualGblink} they induce the same space.
To finish the answer, let's focus again on the top most blink (result
of (C) operation). It is obtained from the central blink by changing
crossing status of the zigzag-edges. In the blink view of the g-blink
this is equivalent to swap the colors of the edges from green to
red and vice-versa. On the BFL view this is just that: change all
the crossings. This has the effect of inverting the writhe of all
components: for example, a component that had writhe 1 becomes one
with writhe -1. So the end effect of this change
is to invert the orientation of the original space. Conclusion:
the g-blinks on the white region of Figure~\ref{fig:blinkWith8Operations}
induce the same space of the gray region g-blinks except for the
orientation that is changed.

\newcommand{\summaryOfInterpretations}{
{ \small \setstretch{1.1}
\setlength{\parskip}{-2pt} \sl
\begin{tabular}{|ccc|} \hline
\multicolumn{3}{|c|}{\textsc{Dual($G$)}} \\
g-blink & blink & BFL \\ \hline
\multicolumn{1}{|p{4.5cm}|}{\protect \centering change parity (P) and swap face-edges and vertex-edges (FV) \\[0.3cm] \includegraphics{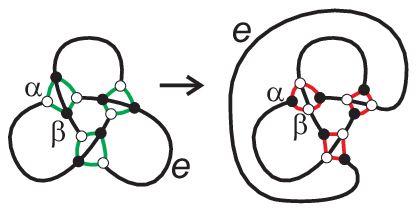}} &
\multicolumn{1}{|p{4.5cm}|}{\protect \centering each face becomes a
vertex, each
edge becomes a dual edge with different color \\[0.3cm] \includegraphics{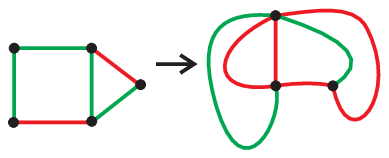}} &
\multicolumn{1}{|p{4.5cm}|}{\protect \centering overpass one external edge \\[0.3cm] \includegraphics{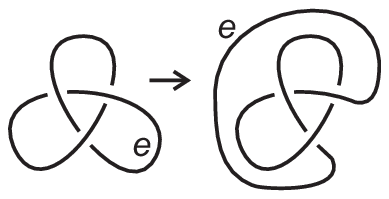}} \\ \hline
\multicolumn{3}{|c|}{\textsc{Reflection($G$)}} \\[0.1cm]
g-blink & blink & BFL \\ \hline \multicolumn{1}{|p{4.5cm}|}{\protect
\centering change parity (P) and change overcross and undercross
status on zigzag-edges (C) \\[0.3cm]
\includegraphics{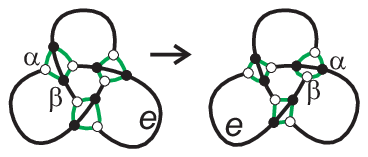}}  &
\multicolumn{1}{|p{4.5cm}|}{\protect \centering
reflect \\[0.3cm] \includegraphics{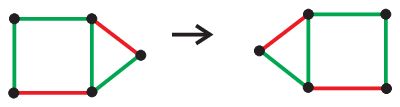}} &
\multicolumn{1}{|p{4.5cm}|}{\protect \centering
reflect and change the crossings  \\[0.3cm]
\includegraphics{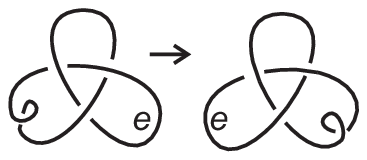}} \\ \hline
\multicolumn{3}{|c|}{\textsc{RefDual($G$)}} \\
g-blink & blink & BFL \\ \hline \multicolumn{1}{|p{4.5cm}|}{\protect
\centering change overcross and undercross
status on zigzag-edges (C) and swap face-edges and vertex-edges (FV) \\[0.3cm]
\includegraphics{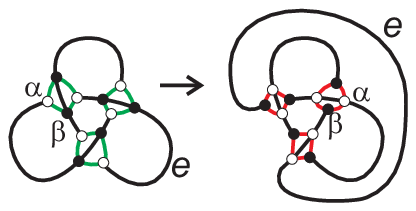}} &
\multicolumn{1}{|p{4.5cm}|}{\protect \centering first make each face
become a vertex and each
edge become a dual edge with the color changed, then reflect
the result \\[0.3cm]
\includegraphics{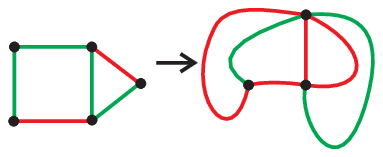}} &
\multicolumn{1}{|p{4.5cm}|}{\protect \centering overpass one
external edge, reflect and change the crossings \\[0.2cm]
\includegraphics{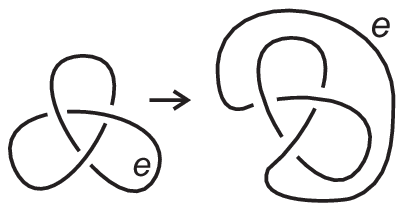}} \\[0.3cm] \hline
\end{tabular} \bigskip \medskip}}

One important consequence of the properties we described here is that we
may search for distinct spaces on only one of the eight possible g-blinks.
The other orientation of the space is trivially obtained from any g-blink.
This saves computational effort of identifying distinct spaces.

We end this section by summarizing the effects on the structures of
a g-blink, blink and BFL when the dual, reflection
and refdual operations are applied.

\summaryOfInterpretations

\newpage

\section{Merging and breaking g-blinks} \label{sec:mergingGBlinks}

\begin{figure}[htp]
   \begin{center}
      \leavevmode
      \includegraphics[width=14cm]{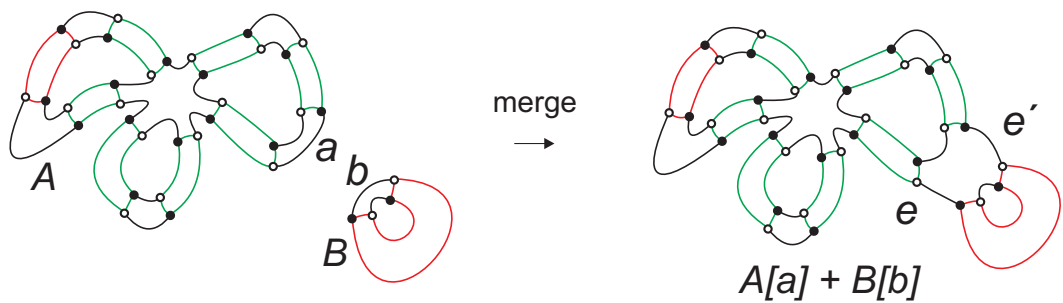}
   \end{center}
   \vspace{-0.7cm}
   \caption{Merging of $A$ and $B$ on {\it basepair} $(a,b)$}
   \label{fig:merging}
\end{figure}

Let $A$ and $B$ be distinct g-blinks. A {\it basepair on $A$ and $B$}
is a pair of angle-edges $(a,b)$ so that $a \in A$ and $b \in B$. The
{\it merging of $A$ and $B$ at basepair $(a,b)$}, denoted by $$ A[a] + B[b] \,\,
,$$ is the g-blink obtained by replacing $a$ and $b$ by new edges
$e$ and $e'$ both connecting $A$ to $B$, having the
same ends as $a$ and $b$ and linking vertices of distinct parity.
See Figure~\ref{fig:merging} for an example.

Observe the result of the merging of Figure~\ref{fig:merging}. The edges $e$
and $e'$ are both incident to the same g-face and g-vertex and we could reverse
the merging by replacing $e$ and $e'$ back with $a$ and $b$. Indeed, any pair of
distinct angle-edges incident to the same g-face and g-vertex on a g-blink
defines a {\it breakpoint}: a point where we can break a g-blink into two
disconnected g-blinks. To {\it break} a g-blink on {\it breakpair} $(e,e')$
is to separate it into two g-blinks by replacing edges $e$ and $e'$
by two new edges incident to same vertices of $e$ and $e'$ obtaining
two disconnected g-blinks. For an example see Figure~\ref{fig:merging}
from right to left.

\begin{Theo}[Theorem on partial dual]
\label{theo:partialDual} Let $A$ and $B$ be arbitrary disjoint
g-blinks and $(a,b)$ a basepair on them. Then
 $A[a] + B[b] \Sequiv A[a] + \textsc{Dual}(B)[b]$.
\end{Theo}
In the language of BFLs the diagrammatic reformulation of Theorem
\ref{theo:partialDual} is given by the diagram below. Note that
$\alpha$ and $\beta$ are the ends of $a$ and $\gamma$ and $\delta$
are the ends of $b$.
$$\includegraphics[width=5cm]{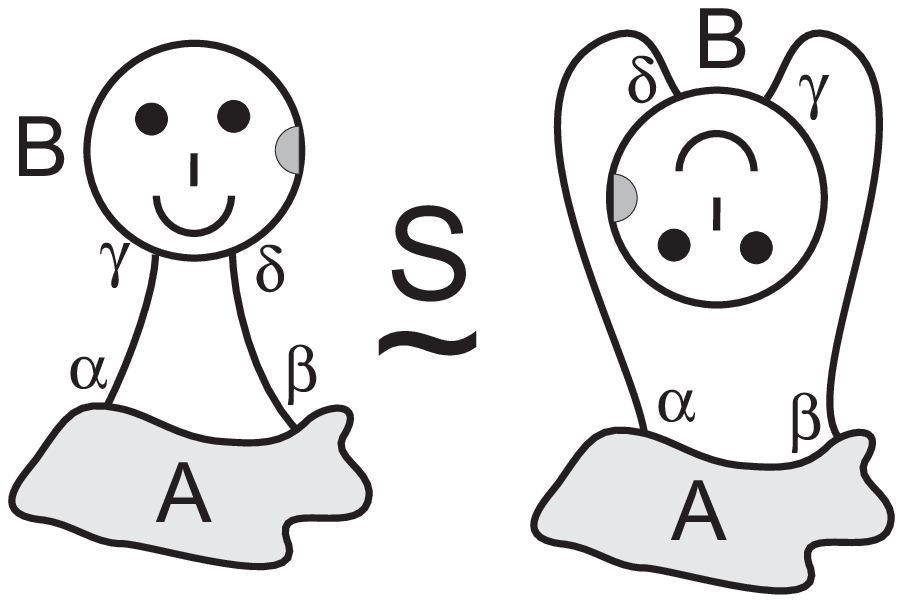}$$ The right
diagram above is obtained by cutting the two wires $\pi$-rotating
$B$ and reconnecting the wires. (Note that the smiling fellow has
only the right rear.) Theorem~\ref{theo:partialDual} was suggested by
computer experiments very early in our research. It is a central
result to curtail the number of relevant blinks: see next section.
Its proof, however, was elusive until October 31, 2006: it has to
wait for the proofs of Theorems \ref{theo:partialReflection} and
\ref{theo:partialRefDual}.

\begin{Theo} [Theorem on partial reflection] \label{theo:partialReflection}
Let $A$ and $B$ be arbitrary disjoint g-blinks, $(a,b)$ a basepair
on them. Then $\, A[a] + B[b] \Sequiv A[a] +
\textsc{Reflection}(B)[b].$
\end{Theo}

In the language of blinks the diagrammatic reformulation of Theorem
\ref{theo:partialReflection} is on the left part of the diagram
below. The right part of it is the reformulation of the same Theorem
in BFL language. Note that the right ear becomes a left ear
indicating a $B$-reflection.
$$\includegraphics[width=10cm]{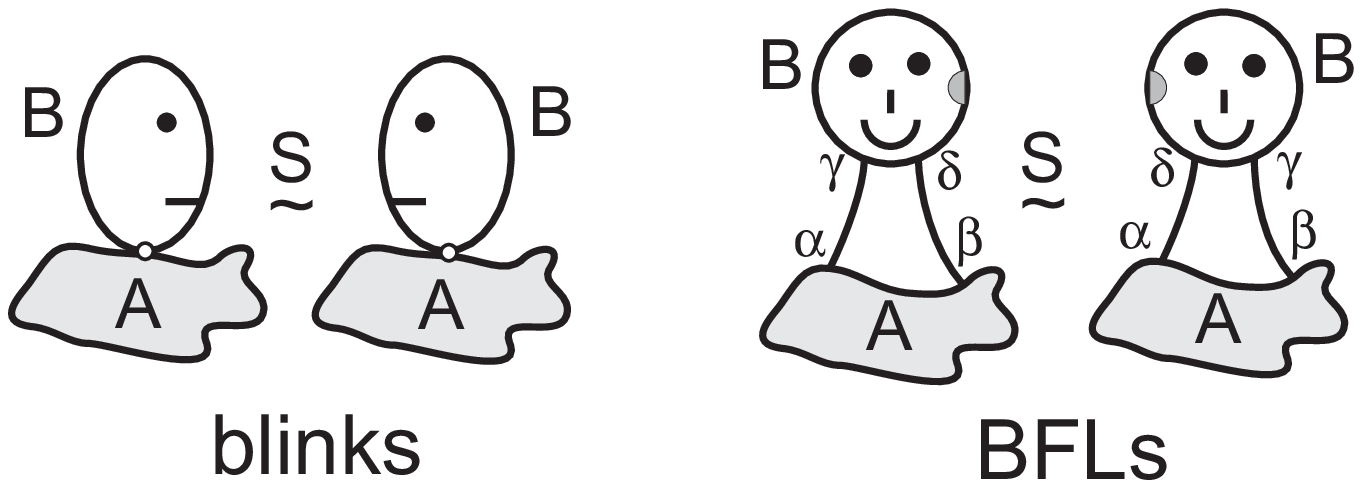}$$
Theorem \ref{theo:partialReflection} is proved by topological
techniques allied to crucial facts on the theory of gems. It is done
in Chapter 4. The proof of Theorem \ref{theo:partialReflection} is
the main theoretical contribution of this thesis.

The diagrammatic reformulation of Theorem
\ref{theo:partialRefDual} in the language of BFLs is the passage
from the first to the third diagram below
$$\includegraphics[width=8cm]{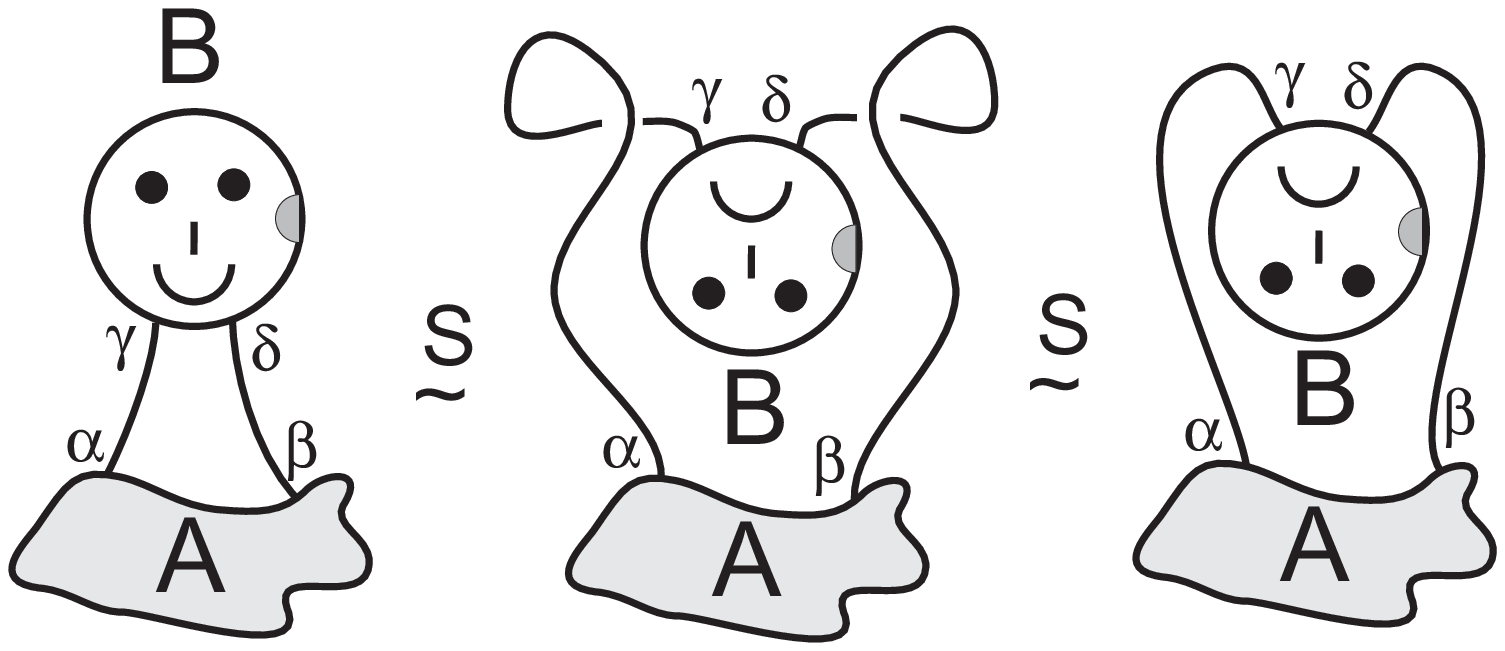}$$ The third
diagram above is obtained by a $3D$-flip on $B$ (getting the central
diagram) followed by a ribbon move, regular isotopies and Whitney
trick. The smiling to frowning change is to indicate that all the
crossings are switched (and that the fellow became angry for being
put upside down and being retracted from the ear). On the contrary
of the previous Theorems we can tackle the proof of Theorem
\ref{theo:partialRefDual} immediately.

\begin{Theo}[Theorem on partial refDual]
\label{theo:partialRefDual} Let $A$ and
$B$ be arbitrary disjoint g-blinks, $(a,b)$ a basepair on them. Then
$A[a] + B[b] \Sequiv A[a] + \textsc{RefDual}(B)[b].$
\end{Theo}

\begin{figure}[htp]
   \begin{center}
      \leavevmode
      \includegraphics[width=8cm]{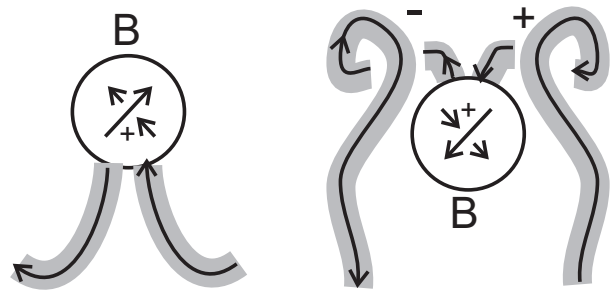}
   \end{center}
   \vspace{-0.7cm}
   \caption{BFLs induce same space because are the same link with the same
   writhe at each component}
   \label{fig:smilingFrowningFellow2}
\end{figure}

\begin{proof}
The proof is easy with the help of the BFL manifestation of the
Theorem. See Figure \ref{fig:smilingFrowningFellow2}. The ambient
isotopy classes of the links corresponding to $A[a] + B[b]$ and to
$A[a] + \textsc{RefDual}(B)[b]$ are the same. It is enough to prove
that the writhe of each component of the BFLs is maintained. Outside
the $B$ there is no change in the crossing numbers. In the interior
of $B$ the crossings are switched and reflected (become upside down)
thus, again, there is no change in the crossing numbers. Finally,
the crossing numbers of the new curls are in the same component and
cancel each other.
\end{proof}

\begin{Lem} For any g-blink $B$,
$$\textsc{Dual(Reflection}(B))=\textsc{Reflection}(\textsc{Dual}(B)).$$
\end{Lem}
\begin{proof} By their combinatorial definitions in g-blinks
the operations of taking the dual and reflecting are seen to be
commuting involutions. Thus, both spaces in the statement of the
lemma are equal to $\textsc{RefDual}(B))$.
\end{proof}

\begin{Lem} \label{lem:23implies1}
Theorem \ref{theo:partialDual} is implied by Theorems
\ref{theo:partialReflection} and \ref{theo:partialRefDual}.
\end{Lem}
\begin{proof}
$A[a] + B[b] \Sequiv A[a] + \textsc{Reflection}(B)[b]$ and $A[a] +
B[b] \Sequiv A[a] + \textsc{RefDual}(B)[b]$ imply by transitivity
that $A[a] + \textsc{Reflection}(B)[b] \Sequiv A[a] +
\textsc{RefDual}(B)[b]$. Note that $b$ is an angle-edge in $C
=\textsc{Reflection}(B)$. Taking $c=b$ we have for any blink $C$,
$A[a] + C[c] \Sequiv A[a] + \textsc{Dual}(C)[c]$, establishing
Theorem \ref{theo:partialDual} for arbitrary disjoint g-blinks
$(A,C)$ and basepairs $(a,c)$.
\end{proof}
From Lemma \ref{lem:23implies1} and Theorem
\ref{theo:partialRefDual}, Theorem \ref{theo:partialDual} will
follow from Theorem \ref{theo:partialReflection}. The proof of this
result is given at the end of Chapter 4.

\newpage

\section{Representative of a g-blink}
\label{sec:gblinkrepresentative}

We learned on Section \ref{sec:gblinks} that a g-blink induces
different blinks. All these blinks induce the same space, which is
defined as the space of the g-blink. We saw also that different
g-blinks may induce the same space: the g-blinks $G$,
$\textsc{Reflection}(G)$, $\textsc{Dual}(G)$ and $\textsc{RefDual}(G)$
dual are different g-blinks but induce the same space.
In this section we define a normalization procedure for g-blinks.
This normalization maps a g-blink into another g-blink that induces
the same space as the first. Our goal with this procedure is to look
for different spaces on fewer g-blinks: we need to look for different spaces
only on g-blinks that are normalized. The normalized version of a
g-blink will be denoted as its {\it representative}.

\begin{figure}[htp]
   \begin{center}
      \leavevmode
      \includegraphics[width=14cm]{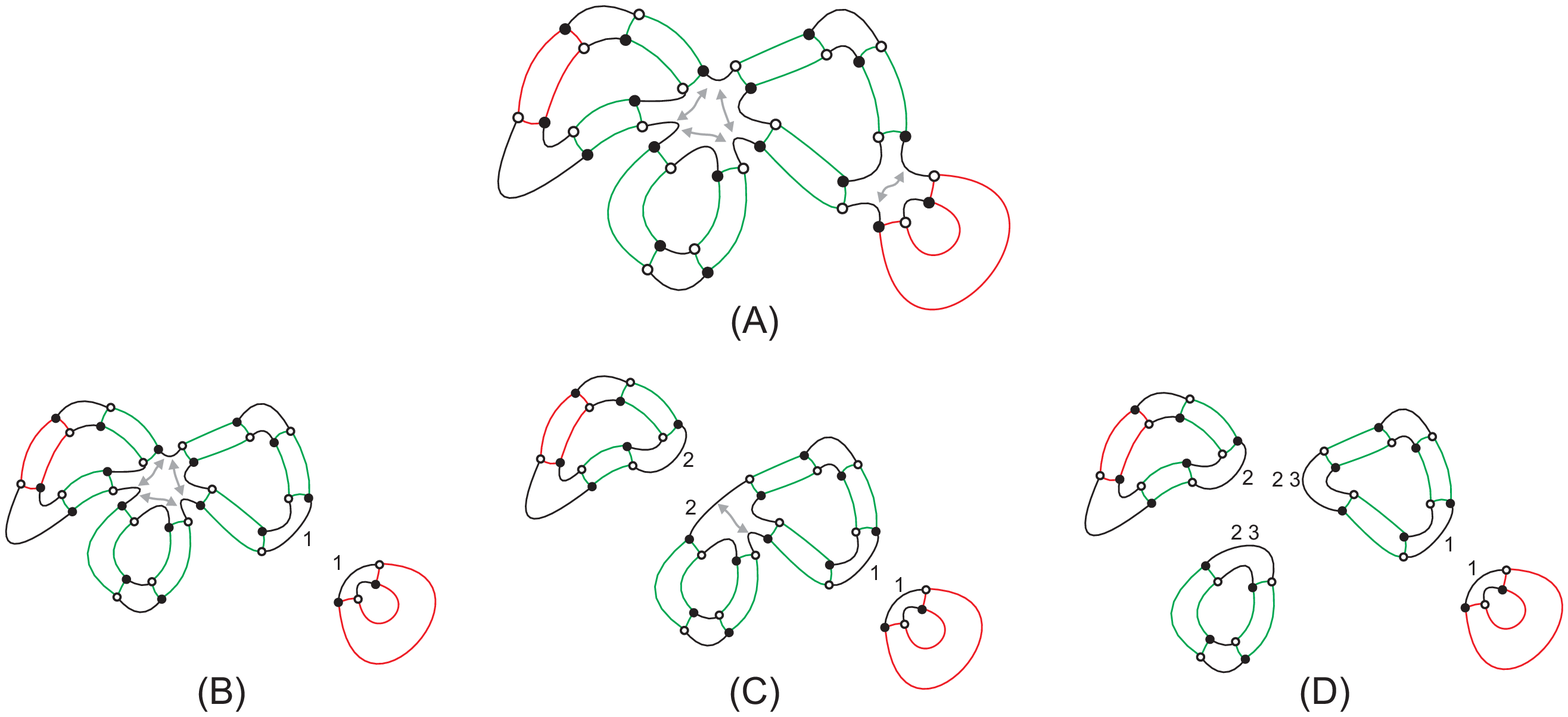}
   \end{center}
   \vspace{-0.7cm}
   \caption{Breakpoints of a g-blink and its breaking}
   \label{fig:breakpair}
\end{figure}

We saw on Section~\ref{sec:mergingGBlinks} that a breakpair on
a g-blink is a pair of angle-edges that are on
the same g-vertex and the same g-face. Given a g-blink and one
breakpair in it we may separate it into two g-blinks.
Figure~\ref{fig:breakpair}A shows a g-blink and its breakpairs: the
gray arrows point to the pair of angle-edges of the breakpair. We
can separate a g-blink in pieces (smaller g-blinks) until there are
no more breakpairs.
Figures~\ref{fig:breakpair}B, \ref{fig:breakpair}C and \ref{fig:breakpair}D
show this separation process. A piece without breakpairs is called a
{\it block}. Figure~\ref{fig:breakpair}D have 4 blocks. No matter
what sequence of breakpairs one uses to separate a g-blink in blocks,
the final blocks are always the same.

Some important properties of a breakpair that we need are related to
the g-zigzags of its g-blink. They are the subject of next two
propositions.

\begin{Prop}
If $p$ is a breakpair on g-blink $G$ and its angle-edges are $e_1$
and $e_2$ then the g-zigzag of $G$ that contains $e_1$ is the same
as the one that contains $e_2$.
\end{Prop}
\begin{proof} Straightforward. \end{proof}

\begin{Prop}
\label{prop:breakpairDisjointGZigzags}
 The only g-zigzag $z$ affected by separating a g-blink $G$ on
a breakpair $p$ is the one that contains both angle-edges of
$p$. If $P_1$ and $P_2$ are the pieces obtained by separating $G$ on
$p$ then the g-zigzags of $P_1$ and $P_2$, except for $z$, were
disjoint in $G$.
\end{Prop}
\begin{proof} Straightforward. \end{proof}

We also saw on Section~\ref{sec:mergingGBlinks} that any pair of angle
edges, each on different g-blinks may be the {\it
basepair} of a g-blink merging operation. Think of the
transition from Figure~\ref{fig:breakpair}B to
Figure~\ref{fig:breakpair}A the basepairs are the angle-edges labeled
1 on Figure~\ref{fig:breakpair}B. So, to merge two g-blinks on a
basepair is to replace the basepair angle-edges by two new edges
connecting the two g-blinks and respecting the parity.

\begin{figure}[htp]
   \begin{center}
      \leavevmode
      \includegraphics[width=14cm]{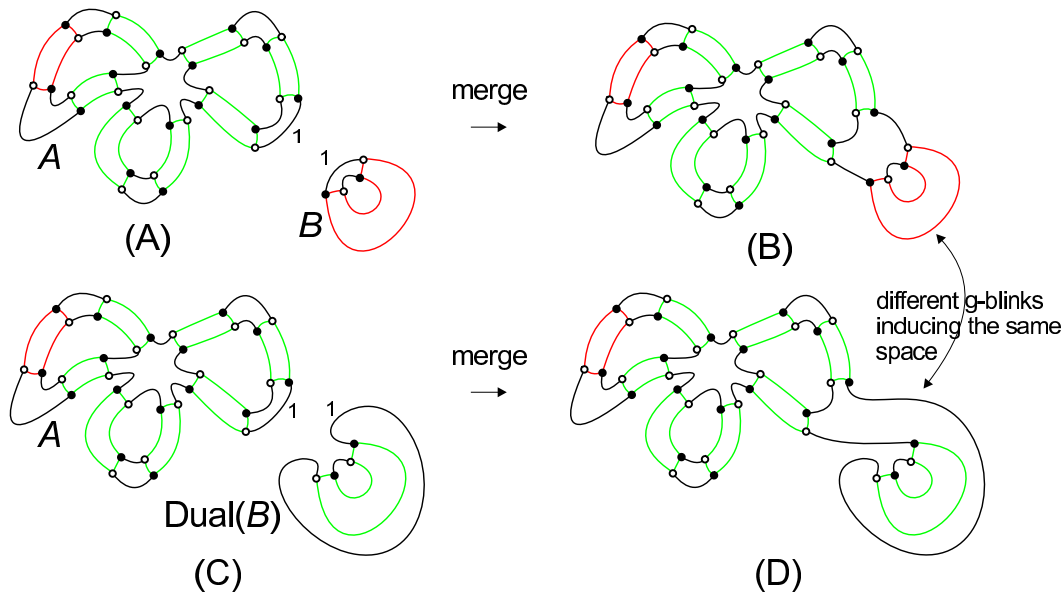}
   \end{center}
   \vspace{-0.7cm}
   \caption{Merging A with B and with \textsc{Dual}(B)}
   \label{fig:mergingWithDual}
\end{figure}

The fact that the two g-blinks at the right of Figure
\ref{fig:mergingWithDual} induce the same space is given by Theorem
\ref{theo:partialDual}. The proof of this Theorem still depends on
the proof of Theorem \ref{theo:partialReflection} which will be
given in Chapter 4. It depends on topological facts and a
reformulation of the BFL in terms of gems. As we have said before,
this is the main theoretical contribution of our Thesis.

Follow Figure~\ref{fig:mergingWithDual} to see an application of
Theorem \ref{theo:partialDual}. The basepair in both rows are the
same: the pair of edges labeled~1. The theorem together with the result
we describe now form the basis of the normalization procedure.

\begin{Prop} \label{prop:sameZigzagsSameSpace}
Let $A$ and $B$ be two g-blinks. Let $p$ be a basepair on them. Let $\alpha$
be the g-zigzag of the angle-edge of $p$ on $A$ and $\beta$ be the
g-zigzag of the angle-edge of $p$ on $B$. Let $p'$ be any other
basepair on g-zigzags $\alpha$ of $A$ and $\beta$ of $B$. The result
of $A$ and $B$ merged on $p$ induces the same space as $A$ and $B$
merged on $p'$.
\end{Prop}

\begin{figure}[htp]
   \begin{center}
      \leavevmode
      \includegraphics[width=12cm]{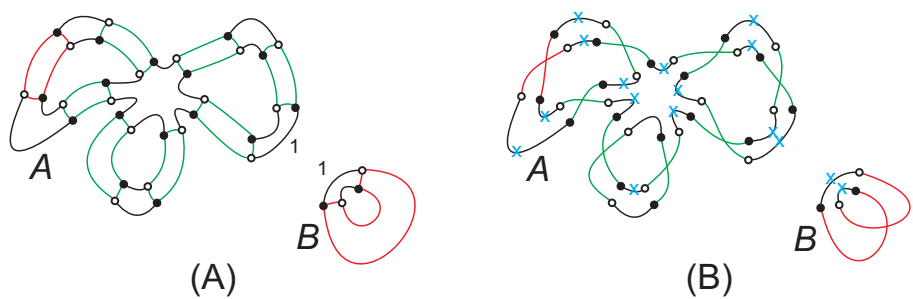}
   \end{center}
   \vspace{-0.7cm}
   \caption{Merging on any angle-edge of the same g-zigzags}
   \label{fig:mergingOnAnyAngleEdge}
\end{figure}

For example, the g-blink resultant of the merge of
Figure~\ref{fig:mergingOnAnyAngleEdge}A on basepair labeled 1
induces the same space as the g-blink resultant of merging any pair
of angle edges tagged with a ``blue X'' one from $A$ and another
from $B$ of Figure~\ref{fig:mergingOnAnyAngleEdge}B. Note that all
angle-edges tagged with a ``blue X'' on
Figure~\ref{fig:mergingOnAnyAngleEdge}B are on the same g-zigzag of
the angle-edges labeled 1 on
Figure~\ref{fig:mergingOnAnyAngleEdge}A.

Having Theorems \ref{theo:partialDual}, Theorem
\ref{theo:partialReflection} and Theorem \ref{theo:partialRefDual}
at our disposal, we are now able to describe the normalization
procedure. The intuitive idea is to separate the g-blink into blocks
and then remount the blocks (or their duals, reflection or refdual,
depending who is ``smaller'') in a canonical way. We divide this
procedure in three phases: separating phase, intermediate phase
and merging phase.

Separating phase. Let $G$ be the g-blink we want to normalize. First
give each g-zigzag of $G$ an unique label. Let $Z$ be the set of
these labels. For each angle-edge $e$ on $G$ record the label of the
g-zigzag that contains $e$ as its {\it zigzag label} $z_e$.
Initialize the ``pieces set'' as ${\cal P} \leftarrow \{G\}$.
Suppose there is a piece $P$ in ${\cal P}$ with a breakpair $p$. Let
$e_1$ and $e_2$ be the two angle edges of the breakpair $p$. Note
that the zigzag labels of $e_1$ and $e_2$ are the same: $z_{e_1} =
z_{e_2}$. Separate $P$ into $P_1$ and $P_2$ and make the two new
edges $e'_1$ on $P_1$ and $e'_2$ on $P_2$ have the same zigzag
labels as $e_1$ and $e_2$: $z_{e'_1} \leftarrow z_{e_1} (= z_{e_2})$
and $z_{e'_2} \leftarrow z_{e_1} (= z_{e_2})$. Replace $P$ with
$P_1$ and $P_2$ on ${\cal P}$. Repeat this until ${\cal P}$ contains
only blocks (g-blinks without breakpairs). The separating phase is
finished.

Intermediate phase. Define a bipartite graph $X$. The vertices of
$X$ are the labels $z \in Z$ and the pieces $P \in {\cal P}$; there
is an edge $(z,P)$ between label $z$ and piece $P$ in $X$ if there
is an angle-edge $e$ in $P$ with $z_e = z$. Note that $X$ is a tree:
no cycles. Let pieces $P_1$ and $P_2$ be neighbors of label $z$ on
$X$. The only common neighbor of $P_1$ and $P_2$ must be $z$
otherwise $P_1$ could not be separated from $P_2$ (See Proposition
\ref{prop:breakpairDisjointGZigzags}). This implies that there
cannot be a cycle in $X$. Remove every vertex $z \in Z$ of $X$ that
has only one neighbor. This asserts that every leaf of $X$ is a
vertex $P$ in ${\cal P}$. A consequence of this is that $X$ has a
single {\it center}. The center of a tree (see
\cite{BondyAndMurty1976}) is obtained by removing all leafs of a
tree in each step until arriving at a pair of vertices or a single
vertex. By applying the tree center algorithm on $X$ the leafs on
each step alternates between $z$ nodes and $P$ nodes. So it must
finish on a single node once there cannot be two adjacent $z$'s or
two adjacent $P$'s. Let $v$ be the center of $X$. We root $X$ at $v$
and $X$ becomes a {\it rooted tree}. The idea now is to organize
this rooted tree in a canonical way. To this aim we must have a way
to compare nodes and subtrees.

Remember from Section \ref{sec:gblinkcode} that every g-blink has a
unique code. So we can compare g-blinks by comparing their codes. We
know that merging a g-blink, or its dual, or its reflection or its refDual
on the same basepair results in a g-blink that induces the same space.
So we normalize ${\cal P}$ by replacing each piece $P$ in it by
$\min\{$ $P$, $\textsc{Dual}(P)$, $\textsc{Reflection}(P)$, $\textsc{RefDual}(P)$ $\}$.
Note that doing this does not affect the zigzag labeling of the angle-edges
because angle-edges are preserved on these operations (\ie dual,
reflection and refdual). Note also that the
$P$ nodes of $X$ are also updated by this criterion. Using the code
of a g-blink we can also organize the rooted tree $X$. To organize
$X$ we mean to define a fixed sequence for the children of a node of
$X$. We do this inductively. The base case is a node without child.
This node is already organized, so we are finished. Consider a node
$u$ with children $w_1,\ldots,w_k$ all of them already organized. To
organize $u$ we need to define a sequence for these children. Using
the code of a g-blink we can define a function
\textsc{CompareTrees}($r_1$,$r_2$) to compare organized rooted trees
that is evaluated to -1 if tree rooted at $r_1$ is smaller than the
tree rooted at $r_2$, to 0 if they are the same and to +1 if tree
rooted at $r_1$ is greater than tree rooted at $r_2$. So to organize
$u$ is a matter of sorting $w_1,\ldots,w_k$ using the
\textsc{CompareTrees} function. With these explanations the problem
of organizing $X$ is solved. This also ends the intermediate phase.
\newcommand{\CompareTreesAlgorithm}{
\parbox[t]{8cm}{
\begin{small}
\begin{algorithmic}[1]
 \Statex
 \Function {\textsc{CompareTrees}}{$r_1,r_2$}
    \State $s_1 \leftarrow$ \textsc{LinearizeTree}($r_1$)
    \State $s_2 \leftarrow$ \textsc{LinearizeTree}($r_2$)
    \State $n_1 \leftarrow$ length($s_1$); $n_2 \leftarrow$ length($s_2$)
    \State $i \leftarrow 1$
    \While {$i \leq {\rm min}(n_1,n_2)$}
       \State $(u_1,{\rm level}_1) \leftarrow s_1[i\,]$
       \State $(u_2,{\rm level}_2) \leftarrow s_2[i\,]$
       \State {\bf if} {$({\rm level}_1 < {\rm level}_2)$} {\bf then return} -1
       \State {\bf else if} {$({\rm level}_1 > {\rm level}_2)$} {\bf then return} +1
       \State {\bf if} {$\kappa(u_1) < \kappa(u_2)$} {\bf then return} -1
       \State {\bf else if} {$\kappa(u_1) > \kappa(u_2)$} {\bf then return} +1
       \State $i \leftarrow i + 1$
    \EndWhile
    \State {\bf if} {$n_1 = n_2$} {\bf then return} 0
    \State {\bf else if} {$n_1 < n_2$} {\bf then return} -1
    \State {\bf else} {\bf return} +1
 \EndFunction
\end{algorithmic}
\end{small}}}
\newcommand{\LinearizeTreeAlgorithm}{
\parbox[t]{7.5cm} {
\begin{small}
\begin{algorithmic}[1]
 \Statex
 \Function {\textsc{LinearizeTree}}{$r$}
    \State $s \leftarrow <>$
    \Procedure {\textsc{LT-DFS}}{$u$,{\rm level}}
       \State $s \leftarrow s \cdot <(u,{\rm level})>$
       \For {every children $v$ of $u$ taken in the ordered sequence}
           \State \textsc{LT-DFS}($v$,level+1)
       \EndFor
    \EndProcedure
    \State \textsc{LT-DFS}($r$,0)
    \State {\bf return} $s$
 \EndFunction
\end{algorithmic}
\end{small}}}

\begin{algorithm}
\caption{\textsc{CompareTrees} Algorithm}
\label{alg:compareTrees}
\begin{center}
\CompareTreesAlgorithm \hspace{0.1cm} \LinearizeTreeAlgorithm

\bigskip

\parbox{12cm}{
\begin{small}
In this algorithm the code $\kappa$ is taken over not only g-blinks
but also on zigzag labels. Consider the code of a zigzag label the
empty word. With this definition two zigzag labels $z_1$ and $z_2$
always satisfy $\kappa(z_1) = \kappa(z_2)$. Note also that a zigzag
node code is smaller than any g-blink code.\end{small}}
\end{center}
\end{algorithm}
Merging phase. The rooted tree $X$ is already organized. The idea
now is to merge the blocks using the order defined on $X$ and using
the code to define a canonical basepair for each merging operation.
Let $z$ be a label in $X$. Let $P_1, \ldots, P_k$ be the neighbors
of $z$. If $z$ is not the root of $X$ then $P_1$ is the parent of
$z$ and $P_2 \ldots P_k$ are the children of $z$ taken in order. If
$z$ is the root of $X$ then $P_1 \ldots P_k$ are the children of $z$
taken in order. We want now to merge the pieces $P_1 \ldots P_k$ in
the order they appear. First $P_1$ with $P_2$, second the result of
$P_1$ and $P_2$ with $P_3$ and so on. The only thing not defined yet
is the base point of each merging. This is solved using the code of
the blocks. For each label $z$ that appears on a block $P$ we define
a {\it canonical basepair angle-edge} $e_P^z$ on the g-zigzag whose
angle-edges are all zigzag labeled with $z$. This is done using the
code labeling of $P$ (see Section \ref{sec:gblinkcode}). Label the
vertices of $P$ with its code labeling. Define $e_P^z$ as the
angle-edge (among all angle-edges with zigzag label $z$ on $P$)
incident to the vertex of $P$ that has the smallest label.
The last thing we need to define is how to update the canonical
basepair angle-edge after a merging operation. In symbols, after
merging $P_1$ and $P_2$ whose canonical basepair angle-edges
were $e_{P_1}^z$ and $e_{P_2}^z$ what will be the canonical basepair
angle-edge $e_{P_1+P_2}^z$ of $P_1 + P_2$? By definition
$e_{P_1+P_2}^z$ will be the new angle-edge incident to the odd
vertex~of~$P_2$. Repeat this merging for all zigzag label $z$
in any order until no more merging may be
done. Finished.

\begin{figure}[htp]
   \begin{center}
      \leavevmode
      \includegraphics[width=12cm]{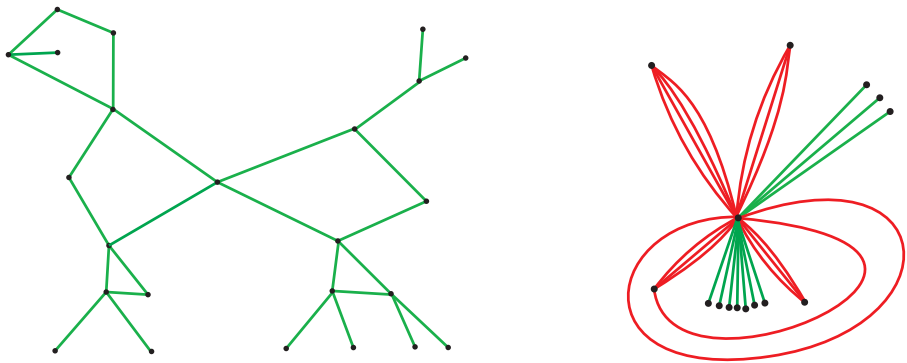}
   \end{center}
   \vspace{-0.7cm}
   \caption{Representative of the dog like blink}
   \label{fig:dog}
\end{figure}

This three-phase procedure results in a unique g-blink $r(G)$, {\it
the representative of $G$}. Both g-blinks, $G$ and $r(G)$ induce the
same space. A g-blink is said to be {\it a representative} if $G =
r(G)$. We finish this section with an example in
Figure~\ref{fig:dog} of the result of the algorithm that we have
implemented for obtaining the representative of a g-blink (or blink).

\section{Towards a census of spaces induced by small blinks}
\label{sec:towardsACensusOfPrimeSpacesInducedBySmallBlinks}

\begin{center}
\it What spaces have a small \underline{blink} presentation?
\end{center}
In this chapter we saw that a g-blink is a family of blinks that
induce the same space and that any blink has an associated g-blink.
Thus, our question is equivalent to
\begin{center}
\it What spaces have a small \underline{g-blink} presentation?
\end{center}
In this form, our original question becomes easier ``to be computed''
once a g-blink is a combinatorial object and the number of g-blinks
with $\,\,\leq k \,\,$ g-edges is finite. Moreover, by being
combinatorial simple objects, \hbox{g-blinks} have a direct way of
going into computers. For instance, a g-blink may be represented
in a computer by its code.

By the fact that every g-blink is associated to a special g-blink
with the same size that induces the same space called its
representative, our question is also equivalent to
\begin{center}
\it What spaces have a small \underline{representative g-blink}
presentation?
\end{center}
In this form, our original question becomes even ``easier'' in the sense
that, with $\,\,\leq k \,\,$ g-edges, there are fewer ``representative
g-blinks'' than ``general g-blinks''. Thus, this form is the one we
use.

Suppose we have generated the set of all representative g-blinks with
$\,\,\leq k \,\,$ g-edges for some $k \geq 1$. We already know that
all spaces that have a presentation with $\,\,\leq k \,\,$ edges are
there. But how to identify them? What elements of this set
(representative g-blinks) induce the same space? In
Sections~\ref{sec:homologyGroup}~and~\ref{sec:quantumInvariant} we
described a way to calculate two space invariants from a g-blink
presentation: the homology group and the quantum invariant.
Calculating these invariants on all g-blinks we can partition this
set in classes where each element of the same class has the same
homology group and same quantum invariant. At this point we are sure
that different classes induce different spaces because there is a
topological invariant that distinguishes them. The remaining problem
is to know if all g-blinks in the same class (same homology group
and same quantum invariant) induce the same space. To prove that two
\hbox{g-blinks} indeed induce the same space we will use a
computational method in 3-Gem Theory that was described in
\cite{Lins1995}. Next chapter will be about 3-Gems and this
computational method.

%% file: chapter4.tex
\chapter{3-Gems}
\label{chap:gems}

\section{Definition}

An {\em $(n+1)$-graph} is a regular graph where all its vertices
have degree $n+1$ and the edges incident to each vertex have
distinct colors $0,1,\ldots,n$. Let $K \subseteq \{0,1,\ldots,n\}$
be a subset of colors and $G$ an $(n+1)$-graph. Define $G[K]$ as the
subgraph of $G$ induced by $K$. We say that each connected component
of $G[K]$ is a {\em $K$-residue} of $G$ (note that $K$ here is a
subset of colors). If $k = |K|$ then a $K$-residue is also said to
be a {\em $k$-residue} of $G$ (note that $k$ here is a number). If
$K$ is a set, we denote by $\overline{K}$ its complement
($\{0,\ldots,n\} \backslash K$). A 2-residue is also called a {\em
bigon} and a 3-residue is also called a {\em triball}. A {\em 3-gem}
(acronym for 3-dimensional graph encoded manifold) is a
$(3+1)$-graph where each of its 3-residues induces the surface of a
sphere, $\IS^2$. Each bipartite gem $G$ corresponds to a unique
space $|G|$. The easiest way to define $|G|$ is to start with
$v_G$ tetrahedra each with its 4 vertices painted each with one color of
$\{0,1,2,3\}=\{h,i,j,k\}$ and glue a pair of tetrahedra $t_u$ and
$t_v$  by identifying its faces opposite to the $i$-colored vertices
so as to match colors $i,j,k$ whenever there is an
$h$-colored edge in $G$ between $u$ and $v$. In this way $G$ is the
dual of the pseudo-triangulation of the pseudo-manifold obtained by
the gluing. In the case of a gem, the pseudo-manifold is a manifold.
Given a $4$-regular properly edge colored graph $G$ denote
$\alpha(G)=b_G-v_G-t_G$ the {\em agemality} of $G$,
where $b_G$ is the number of 2-residues of $G$,
$v_G$ is the number of vertices of $G$ and
$t_G$ is the number of 3-residues of $G$.
The agemality is non-negative and it is $0$ if and only if $G$ is a gem. Indeed we
have (\cite{Lins1995}):

\begin{Prop}
Let $G$ be a $(3+1)$-graph with $b_G$ 2-residues, $t_G$ 3-residues
and $v_G$ vertices, then $G$ is a 3-gem if and only if its agemality
is zero, that is, $$v_G + t_G = b_G.$$
\end{Prop}

\section{Moves on gems}
Let $G$ be a 3-gem. An $i$-colored edge $\alpha$ of $G$ is a {\em
1-dipole} if the vertices incident to $\alpha$ are in different
$\overline{\{ i \}}$-residues. A pair of edges of $G$ one with color
$i$ and the other with color $j$ and with equal ends  is a {\em
2-dipole} if these ends are in different $\overline{\{ i,j
\}}$-residues. The creation and cancelation of a
$k$-dipole ($k=1,2$) does not change the induced space. A 3-gem free
of 1-dipoles is said to be a {\em 3-crystallization}.

A {\em $\rho$-pair} in a $(3+1)$-graph is a pair of edges of the
same color that are incident to 2 or 3 common bigons (the two edges
are both contained in 2 or 3 bigons of $G$). If the edges of the
pair are incident to only two common bigons then the pair is said to
be a {\em $\rho_2$-pair}. If the edges of the pair are incident to
three common bigons then the pair is said to be a {\em
$\rho_3$-pair}. If a $\rho$-pair is found in a gem we can get a
smaller gem inducing the same space.

\section{Simplifying dynamics}
\begin{figure}[htp]
   \begin{center}
      \leavevmode
      \includegraphics{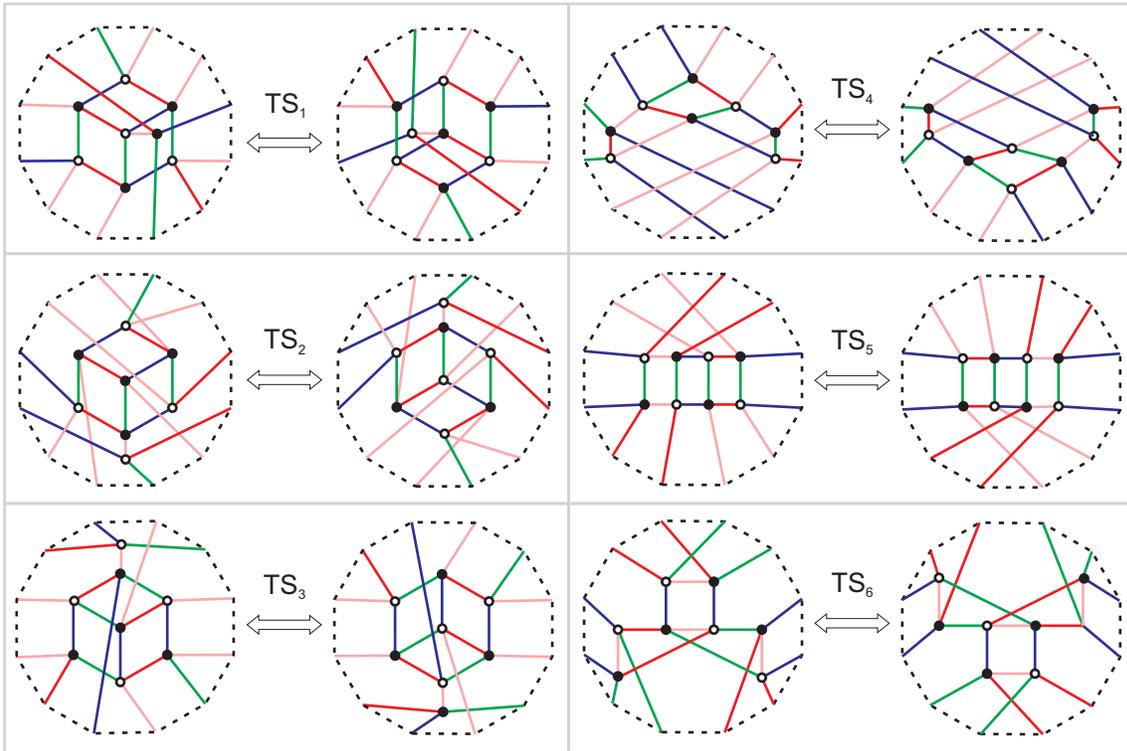}\\
   \end{center}
   \vspace{-0.7cm}
   \caption{The six TS-moves}
   \label{fig:tsmoves}
   \vspace{-0.5cm}
\end{figure}

In this section we briefly review the simplifying dynamics on gems.
This technique is developed in~\cite{Lins1995} and it uses the so
called $TS$-moves and  $U$-move which maintain the induced
3-manifold. The relevant algorithm to simplify gems and get to an
attractor for the spaces induced by a gem is named the $TS_\rho
U$-algorithm (\cite{Lins1995}). We have re-implemented this
algorithm which is the basis for the proof that the blinks with the
same homology and the same quantum invariants up to $r=12$ indeed
induce the same spaces. The six TS-moves on gems are defined in
Figure~\ref{fig:tsmoves}.

\begin{figure}[htp]
   \begin{center}
      \leavevmode
      \includegraphics[width=10cm]{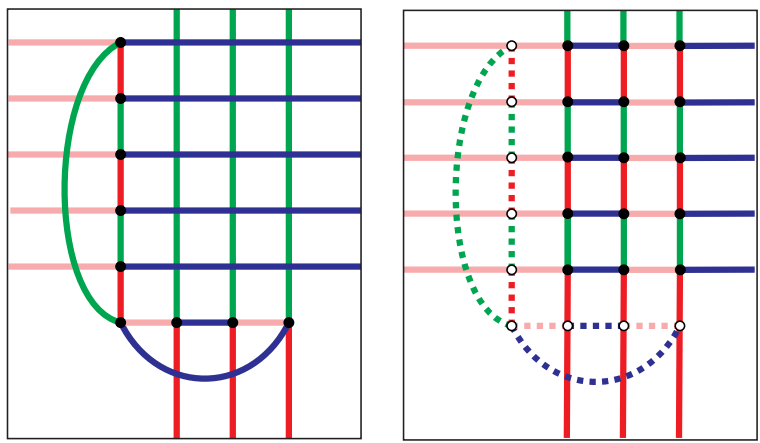}\\
   \end{center}
   \vspace{-0.7cm}

  \caption{$U_{2,3}$-move applied to a
  1-monopole of type (2,3)}
  \label{fig:umove}
\end{figure}

A {\em monopole}\index{monopole} in a $(3+1)$-graph is a vextex
which is the only intersection of an $hi$-gon and a $jk$-gon,
$(h,i,j,k)$ a permutation of $(0,1,2,3)$. This defines a
configuration which induces a fundamental move in the classification
of gems. A $U_{mn}$-move is defined on a monopole, by making the
$hi$-gon of size $2m$ and the $jk$-gon of size $2n$ (whose union has
$2m+2n-1$ vertices) disappear, being replaced by a cluster of
squares with $(2m-1) \times (2n-1)$ vertices. A $U_{mn}$-move does
not change the induced space of the gem. We give an example in
Figure~\ref{fig:umove} of $U_{23}$ move. In general the
$U_{mn}$-move increases the number of vertices of a gem. However, in
conjuntion with the TS-moves and $\rho$-pairs the $U_{mn}$-moves
have been so far sufficient to classify gems up to 30 vertices.

\enlargethispage{2cm}

\newpage

\section{From g-blink to 3-gem}
\label{sec:fromGBlinkTo3Gem}

Assume that $G$ is a g-blink with no g-zigzags with g-edges
alternating red and green. This kind of g-zigzag corresponds in the
BFL to a component that goes totally over or totally under and can
be separated from the rest of the BFL by Reidemeister moves $II$ and
$III$. Assume the following convention on the colors of a gem: $0
\equiv$ pink, $1 \equiv$ blue, $2 \equiv$ red, $3 \equiv$ green. We
begin by proving a result which simplifies considerably the passage
``blink $\rightarrow$ gem'' first given in
\cite{KauffmanAndLins1994}.

\begin{figure}[htp]
   \begin{center}
      \leavevmode
      \includegraphics[width=10cm]{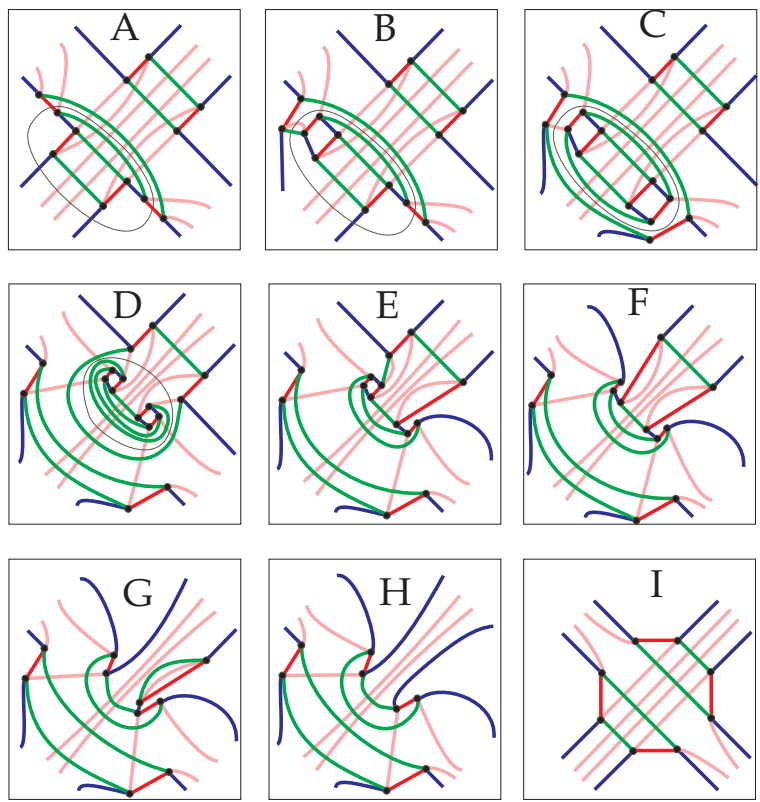}\\
   \end{center}
   \vspace{-0.7cm}
   \caption{Simplifying the gem of a blink: from 12 to 8 vertices by crossing}
   \label{fig:The12x8Move}
\end{figure}

\begin{Theo} \label{theo:blinkToGem}
Given a g-blink $B$ with no alternating g-zigzags it is possible to
obtain a gem $J^\downarrow$ where each edge of the blink (which
corresponds to a crossing of the associated BFL) becomes the
sub-configuration of $8$ vertices shown in Figure~\ref{fig:The12x8Move}I
so that $B$ and $G(B)$ induce the same
space.
\end{Theo}
\begin{proof}
It is proved in \cite{KauffmanAndLins1994} that replacing each
crossing of the BFL associated to the blink by the configuration of
Figure~\ref{fig:The12x8Move}A the final gem $J^\downarrow$ will have
the desired property. The rest of the proof consists in effecting
dipole moves in $J'$ so as to arrive at $J$. The sequence of dipole
moves are depicted in Figure \ref{fig:The12x8Move}. The dipole
moves are local and should be made in the neighborhood of each
original crossing of the BFL. The resulting gem is $J^\downarrow$.
\end{proof}

The gem obtained from a blink by replacing each color of the BFL by
the configuration of Figure \ref{fig:The12x8Move}I is called the
{\em reduced canonical gem} of the blink.

We introduce the following notation to represent both a crossing and
its switched form. The {\em octagon} of a crossing corresponds to an
unidentified crossing. This is indicated by light green edges in the
place of the normal green ones.
\begin{center}
      \includegraphics[width=10cm]{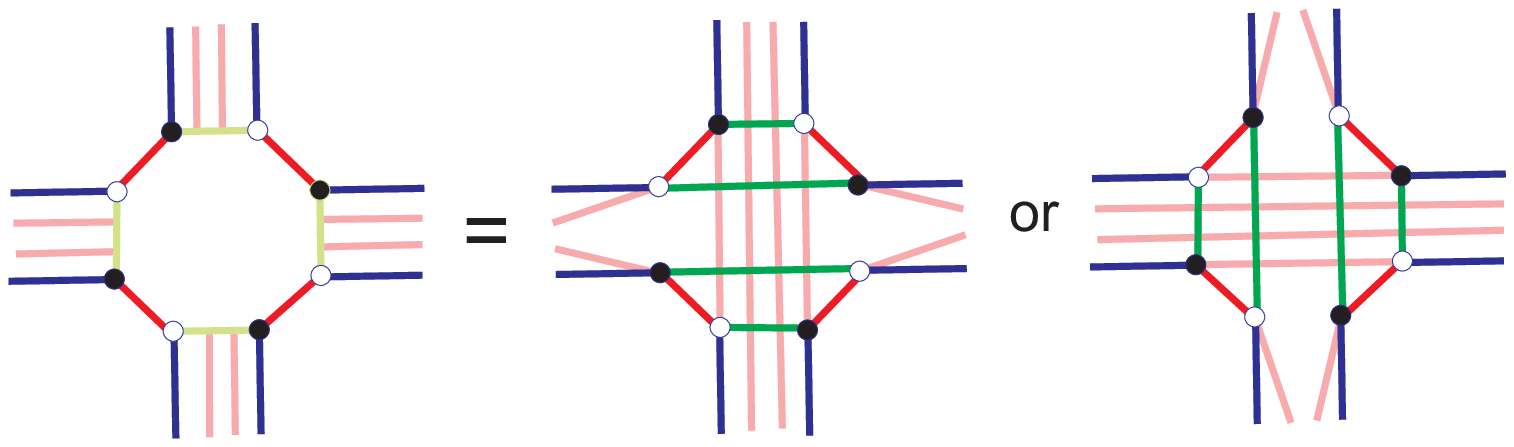}
\end{center}

\begin{figure}[htp]
   \begin{center}
      \leavevmode
      \includegraphics[width=15cm]{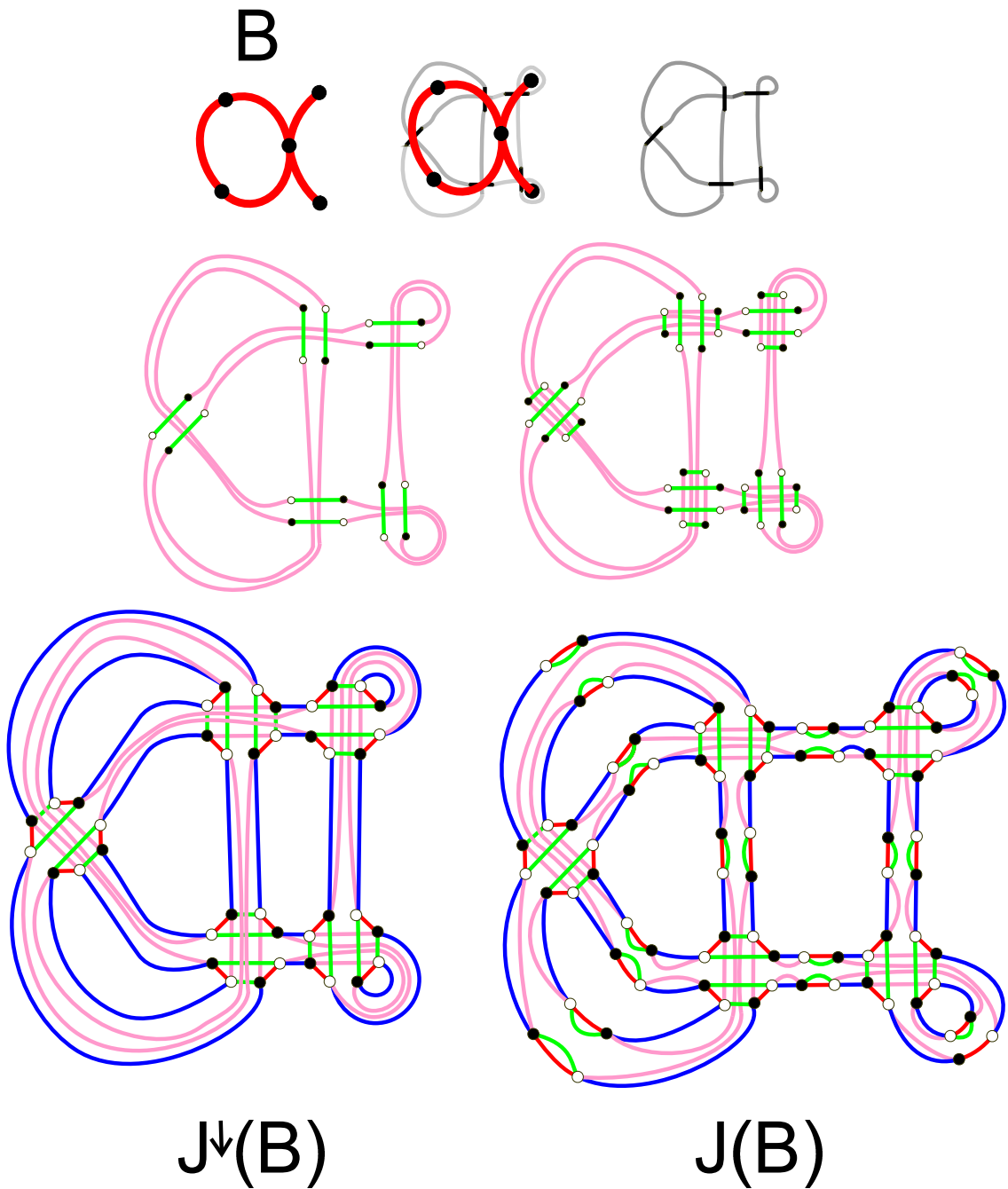}
   \end{center}
   \vspace{-0.7cm}
   \caption{ Obtaining the canonical gems $J^\downarrow(B)$ and
   $J(B)$ from a blink $B$}
   \label{fig:poincareGemBlink}
\end{figure}

In Figure \ref{fig:poincareGemBlink} we display a complete example
of the above algorithm to go from a blink $B$ to its canonical gem
$J(B)$ inducing the same space. This example corresponds to
Poincaré's homology sphere. Observe that an immersion of the gem in
the plane is directly obtained from the embedding of the BFL. The
gem obtained is bipartite. In going clockwise along the
(blue,red)-gons (which corresponds to the faces of the BFL) the red
edges go from a black to a white vertex. Observe that the pink-green
gons form a neighborhood of the original blink. The convention here
is that the green edges are the overpasses, while the pink edges the
underpasses.

The difference between the {\em canonical reduced} gem of the blink
$J^\downarrow(B)$ and the {\em canonical} gem of the blink $J(B)$ is
that we introduce in the latter two 2-dipoles (red-green digons) at
each site corresponding to a g-edge of the original g-blink. While
redundant these $4|E(B)|$ vertices are convenient for our purposes
as we show next. In the computer implementation we use only
$J^\downarrow(B)$. The construction of
Figure~\ref{fig:poincareGemBlink} emphasizes the geometric
simplicity of the algorithm.

The main point in using the auxiliary red-green digons in defining
$J=J(B)$ is that they induce, for each directed g-edge $e$ between
crossings $\alpha$ and $\gamma$, three cylinders $\C_e, \C_\alpha,
\C_\gamma$. For a color $i$ and a vertex $a$ of a gem, denote by
$a_i$ the $i$-colored edge incident to $a$. Let $a_i^\star$ denote
the $2$-simplex in the dual pseudo-complex $J^\star$ corresponding
to the edge $a_i$ of the gem $J$.

\begin{figure}[htp]
   \begin{center}
      \leavevmode
      \includegraphics[width=12cm]{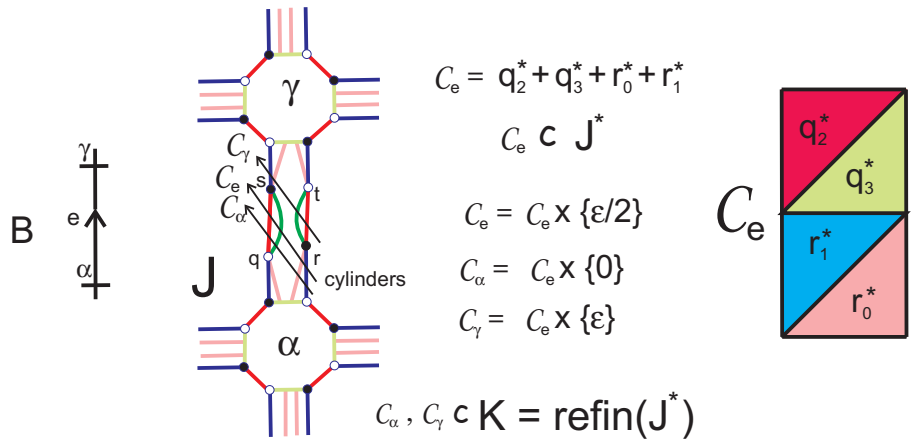}
   \end{center}
   \vspace{-0.7cm}
  \caption{Some cylinders induced in $K(B)$ by a directed
  edge of the BFL between crossings $\alpha$ and $\gamma$}
  \label{fig:cylinderInEdge}
\end{figure}

\begin{Lem} Let $q,r,s,t$ be the ends of the two
red-green digons induced in $J$ by the directed edge $e$ of $B$, as
shown in Figure~\ref{fig:cylinderInEdge}. Then the sub-complex
$\C_e=q_2^\star + q_3^\star + r_1^\star + r_0^\star$ is a
non-singular cylinder in $J^\star$.
\end{Lem}
\begin{proof}
Two $2$-simplexes of $J^\star$ in colors $i$ and $j$ have a common
$1$-simplex if and only if the dual edges are in the same
$(i,j)$-gon. Note that $q_2$ and $q_3$ are in the same $(2,3)$-gon,
$q_3$ and $r_1$ are in the same $(3,1)$-gon, $r_1$ and $r_0$ are in
the same $(1,0)$-gon and $r_0$ and $q_2$ are in the same
$(0,2)$-gon. To complete the proof just note that there are $4$
distinct vertices in the subcomplex $\C_e$, that $q_2$ and $r_1$ are
not in the same $(2,1)$-gon and finally, that $q_3$ and $r_0$ are
not in the same $(3,0)$-gon.
\end{proof}

The cylinder $\C_e$ is contained in the dual pseudo-complex
$J^\star(B)$. Take a neighborhood $\C_e \times [0,\epsilon]$ in
$|J|$ identify $\C_e \times \{\epsilon/2\}\equiv \C_e$ and define
$\C_\alpha = \C_e \times \{0\}$ and $\C_\gamma = \C_e \times
\{\epsilon\}$. Let $K$ be a simplicial complex which is a refinement
of $J^\star$ containing both $\C_\alpha$ and $\C_\gamma$ as
sub-complexes. We observe that $|J|=|J^\star|=|K|$ and that
vertices $r$ and $s$ of gem $J$ are in $\C_e \times [0,\epsilon]$.

\begin{figure}[htp]
   \begin{center}
      \leavevmode
      \includegraphics[width=9cm]{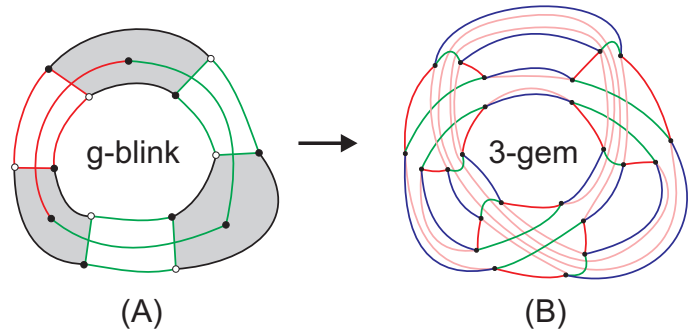}
   \end{center}
   \vspace{-0.7cm}
   \caption{ g-blink $G$ and its reduced canonical
   3-gem $J^\downarrow(G)$ defined by \textsc{GBlink2Gem}}
   \label{fig:gBlinkAndGem}
\end{figure}

The procedure \textsc{GBlink2Gem} that follows apply to g-blinks
without alternating zigzags. The procedure is purely combinatorial
and it teaches the computer to go from a g-blink $G$ to the 3-gem
$J^\downarrow$ given in Theorem \ref{theo:blinkToGem}. For each
vertex $v$ of $G$ we define two vertices $v_i$ and $v_o$ in $j$.
\begin{figure}[htp]
   \begin{center}
      \leavevmode
      \includegraphics{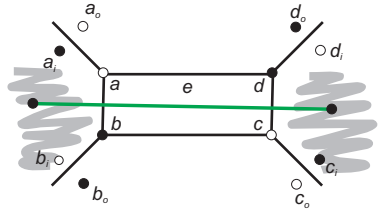}
   \end{center}
   \vspace{-0.7cm}
   \caption{ Scheme to define a 3-gem from a g-blink}
   \label{fig:gBlink2gem}
\end{figure}
Let $e$ be a g-edge on $G$ with vertices $a, b, c, d$ as shown in
Figure~\ref{fig:gBlink2gem}. Note that $(a,b)$ and $(c,d)$ are
face-edges and $(a,d)$ and $(b,c)$ are vertex-edges. The vertices
$a_i$, $a_o$, $b_i$, $b_o$, $c_i$, $c_o$, $d_i$, $d_o$ of
$J^\downarrow$ corresponding to $a$, $b$, $c$, $d$ are also shown on
Figure~\ref{fig:gBlink2gem}. The $i$ (in) index indicates that the
vertex is drawn inside a g-vertex and the $o$ (out) index indicates
that the vertex must be drawn inside a g-face (or outside the
g-vertex). The edges of $J^\downarrow$ are defined according to the
following procedure:
\begin{enumerate}
\item Each g-edge $e$ of $G$ aligned like the scheme of
Figure~\ref{fig:gBlink2gem} induce the following
edges~on~$J^\downarrow:$
\begin{center}
\begin{small}
\begin{tabular}{l|p{2.7cm}|p{2.7cm}|p{2.7cm}|p{2.7cm}} \hline
& color 0 & color 1 & color 2 & color 3 \\[0.1cm] \hline
   $e$ is green &
   $(a_i,c_o)$, $\,\,(a_o,c_i)$ &
   &
   $(a_i,b_i)$, $\,(c_o,b_o)$, $\,(c_i,d_i)$, $\,(d_o,a_o)$ &
   $(a_i,a_o)$, $\,(b_i,d_o)$, $\,(b_o,d_i)$, $\,(c_i,c_o)$
   \\[0.1cm] \hline
   $e$ is red &
   $(b_i,d_o)$, $\,(b_o,d_i)$ &
   &
   $(a_i,b_i)$, $\,(c_o,b_o)$, $\,(c_i,d_i)$, $\,(d_o,a_o)$ &
   $(a_i,c_o)$, $\,(b_i,b_o)$, $\,(a_o,c_i)$, $\,(d_i,d_o)$
   \\[0.1cm] \hline
\end{tabular}
\end{small}
\end{center}
\item For every angle-edge $\hat{e} = (u,v)$ in $G$ edges
$(u_i,v_i)$ and $(u_o,v_o)$, both with color 1, are added to
$J^\downarrow$.
\newpage

\item At this point, some vertices in $J^\downarrow$
do not have a color 0 incident edge. Let $u$ be a vertex in
$J^\downarrow$ without a neighbor of color 0. We add the edge
$(u,v)$ with color 0 in $J^\downarrow$, where $v$ is the result of
$$
\begin{array}[t]{l}
x \leftarrow {\rm neighbor}(v,1) \\
c \leftarrow 0 \\
\hbox{\bf while } {\rm neighbor}(x,c) \hbox{ is defined} \\
\hspace{1cm} x \leftarrow {\rm neighbor}(x,c) \\
\hspace{1cm} c \leftarrow (c + 1) \,\, {\rm mod} \,\, 2\\
v \leftarrow x \\
\end{array}.
$$
The expression ${\rm neighbor}(x,c)$ denotes the vertex adjacent to
$x$ by color $c$ in $J^\downarrow$. We do this until every vertex
has an incident color 0 edge.
\end{enumerate}
According to Theorem~\ref{theo:blinkToGem}, $J^\downarrow$ defined
this way is a 3-gem and it induces the same space as $G$ does. We
denote this procedure described here as \textsc{GBlink2Gem}. A
complete example of a g-blink and the 3-gem defined by
\textsc{GBlink2Gem} is depicted on Figure~\ref{fig:gBlinkAndGem}.

\newpage
\section{A proof of the partial reflection theorem}
\label{sec:proofOfThePartialReflectionTheorem}
We first show that a breakpair $\{e,f\}$ in a g-blink $C$
corresponds in $J^\star=J^\star(C)$ to a separating non-singular
$2$-torus $T^2_{ef}$.

\begin{figure}[htp]
   \begin{center}
      \leavevmode
      \includegraphics[width=14cm]{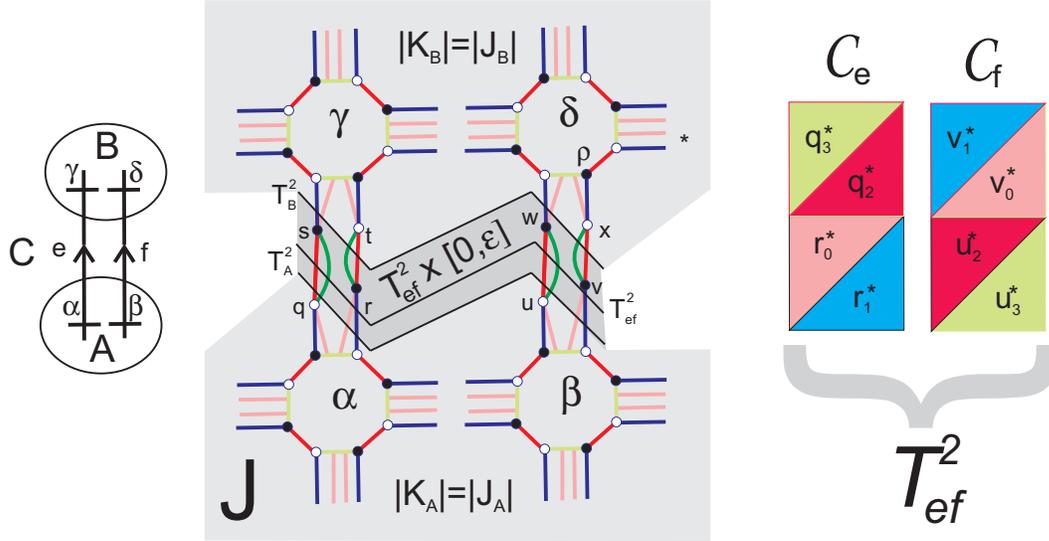}\\
   \end{center}
   \vspace{-0.7cm}
  \caption{$1-1$ correspondence: breakpair $\{e,f\}$
  in g-blink $C$ $\leftrightarrow$
  separating 2-torus $T_{ef}^2$ in $J^\star=J^\star(C)$ }
  \label{fig:breakpairAndTorus}
\end{figure}

\begin{Lem} \label{lem:torus}
The subcomplex $T^2_{ef} = \C_e + \C_f$ is a non-singular separating
torus in the dual pseudo-complex $J^\star$.
\end{Lem}

\begin{proof} We have seen already that $\C_e$ and $\C_f$ are
cylinders in $J^\star$. It remains to show that these cylinders have
the same boundary and that $\C_e + \C_f$ is a torus. We refer to
Figure~\ref{fig:breakpairAndTorus}. Note that $q_2$ and $v_1$ are in
the same $(2,1)$-gon, $r_1$ and $u_2$ are in the same $(1,2)$-gon,
$q_3$ and $v_0$ are in the same $(3,0)$-gon and that $r_0$ and $u_3$
are in the same $(0,3)$-gon. So, $T^2_{ef}=\C_e + \C_f$ is a
non-singular torus. It clearly separates. This completes the proof.
\end{proof}

To simplify the notation henceforth we write $T^2$ in place of
$T^2_{ef}$. Consider an $\epsilon$-neighborhood $T^2 \times
[0,\epsilon]$ of $T^2 \subset K$ so that $T^2 \equiv T^2 \times
\{\epsilon/2\}$. If we now remove $T^2 \equiv T^2 \times
\{0,\epsilon\}$ from $|K|$ then we have two disjoint spaces $|K_A|$
with boundary $\C_\alpha + \C_\beta = T_{\alpha\beta}^2 \equiv T^2
\times \{0\}$ and  $|K_B|$ with boundary $\C_\gamma + \C_\delta =
T_{\gamma\delta}^2 \equiv T^2 \times \{\epsilon\}$ as shown in
Figure~\ref{fig:breakpairAndTorus}. It follows that
\begin{equation}\label{eq:tripartition1} |K| = |K_A| \cup \left(T^2
\times [0,\epsilon]\right)\ \cup\ |K_B|,\end{equation} with
\begin{equation}\label{eq:tripartition2} |K_A| \cap \left(T^2 \times
[0,\epsilon]\right) = T_{\alpha\beta}^2 , \hspace {8mm} \left(T^2
\times [0,\epsilon] \cap |K_B| \right) = T^2_{\gamma\delta}, \hspace
{8mm} |K_B| \cap |K_A| = \emptyset.\end{equation}

For the proof of the Partial Reflection Theorem we present the
$2$-torus as the quotient space of $\IR^2$ by the lattice of integer
points: $T^2={\frac{\IR\times \IR}{\IZ\times \IZ}}$. Seeing $T^2$ in
this way, the $\pi$-rotational symmetry that we will need becomes
simply $(x,y) \mapsto (-x,-y)$. Let $F = \IR\times \IR \times
[0,\pi]$. Consider the auto-homeomorphism $\mu$ of $F$ given by
$$\mu(x,y,\theta)=\left(x\, \cos\, \theta
+y\, \sin\, \theta, -x\, \sin\, \theta+y\, \cos\, \theta,\,
\theta\right).$$

Define $\equiv'$ as the equivalence relation on $F$: $(x,y,\theta)
\equiv' (x',y',\theta')$ if $\theta'=\theta$, $x-x' \in \IZ$ and
$y-y' \in \IZ$. The quotient space $F / \equiv'$ is denoted by $F'$.
The image under $\mu$ of a the vertical segment linking $(x,y,0)$ to
$(x,y,\pi)$ is a helicoidal curve that starts at $(x,y,0)$, and
finishes at $(-x,-y,\epsilon)$. Note that any two vertical segments
in $F$ whose distance is an integer are identified. Clearly $F'
\approx \IS^1 \times \IS^1 \times [0,\epsilon]$.

\begin{figure}[htp]
   \begin{center}
      \leavevmode
      \includegraphics[width=15cm]{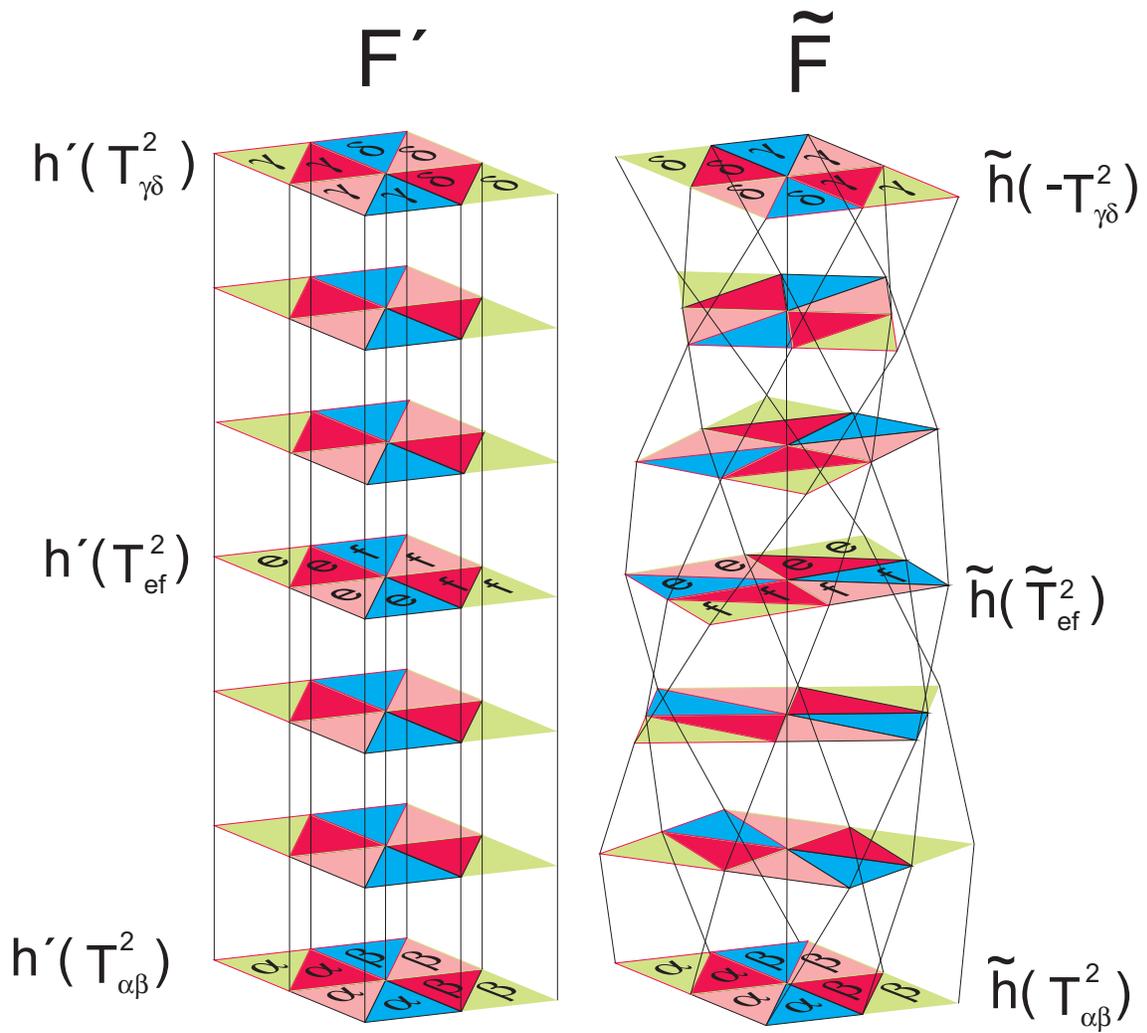}\\
   \end{center}
   \vspace{-0.7cm}
  \caption{The action of $\mu$ on the fundamental domain
  centered at the origin mapping $F'$ onto $\widetilde{F}$}
  \label{fig:filme}
\end{figure}

Define $\equiv_\mu$ as another equivalence relation on $F$ given by
$$(x,y,\theta) \equiv_\mu (x',y',\theta') {\hspace{5mm}\rm
if \hspace{5mm} } \mu^{-1}(x,y,\theta) \equiv'
\mu^{-1}(x',y',\theta'), \rm {\ that\ is,\ if \ } \theta=\theta'
\rm{\ and} $$  $$ x\,\cos\, \theta -y\, \sin\, \theta - x'\, \cos\,
\theta'+ y'\, \sin\, \theta' \in \IZ$$
$$x\, \sin\, \theta + y\,
\cos\, \theta - x'\, \sin\, \theta' - y'\, \cos\, \theta' \in \IZ.$$
The space $F/\equiv_\mu$ is denoted by $\widetilde{F}$. Observe the
simple fact that $\mu$ induces a homeomorphism sending $F'$ onto
$\widetilde{F}$, also named $\mu$ by abuse of language. As a
consequence of our definitions, two helicoidal curves in $F$ are
identified in $\widetilde{F}$ if their pre-images under $\mu$ are
two vertical segments identified in $F'$. The action of $\mu$ in the
fundamental domain centered at the origin from $F'$ to
$\widetilde{F}$ is shown in Figure~\ref{fig:filme}.

We are now ready to prove Theorem~\ref{theo:partialReflection}

\begin{figure}[htp]
   \begin{center}
      \leavevmode
      \includegraphics[width=15cm]{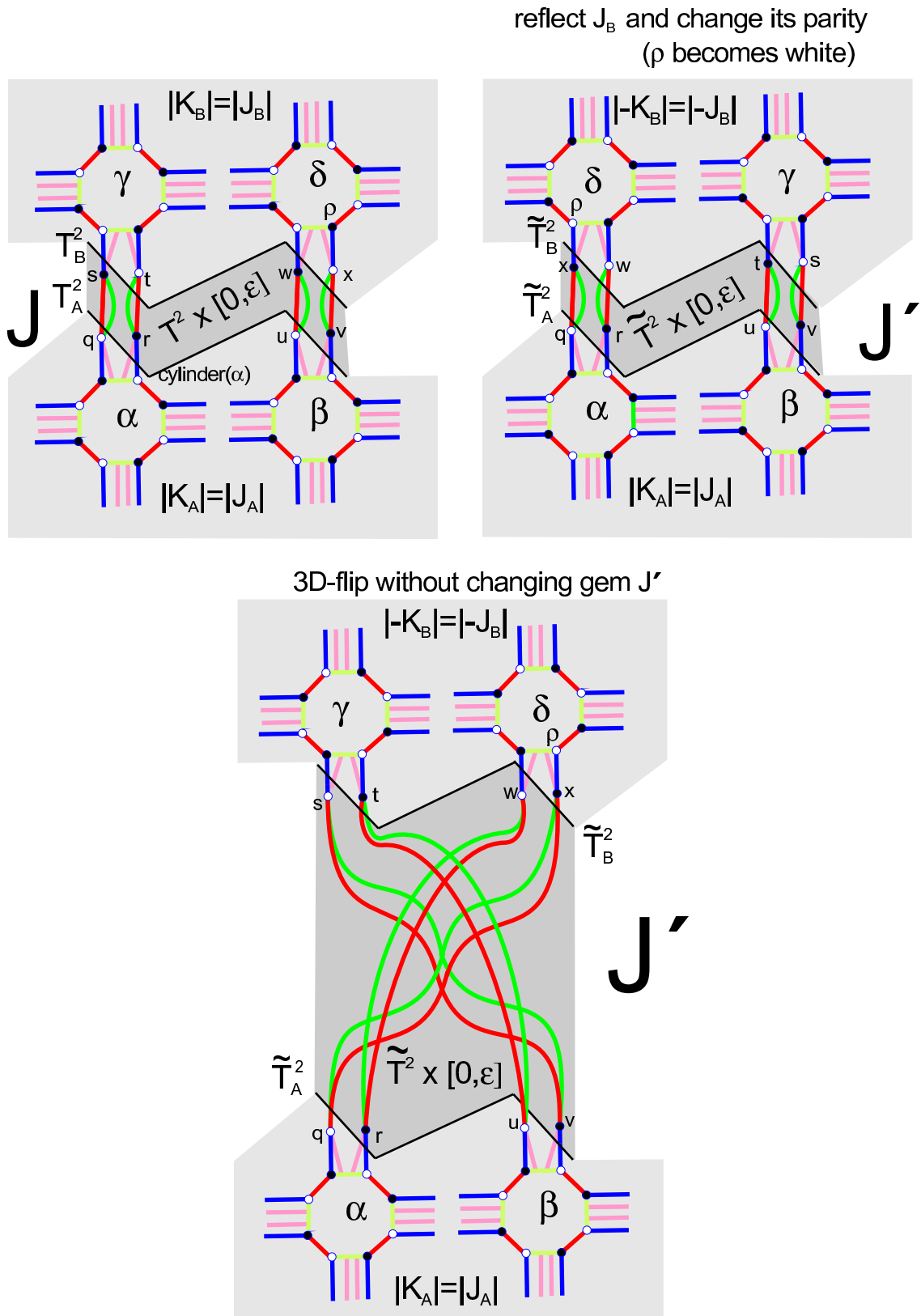}\\
   \end{center}
   \vspace{-0.7cm}
  \caption{For the proof of the Partial Reflection Theorem:
  gems $J$ and $J\, '$ induce the same space}
  \label{fig:inThreeParts}
\end{figure}

\noindent {\bf Proof of Theorem~\ref{theo:partialReflection}:\ } Let
$A$ and $B$ be arbitrary disjoint g-blinks, $(a,b)$ a basepair on
them. Then $A[a] + B[b] \Sequiv A[a] + \textsc{Reflection}(B)[b].$
\begin{proof}
Let $J$ be the canonical gem of the g-blink $A[a] + B[b]$ and $J\,
'$ be the canonical gem of g-blink $A[a] +
\textsc{Reflection}(B)[b]$. Let $K$ be a simplicial refinement of
$J^\star$ containing the $2$-torus $T^2$ (given in
Lemma~\ref{lem:torus}) as a subcomplex. Let $\widetilde{K}$ be a
simplicial refinement of $(J')^\star$ containing the $2$-torus
$\widetilde{T^2}$ (which plays in $J'$ the same role of $T^2$ in
$J$) as a subcomplex. Let $h'$ be a fixed homeomorphism which maps
$T^2 \times [0,\epsilon]$ onto $F'$ and $\widetilde{h}$ be a fixed
homeomorphism which maps $\widetilde{T^2} \times [0,\epsilon]$ onto
$\widetilde{F}$. By applying Equations~\ref{eq:tripartition1} and
\ref{eq:tripartition2} to the tori $T^2$ and $\widetilde{T^2}$ we
have ($-K_B$ is $K_B$ with orientation reversed: they are oriented
simplicial complexes):
$$ |K| = |K_A| \cup \left(T^2 \times
[0,\epsilon]\right)\ \cup\ |K_B|, \hspace{3mm} |\widetilde{K}| =
|K_A| \cup \left(\widetilde{T_2} \times [0,\epsilon]\right)\ \cup\
|-K_B|, \hspace{3mm} |K_A| \cap |K_B| = \emptyset.$$ with
\begin{equation}
\label{eq:tripartition4} |K_A| \cap \left( T^2 \times
[0,\epsilon]\right) = T_{\alpha\beta}^2, \hspace {8mm} \left({T^2}
\times [0,\epsilon]\right) \cap |K_B| = T_{\gamma\delta}^2,
\end{equation}
$$ |K_A| \cap \left( \widetilde{T^2} \times
[0,\epsilon]\right) = T_{\alpha\beta}^2, \hspace {8mm} \left(
\widetilde{T^2} \times [0,\epsilon]\right) \cap |-K_B| =
-T_{\gamma\delta}^2.$$ Define the map $\rho$ from $|K|$ to
$|\widetilde{K}|$ to be the identity in $|K_A| \cup |K_B|$. For $x
\in T^2 \times [0,\epsilon]$, define $\rho(x) = \left[(\,
\widetilde{h}\, )^{-1}\ {\circ} \ \mu\ {\circ} \ h'
\right](x) \in \widetilde{T^2} \times [0,\epsilon]$. Map $\rho$ is
the desired homeomorphism taking $|K|$ onto $|\widetilde{K}|$.
\end{proof}

\begin{figure}[htp]
   \begin{center}
      \leavevmode
      \includegraphics[width=12cm]{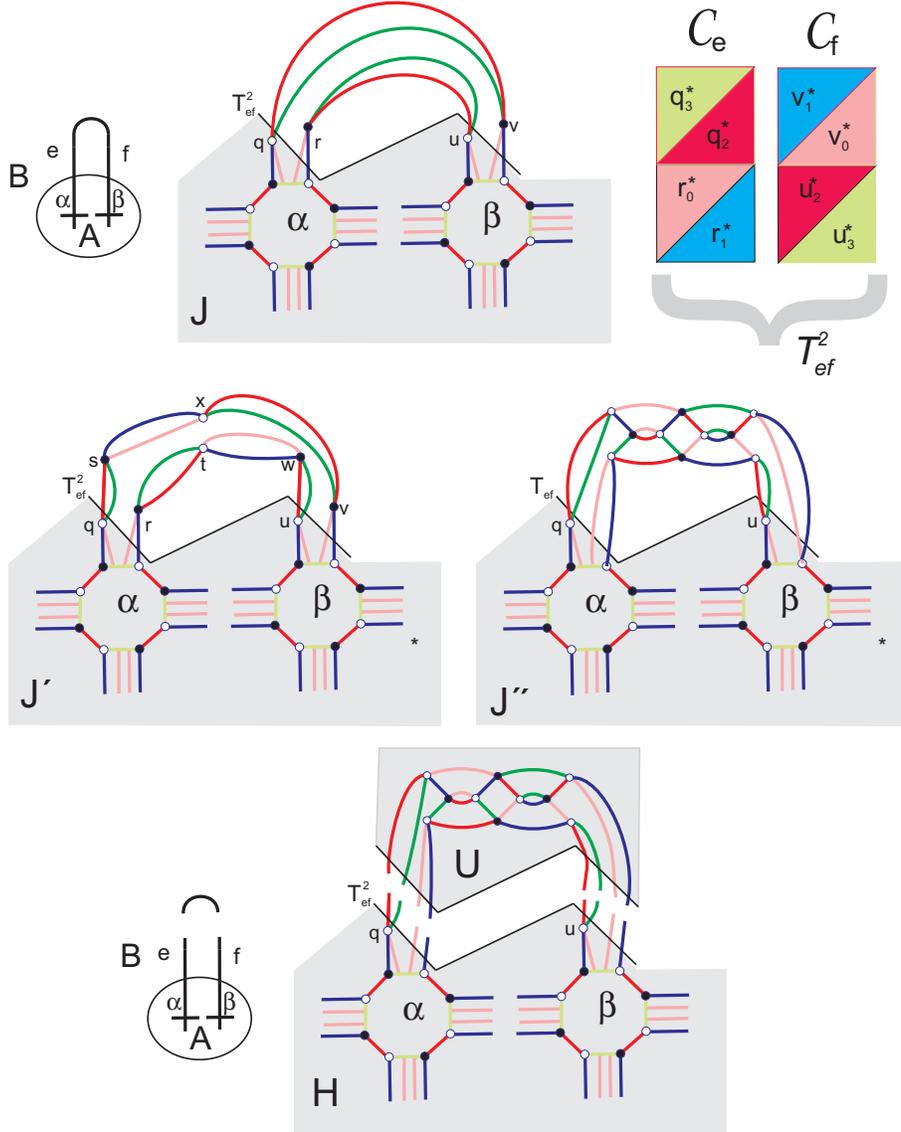}\\
   \end{center}
   \vspace{-0.7cm}
  \caption{Space $|H|$, $\partial(|H|)=T^2_{ef}$, canonical way
  to close it: $|B|=|H| \cup_{T^2_{H}\equiv T^2_U} |U|$, $|U|$ solid torus}
  \label{fig:spaceToroidalBoundary}
\end{figure}

We finish this chapter by proving the following Theorem about BFLs
with a segment between crossings removed:
\begin{Theo} \label{BFLwithSegmentRemoved}
Let $B^\circ$ be a BFL $B$ with a segment between crossings removed.
There exists a well defined 3-manifold with toroidal boundary
$S^\circ$ which can be associated to $B^\circ$. Moreover, there
exists a canonical way to close $S^\circ$ by attaching a solid torus
to its boundary to get a space $S$ such that $|B|=S$.
\end{Theo}
\begin{proof} The proof should be followed in Figure~\ref
{fig:spaceToroidalBoundary}. Gem $J$ is the canonical gem of the BFL
$B$. Gems $J\, '$ and $J\, '\, '$ are obtained from $J$ by 2-dipole
creations. The last gem is subdivided into two gems with boundary
$H$ and $U$. The boundaries of these gems are homeomorphic to the
$2$-torus $T^2_{ef} = \C_e + \C_f$ , given in Lemma~\ref{lem:torus}.
It can be shown that gem with toroidal boundary $U$ induces a solid
torus, thus completing the proof.
\end{proof}

%% file: chapter5.tex
\chapter{Computational experiments and results}
\label{chap:census}

\section{A census of prime spaces induced by small g-blinks}
\label{sec:census}

\newcommand{\kpu}{\,$k$-prime-unavoidable }
\newcommand{\npu}[1]{\,$#1$-prime-unavoidable }
\newcommand{\Rk}{{\cal R}_k }
\newcommand{\Rkl}{{\cal R}'_k}

In Section~\ref{sec:towardsACensusOfPrimeSpacesInducedBySmallBlinks}
we saw that if we have a set with all representative g-blinks with
$\,\,\leq k \,\,$ g-edges, then all spaces induced by blinks with
$\,\,\leq k \,\,$ edges are there. Even though, for any fixed $k \geq
1$ this set is finite, we would like to search for spaces in an
even smaller set and not lose any space. This may be done if we
believe in the following reasonable conjecture.

\begin{Conj} \label{conj:noCompositeHasFewerEdgesThanItsMinimalPrimeSum}
Let a space $S$ be the connected sum of prime spaces $A$ and $B$.
If the minimal blink presentations for $A$ and $B$ have respectively $n_A$ and
$n_B$ edges, then $n_S$, the number of edges for the minimal blink presentation
of $S$, satisfies  $n_S = n_A + n_B$.
\end{Conj}

A {\it minimal blink for a space} is a blink that
has the same number of edges or fewer edges than any other blink
presentation for that space. A {\it prime space} is one that cannot
be expressed by a connected sum of two or more spaces
different from $\IS^3$. A {\it composite space} is a space
that is not prime, \ie one space that can be expressed by a connected sum of
two or more spaces different from $\IS^3$. A blink presentation for any composite space
is obtained by drawing the blinks of each of its prime pieces
separately in the same drawing: a red-green plane graph with more than one
connected component. This construction clearly defines an upper bound for
the minimal number of edges of the composite space: the sum of
the number of edges of a minimal blink for each of its
prime pieces. Suppose space $S$ is the connected sum of $n$ prime
spaces, a minimal blink for $S$ with $n$ connected
components satisfies this upper bound, otherwise there would be
a blink for one of its prime pieces with fewer edges than its minimum
number of edges which is a contradiction. The only
possibility that remains to the conjecture to be false is
a blink presentation with less than $n$ components. On the other
hand we have,

\begin{Prop} \label{prop:primeIsSufficient}
If Conjecture~\ref{conj:noCompositeHasFewerEdgesThanItsMinimalPrimeSum} is true
then, knowing one minimal blink for each prime space that has a blink presentation
with $\leq k$ edges, it is possible to exhibit a minimal blink
for all spaces (composite or prime) that have a blink presentation with
$\leq k$ edges.
\end{Prop}
\begin{proof}
If $S$ is a prime space that has a blink with fewer than $k$ edges then,
by hypothesis, we already know one minimal blink for it. If $S$ is composite,
then exhibit together in a single drawing the minimal blink of all its prime
pieces. By the Conjecture~\ref{conj:noCompositeHasFewerEdgesThanItsMinimalPrimeSum}
this is a minimal blink for $S$.
\end{proof}

A consequence of Proposition~\ref{prop:primeIsSufficient} is that, if
Conjecture~\ref{conj:noCompositeHasFewerEdgesThanItsMinimalPrimeSum} is true, then
to identify all spaces that have a blink presentation with fewer than $k$ edges
it is sufficient to know only the prime spaces that have a blink
presentation with fewer than $k$ edges. So, in our experiment,
this is what we do. We focus only on prime spaces. We simplify our
computational effort by answering not what are all the spaces with
a blink presentation with fewer than $k$ edges, but what are all the
\underline{prime} spaces with a blink presentation with fewer than $k$
edges. So, we may reduce the set of representative g-blinks to
representative g-blinks that are prime or, in practice, that are
not easily shown composite.

Spaces have an orientation. If we swap the red-green edges of a
blink the effect on the induced space is its change of orientation.
So, any set of blinks $B$ that induces spaces $S$ may be easily
extended to a set of blinks $B'$ that induces $S$ and also the
changed orientation version of the spaces in $S$. The set $B'$ is
just $B$ plus the blinks of $B$ with the red-green edges swapped.
This leads us to the following definition: a {\it set of blinks} is
said to be {\it \kpu} if every space with a blink presentation with
$\,\,\leq k \,\,$ edges is induced by some blink in this set or a
red-green swapped version of some blink in this set. These notions
are analogously extended to g-blinks. A {\it set of g-blinks} is
said to be {\it \kpu} if every space with a g-blink presentation with
$\,\,\leq k \,\,$ g-edges is induced by some g-blink in this set or
a red-green g-edges swapped version of some g-blink in this set.

\newpage

What we mean concretely by the title of this section:
{\it a census of prime spaces induced by
g-blinks} is a triple \vspace{-0.3cm}$$(k,{\cal B},f:{\cal B}
\rightarrow \{1, \ldots, n\}),\vspace{-0.3cm}$$ where $k$ is a
positive integer, ${\cal B}$ is a \kpu set of g-blinks and $f$ is a
surjective function that maps each g-blink in ${\cal B}$ to an integer
in $\{1, \ldots, n\}$ satisfying the constraints: if $B_1,B_2 \in
{\cal B}$ induce the same space or induce the same space with
swapped orientations then $f(B_1) = f(B_2)$, else
$f(B_1) \neq f(B_2)$. Note that $f$ defines a partition of ${\cal B}$
into $n$ classes where the g-blinks in each class induce the same space
modulo orientation. For this reason we call $f$ the
{\it partition function} of the census. In view of this definition of a
census of prime spaces, the steps to build one are: (1) define $k$;
(2) define a \kpu set of g-blinks; (3) define the partition function $f$.

\newpage
\section{A prime-unavoidable set of g-blinks: $U$}

To obtain a census of prime spaces induced by g-blinks with $\,\,\leq
k \,\,$ edges, a set of \kpu g-blinks is needed. Before defining the
specific way we did this generation it is good to say that in theory
what is needed to do is simple: enumerate all representative g-blinks up to size
$k$ and discard those that you can show that are composite or that are
not minimal. In practice we use some shortcuts to avoid a full
enumeration of all representative g-blinks.

\begin{figure}[htp]
   \begin{center}
      \leavevmode
      \includegraphics[width=12cm]{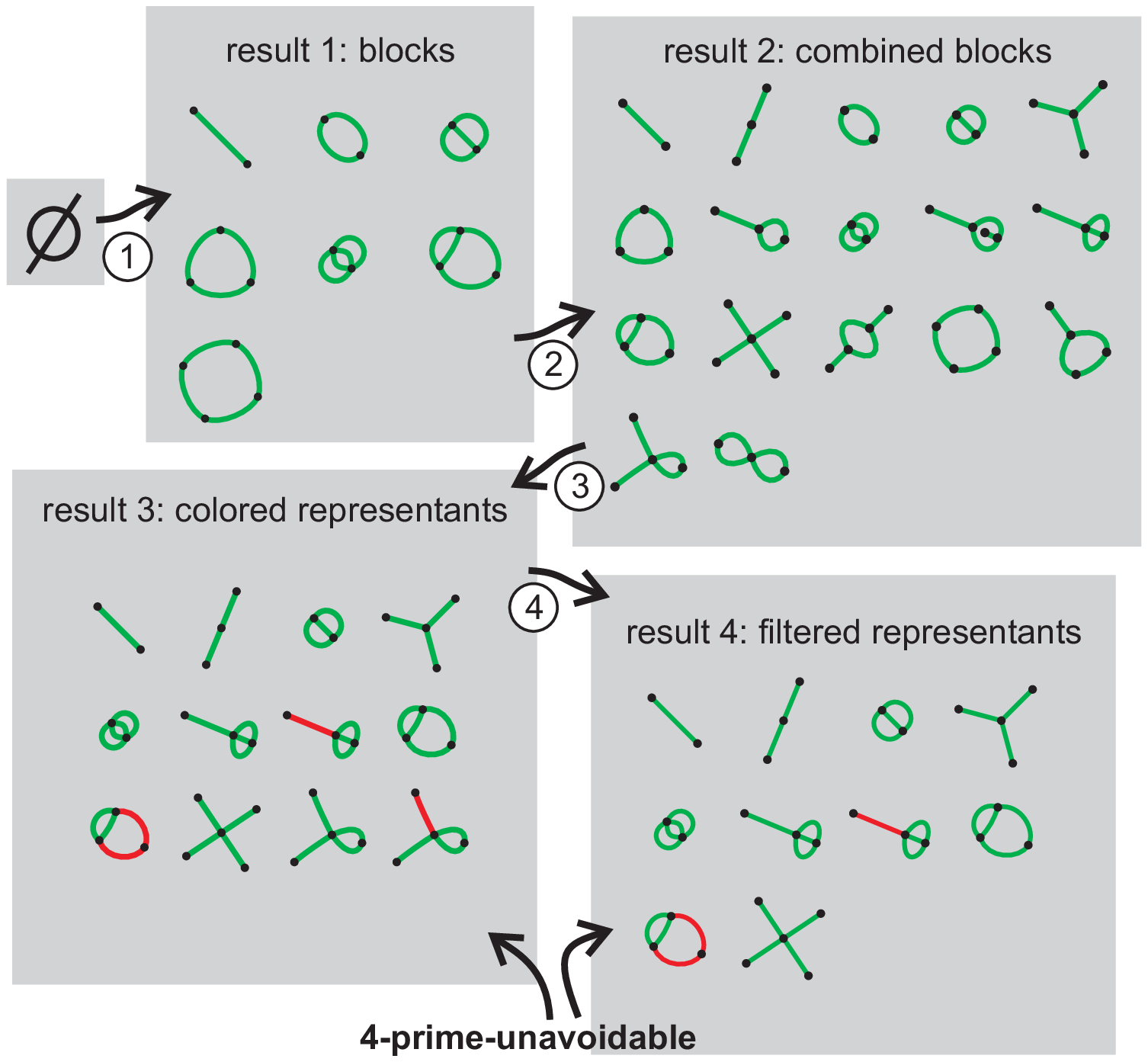}
   \end{center}
   \vspace{-0.7cm}
   \caption{Pipeline of the \kpu set generation for $k=4$}
   \label{fig:unavoidableSetPipeline}
\end{figure}

The procedure we defined to obtain a \kpu set was a pipeline with 4
steps. The output of each step was the input to the next one. The
final and intermediate results of this procedure for $k=4$ is shown
on Figure~\ref{fig:unavoidableSetPipeline}. The steps on this pipeline,
presented on this figure by an arrow and a number, are named:
(1) \textsc{BlockGeneration}, (2) \textsc{BlockCombination}, (3) \textsc{Coloring}
and (4) \textsc{Filtering}.

\begin{figure}[htp]
   \begin{center}
      \leavevmode
      \includegraphics[width=12cm]{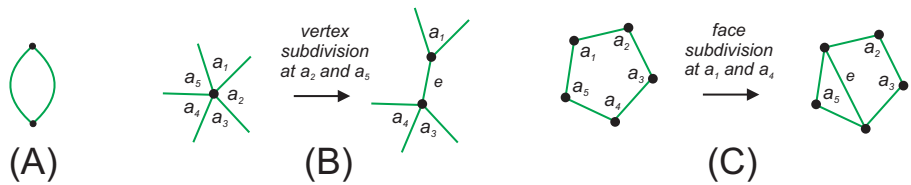}
   \end{center}
   \vspace{-0.7cm}
   \caption{Block generation}
   \label{fig:blockGeneration}
\end{figure}

The \textsc{BlockGeneration} step has as input one positive integer $k$:
the maximum number of edges. Its output is all possible 2-connected
plane graphs with number of g-edges not exceeding $k$ plus the single-edge
plane graph. By convention we see all these plane graphs as
green-edged blinks. The blocks, besides the single-edged one, are
obtained from the plane graph with two parallel green edges shown on
Figure~\ref{fig:blockGeneration}A by applying inductively
and in all possible ways {\it vertex subdivisions} and
{\it face subdivisions}. An example of vertex subdivision may be seen on
Figure~\ref{fig:blockGeneration}B. Any two distinct angles on a
vertex are a base for this operation.  An example of face
subdivision may be seen on Figure~\ref{fig:blockGeneration}C.
Any two distinct angles on a face are a base for
this operation. For $k=4$, the number of resulting
blocks is 7 as it is shown on Figure~\ref{fig:unavoidableSetPipeline}.
The block term used here is also aligned to the fact that the
g-blinks induced from these resulting green-edged blinks do not
have breakpairs.

The \textsc{BlockCombination} step has as input $k$, the maximum number
of edges or g-edges, and the resulting blocks from the \textsc{BlockGeneration}
procedure. Here we see this input as green g-blinks and apply the
following algorithm. Let $B$ the set with these
input g-blinks or blocks. Make $A_1 = B$. For $i$ from 2 to $k$ make
$A_i$ the result of {\it combining} every g-blink at $A_{i-1}$ with each
g-blink in $B$. Combining a g-blink $G$ with $n_G$ g-zigzags to a
g-blink $G'$ with $n_{G'}$ g-zigzags results in $n_G \times n_{G'}$
g-blinks. This is the result of merging $G$ and $G'$ on
basepairs coming from all distinct combinations of g-zigzags. This
includes all possible spaces obtainable from merging $G$ and $G'$
as asserts Proposition~\ref{prop:sameZigzagsSameSpace}. The g-blinks
that overflows the maximum number of g-edges $k$ are discarded. For $k=4$,
the number of resulting combinations is 17 as it is shown on
Figure~\ref{fig:unavoidableSetPipeline}. Now we do an important
observation. Merging two g-blinks and then assigning a color to
each of its g-edges is the same as assigning the right colors in each
of the two g-blinks before merging and then merging them. This implies
that coloring the all-green g-blinks resulting from this step in all
possible ways really spans all possible spaces.

The \textsc{Coloring} step has as input the  all-green
g-blinks from the \textsc{BlockCombination} procedure.
The idea is to assign all possible g-edge color combination
to each of the given g-blinks. For each g-blink assigned
with a coloring some tests are made and this g-blink may
be discarded if it is asserted that, by doing this, we are
not losing a minimal g-blink to that same space (or its swapped
orientation version). Let $G$ be a g-blink already assigned
a coloring, these tests are the following:
\begin{enumerate}
\item If the number of red g-edges on $G$ is greater than the
number of green g-edges then it is discarded. This is
justified by the fact the red-green g-edges swapped version of $G$
will not be discarded by this rule (green g-edges is greater than
red g-edges) and it induces the same space as $G$ with orientation
changed.

\item If $G$ contains the structure shown
on the left side of Figure~\ref{fig:unecessaryGBlinkStructure}A it
is possible to apply a Reidemeister move of type II reducing by 2
the number of crossings and preserving the space. So $G$ is
unnecessary once its induced space was already considered by some
g-blink with fewer g-edges.

\begin{figure}[htp]
   \begin{center}
      \leavevmode
      \includegraphics[width=12cm]{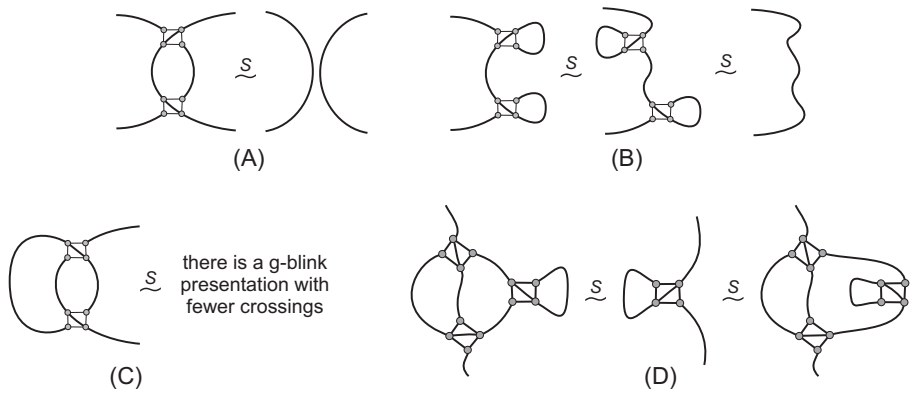}
   \end{center}
   \vspace{-0.7cm}
   \caption{ Structures to identify g-blinks that may be discarded}
   \label{fig:unecessaryGBlinkStructure}
\end{figure}

\item If $G$ contains the structure shown
on the left side pattern or the middle pattern of
Figure~\ref{fig:unecessaryGBlinkStructure}B then it is unnecessary.
If it contains the left side pattern of
Figure~\ref{fig:unecessaryGBlinkStructure}B, by a ribbon move it is
converted to the pattern on the middle of
Figure~\ref{fig:unecessaryGBlinkStructure}B that may be converted by
Whitney Trick to the right pattern of
Figure~\ref{fig:unecessaryGBlinkStructure}B. All of them induce the
same space, and the right pattern has fewer crossings.

\item If $G$ contains the pattern on
Figure~\ref{fig:unecessaryGBlinkStructure}C then it has a {\it
circumcised component} and may be simplified by moves in the BFL
calculus to a g-blink with fewer crossings as it is explained on
pages 138--140 of \cite{KauffmanAndLins1994}. So, it may be discarded.

\item If $G$ contains the left or the right pattern on
Figure~\ref{fig:unecessaryGBlinkStructure}D then it may be simplified
by the move $K_4(1)$ of BFL calculus to the middle pattern with only
one crossing. So, it may be discarded.

\item If $G$ seen as a BFL contains more than one
component (more than one g-zigzag) and one of the components is completely
overcrossing the others or completely undercrossing the others than it may
be separated by Reidemeister moves and is not a minimal presentation.
\end{enumerate}
If g-blink $G$ passes all tests then one last transformation is done.
The g-blink included in the result set of this step is
actually $\min\{r(G), r(-G)\}$: the smallest g-blink
between the representative of $G$ or the representative of $-G$
\,\,(\ie g-blink $G$ with all crossings changed (C) or, equivalently
all g-edge colors swapped). This resulting g-blink is asserted to
induce the same space as $G$ or its changed orientation version. Observe
that the g-blinks resulting from this step are all representatives.  For $k=4$,
the number of g-blinks resulting from this step is 12 (see
Figure~\ref{fig:unavoidableSetPipeline}).

The \textsc{Filtering} step has as input the representative
g-blinks resulting from the \textsc{Coloring} step. Let $B$ be this
input set and $R = \{\}$, initially empty, be the result set.
The filtering algorithm flows like this:

{\setstretch{1.25}
\begin{center}
\hbox{
\parbox{15cm}{
\begin{small}
\begin{algorithmic}[1]
 \While{$B$ \, is not empty}
    \State $G \leftarrow$ an element of $B$
    \State $RM3 \leftarrow$ closure of $G$ by Reidemeiter III Move
    \If{no element of $RM3$ may be discarded by rules 2 to 6 of the \textsc{Coloring} step}
       \State $R \leftarrow R \cup \{G\}$
    \EndIf
    \State $B \leftarrow B \backslash RM3$
 \EndWhile
\end{algorithmic}
\end{small}}}\end{center}
}
The idea of this step is to use the Reidemeister III move, that preserves
the number of crossings and the space, to find some blink version of the
space that may be simplified by the rules 2 to 6 explained on the
\textsc{Coloring} step. If there exists such a version then that
g-blink and the whole closure of g-blinks obtained from it by
Reidemeister III may be discarded. Otherwise only the given
g-blink on its Reidemeister III closure may be preserved (that
is why we remove all RM3 set on line 7 of the above algorithm).
The result of this step for $k = 4$ is a set with 10 g-blinks
(see Figure~\ref{fig:unavoidableSetPipeline}).

We name $U$ the set resulting from this pipeline for $k=9$. This set
is, as we saw in its construction, a \hbox{\npu{9}} set and has 3437 g-blinks
divided in 1 g-blink with 1 g-edge, 1 g-blink with 2 g-edges, 2 g-blinks
with 3 g-edges, 6 g-blink with 4 g-edges, 12 g-blinks with 5 g-edges, 43
g-blinks with 6 g-edges, 133 g-blinks with 7 g-edges, 585 g-blinks with
8 g-edges and 2654 g-blinks with 9 g-edges. We denote by $U[1]$ the
smallest g-blink in $U$, $U[2]$ the second smallest g-blink in $U$,
up to $U[3437]$ the greatest g-blink in $U$. The time elapsed to
generate the set $U$ was less than twelve hours. At this point we have
the first two ingredients to a census of prime spaces up to 9 g-edges:
$(9, U, ?)$. The only missing part is the third ingredient of a census:
the partition function. This is the subject of next section.

\section{Topological classification of g-blinks in $U$}
\label{sec:topologicalClassificationOfU}

The set $U$ is a \npu{9} set of g-blinks. The next step to define
a census of prime spaces with a g-blink (or blinks) presentation
with up to 9 g-edges (edges) is to identify what g-blinks induce
different spaces and what g-blinks induce the same space
(modulo the orientation). To reach
this goal the first thing we did was to
calculate, for each g-blink in $U$, the homology group and the
Witten-Reshetikhin-Turaev quantum invariant (for
$r \in \{3,4,5,6,7,8\}$) of its induced space. In
Sections~\ref{sec:homologyGroup}~and~\ref{sec:quantumInvariant}
we show how to do this calculation from a g-blink
presentation of a space. To help on this exposition
we will use HG, QI and HGQI when we want to refer
to the homology group, quantum invariant, respectively.
The time elapsed to calculate the HG and QI of all
g-blinks in $U$ was less than half an hour.

The effect on the quantum invariant of changing the
orientation of a space is that each complex number in
its sequence becomes its conjugate (remember that the
quantum invariant is a sequence of complex numbers). So
when two g-blinks have their QIs differing by, for each
$r$, one being the conjugate of the other, then these
g-blinks may induce the same space in different
orientations. As it was defined for a census, g-blinks
that induce distinct orientations of the same space are
mapped, by the partition function, to the same value. We
are interested in spaces modulo orientations. For this
reason, we mounted from the HG and QI data of each
g-blink in $U$ the information named HGnQI (HG and
normalized QI). It is just the pair HG and nQI where
nQI is the normalized version of QI: if the first complex
entry with imaginary part in QI is negative then nQI
entries are the conjugate of QI entries, otherwise nQI
is equal to QI.

Using the HGnQI information of each g-blink we partitioned the
set $U$ into 501 classes. The 3437 g-blinks of $U$ induced 501
distinct HGnQIs. One consequence of this fact is that $U$ induces, at
least, 501 different (modulo orientation) spaces. This HGnQI partition
of the set $U$ is a first candidate for the partition function
to the census we want. If the homology group together with
the quantum invariant is a strong enough invariant of space,
then we already have the exact partition function we want. To
prove this, it remains to show that all entries in the
same HGnQI class indeed induce the same space. To do
this, we need another tool. Before entering into this topic,
we want to make some comments about the HGnQI partition of $U$.

\begin{figure}[htp]
   \begin{center}
      \leavevmode
      \includegraphics[width=15cm]{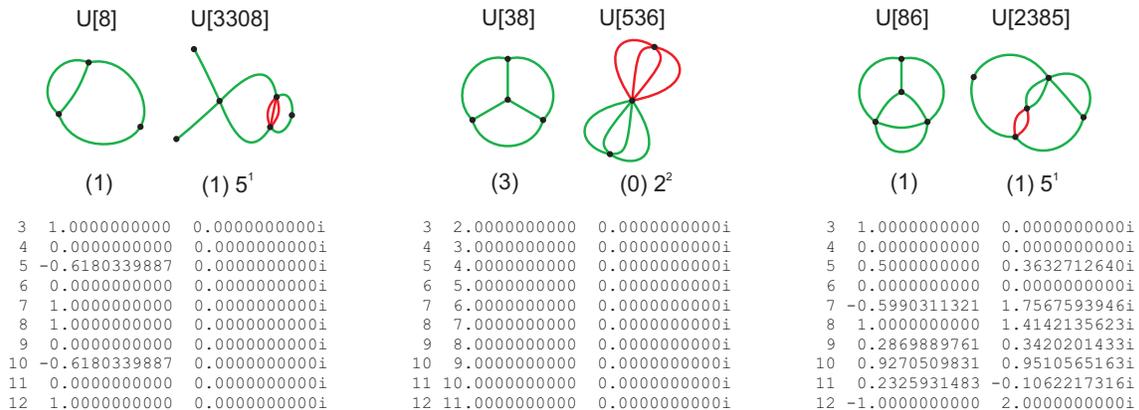}
   \end{center}
   \vspace{-0.7cm}
   \caption{ The 3 cases in $U$ where HG helped QI to distinguish spaces}
   \label{fig:qiFailure}
\end{figure}

After partitioning the set $U$ in HGnQI classes, a very apparent
fact was that the quantum invariant was almost perfect in
identifying the 501 classes. It, alone, separated $U$ into 498
classes. In only 3 cases the homology group was important to
distinguish spaces that the quantum invariant did not. In
Figure~\ref{fig:qiFailure} we show a blink presentation for
these 3 cases.
In the first case, $U[8]$ has HG $(1)$
while $U[3308]$ has HG $(1) 5^1$ and the 12 first
entries of the quantum invariants, as is shown,
are all real numbers. Indeed, all entries in the QI of $U[8]$
sequence are real once it has only one orientation.
For a proof of this fact note that
$\textsc{Dual}(U[8])$ is $U[8]$ with
all edges being red which is also $U[8]$ after applying (C)
(change crossings). As these two g-blinks are the same
g-blink and they are the two possible orientations
for the space, we can conclude that this space has only
one orientation. In the second case, $U[38]$ has HG $(3)$
while $U[536]$ has HG $(0) 2^2$ and the 12 first
entries of the quantum invariants, as is shown,
are all integer numbers. In the third case, $U[86]$ has HG $(1)$
while $U[2385]$ has HG $(1) 5^1$ and the 12 first
entries of the quantum invariants, as is shown,
are all complex numbers. It might be the case that
the quantum invariant in some point distinguishes these
spaces as the homology group did. We did not check this.

Now let's return to our open problem. Are the 501 HGnQI
classes really inducing the same space or some of them
induce more than one space? To answer this question we
used 3-gem theory. We saw in Section~\ref{sec:fromGBlinkTo3Gem}
that from a g-blink we can obtain a 3-gem inducing the
same space as it does. This fact enables us to change
our question in g-blink language into a question in 3-gem language. The
idea is to take, for each of the 501 HGnQI classes, all g-blinks
in the same HGnQI class, calculate a 3-gem version for it
and then try to find a proof that they are the same space, \ie
a path of ``moves'' in 3-gems that preserve the induced space
connecting all these 3-gems.

\newcommand{\GemOfGBlink}{
\parbox{7.5cm}{
\begin{footnotesize}
\setstretch{1.25}
\begin{algorithmic}[1]
   \Function{GemOfGBlink}{$G$}
   \State $J \leftarrow \textsc{GBlink2Gem}(G)$
   \While{true}
      \State $\textsc{SimplifyGem}(J)$
      \State $\textsc{SearchInTSClass}(J,12\, {\rm seconds})$
      \If{$J$ has no dipole, no $\rho_2$-pair and no $\rho_3$-pair}
         \State break
      \EndIf
   \EndWhile
   \State {\bf return} \,\, $J$
   \EndFunction
\end{algorithmic}
\end{footnotesize}
}}

\newcommand{\SearchInTSClass}{
\parbox{7.5cm}{
\begin{footnotesize}
\setstretch{1.25}
\begin{algorithmic}[1]
   \Procedure{SearchInTSClass}{$J, maxtime$}
      \State $C \gets \{ J \}$ \,\, $U \gets \{ J \}$ \Comment{$C$ is the current TS-class of $J$ and $U$ are the unprocessed gems}
      \While{$U$ is not empty and elapsed time < $maxtime$}
         \State $J' \leftarrow$ a gem in $U$
         \State $U \leftarrow U \backslash \{J'\}$
         \For{all possible TS-moves $m$ in $J'$ }
            \State $J'' \leftarrow J'$ with TS-move $m$ applied
            \If{$J'' \notin C$}
                \If{there is a dipole or $\rho_2$-move or $\rho_3$-move in $J''$}
                    \State $J \gets J''$ and exit
                \Else
                    \State $U \leftarrow U \cup \{J''\}$
                    \State $C \leftarrow C \cup \{J''\}$
                \EndIf
            \EndIf
         \EndFor
      \EndWhile
      \State $J \gets \min\{J'\in C\}$
   \EndProcedure
\end{algorithmic}
\end{footnotesize}
}}

\newcommand{\SimplifyGem}{
\parbox{7.5cm}{
\begin{footnotesize}
\setstretch{1.25}
\begin{algorithmic}[1]
   \Procedure{SimplifyGem}{$J$} \Comment{$J$ becomes its simplified version}
   \While{true}
      \If{there is a dipole in $J$}
         \State apply dipole cancelation in $J$
      \ElsIf{there is a $\rho_3$-pair or $\rho_2$-pair in $J$}
         \State apply $\rho_3$-pair or apply $\rho_2$-pair in $J$
      \Else  \State break
      \EndIf
   \EndWhile
   \EndProcedure
\end{algorithmic}
\end{footnotesize}
}}

The 3-gem that we associated to each g-blink in $U$ was given by
the function \textsc{GemOfGBlink} shown in Algorithm \ref{alg:algorithmFor3Gems}.
The idea of this function is to simplify the initial 3-gem of
the g-blink given by the \textsc{GBlink2Gem} procedure explained in
Section~\ref{sec:fromGBlinkTo3Gem} using dipole cancelations,
$\rho_2$-moves, $\rho_3$-moves and TS-moves until it cannot be simplified
anymore or until a certain timeout occurs. This step resulted in 999
distinct gems for the 3437 g-blinks of $U$. We used a timeout of 12 seconds.
From these 999 3-gems, 657 (or 65\%) gems were proven to be
TS-class representatives (minimum 3-gem in the class) such that the entire class
had no simplifications of the types: dipole cancelation, $\rho_2$-move and $\rho_3$-move.
The remaining 342 3-gems were the minimum 3-gem obtained before the timeout
occurred. The 3-gem also encodes the orientation of the space, but, in
this case we used 3-gems modulo orientation. In other words, the 3-gem we
associated to each g-blink could be exactly the same space, or the same
space with orientation changed. As it might be clear now, this is enough
here: spaces modulo orientation.
\begin{algorithm}
\caption{Algorithms for 3-Gems}
\label{alg:algorithmFor3Gems}
\smallskip
\begin{tabular}{c|c}
\begin{tabular}{c}
\GemOfGBlink\\[2.9cm] \hline
\\[-0.05cm]
\SimplifyGem
\end{tabular} &
\SearchInTSClass
\end{tabular}
\end{algorithm}

The remaining challenge at this point was to find whether these 999
3-gems, seen as nodes of a graph, could be connected in 501
connected components, where a connected component means that all
gems in the same component induce the same space (modulo
orientation). So, we started to insert edges in this graph of 999
nodes and initially no edge. This was done by ``perturbing'' the
gems on the nodes by using U-moves and then applying the same
simplification procedure used in the function \textsc{GemOfGBlink}
until a gem with no simplification or a timeout occurred. This final
gem, if not yet in our graph, was added as a new node. An edge, if
not existent, from the perturbed 3-gem node to this new, or already
existent, node, was also added to the graph. This procedure was
oriented by the HGnQI classes, so if a HGnQI class was already a
single connected component then nothing more was needed to be done
there: the HGnQI class was proved to be a single space (modulo
orientation). This procedure of connecting the gems of a HGnQI class
on this graph took about 3 days with manual interference being
important: by looking at the graph we perturbed the most promising
nodes. The final result was: 499 of the 501 HGnQI classes were
proven to induce a single space (modulo orientation). In only two
HGnQI classes we could not find a single connected component.

\begin{figure}[htp]
   \begin{center}
      \leavevmode
      \includegraphics[width=15cm]{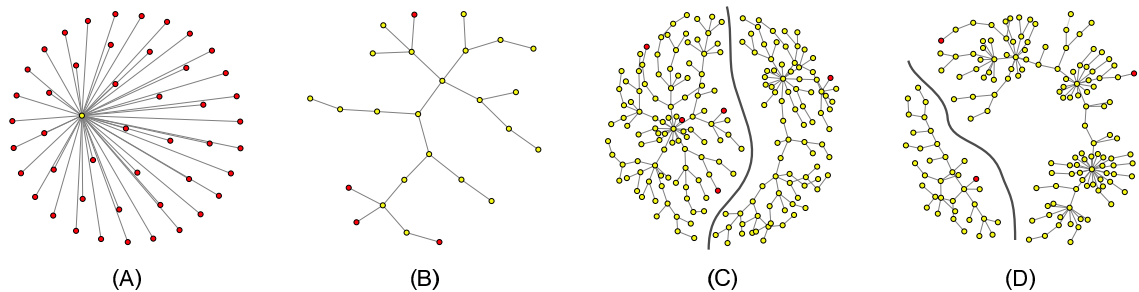}
   \end{center}
   \vspace{-0.7cm}
   \caption{Graphs of g-blinks (red nodes) and gems (yellow nodes). The first two are trees
   and the last two are forests with two components (the two uncertainties)}
   \label{fig:subgraphs}
\end{figure}

Figure~\ref{fig:subgraphs} shows subgraphs (trees) for 4 HGnQI
classes on the final graph. The red nodes are g-blinks from $U$. The
yellow nodes are the 3-gems. Note that every red node is connected
to a single yellow node: this yellow node is the result of the
\textsc{GemOfGBlink} applied to this g-blink.
Figure~\ref{fig:subgraphs}A was an easy case where all g-blinks
of the same HGnQI class were pointing right to the same 3-gem.
Nothing was needed to do in this case. Figure~\ref{fig:subgraphs}B
was one of the difficult cases: many redundant edges (not shown) and
different 3-gems were generated before all g-blinks were connected.
\begin{figure}[htp]
   \begin{center}
      \leavevmode
      \includegraphics[width=16cm]{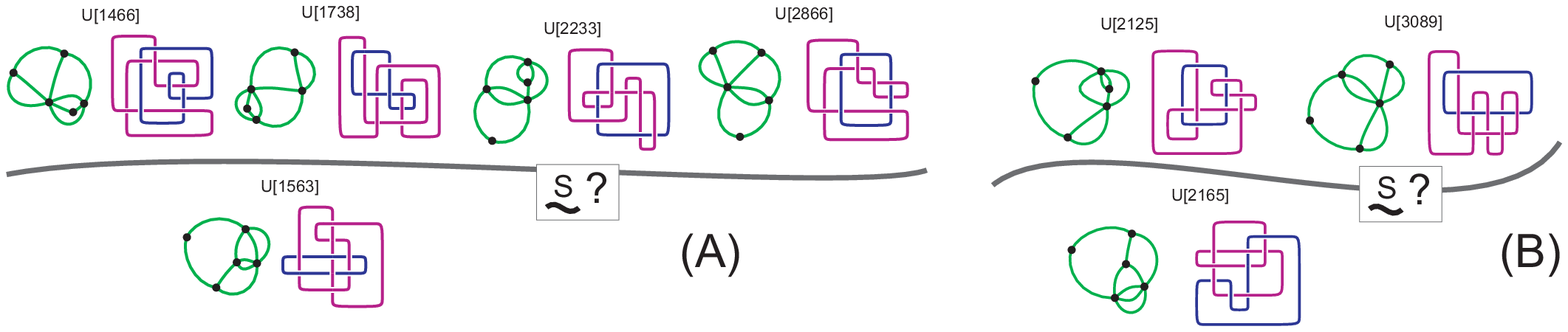}
   \end{center}
   \vspace{-0.7cm}
   \caption{ The only 2 classes with same HGnQI where a proof of the homeomorphism was not found}
   \label{fig:doubts}
\end{figure}

Figures~\ref{fig:subgraphs}C~and~\ref{fig:subgraphs}D presents the
two cases where one doubt was left. In each of these cases, two
connected components remained: they are shown with the dark line
separating them. In the first case, Figure~\ref{fig:subgraphs}C, the
HGnQI class had 5 g-blinks where 4 of them were proven to be the
same space. In the second case, the HGnQI class had 3 g-blinks where
2 of them were proven to induce the same space. A blink and BFL
presentation for the g-blinks involved in these doubts are shown in
Figure~\ref{fig:doubts}. In the first doubt, the g-blinks involved
are $U[1466]$, $U[1563]$, $U[1738]$, $U[2233]$, $U[2866]$ and
$U[1563]$ is the only g-blink we did not find a proof as being the
same space (modulo orientation) of the others. In the second doubt,
the g-blinks involved are $U[2125]$, $U[2165]$, $U[3089]$ and
 $U[2165]$ is the only g-blink we did not find a proof as
being the same space of the others. Is there a proof for these
two cases and we just could not find them or are these the only
weak points of the HGnQI invariant on the set $U$? We leave
this question open and register it as the following conjectures.

\begin{Conj} \label{conj:conjecture1}
The spaces induced by all 5 blinks or BFLs on Figure~\ref{fig:doubts}A
are the same.
\end{Conj}

\begin{Conj} \label{conj:conjecture2}
The spaces induced by all 3 blinks or BFLs on Figure~\ref{fig:doubts}B
are the same.
\end{Conj}

The only reason we conjecture these stems from the fact that
HGnQI have not failed in all other 499 cases. But,
the fact that we had no success, after various days of
computational effort trying to prove these conjectures
using the simplification combinatorial dynamics of 3-gems,
suggests the contrary: these conjectures are false.
Figures~\ref{fig:subgraphs}C~and~\ref{fig:subgraphs}D show the
two trees that could not be connected for
each case after all the computational effort.

All data involved in all the experiments we explained here are in a computer program named
\textsc{Blink}. So a proof that all HGnQI classes indeed induce the same space,
except for the two cases explained, can be exhibited by this program.

\enlargethispage{0.5cm}
In the 3-gem presentation it is sometimes possible to identify that
its induced space is composite. For example the space induced by g-blink $U[31]$
is also induced by a 3-gem ($r_4^{18}$ in the 3-gems catalogue of \cite{Lins1995}) that
contains a {\it disconnecting quartet}, \ie four edges with distinct colors
that disconnected the 3-gem. The existence of this structure in a 3-gem
or the existence of handles, \ie connected sums with $\IS^2 \times \IS^1$,
is a proof that the induced space is composite. In the 501 HGnQI classes, using
this kind of 3-gem information, we could prove that 14 of them were composite. The
rules that we used in the construction of $U$ were not able to identify that
some g-blinks induced composite spaces, but, anyway, that was not the goal there. The
goal there was to create a small set of g-blinks that did not lose a minimal
presentation by g-blink of a prime space. This is the important property of $U$:
all prime spaces have a minimal g-blink presentation in $U$. Using this information
of the 14 composite classes, we named each of the 501 HGnQI classes like this:
the 487 classes that were not proven composite gained names 1.1, 2.1, 3.1  $\ldots$  3.2,
4.1 $\ldots$ 4.5, 5.1 $\ldots$ 5.6, 6.1 $\ldots$ 6.19, 7.1 $\ldots$ 7.38, 8.1 $\ldots$ 8.119
and 9.1 $\ldots$ 9.296; the 14 classes that were proven composite gained names 6.1c \ldots 6.3c,
8.1c \ldots 8.5c and 9.1c $\ldots$ 9.6c. The number before the point stands for
the number of g-edges of the minimal g-blink in $U$ found for that space.
Let $U[n.i]$ denote the smallest g-blink (\ie smallest code) in class $n.i$, \ie
$U[n.i] = \min\{G \in n.i\}$. For instance $U[5.1]$ is $U[11]$ and $U[6.1{\rm c}]$ is $U[31]$.
The number after the point stands for the following: $n.1$ is the class
where g-blink $U[n.1]$ has $n$ g-edges and is the smallest g-blink among all
classes $U[n.j]$, for any $j$ that defines a valid class name; $n.2$ is the class
where g-blink $U[n.2]$ has $n$ g-edges and is the second smallest g-blink among all
classes $U[n.j]$, for any $j$ that defines a valid class name; and so on. The two
classes that we do not know whether they induce a single space or two spaces are 9.126
(Figure~\ref{fig:doubts}A) and 9.199 (Figure~\ref{fig:doubts}B).

\begin{figure}[h!tp]
   \begin{center}
      \leavevmode
      \includegraphics[width=16cm]{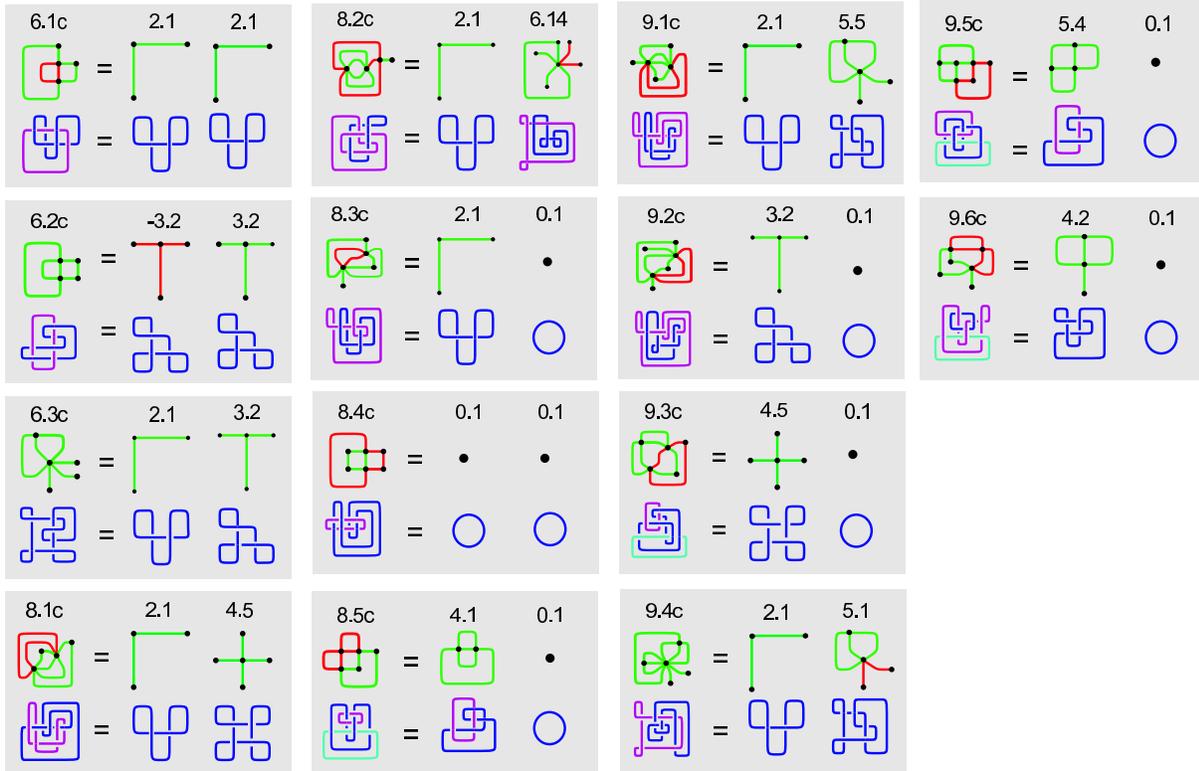}
   \end{center}
   \vspace{-0.7cm}
   \caption{ The 14 composite spaces in $U$}
   \label{fig:compositeSpaces}
\end{figure}

The 14 composite spaces in $U$ are in Appendix~\ref{chap:compositeCatalogue}.
The quantum invariant at level $r$ of the connected sum of spaces $A_1 \ldots A_n$ is the product
of their quantum invariants at the same level divided by the $r$-th quantum invariant of $\IS^3$ to
the power $n-1$. Using this we could align the orientations of the prime spaces that
produced these composite spaces. Figure~\ref{fig:compositeSpaces} shows explicitly these
14 spaces as a prime space composition with the correct orientation.

\begin{figure}[h!tp]
   \begin{center}
      \leavevmode
      \includegraphics[width=14cm]{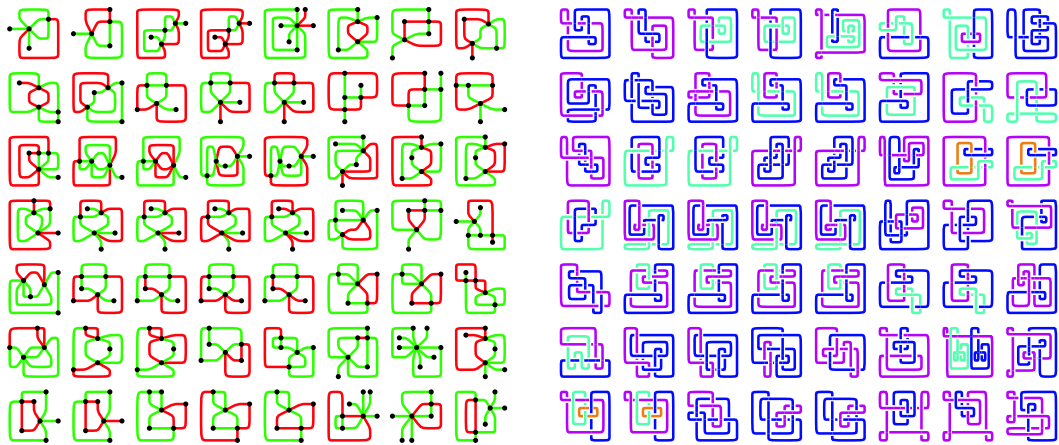}
   \end{center}
   \vspace{-0.7cm}
   \caption{ Blink and BFL presentations for g-blinks in 6.5: space $\IS^2 \times \IS^1$}
   \label{fig:gblinksForS1xS2inU}
\end{figure}

The space $\IS^2 \times \IS^1$ has a blink presentation that is just a vertex
and no edges, \ie a BFL that has no crossings and is just a closed loop. This
space is a special one as it is the only prime space that has a blink presentation
without edges. By the rules we used on the construction of set $U$ this space
needed not to appear once: (1) we did not include blinks without edges and
(2) only one minimal presentation of a space was asserted to appear. This
space should be included artificially after. In spite of that, $\IS^2 \times \IS^1$
appeared as class 6.5. Figure~\ref{fig:gblinksForS1xS2inU} shows a blink
and a BFL for the 36 g-blinks in class 6.5. In a strict sense, this class
could be named 0.1 and spaces 6.6 to 6.19 would be decreased by one to 6.5
to 6.18, but we do not do this.

\begin{figure}[h!tp]
   \begin{center}
      \leavevmode
      \includegraphics[width=15cm]{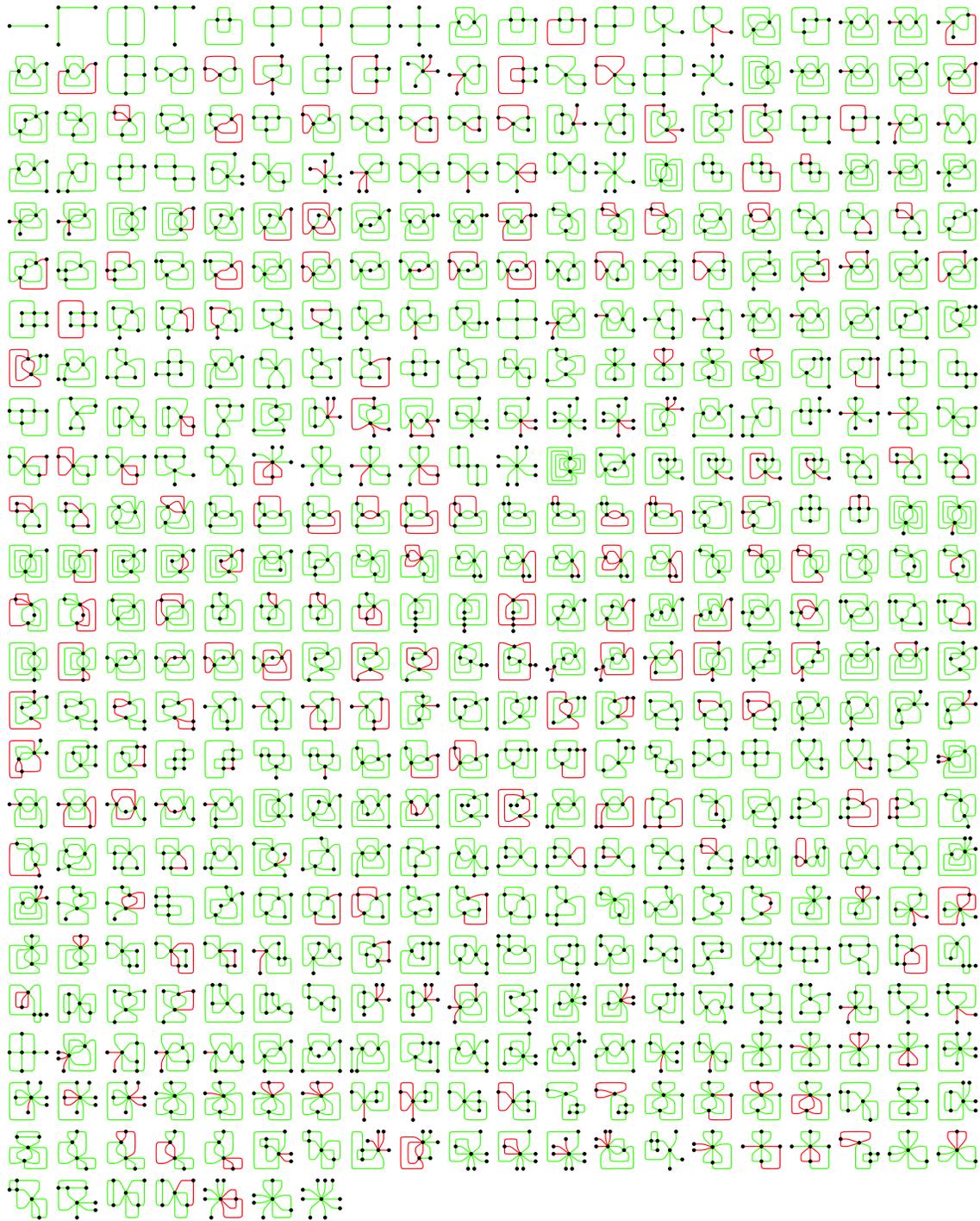}
   \end{center}
   \vspace{-0.7cm}
   \caption{ List of 487 blinks that induce once any prime space (modulo orientation)
   that has a blink presentation with $\leq$ 9 edges}
   \label{fig:primeSpace487Representatives}
\end{figure}

\begin{Theo}
\label{theo:primeSpacesUpTo9Edges}
Any prime space that has a blink presentation with $\leq$ 9 edges induces the
same space (modulo orientation) as one and only one of the 487 blinks in
Figure~\ref{fig:primeSpace487Representatives} or the 487 BFLs in
Figure~\ref{fig:primeSpace487RepresentativeLinks} or the 487 spaces shown
in Appendix~\ref{chap:primeCatalogue}.
\end{Theo}

\enlargethispage{2cm}

\begin{proof}
The construction of set $U$ asserts that it contains at least one minimal
g-blink for each prime space except for space $\IS^2 \times \IS^1$,
which is a special case where its minimal blink presentation
has no edges: only a single vertex. In spite of that $\IS^2 \times \IS^1$
appears in $U$ as class 6.5 so any prime space is included. The proof
that there are only 487 (with 2 doubts) is in the program \textsc{Blink}.
\end{proof}

\newpage

\newpage

\begin{figure}[h!tp]
   \begin{center}
      \leavevmode
      \includegraphics[width=15cm]{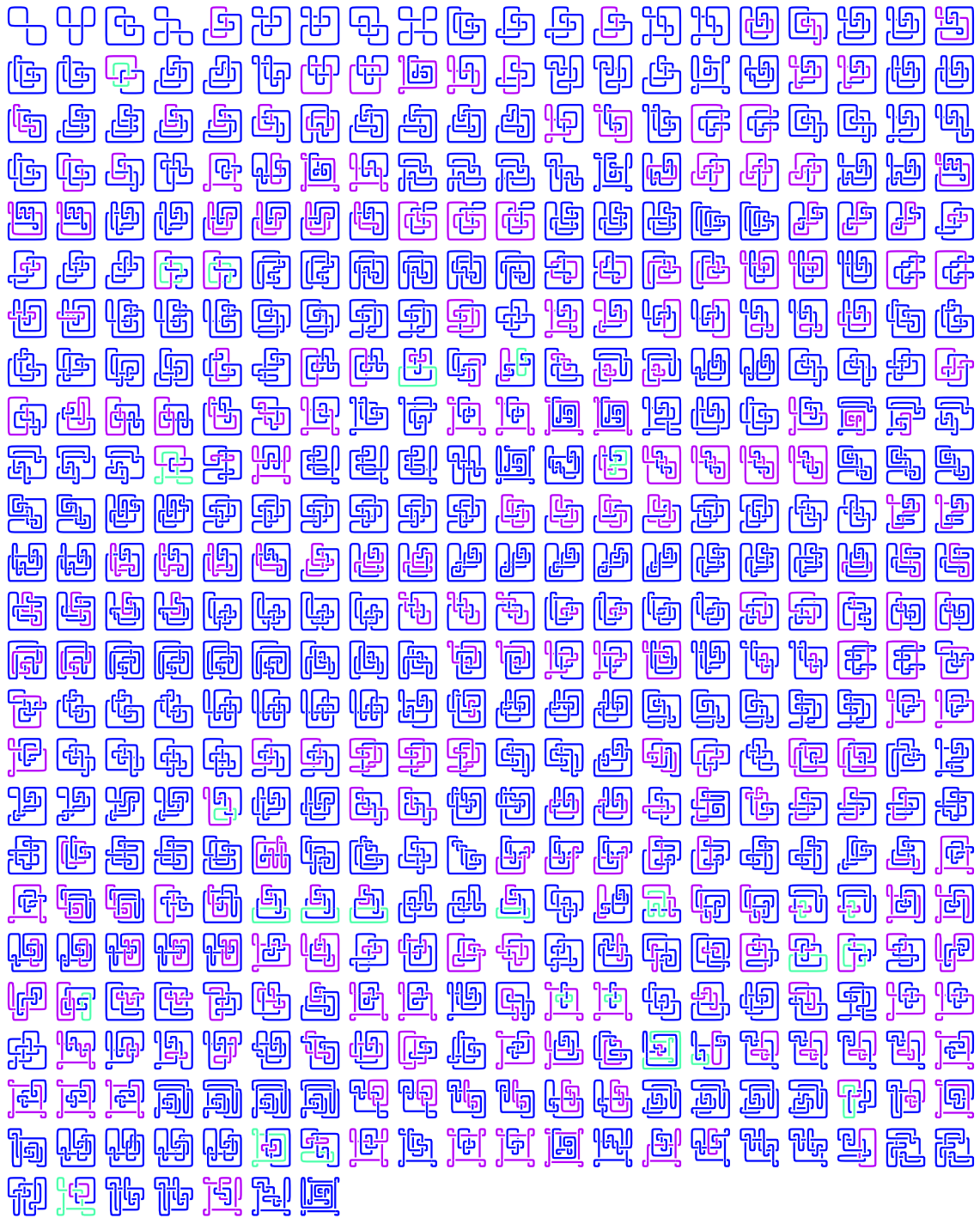}
   \end{center}
   \vspace{-0.7cm}
   \caption{ List of 487 BFLs that induce once any prime space (modulo orientation)
   that has a BFL presentation with $\leq$ 9 crossings}
   \label{fig:primeSpace487RepresentativeLinks}
\end{figure}

\newpage

\section{Spaces induced by simple 3-connected monochromatic blinks}

Blinks bring to the stage a very interesting connection: spaces and
plane graphs. Can concepts of graph theory when interpreted in space
language bring light to some unknown aspect of spaces? Some new
invariant for spaces?

With this spirit, what can we say about the space of a blink that is
$k$-connected? In Chapter~\ref{chap:blinks} we saw that the blocks
(2-connected pieces) of a blink may be recombined in different
ways leading to the same space. What are these blocks?
In this crude form, this concept of block or more general
$k$-connected blink does not mean something useful in
the language of spaces because of the following
observation: using the $B_2$ move of blink calculus
(\ie $RM_2$ in BFL calculus) explained
on Section~\ref{sec:blinkCalculus} one
may obtain blinks with higher connectivity inducing
the same space. But this comes
at a price, these equivalent versions with higher
connectivity contains local simplifications (moves
that reduce the number of edges) that leads back to
the first blink we started. A family of blinks
that do not contain these local simplifications are the
monochromatic blinks. Note that all simplification moves
on the blink calculus shown in Section~\ref{sec:blinkCalculus}
are, except for $B_4(1)$, from pieces with two colors. When
talking about blocks or higher connected monochromatic
blinks there is no local simplification at all.
So, the connectivity issue on monochromatic blinks
might mean something on spaces.

\begin{figure}[htp]
   \begin{center}
      \leavevmode
      \includegraphics{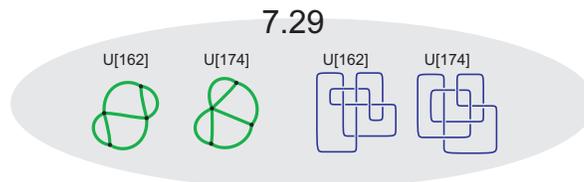}
   \end{center}
   \vspace{-0.7cm}
   \caption{ Non-trivial pair of green blocks (2-connected blinks) inducing the same space}
   \label{fig:nonTrivialGreenBlocks}
\end{figure}

Let $B$ be a green (all edges green) blink, $B'$ be a green blink
whose map (plane graph) is the dual map of $B$, $B''$ be $B$
reflected on the plane and $B'''$ be $B'$ reflected on the plane.
As we saw in Chapter~\ref{chap:blinks} all these blinks induce the
same space modulo orientation. Let's denote as {\it trivial}
a pair of blinks $A$ and $B$ if they induce the same space modulo
orientation and $A \in \{B,B',B'',B'''\}$. A pair of blinks inducing
the same space modulo orientation that is not trivial is called
{\it non-trivial}. Are all pair of green blocks (2-connected blinks)
that induce the same space modulo orientation trivial?
No. Space $7.29$ in Appendix~\ref{chap:primeCatalogue}
has a counterexample. Figure~\ref{fig:nonTrivialGreenBlocks} shows
a pair of non-trivial green 2-connected blinks in space 7.29 that induces
the same space. By the fact that they  all induce the same space, there must exist
paths connecting these blinks (or BFLs) using the moves
on the blink calculus (or BFL
calculus). Can you find such a path? We found the path via gem theory.

\begin{figure}[htp]
   \begin{center}
      \leavevmode
      \includegraphics{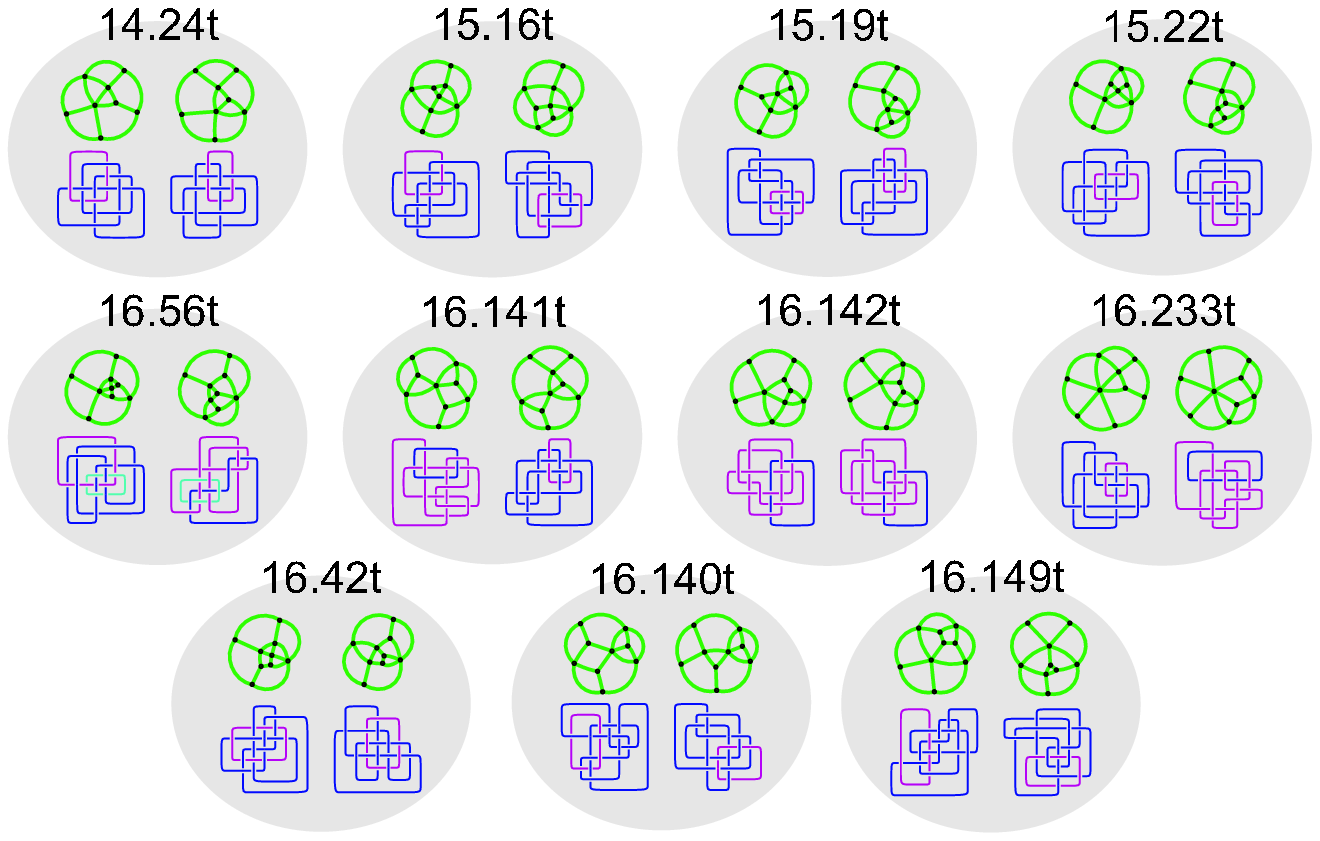}
   \end{center}
   \vspace{-0.7cm}
   \caption{ Doubts on simple 3-connected all green blinks}
   \label{fig:doubts3ConnectedIsolated}
\end{figure}

What about simple 3-connected monochromatic blinks?  Are there
non-trivial pairs of simple 3-connected monochromatic blinks?
To answer this question we generated a set named $T$ with all simple
3-connected green blinks up to 16 edges\footnote{To generate the simple
3-connected maps we started from the wheel maps (maps that are a polygons
with its vertices connected to a central vertex) and then, inductively, we
subdivided the faces and vertices in all possible ways preserving the
3-connectivity property.}
and calculated their HGnQI invariants (QI up to level 8).
The result was interesting. There are
708 simple 3-connected monochromatic g-blinks and they are divided in 381
classes HGnQI. These classes were named: $6.1$t, $8.1$t, $9.1$t,
10.1t$\ldots$ 10.2t, 11.1t $\ldots$ 11.2t, 12.1t $\ldots$ 12.9t,
13.1t $\ldots$ 13.11t, 14.1t $\ldots$ 14.36t, 15.1t $\ldots$ 15.76t
and 16.1t $\ldots$ 16.242t. This name convention is analogous to
the convention of the HGnQI classes in $U$ except for the
letter ``t'' at the end. These classes are presented with details in
Appendix~\ref{chap:catalolgue3con}. In these 381 classes there are
only 11 classes with exactly one non-trivial pair candidate. They
are shown in Figure~\ref{fig:doubts3ConnectedIsolated}.

\enlargethispage{2cm}

Is any of these a non-trivial pair or all of them induce different spaces
modulo orientation that the HGnQI could not capture? We leave this question
open.

%% file: chapter6.tex
\chapter{Conclusions and future work}
\label{chap:conclusion}

\section{Results, uncertainties and the need of new invariants}

\begin{figure}[htp]
   \begin{center}
      \leavevmode
      \includegraphics[width=10cm]{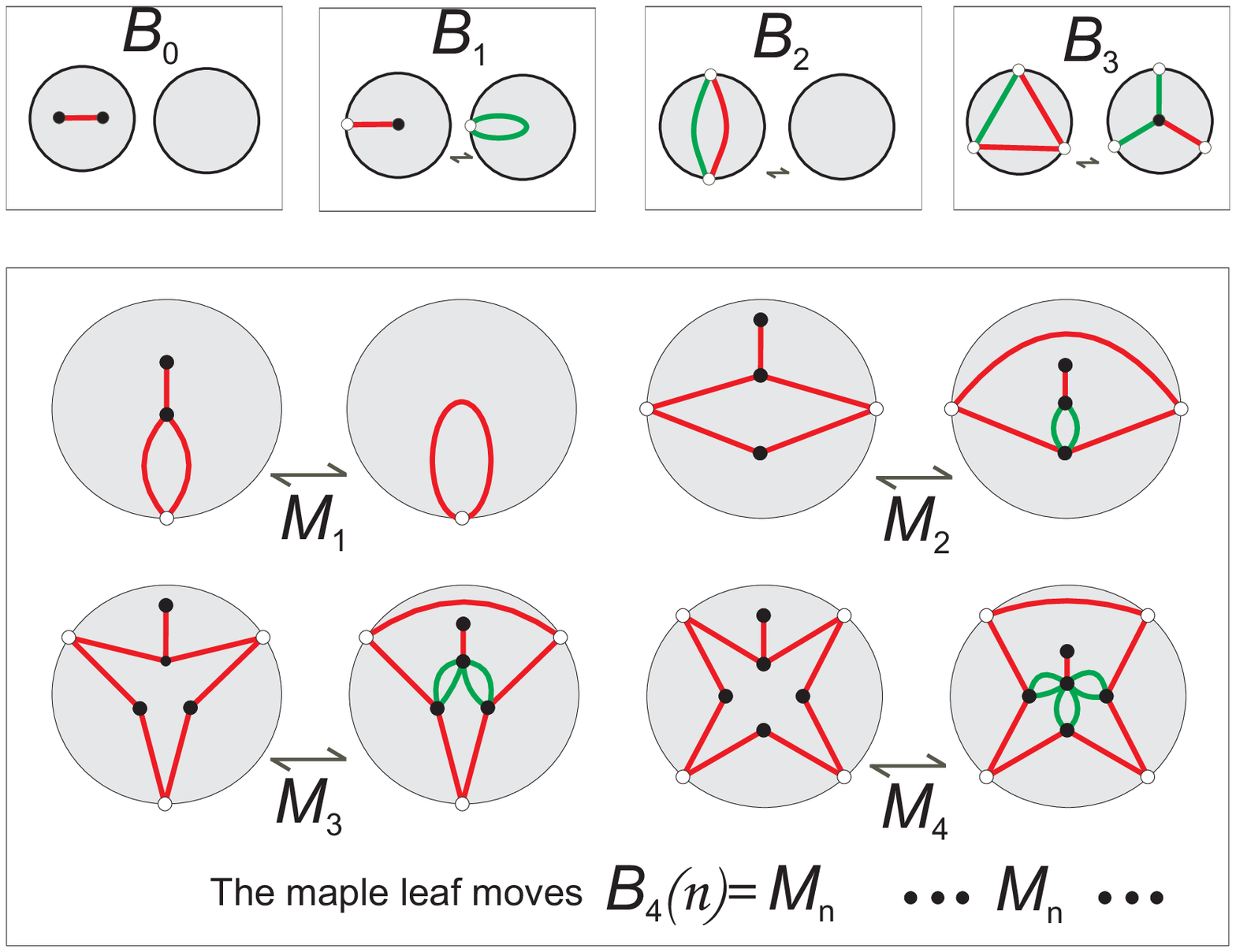}
   \end{center}
   \vspace{-0.7cm}
   \caption{Blink formal calculus by local coins replacements}
   \label{fig:blinkCalculusOnCoins2}
\end{figure}

\bigskip \centerline{\bf \textsc{Blink Calculus}} \bigskip
The first contribution of this thesis that we want to stress
here was given in Section~\ref{sec:blinkCalculus}. Based on the BFL
calculus (\ie Kirby's calculus reformulated in BFL language) we obtained
a purely blink calculus. This calculus is a formal language
which is a counterpart for homeomorphism of spaces. \hbox{Figure~\ref{fig:blinkCalculusOnCoins2}}
presents again our blink calculus. Although theoretically complete (it is supported by Kirby's Calculus) this calculus was not
used in our computational experiments as a tool to prove homeomorphisms.
For this task we used the combinatorial simplification dynamics of
3-gems. In spite of that, we think that the blink calculus can help
in the search for new space invariants.

\bigskip \centerline{\bf \textsc{Decomposition/Composition Theory}}
A second contribution of this work that was important to the
computational results are the following propositions and theorems
in g-blink language:

\noindent {{\bf (Theorem on partial dual \ref{theo:partialDual})} \,\, \it
Let $A$ and $B$ be arbitrary disjoint
g-blinks, $(a,b)$ a basepair on them. Then
 $A[a] + B[b] \Sequiv A[a] + \textsc{Dual}(B)[b]$.
}

\noindent {{\bf (Theorem on partial reflection \ref{theo:partialReflection})} \,\, \it
Let $A$ and $B$ be arbitrary disjoint g-blinks, $(a,b)$ a basepair
on them. Then($A[a] + B[b] \Sequiv A[a] +
\textsc{Reflection}(B)[b].$
}

\noindent {{\bf (Theorem on partial refDual \ref{theo:partialRefDual})} \,\, \it
Let $A$ and
$B$ be arbitrary disjoint g-blinks, $(a,b)$ a basepair on them. Then
$A[a] + B[b] \Sequiv A[a] + \textsc{RefDual}(B)[b].$
}

The third theorem on partial refDuals is obtained directly from the framed link theory.
Given this third theorem, the first two theorems are equivalent: given one we have the other.
The theorem on partial reflection was tricky to obtain. We
used both the theory of gems and the \textsc{Blink2Gem} algorithm as well as topological machinery
to exhibit an explicit homeomorphism.
Section~\ref{sec:proofOfThePartialReflectionTheorem} contains
this proof. These theorems yield a block decomposition/composition
theory which leads to the representative concept and curtailed search
spaces of our computational experiments.

\bigskip \centerline{\bf \textsc{An unavoidable set of blinks up to 9 edges}} \bigskip

We achieved our initial main objective which was to classify
spaces presentable by blinks with small number $n$ of edges.
At the level of $n \leq 9$ the combination of tools
\begin{itemize}
\item theory of decomposition/composition leading to representative g-blinks --- which reduces the search space
\item quantum invariants and homology --- which provide distinctiveness
\item combinatorial simplification dynamics of 3-gem theory --- which provides similarity
\end{itemize}
was as effective as leaving only two uncertainties in more than 500 spaces.
These uncertainties, as we saw in Section~\ref{sec:topologicalClassificationOfU}, were registered
as Conjecture\ref{conj:conjecture1}~and~Conjecture~\ref{conj:conjecture2}.
To be honest, these conjectures are actually doubts and
seem an interesting research problem. It could be answered by a new invariant
which complements the HGnQI invariant. In any case (\ie the two
conjectures are false, or one is true and the other is false, or both are true) the
relevant fact is that any space that has a blink presentation with up to 9
edges is induced by only one of the classes in Appendix~\ref{chap:primeCatalogue},
where classes 9.126 and 9.199 may be broken into two classes each. A space
that is not prime and has a blink presentation with $\leq 9$ edges is
just a blink with more than one prime component which is in the catalogue
(Section~\ref{sec:census}).

\begin{figure}[h!tp]
   \begin{center}
      \leavevmode
      \includegraphics[width=16cm]{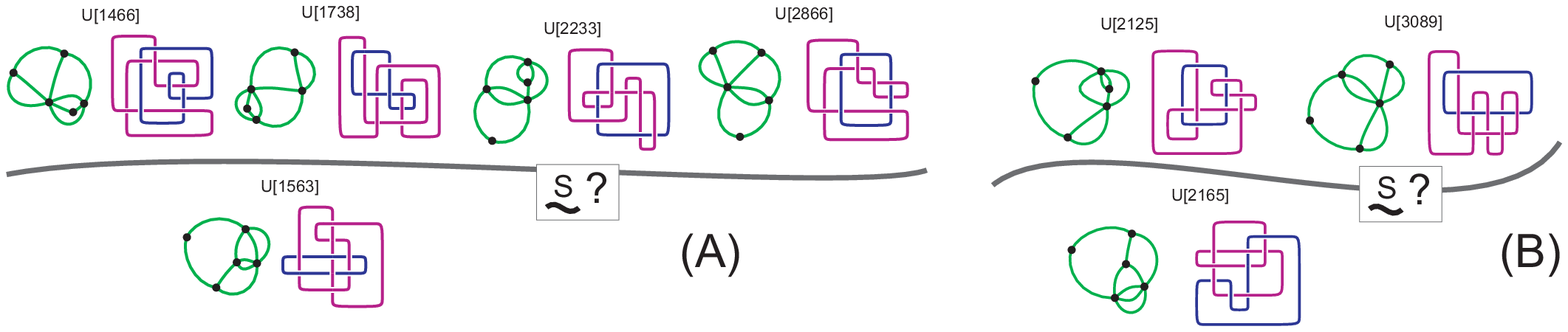}
   \end{center}
   \vspace{-0.7cm}
   \caption{ The only 2 classes with same HGnQI where a proof of the homeomorphism was not found}
\end{figure}

\begin{figure}[h!tp]
   \begin{center}
      \leavevmode
      \includegraphics{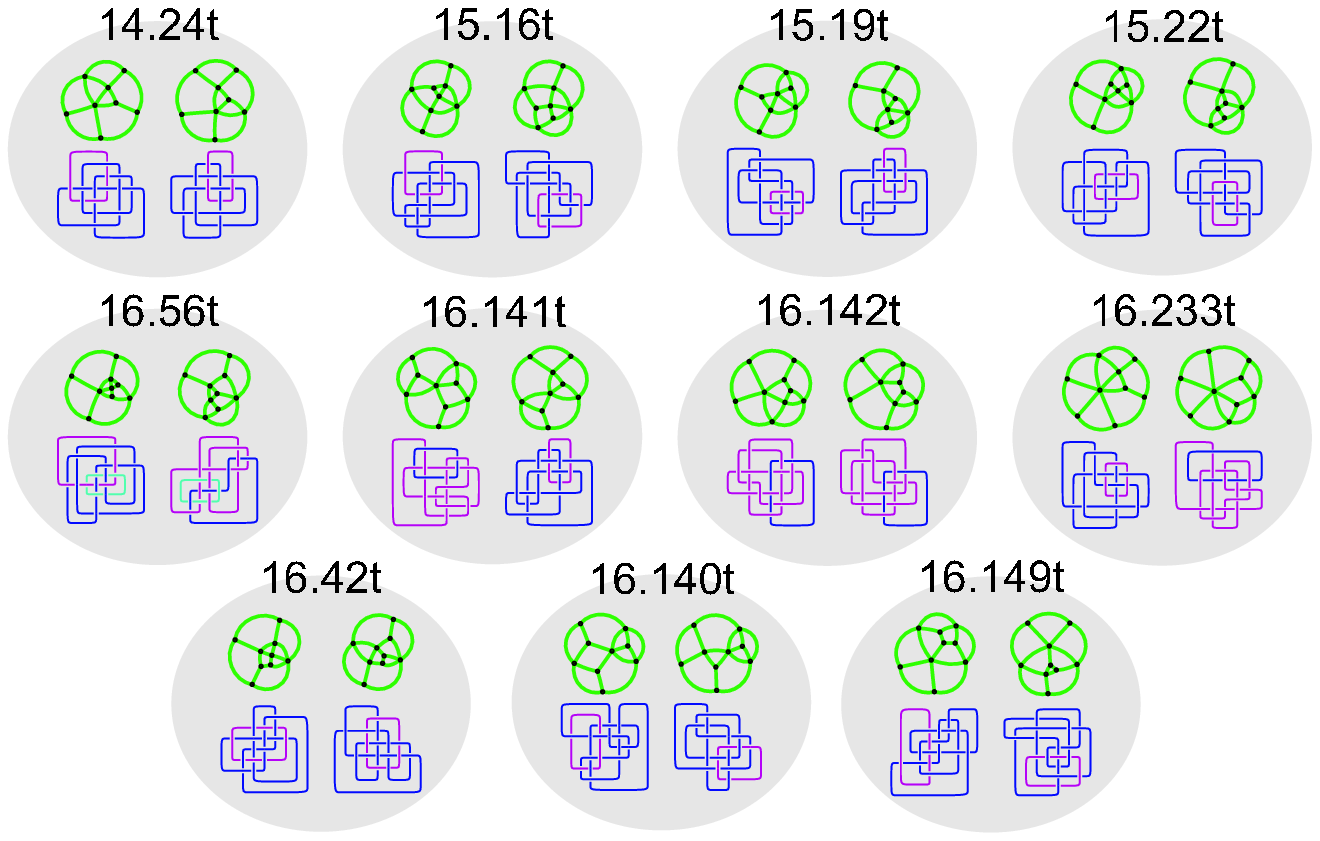}
   \end{center}
   \vspace{-0.7cm}
   \caption{ Doubts on simple 3-connected all green blinks}
   \label{fig:doubts3ConnectedIsolated2}
\end{figure}

In the case of simple 3-connected monochromatic blinks with up to 16
edges there are only the 11 uncertainties shown in Figure~\ref{fig:doubts3ConnectedIsolated2}.
We did not use the simplification combinatorial dynamics of 3-gems
to deal with these cases.

Putting together blinks that induce the same space in a non-trivial
way (Appendix~\ref{chap:primeCatalogue},
Appendix~\ref{chap:compositeCatalogue} and
Appendix~\ref{chap:catalolgue3con}) we hope to
be contributing with non-trivial examples that can
motivate and help the search for new effective
subtle invariants of spaces to complement the
HGnQI invariant.

\newpage

\section{The inverse algorithm: from gem to blink}

A rather frustrating fact up to now is that we could not find a blink
for the space ${\rm EUCLID}_1$. This space is generated by the rigid
gem $r_{5}^{24}$ (notation of 3-gems catalog of \cite{Lins1995}). Blinks
and BFLs for the other euclidean spaces are given below.
\begin{center}
\includegraphics[width=14cm]{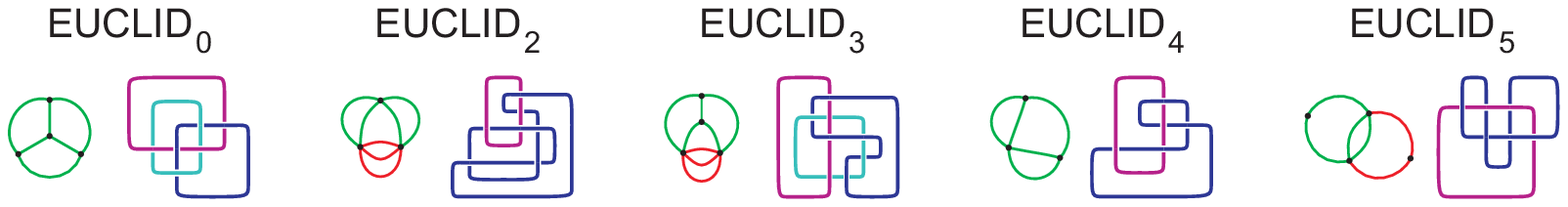}
\end{center}
They correspond, respectively, to spaces 6.8, 7.10, 8.32, 5.4 and 6.13. By looking at quantum
invariants of these spaces (see Appendix~\ref{chap:primeCatalogue}) we are
led to the following conjecture.
\begin{Conj}
The absolute value of the quantum invariants of the euclidean
spaces are non-negative integers for all levels $r$.
\end{Conj}
The missing ${\rm EUCLID}_1$ space motivates the
following discussion.

\begin{figure}[htp]
   \begin{center}
      \leavevmode
      \includegraphics[width=8cm]{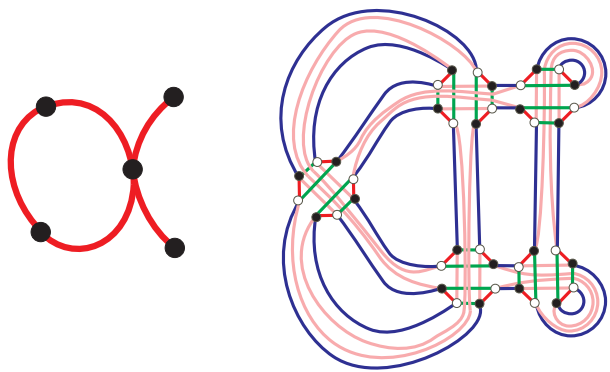}
   \end{center}
   \vspace{-0.7cm}
   \caption{ Blink to gem algorithm: indispensable to prove homeomorphisms of blinks}
   \label{fig:blink2GemAlgorithm}
\end{figure}

There exists a rather simple algorithm to go from a framed link
inducing a space to a triangulation of the same space. This was
first done in chapter 11 of~\cite{KauffmanAndLins1994} via {\em
graph encoded 3-manifolds} or {\em gems}. This algorithm was
improved here in Section~\ref{sec:fromGBlinkTo3Gem} and it is
a central tool in the \textsc{Blink} program to prove spaces
in $U$ are homeomorphic. Figure~\ref{fig:blink2GemAlgorithm}
shows this algorithm. Thus to get a gem from a blackboard
framed link is a direct task.

\begin{figure}[htp]
   \begin{center}
      \leavevmode
      \includegraphics[width=8cm]{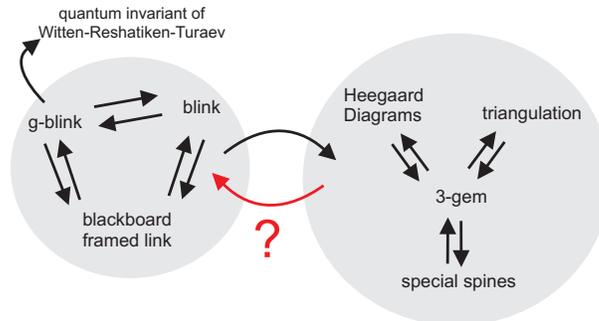}
   \end{center}
   \vspace{-0.7cm}
   \caption{ Blink based presentation and 3-Gem based presentation}
   \label{fig:spacePresentationsSchema}
\end{figure}

However, the contrary, given a gem to find by
a polynomial algorithm a blackboard framed link inducing the same
3D-space is, as far as we know, an untouched problem in the
literature. Figure~\ref{fig:spacePresentationsSchema} shows this computational gap as a red arrow.
The reason why it is desirable to have this arrow in black
stems from the fact that the quantum invariants are not computable
from a triangulation or gem based presentation of 3D-spaces. The two
languages, triangulations and blackboard framed links have at
present only a one way translation.

Trying to get this converse algorithm took a long a time of our research for this thesis.
Only recently we got confident that we have succeeded. A first step in this direction was given
in the paper~\cite{Lins2007}, where a linear algorithm to prove the
Lickorish-Wallace Theorem is provided. The second part, which actually
presents the blink from the gem
is a joint work with S. Lins \cite{LinsLins2007} and
awaits a proper computer implementation. The first test of
this implementation will be to get a blink for ${\rm EUCLID}_1$.

\section{The \textsc{Blink} computer program}

A computer program to manipulate spaces through its many
possible presentations was one of our goals in this work.
Indeed, a great effort was made to bring \textsc{Blink}
to life: a program written in Java that, at this moment,
has more than 800 hundred classes and more than 70000
lines of code. Today, \textsc{Blink} supports
blinks, g-blink, BFLs and 3-gems. The idea is, in
the future, to bring other possible space presentations,
like special spines, into it.

To make \textsc{Blink} a flexible program we decided that its
interface would be a {\it Command Line Interface}. It displays
a prompt and the user enters a command or a script written in
a {\it language} that we also name \textsc{Blink}. Once a
command or script has been entered, the program calculates
the script result and shows it to the user. The flexibility we
get in this type of design is good; for example we can
combine functions and easily express more complex functions.

Besides the calculation of invariants, the identification of certain
structures into 3-gems (\eg disconnecting quartets, dipoles) or
into g-blinks (\eg simplification points) one of the main characteristics of
\textsc{Blink} is its capability of presenting drawings or
diagrams for blinks, g-blink, BFLs and 3-gems. Almost all
drawings on this thesis came from \textsc{Blink}. To get
good looking and correct drawings for blinks,
g-blinks and BFLs took us a long time once we didn't
know a good way of doing it. But finally we found a
great solution: {\em Tamassia's Algorithm \cite{Tamassia1987}}.

We have implemented the following four algorithms to
deal with the drawing issue. Except for the first algorithm,
the other three are further fine-tuned with {\em Bézier curves and
splines techniques} (\cite{FoDaFeHu1990}) to produce
rounded-drawings with curved edges.

\begin{enumerate}
\item {\em Coin-drawing Algorithm:} this was our own first original
algorithm which we implemented to correctly draw in a visible scale
the whole of any plane graph. The drawing is in the interior of a
disk named a {\em coin}. The coin-drawing algorithm chooses and draw
a spanning tree of the graph with appropriate lengths and angles.
These ensure that the remaining edges can be displayed as a path
which is a line segment, an arc of circle and another line segment.
This is the simplest algorithm producing the less pleasing
aesthetical effect. Nevertheless, these {\em coin-drawings} are
important because they were, for a long time in our work, the only
general method with total visibility. Figure~\ref{fig:coinDrawing} presents
an example of our coin drawing algorithm.

\begin{figure}[htp]
   \begin{center}
      \leavevmode
      \includegraphics{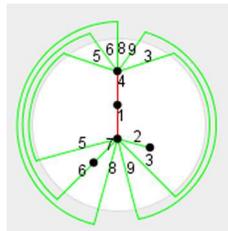}
   \end{center}
   \vspace{-0.7cm}
   \caption{ Coin drawing of $U[1078]$}
   \label{fig:coinDrawing}
\end{figure}

\item {\em Tutte's Barycentric Algorithm \cite{Tut1967}:}
we have implemented this well known algorithm that draws a
$3$-connected plane graph by choosing the external face and
extending the drawing so that every interior vertex is in the
barycenter of its neighbors. Frequently, it produces pleasant
drawings. However it does not treat the less connected graphs which
are central for our work: loops, pendant vertices and cut-vertices
are of fundamental importance in blink theory. Another problem that
occurs with Tutte's based algorithm is the one of discrepant scales:
some parts of the drawing are exponentially smaller that others, and
simply disappear from the drawings. Despite of these
disadvantages, Tutte's algorithm works well for the majority of
blinks in the set $U$.

\item {\em Koebe, Andreev
and Thurston's Theorem on circle packing in the hyperbolic plane:}
beautiful drawings of plane graphs are possible to obtain from the
geometry of the hyperbolic plane. Given a $3$-connected plane graph,
there exist circles centered at the vertices of the graph so that
the edges are defined by the contact points of two circles. See the
articles of Smith~\cite{Smith1994}, Stephenson~\cite{Ste2003} of
Collins and Stephenson~\cite{CoSte2003} where algorithms are
outlined for the case of triangulations.  The Theorem yielding the
circle packing was proved independently by Koebe \cite{Koebe1936},
Andreev \cite{Andreev1970} and by Thurston \cite{Thurston1982}.
We have implemented our own version of the algorithm which works in
the case of $3$-connected graphs. However, it suffers the same
disadvantages as Tutte's Algorithm. Nevertheless, when it works it
produces the nicest results. Figure~\ref{fig:circlePacking} presents a
blink, its BFL, the circle packing that defined the first two drawings,
and the first three drawings together.

\begin{figure}[h!tp]
   \begin{center}
      \leavevmode
      \includegraphics[width=14cm]{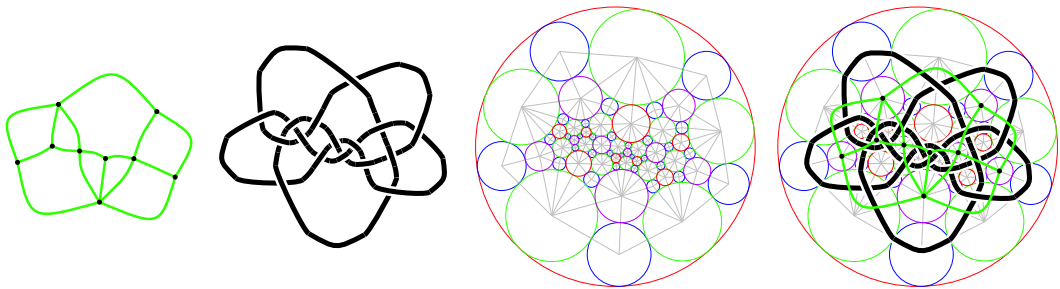}
   \end{center}
   \vspace{-0.7cm}
   \caption{ Circle packing of a 3-connected blink}
   \label{fig:circlePacking}
\end{figure}

\item {\em Tamassia's Algorithm \cite{Tamassia1987}: to embed an arbitrary
plane graph with valency at most 4 in the rectilinear grid so as to
minimize the number of bends.} This algorithm came to our attention
only at latter phase of our research. Fine-tunings of it has all the
properties we needed: it correctly draws any plane graph and it does
not suffer from the undesirable phenomenon of discrepant scales:
all the vertices and the edges are entirely well visible. The
objective of the method is to minimize (in a precisely defined
mathematical way) the number of bends. The algorithm depends three
times on the algorithm to compute a minimum cost-flow in a network.
The essence of this algorithm is in the design of the first network.
Each feasible flow in this network encodes valid ``shapes'' for the edges
of the graph. This encoding tells, for example, that an edge has no
bends or that it has one bend to the right then one bend to the left.
When a minimum flow is found in this network the minimum number of
bends for a valid rectilinear embedding of the given graph is found. The
other two networks are used to find the lengths of the horizontal
end vertical segments of the edges. In order to have Tamassia's
algorithm available, we had first to implement the minimum cost flow
algorithm in its full generality via the network simplex method. We
have based our implementation on the lucid exposition Chapter 19
of \cite{Chvatal1983} and in the
(editor's categorization) ``Exceptional Paper'' \cite{BrBrGr1977},
which is the original source of the network simplex method.
Tamassia's algorithm is an unexpected application of network flow
theory in its full strength. Since its publication in 1987, it has
become a theoretically beautiful at the same time a practical
device used on dozens of applications. Having implemented it from
scratch, we had the opportunity of tailoring it to fulfill our
expectations on drawing of general plane graphs. In particular, the
restriction about the maximum degree 4 is easy to overcome. Finally,
the use of Bézier curves and splines (\cite{FoDaFeHu1990}) makes the
drawings more pleasing aesthetically, with smaller perceptual
complexity. All drawings in the Appendices are based in this algorithm.
\end{enumerate} \enlargethispage{1cm}

We want to make \textsc{Blink} an open source project on the
internet but this wasn't done yet.

\newpage

\centerline{\bf \textsc{An example of Blink usage}} \bigskip

We finish this section with an example of the usage of
\textsc{Blink}. Here is the code:

\begin{center}
\begin{minipage}{14cm}
\begin{footnotesize}
\setstretch{1.1}
\begin{verbatim}
// associates B to the g-blink with 4 parallel green edges: U[5]
B = gblink(5)
// all possible toroidal sums with B up to 24 edges
C = combineGBlinks({B},24)
// calculate representative
C = rep(C)
// remove duplicates
C = set(C)
// homology groups on all g-blink of C
HGs = hg(C)
// calculate quantum invariant on all g-blink of C up to level 4
QIs = qi(C,4)
// produces the blink drawings
db(C,cols=10,rows=4,eps="blinks.eps")
// produces the link drawings
dl(C,cols=10,rows=4,eps="links.eps")
\end{verbatim}
\end{footnotesize}
\end{minipage}
\end{center}

\begin{figure}[h!tp]
   \begin{center}
      \leavevmode
      \includegraphics{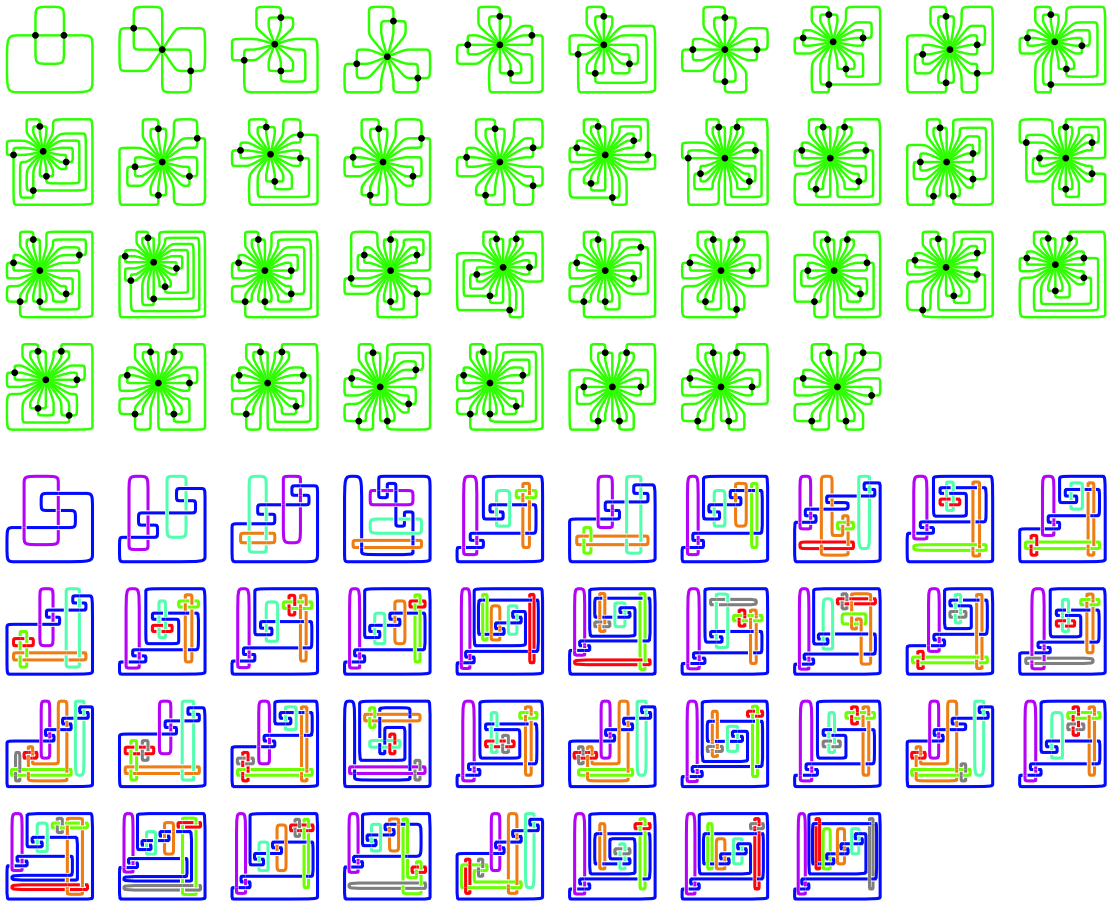}
   \end{center}
   \vspace{-0.7cm}
   \caption{Toroidal sums or g-blink merges up to six copies of the quaternionic space}
   \label{fig:blinkLanguageExample}
\end{figure}

The BFL presentation of $U[5]$ is the first in $U$
to have two components. It is a blink with four parallel
green edges. By merging $U[5]$ with itself in all possible
ways we obtain 38 representative blinks with $\leq$
24 edges. These 38 blinks and their associated
BFLs are shown in Figure~\ref{fig:blinkLanguageExample}.
We calculated homology group and  the quantum invariant
up to level 4 of these 38 blinks. We could distinguish
24 spaces. The blinks we cannot distinguish with this
experiment are:
$\{ 9, 12\}$,
$\{10, 13\}$,
$\{16, 27, 36\}$,
$\{17, 19, 25, 28\}$,
$\{18, 21, 26, 29, 31\}$,
$\{20, 24\}$,
$\{23, 30\}$,
$\{32, 34\}$.

\def\rsep{-0.15cm}

\begin{center}
\begin{scriptsize}
\begin{tabular}{cccc}
\# & HG & r=3 & r=4 \\[\rsep]
01 & $2^2    $ & $1.41421 + 0.00000i$ & $ 1.50000 -  0.50000i$ \\[\rsep]
02 & $(1)\, 2^2$ & $2.00000 + 0.00000i$ & $ 2.00000 -  1.00000i$ \\[\rsep]
03 & $2^4    $ & $2.82842 + 0.00000i$ & $ 3.00000 -  2.00000i$ \\[\rsep]
04 & $(2)\, 2^2$ & $2.82842 + 0.00000i$ & $ 3.00000 -  1.00000i$ \\[\rsep]
05 & $(1)\, 2^4$ & $4.00000 + 0.00000i$ & $ 5.00000 -  3.00000i$ \\[\rsep]
06 & $(1)\, 2^4$ & $4.00000 + 0.00000i$ & $ 4.00000 -  4.00000i$ \\[\rsep]
07 & $(3)\, 2^2$ & $4.00000 + 0.00000i$ & $ 6.00000 +  0.00000i$ \\[\rsep]
08 & $(2)\, 2^4$ & $5.65685 + 0.00000i$ & $ 8.00000 -  6.00000i$ \\[\rsep]
09 & $2^6    $ & $5.65685 + 0.00000i$ & $ 7.00000 -  7.00000i$ \\[\rsep]
10 & $(2)\, 2^4$ & $5.65685 + 0.00000i$ & $ 6.00000 -  6.00000i$ \\[\rsep]
11 & $2^6    $ & $5.65685 + 0.00000i$ & $ 5.00000 -  7.00000i$ \\[\rsep]
12 & $2^6    $ & $5.65685 + 0.00000i$ & $ 7.00000 -  7.00000i$ \\[\rsep]
13 & $(2)\, 2^4$ & $5.65685 + 0.00000i$ & $ 6.00000 -  6.00000i$ \\[\rsep]
14 & $(2)\, 2^4$ & $5.65685 + 0.00000i$ & $10.00000 -  4.00000i$ \\[\rsep]
15 & $(4)\, 2^2$ & $5.65685 + 0.00000i$ & $14.00000 +  2.00000i$ \\[\rsep]
16 & $(1)\, 2^6$ & $8.00000 + 0.00000i$ & $14.00000 - 12.00000i$ \\[\rsep]
17 & $(1)\, 2^6$ & $8.00000 + 0.00000i$ & $ 8.00000 - 12.00000i$ \\[\rsep]
18 & $(1)\, 2^6$ & $8.00000 + 0.00000i$ & $10.00000 - 12.00000i$ \\[\rsep]
19 & $(1)\, 2^6$ & $8.00000 + 0.00000i$ & $ 8.00000 - 12.00000i$ \\[\rsep]
\end{tabular}
\begin{tabular}{cccc}
\# & HG & r=3 & r=4 \\[\rsep]
20 & $(1)\, 2^6$ & $8.00000 + 0.00000i$ & $ 8.00000 - 14.00000i$ \\[\rsep]
21 & $(1)\, 2^6$ & $8.00000 + 0.00000i$ & $10.00000 - 12.00000i$ \\[\rsep]
22 & $(1)\, 2^6$ & $8.00000 + 0.00000i$ & $ 6.00000 - 12.00000i$ \\[\rsep]
23 & $(1)\, 2^6$ & $8.00000 + 0.00000i$ & $ 8.00000 - 10.00000i$ \\[\rsep]
24 & $(1)\, 2^6$ & $8.00000 + 0.00000i$ & $ 8.00000 - 14.00000i$ \\[\rsep]
25 & $(1)\, 2^6$ & $8.00000 + 0.00000i$ & $ 8.00000 - 12.00000i$ \\[\rsep]
26 & $(1)\, 2^6$ & $8.00000 + 0.00000i$ & $10.00000 - 12.00000i$ \\[\rsep]
27 & $(1)\, 2^6$ & $8.00000 + 0.00000i$ & $14.00000 - 12.00000i$ \\[\rsep]
28 & $(1)\, 2^6$ & $8.00000 + 0.00000i$ & $ 8.00000 - 12.00000i$ \\[\rsep]
29 & $(1)\, 2^6$ & $8.00000 + 0.00000i$ & $10.00000 - 12.00000i$ \\[\rsep]
30 & $(1)\, 2^6$ & $8.00000 + 0.00000i$ & $ 8.00000 - 10.00000i$ \\[\rsep]
31 & $(1)\, 2^6$ & $8.00000 + 0.00000i$ & $10.00000 - 12.00000i$ \\[\rsep]
32 & $(3)\, 2^4$ & $8.00000 + 0.00000i$ & $14.00000 - 10.00000i$ \\[\rsep]
33 & $(3)\, 2^4$ & $8.00000 + 0.00000i$ & $12.00000 -  8.00000i$ \\[\rsep]
34 & $(3)\, 2^4$ & $8.00000 + 0.00000i$ & $14.00000 - 10.00000i$ \\[\rsep]
35 & $(3)\, 2^4$ & $8.00000 + 0.00000i$ & $ 8.00000 -  8.00000i$ \\[\rsep]
36 & $(1)\, 2^6$ & $8.00000 + 0.00000i$ & $14.00000 - 12.00000i$ \\[\rsep]
37 & $(3)\, 2^4$ & $8.00000 + 0.00000i$ & $22.00000 -  6.00000i$ \\[\rsep]
38 & $(5)\, 2^2$ & $8.00000 + 0.00000i$ & $32.00000 +  4.00000i$ \\[\rsep]
\end{tabular}
\end{scriptsize}
\end{center}
\enlargethispage{2cm}
Observe the curious fact that, except for the first blink, the quantum
invariants at level $r=4$ are Gauss integers (\ie $a + bi$ with $a,
b$ integers). This type of experiment is very easy to do with \textsc{Blink}.

\newpage

\section{Two final remarks}

\enlargethispage{1cm}

First. Recently we have extended the $U$ set to blinks with up to 10
edges. The number of blinks was increased from 3437 to 17948. The
number of potentially prime classes increased from 487 to 1025. The
number of composite classes increased from 14 to 40. We did not
attempt the topological classification of the classes 10.\_\_ using
3-gems.

Second. We have a contract with World Scientific Publisher to write
a book to be co-authored by S. Lins based on the material of this
thesis. The tentative title of this book: {\it All Shapes of Spaces:
a Genealogy of Closed Oriented 3-Manifolds} and it should be
finished by the year 2008.

%% file: appendix1.tex
\chapter{The 487 potentially prime spaces in $U$}
\label{chap:primeCatalogue}

We here present the 487 spaces that are ``potentially prime'' once we could
not prove them composite in our tests. One thing is certain, as stated in
Theorem~\ref{theo:primeSpacesUpTo9Edges}: any prime space that can be
presented as a blink with $\leq$ 9 edges induces the same space (modulo
orientation) as one and only one of these 487 spaces. Actually there
are two points where this last statement may fail: space 9.126 and space 9.199 (although
they have the same HGnQI we could not find a proof of homeomorphism
between g-blink $U[1563]$ and the other g-blinks in 9.126 and g-blink $U[2165]$
and the others in 9.199). All 3437 g-blinks in $U$ appears in this Appendix or
in Appendix~\ref{chap:compositeCatalogue}.

\begin{figure}[htp]
   \begin{center}
      \leavevmode
      \includegraphics[width=14.5cm]{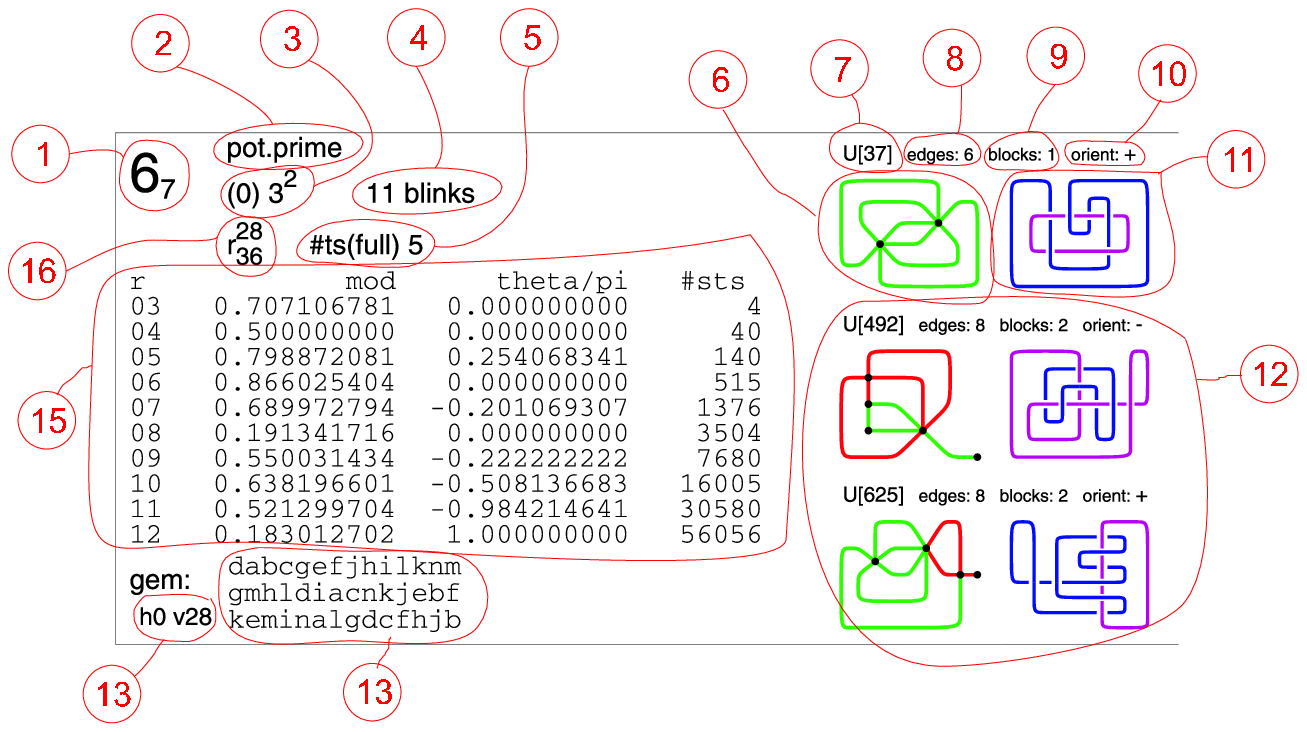}
   \end{center}
   \vspace{-0.7cm}
   \caption{ Elements of catalogue}
   \label{fig:catalogueExplanation}
\end{figure}

The elements of this catalogue are: (1) the space name: $6_7$ is
a synonym for $6.7$; (2) the primality test outcome; (3) the homology
group; (4) the number of g-blinks in $U$ that induces this space; (5) number
of 3-gems identified in the same ts-class of the minimum 3-gem
found for this space: {\it full} means that all ts-class
was identified, {\it partial} means that we do not know if
all ts-class was identified; (6) the minimum blink
presentation for this space in set $U$ and also a
minimal presentation for this space (this is always true,
except for class $6.5$ that should be $0.1$); (7) the name
of the g-blink in $U$; (8) its number of edges; (9) its number
of blocks in the blink presentation (2-connected components);
(10) its orientation compared to the orientation of the QI shown:
+ sign means the same and - sign means different; (11) the
corresponding BFL presentation; (12) other g-blinks in the
same space; (13) the code of the minimal 3-gem found for
this space the code convention is defined in \cite{Lins1995};
 (14) the number of handles (composition with
$\IS^2 \times \IS^1$ and the number of vertices of this 3-gem); (15)
the quantum invariant of this space in polar form where the
angle is divided by $\pi$; (16) the name of this minimal
3-gem in the catalogue of \cite{Lins1995} when it is present
in this catalogue.

The spaces that have integral quantum invariants up to level 12 are:
6.5 $(\IS^2 \times \IS^1)$, 6.8, 6.18 and 8.32. The spaces that have real but not integral
quantum invariant up to level 12 are 1.1, 2.1, 4.4, 6.14, 6.19, 8.58,
8.70, 8.75, 8.76, 8.81, 8.86, 8.87, 8.89, 8.100, 8.102, 8.103, 8.117,
9.23, 9.183. The remaining classes have entries with non-zero imaginary
part (\ie $\theta/\pi \notin \{0, 1\}$).

\newcount\ii \newcount\jj   
\def\producePages#1#2{
\ii=#1                      
\jj=#2                      
\advance\jj by 1            
\loop   
   \ifnum\ii<\jj
{
   \hspace{-1.8cm}
   \enlargethispage{5cm}
   {\centering
   \includegraphics[width=18cm]{fig/catalog\ifnum\ii<100 0\fi\ifnum\ii<10 0\fi\number\ii.eps}
   }
   \newpage}
      \advance\ii by 1
   \repeat
}

\newpage
\setlength{\topmargin}{-1.2cm}

\producePages{1}{101}

%% file: appendix2.tex
\chapter{The 14 composite spaces in $U$}
\label{chap:compositeCatalogue}

We here present the 14 spaces induced from g-blinks in $U$.
Their ``connected sum'' details: what prime spaces compose to
them are shown in Chapter~\ref{chap:census}. The elements of
this presentation are the same as the explained in
Appendix~\ref{chap:primeCatalogue}.

\newcount\ii \newcount\jj   
\def\producePagesTwo#1#2{
\ii=#1                      
\jj=#2                      
\advance\jj by 1            
\loop   
   \ifnum\ii<\jj
{
   \hspace{-1.8cm}
   \enlargethispage{5cm}
   {\centering
   \includegraphics[width=18cm]{fig/compcatalog\ifnum\ii<100 0\fi\ifnum\ii<10 0\fi\number\ii.eps}
   }
   \newpage}
      \advance\ii by 1
   \repeat
}

\newpage
\setlength{\topmargin}{-1.2cm}

\producePagesTwo{1}{3} 

%% file: appendix3.tex
\chapter{Simple 3-connected monochromatic blinks up to 16 edges}
\label{chap:catalolgue3con}

We here present all simple 3-connected green blinks with
$\leq 16$ edges divided in 381 HGnQI classes. The quantum
invariant was calculated up to level $r=8$ for each of
these blinks. There are left 11 uncertainties: 14.24t,
15.16t, 15.19t, 15.22t, 16.42t, 16.56t, 16.141t, 16.142t,
16.149t, 16.233t. Except for these classes the other 370
consisted of only one blink (or the two orientations of
the same space). This fact suggests that if $A$ and $B$
are two different simple 3-connected monochromatic blinks,
that do not form a trivial pair (trivially induce the
same space), then they probably induce different spaces.
Are the 11 uncertainties examples of non-trivial pairs?

\begin{center}
\includegraphics{fig/doubts3connectedisolated.eps}
\end{center}

\newcount\ii \newcount\jj   
\def\producePagesThree#1#2{
\ii=#1                      
\jj=#2                      
\advance\jj by 1            
\loop   
   \ifnum\ii<\jj
{
   \hspace{-1.8cm}
   \enlargethispage{5cm}
   {\centering
   \includegraphics[width=18cm]{fig/con3catalog\ifnum\ii<100 0\fi\ifnum\ii<10 0\fi\number\ii.eps}
   }
   \newpage}
      \advance\ii by 1
   \repeat
}

\newpage

\setlength{\topmargin}{-1.2cm}

\producePagesThree{1}{39}
